\documentclass[a4paper,12pt]{article}

\usepackage[divide={1in,*,1in}]{geometry}
\usepackage{tabularx}
\usepackage{hhline}
\usepackage{enumitem}
\usepackage{slashbox}

\setlist{noitemsep,topsep=0pt,parsep=0pt,partopsep=0pt}
\setenumerate{label={\arabic*)}}

\usepackage{amsmath,amsthm}
\usepackage{amsfonts,amssymb}
\usepackage{mathtools}
\usepackage{stmaryrd}
\usepackage{BOONDOX-calo}

\usepackage{centerpict}
\usepackage{fntsymbols}
\usepackage{pst-3d}
\usepackage{diag-arc}
\usepackage{movies}
\usepackage{mywrap}
\usepackage{calc}

\usepackage{diagps}						% PS version of the previous one - if needed
\usepackage{diagps-arrows}		% To make in-text arrows compatible with those in diagrams

\usepackage[
		colorlinks,
		pdfauthor={Anna Beliakova, Krzysztof K. Putyra, Stephan M. Wehrli},
		pdftitle={Quantum Link Homology via Trace Functor I}]{hyperref}
\usepackage[all]{hypcap}		%	Fix hyperlinks to figures. Use \capstart for wrapfigure, etc.
\usepackage[savepos]{zref}

\newFntSymbol{big}{Circle}:|:{\pscircle[linewidth=0.5pt](0,0.05){0.45}}

\newFntSymbol{big}{NWSECr}[r]{%
		\psline(-0.5,-0.5)(0.5,0.5)
		\psline[border=4\pslinewidth](-0.5,0.5)(0.5,-0.5)
}%

\newFntSymbol{big}{NESWCr}[r]{%
		\psline(-0.5,0.5)(0.5,-0.5)
		\psline[border=4\pslinewidth](-0.5,-0.5)(0.5,0.5)
}%

\newFntSymbol{big}{VertRes}[r]{%
		\psbezier(-0.5,-0.5)(-0.08,0)(-0.08,0)(-0.5,0.5)
		\psbezier( 0.5,-0.5)( 0.08,0)( 0.08,0)( 0.5,0.5)
}%

\newFntSymbol{big}{HorRes}[r]{%
		\psbezier(-0.5,-0.5)(0,-0.08)(0,-0.08)(0.5,-0.5)
		\psbezier(-0.5, 0.5)(0, 0.08)(0, 0.08)(0.5, 0.5)
}%

\newFntSymbol{big}{HVCob}[r]{%
		\psbezier(-0.5,-0.5)(0,-0.08)(0,-0.08)(0.5,-0.5)
		\psbezier(-0.5, 0.5)(0, 0.08)(0, 0.08)(0.5, 0.5)
		\psline[linewidth=5\pslinewidth](0,-0.18)(0,0.18)
}%

\newFntSymbol{big}{VHCob}[r]{%
		\psbezier(-0.5,-0.5)(-0.08, 0)(-0.08, 0)(-0.5, 0.5)
		\psbezier( 0.5,-0.5)( 0.08, 0)( 0.08, 0)( 0.5, 0.5)
		\psline[linewidth=5\pslinewidth](-0.18, 0)(0.18, 0)
}%

\newFntSymbol{big}{PosCr}[r]{%
		\psline{->}(0.5,-0.5)(-0.55,0.55)
		\psline[border=4\pslinewidth]{->}(-0.5,-0.5)(0.55,0.55)
}%

\newFntSymbol{big}{NegCr}[r]{%
		\psline{->}(-0.5,-0.5)(0.55,0.55)
		\psline[border=4\pslinewidth]{->}(0.5,-0.5)(-0.55,0.55)
}%

\newFntSymbol{big}{NoCr}[r]{%
		\psbezier{->}(-0.5,-0.5)(-0.08,0)(-0.08,0)(-0.55,0.55)
		\psbezier{->}( 0.5,-0.5)( 0.08,0)( 0.08,0)( 0.55,0.55)
}%

\newFntSymbol{big}{PosCrArrowUp}[r]{%
		\psline(0.5,-0.5)(-0.55,0.55)
		\psline[border=5\pslinewidth](-0.5,-0.5)(0.55,0.55)
		\psline[linewidth=0.5\pslinewidth,border=1\pslinewidth]{->}(0,-0.5)(0,0.7)
}%

\newFntSymbol{big}{PosCrArrowDown}[r]{%
		\psline(0.5,-0.5)(-0.55,0.55)
		\psline[border=5\pslinewidth](-0.5,-0.5)(0.55,0.55)
		\psline[linewidth=0.5\pslinewidth,border=1\pslinewidth]{<-}(0,-0.7)(0,0.5)
}%

\def\fntDrawCupDiagram#1#2{%
	\psline{#1-c}(-0.45,0.35)(-0.45,0)
	\psarc[dimen=middle]{c-c}(0,0){0.45}{180}{360}
	\psline{c-#2}(0.45,0)(0.45,0.35)
}
\def\fntDrawCapDiagram#1#2{%
	\psline{#1-c}(0.45,-0.35)(0.45,0)
	\psarc[dimen=middle]{c-c}(0,0){0.45}{0}{180}
	\psline{c-#2}(-0.45,0)(-0.45,-0.35)
}

\newFntSymbol{big}{Cup}[r]{\fntDrawCupDiagram{}{}}
\newFntSymbol{big}{LCup}[r]{\fntDrawCupDiagram{<}{}}
\newFntSymbol{big}{RCup}[r]{\fntDrawCupDiagram{}{>}}

\newFntSymbol{big}{Cap}[r]{\fntDrawCapDiagram{}{}}
\newFntSymbol{big}{LCap}[r]{\fntDrawCapDiagram{}{>}}
\newFntSymbol{big}{RCap}[r]{\fntDrawCapDiagram{<}{}}

\newFntSymbol{tangle}{FullTwist}:|:<1.25>[c]{%
	\psclip{\psframe[linestyle=none](0,-1)(1,1)}
		\psbezier(-0.433, 0.25)(-0.15,-0.33)(0.15,-0.33)(0.433, 0.25)
	\endpsclip
	\psbezier[border=3\pslinewidth](-0.433,-0.25)(-0.15, 0.33)(0.15, 0.33)(0.433,-0.25)
	\psclip{\psframe[linestyle=none](0.1,-1)(-1,1)}
		\psbezier[border=3\pslinewidth](-0.433, 0.25)(-0.15,-0.33)(0.15,-0.33)(0.433, 0.25)
	\endpsclip
}

\newFntSymbol{tangle}{PosTwist}(-1,-0.5)(1,0.5)[r]{%
	\psbezier{C-C}( 0.5,-0.5)( 0.2,-0.5)(-0.2,0.5)(-0.5,0.5)
	\psbezier[border=3\pslinewidth]{C-C}(-0.5,-0.5)(-0.2,-0.5)( 0.2,0.5)( 0.5,0.5)
	\psarc{C-C}(-0.5,0){0.5}{90}{270}%
	\psarc{C-C}(0.5,0){0.5}{-90}{90}%
}

\newFntSymbol{tangle}{PosTwistHor}(-1,-0.5)(1,0.5)[r]{\pscustom{%
	\moveto(-0.5,0.5)
	\psarc(-0.5,0){0.5}{90}{270}\curveto(-0.3,-0.5)(-0.2,-0.2)(0,-0.2)\curveto( 0.2,-0.2)( 0.3,-0.5)( 0.5,-0.5)
	\psarc( 0.5,0){0.5}{-90}{90}\curveto( 0.3, 0.5)( 0.2, 0.2)(0, 0.2)\curveto(-0.2, 0.2)(-0.3, 0.5)(-0.5, 0.5)
	\stroke
}}

\newFntSymbol{tangle}{PosTwistVert}(-1,-0.5)(1,0.5)[r]{%
	\psarc{C-C}( 0.5,0){0.5}{-90}{90}\psellipticarc{C-C}( 0.5,0)(0.4,0.5){90}{270}%
	\psarc{C-C}(-0.5,0){0.5}{90}{270}\psellipticarc{C-C}(-0.5,0)(0.4,0.5){-90}{90}%
}

\newFntSymbol{tangle}{RIline}:|:<1.25>[c]{%
	\psbezier(-0.35355, 0.35355)(-0.25355, 0.15355)(0, 0.2)(0,0)
	\psbezier(-0.35355,-0.35355)(-0.25355,-0.15355)(0,-0.2)(0,0)
}

\newFntSymbol{tangle}{RIx}:|:<1.25>[c]{%
	\psbezier(-0.35355, 0.35355)(-0.25355, 0.15355)(-0.1,-0.2)(0.1,-0.2)
	\psbezier[border=3\pslinewidth](-0.35355,-0.35355)(-0.25355,-0.15355)(-0.1, 0.2)(0.1, 0.2)
	\psbezier(0.1,0.2)(0.4,0.2)(0.4,-0.2)(0.1,-0.2)
}
	
\newFntSymbol{tangle}{RIh}:|:<1.25>[c]{%
	\psbezier(-0.35355, 0.35355)(-0.2, 0)(-0.1, 0.2)(0.1, 0.2)
	\psbezier(0.1,0.2)(0.4,0.2)(0.4,-0.2)(0.1,-0.2)
	\psbezier(-0.35355,-0.35355)(-0.2, 0)(-0.1,-0.2)(0.1,-0.2)
}

\newFntSymbol{tangle}{RIv}:|:<1.25>[c]{%
	\psbezier(-0.35355, 0.35355)(-0.2,0)(-0.2,0)(-0.35355,-0.35355)
	\pscircle(0.1,0){0.2}
}

\newFntSymbol{tangle}{RIcob}:|:<1.25>[c]{%
	\psbezier(-0.35355, 0.35355)(-0.2,0)(-0.2,0)(-0.35355,-0.35355)
	\pscircle(0.1,0){0.2}
	\psline[linewidth=3\pslinewidth](-0.24,0)(-0.1,0)
}

\newFntSymbol{tangle}{RIIxx}:|:<1.25>[c]{%
	\psbezier(-0.433, 0.25)(-0.15,-0.33)(0.15,-0.33)(0.433, 0.25)
	\psbezier[border=3\pslinewidth](-0.433,-0.25)(-0.15, 0.33)(0.15, 0.33)(0.433,-0.25)
}

\newFntSymbol{tangle}{RIIhh}:|:<1.25>[c]{%
	\psbezier(-0.433, 0.25)(-0.3,-0.1)(-0.2, 0.2)(0, 0.2)
	\psbezier(-0.433,-0.25)(-0.3, 0.1)(-0.2,-0.2)(0,-0.2)
	\psbezier(0.433, 0.25)(0.3,-0.1)(0.2, 0.2)(0, 0.2)
	\psbezier(0.433,-0.25)(0.3, 0.1)(0.2,-0.2)(0,-0.2)
}

\newFntSymbol{tangle}{RIIhv}:|:<1.25>[c]{%
	\psbezier(0.433,0.25)(0.3,0)(0.3,0)(0.433,-0.25)
	\psarc(0,0){0.2}{-90}{90}
	\psbezier(-0.433, 0.25)(-0.3,-0.1)(-0.2, 0.2)(0, 0.2)
	\psbezier(-0.433,-0.25)(-0.3, 0.1)(-0.2,-0.2)(0,-0.2)
}

\newFntSymbol{tangle}{RIIvh}:|:<1.25>[c]{%
	\psbezier(-0.433,0.25)(-0.3,0)(-0.3,0)(-0.433,-0.25)
	\psarc(0,0){0.2}{90}{270}
	\psbezier(0.433, 0.25)(0.3,-0.1)(0.2, 0.2)(0, 0.2)
	\psbezier(0.433,-0.25)(0.3, 0.1)(0.2,-0.2)(0,-0.2)
}

\newFntSymbol{tangle}{RIIvv}:|:<1.25>[c]{%
	\psbezier(-0.433,0.25)(-0.3,0)(-0.3,0)(-0.433,-0.25)
	\pscircle(0,0){0.2}
	\psbezier(0.433,0.25)(0.3,0)(0.3,0)(0.433,-0.25)
}

\newFntSymbol{tangle}{RIIhhs}:|:<1.25>[c]{%
	\psbezier(-0.433, 0.25)(-0.2,0)(0.2,0)(0.433, 0.25)
	\psbezier(-0.433,-0.25)(-0.2,0)(0.2,0)(0.433,-0.25)
}

\newFntSymbol{tangle}{RIIIax}:|:<1.25>[c]{%
	\psline(-0.35355,-0.35355)(0.35355,0.35355)
	\psline[border=2\pslinewidth](-0.35355,0.35355)(0.35355,-0.35355)
	\psbezier[border=2\pslinewidth](-0.5,0)(-0.25,0)(-0.2,0.35)(0,0.35)
	\psbezier[border=2\pslinewidth]( 0.5,0)( 0.25,0)( 0.2,0.35)(0,0.35)
}

\newFntSymbol{tangle}{RIIIah}:|:<1.25>[c]{%
	\psbezier(-0.35355,-0.35355)(0,0)(0,0)( 0.35355,-0.35355)
	\psbezier(-0.35355, 0.35355)(0,0)(0,0)( 0.35355, 0.35355)
	\psbezier[border=2\pslinewidth](-0.5,0)(-0.25,0)(-0.2,0.35)(0,0.35)
	\psbezier[border=2\pslinewidth]( 0.5,0)( 0.25,0)( 0.2,0.35)(0,0.35)
}

\newFntSymbol{tangle}{RIIIav}:|:<1.25>[c]{%
	\psbezier(-0.35355,-0.35355)(0,0)(0,0)(-0.35355, 0.35355)
	\psbezier( 0.35355, 0.35355)(0,0)(0,0)( 0.35355,-0.35355)
	\psbezier[border=2\pslinewidth](-0.5,0)(-0.25,0)(-0.2,0.35)(0,0.35)
	\psbezier[border=2\pslinewidth]( 0.5,0)( 0.25,0)( 0.2,0.35)(0,0.35)
}

\newFntSymbol{tangle}{RIIIac}:|:<1.25>[c]{%
	\psbezier(-0.35355,-0.35355)(0,0)(0,0)( 0.35355,-0.35355)
	\psbezier(-0.35355, 0.35355)(0,0)(0,0)( 0.35355, 0.35355)
	\psline[linewidth=5\pslinewidth](0,-0.08839)(0,0.08839)
	\psbezier[border=2\pslinewidth](-0.5,0)(-0.25,0)(-0.2,0.35)(0,0.35)
	\psbezier[border=2\pslinewidth]( 0.5,0)( 0.25,0)( 0.2,0.35)(0,0.35)
}

\newFntSymbol{tangle}{RIIIahhh}:|:<1.25>[c]{%
	\psbezier(-0.35355,-0.35355)(0.05, 0.05)(-0.05, 0.05)( 0.35355,-0.35355)
	\psbezier(-0.35355, 0.35355)(-0.2,0.2)(-0.25,0)(-0.5,0)
	\psarc(0,0.2){0.15}{45}{315}
	\psbezier(0.10607,0.30607)(0.15607,0.25607)(0.2,0.2)(0.35355,0.35355)
	\psbezier(0.10607,0.09393)(0.25607,0.24393)(0.3,0)( 0.5,0)
}

\newFntSymbol{tangle}{RIIIahhv}:|:<1.25>[c]{%
	\psbezier(-0.35355,-0.35355)(0.05, 0.05)(-0.05, 0.05)( 0.35355,-0.35355)
	\pscircle(0,0.2){0.15}
	\psbezier(-0.35355, 0.35355)(-0.2,0.2)(-0.25,0)(-0.5,0)
	\psbezier( 0.35355, 0.35355)( 0.2,0.2)( 0.25,0)( 0.5,0)
}

\newFntSymbol{tangle}{RIIIahvh}:|:<1.25>[c]{%
	\psbezier(-0.35355,-0.35355)(0.05, 0.05)(-0.05, 0.05)( 0.35355,-0.35355)
	\psarc(0,0.2){0.15}{45}{135}
	\psarc(0,0.2){0.15}{225}{315}
	\psbezier(-0.10607,0.30607)(-0.15607,0.25607)(-0.2,0.2)(-0.35355,0.35355)
	\psbezier(-0.10607,0.09393)(-0.25607,0.24393)(-0.3,0)(-0.5,0)
	\psbezier(0.10607,0.30607)(0.15607,0.25607)(0.2,0.2)(0.35355,0.35355)
	\psbezier(0.10607,0.09393)(0.25607,0.24393)(0.3,0)( 0.5,0)
}

\newFntSymbol{tangle}{RIIIahvv}:|:<1.25>[c]{%
	\psbezier(-0.35355,-0.35355)(0.05, 0.05)(-0.05, 0.05)( 0.35355,-0.35355)
	\psbezier( 0.35355, 0.35355)(0.2,0.2)(0.25,0)(0.5,0)
	\psarc(0,0.2){0.15}{-135}{135}
	\psbezier(-0.10607,0.30607)(-0.15607,0.25607)(-0.2,0.2)(-0.35355,0.35355)
	\psbezier(-0.10607,0.09393)(-0.25607,0.24393)(-0.3,0)(-0.5,0)
}

\newFntSymbol{tangle}{RIIIavhh}:|:<1.25>[c]{%
	\psline(-0.35355,-0.35355)(-0.10607,-0.10607)
	\psline( 0.35355,-0.35355)( 0.10607,-0.10607)
	\psarc(-0.20607,-0.00607){0.14142}{-45}{45}
	\psarc( 0.20607,-0.00607){0.14142}{135}{225}
	\psarc(0,0.2){0.15}{45}{225}
	\psbezier(0.10607,0.30607)(0.15607,0.25607)(0.2,0.2)(0.35355,0.35355)
	\psbezier(0.10607,0.09393)(0.25607,0.24393)(0.3,0)( 0.5,0)
	\psbezier(-0.35355, 0.35355)(-0.2,0.2)(-0.25,0)(-0.5,0)
}

\newFntSymbol{tangle}{RIIIavhv}:|:<1.25>[c]{%
	\psline(-0.35355,-0.35355)(-0.10607,-0.10607)
	\psline( 0.35355,-0.35355)( 0.10607,-0.10607)
	\psarc(-0.20607,-0.00607){0.14142}{-45}{45}
	\psarc( 0.20607,-0.00607){0.14142}{135}{225}
	\psarc(0,0.2){0.15}{-45}{225}
	\psbezier(-0.35355, 0.35355)(-0.2,0.2)(-0.25,0)(-0.5,0)
	\psbezier( 0.35355, 0.35355)( 0.2,0.2)( 0.25,0)( 0.5,0)
}

\newFntSymbol{tangle}{RIIIavvh}:|:<1.25>[c]{%
	\psline(-0.35355,-0.35355)(-0.10607,-0.10607)
	\psline( 0.35355,-0.35355)( 0.10607,-0.10607)
	\psarc(-0.20607,-0.00607){0.14142}{-45}{45}
	\psarc( 0.20607,-0.00607){0.14142}{135}{225}
	\psarc(0,0.2){0.15}{45}{135}
	\psbezier(-0.10607,0.30607)(-0.15607,0.25607)(-0.2,0.2)(-0.35355,0.35355)
	\psbezier(-0.10607,0.09393)(-0.25607,0.24393)(-0.3,0)(-0.5,0)
	\psbezier(0.10607,0.30607)(0.15607,0.25607)(0.2,0.2)(0.35355,0.35355)
	\psbezier(0.10607,0.09393)(0.25607,0.24393)(0.3,0)( 0.5,0)
}

\newFntSymbol{tangle}{RIIIavvv}:|:<1.25>[c]{%
	\psline(-0.35355,-0.35355)(-0.10607,-0.10607)
	\psline( 0.35355,-0.35355)( 0.10607,-0.10607)
	\psarc(-0.20607,-0.00607){0.14142}{-45}{45}
	\psarc( 0.20607,-0.00607){0.14142}{135}{225}
	\psarc(0,0.2){0.15}{-45}{135}
	\psbezier( 0.35355, 0.35355)(0.2,0.2)(0.25,0)(0.5,0)
	\psbezier(-0.10607,0.30607)(-0.15607,0.25607)(-0.2,0.2)(-0.35355,0.35355)
	\psbezier(-0.10607,0.09393)(-0.25607,0.24393)(-0.3,0)(-0.5,0)
}

\newFntSymbol{tangle}{RIIIbx}:|:<1.25>[c]{%
	\psline(-0.35355,-0.35355)(0.35355,0.35355)
	\psline[border=2\pslinewidth](-0.35355,0.35355)(0.35355,-0.35355)
	\psbezier[border=2\pslinewidth](-0.5,0)(-0.25,0)(-0.2,-0.35)(0,-0.35)
	\psbezier[border=2\pslinewidth]( 0.5,0)( 0.25,0)( 0.2,-0.35)(0,-0.35)
}

\newFntSymbol{tangle}{RIIIbh}:|:<1.25>[c]{%
	\psbezier(-0.35355,-0.35355)(0,0)(0,0)( 0.35355,-0.35355)
	\psbezier(-0.35355, 0.35355)(0,0)(0,0)( 0.35355, 0.35355)
	\psbezier[border=2\pslinewidth](-0.5,0)(-0.25,0)(-0.2,-0.35)(0,-0.35)
	\psbezier[border=2\pslinewidth]( 0.5,0)( 0.25,0)( 0.2,-0.35)(0,-0.35)
}

\newFntSymbol{tangle}{RIIIbv}:|:<1.25>[c]{%
	\psbezier(-0.35355,-0.35355)(0,0)(0,0)(-0.35355, 0.35355)
	\psbezier( 0.35355, 0.35355)(0,0)(0,0)( 0.35355,-0.35355)
	\psbezier[border=2\pslinewidth](-0.5,0)(-0.25,0)(-0.2,-0.35)(0,-0.35)
	\psbezier[border=2\pslinewidth]( 0.5,0)( 0.25,0)( 0.2,-0.35)(0,-0.35)
}

\newFntSymbol{tangle}{RIIIbc}:|:<1.25>[c]{%
	\psbezier(-0.35355,-0.35355)(0,0)(0,0)( 0.35355,-0.35355)
	\psbezier(-0.35355, 0.35355)(0,0)(0,0)( 0.35355, 0.35355)
	\psline[linewidth=5\pslinewidth](0,-0.08839)(0,0.08839)
	\psbezier[border=2\pslinewidth](-0.5,0)(-0.25,0)(-0.2,-0.35)(0,-0.35)
	\psbezier[border=2\pslinewidth]( 0.5,0)( 0.25,0)( 0.2,-0.35)(0,-0.35)
}

\newFntSymbol{tangle}{IIIVert}:|:<1.25>[c]{%
	\psbezier(-0.35355,0.35355)(-0.08,0.08)(0.08,0.08)(0.35355,0.35355)
	\psline(-0.5,0)(0.5,0)
	\psbezier(-0.35355,-0.35355)(-0.08,-0.08)(0.08,-0.08)(0.35355,-0.35355)
}

\newFntSymbol{tangle}{IIIVertT}:|:<1.25>[c]{%
	\psbezier(-0.35355,0.35355)(-0.08,0.08)(0.08,0.08)(0.35355,0.35355)
	\psline(-0.5,0)(0.5,0)
	\psbezier(-0.35355,-0.35355)(-0.08,-0.08)(0.08,-0.08)(0.35355,-0.35355)
	\psline[linewidth=3\pslinewidth](0,0)(0,0.14839)
}

\newFntSymbol{tangle}{IIIVertB}:|:<1.25>[c]{%
	\psbezier(-0.35355,0.35355)(-0.08,0.08)(0.08,0.08)(0.35355,0.35355)
	\psline(-0.5,0)(0.5,0)
	\psbezier(-0.35355,-0.35355)(-0.08,-0.08)(0.08,-0.08)(0.35355,-0.35355)
	\psline[linewidth=3\pslinewidth](0,0)(0,-0.14839)
}

\newFntSymbol{tangle}{CircLeft}(-0.375,-0.5)(0.5,0.5):|:<1.25>[r]{%
	\psline(0.375,-0.5)(0.375,0.5)
	\pscircle(0,0){0.25}
}
\newFntSymbol{tangle}{CircRight}(-0.5,-0.5)(0.375,0.5):|:<1.25>[r]{%
	\psline(-0.375,-0.5)(-0.375,0.5)
	\pscircle(0,0){0.25}
}
\newFntSymbol{tangle}{CircOver}(-0.375,-0.5)(0.375,0.5):|:<1.25>[r]{%
	\psline(0,-0.5)(0,0.5)
	\pscircle[border=3\pslinewidth](0,0){0.25}
}
\newFntSymbol{tangle}{CircUnder}(-0.375,-0.5)(0.375,0.5):|:<1.25>[r]{%
	\pscircle(0,0){0.25}
	\psline[border=3\pslinewidth](0,-0.5)(0,0.5)
}

\newFntSymbol{tangle}{CircOverhh}(-0.375,-0.5)(0.375,0.5):|:<1.25>[r]{%
	\pscustom{\moveto(0,-0.5)
		\curveto(0,-0.125)(0.25,-0.3)(0.25,0)
		\curveto(0.25,0.15)(0,0.35)(0,0)
		\curveto(0,-0.35)(-0.25,-0.15)(-0.25,0)
		\curveto(-0.25,0.3)(0,0.125)(0,0.5)
		\stroke
	}
}
\newFntSymbol{tangle}{CircOverhv}(-0.375,-0.5)(0.375,0.5):|:<1.25>[r]{%
	\pscustom{\moveto(0,-0.5)\curveto(0,-0.125)(0.25,-0.3)(0.25,0)\curveto(0.25,0.3)(0,0.125)(0,0.5)\stroke}%
	\pscustom{\moveto(0,0)\curveto(0,-0.35)(-0.25,-0.15)(-0.25,0)\curveto(-0.25,0.15)(0,0.35)(0,0)\stroke}%
}
\newFntSymbol{tangle}{CircOvervh}(-0.375,-0.5)(0.375,0.5):|:<1.25>[r]{%
	\pscustom{\moveto(0,-0.5)\curveto(0,-0.125)(-0.25,-0.3)(-0.25,0)\curveto(-0.25,0.3)(0,0.125)(0,0.5)\stroke}%
	\pscustom{\moveto(0,0)\curveto(0,-0.35)(0.25,-0.15)(0.25,0)\curveto(0.25,0.15)(0,0.35)(0,0)\stroke}%
}
\newFntSymbol{tangle}{CircOvervv}(-0.375,-0.5)(0.375,0.5):|:<1.25>[r]{%
	\pscustom{\moveto(0,-0.5)
		\curveto(0,-0.125)(-0.25,-0.3)(-0.25,0)
		\curveto(-0.25,0.15)(0,0.35)(0,0)
		\curveto(0,-0.35)(0.25,-0.15)(0.25,0)
		\curveto(0.25,0.3)(0,0.125)(0,0.5)
		\stroke
	}
}

\def\pictDeloopTR{\begin{pspicture}(0,0)(0.4,0.6)
	\psellipse(0.2,0.5)(0.2,0.07)
	\psbezier(0,0.5)(0,0)(0.4,0)(0.4,0.5)
	\psdot(0.2,0.3)
\end{pspicture}}

\def\pictDeloopTL{\begin{pspicture}(0,0)(0.4,0.6)
	\psellipticarc(0.2,0.1)(0.2,0.07){-180}{0}
	\psdash\psellipticarc(0.2,0.1)(0.2,0.07){0}{180}
	\psbezier(0.4,0.1)(0.4,0.6)(0,0.6)(0,0.1)
\end{pspicture}}

\def\pictDeloopBR{\begin{pspicture}(0,0)(0.4,0.6)
	\psellipse(0.2,0.5)(0.2,0.07)
	\psbezier(0,0.5)(0,0)(0.4,0)(0.4,0.5)
\end{pspicture}}

\def\pictDeloopBLa{\begin{pspicture}[shift=-0.25](-0.1,0)(0.5,0.5)
	\psellipticarc(0.2,0.1)(0.2,0.07){-180}{0}
	\psdash\psellipticarc(0.2,0.1)(0.2,0.07){0}{180}
	\psbezier(0.4,0.1)(0.4,0.6)(0,0.6)(0,0.1)
	\psdot(0.2,0.3)
\end{pspicture}}

\def\pictDeloopBLb{\begin{pspicture}[shift=-0.25](-0.1,0)(0.5,0.5)
	\psellipticarc(0.2,0.1)(0.2,0.07){-180}{0}
	\psdash\psellipticarc(0.2,0.1)(0.2,0.07){0}{180}
	\psbezier(0.4,0.1)(0.4,0.4)(0,0.4)(0,0.1)
	\pscircle(0.2,0.6){0.2}
	\psdot(0.2,0.52)
	\psdot(0.2,0.68)
\end{pspicture}}

\def\psdash#1{#1[linestyle=dashed,linewidth=0.15pt]}

\def\pictRelS{\begin{centerpict}(-0.6,-0.5)(0.6,0.5)
	\psset{linewidth=0.5pt,dash=1pt 1.5pt,dimen=middle}%
	\pscircle(0,0){0.5}
	\psellipticarc(0,0)(0.5,0.15){-180}{0}
	\psdash\psellipticarc(0,0)(0.5,0.15){0}{180}
\end{centerpict}}

\def\pictRelT{\begin{centerpict}(-0.1,-0.7)(1.3,0.7)
	\psset{linewidth=0.5pt,dash=1pt 1.5pt,dimen=middle}
	\psdash\psellipticarc(0.2,0)(0.2,0.1){0}{180}
	\psdash\psellipticarc(1,0)(0.2,0.1){0}{180}
	\psellipticarc(0.2,0)(0.2,0.1){-180}{0}
	\psellipticarc(1,0)(0.2,0.1){-180}{0}
	\psbezier(0.0,0)(0.00,0.9)(1.20,0.9)(1.2,0)
	\psbezier(0.4,0)(0.35,0.4)(0.85,0.4)(0.8,0)
	\psbezier(0.0,0)(0.00,-0.9)(1.20,-0.9)(1.2,0)
	\psbezier(0.4,0)(0.35,-0.4)(0.85,-0.4)(0.8,0)
\end{centerpict}}

\def\pictRelTuCircles{%
	\psset{linewidth=0.5pt,dash=1pt 1.5pt,dimen=middle}
	\psellipticarc(0.2,0.1)(0.2,0.07\psxunit){-180}{0}
	\psellipticarc(1.0,0.1)(0.2,0.07\psxunit){-180}{0}
	\psdash\psellipticarc(0.2,0.1)(0.2,0.07\psxunit){0}{180}
	\psdash\psellipticarc(1.0,0.1)(0.2,0.07\psxunit){0}{180}
	\psellipse(0.2,1.4)(0.2,0.07\psxunit)
	\psellipse(1.0,1.4)(0.2,0.07\psxunit)
}

\def\pictRelTuL{\begin{centerpict}(-0.1,0)(1.3,1.5)
	\pictRelTuCircles
	\psline(0.0,0.1)(0.0,1.4)\psline(0.4,0.1)(0.4,1.4)
	\psbezier(0.8,0.1)(0.8,0.7)(1.2,0.7)(1.2,0.1)
	\psbezier(0.8,1.4)(0.8,0.8)(1.2,0.8)(1.2,1.4)
\end{centerpict}}

\def\pictRelTuR{\begin{centerpict}(-0.1,0)(1.3,1.5)
	\pictRelTuCircles
	\psbezier(0.0,0.1)(0.0,0.7)(0.4,0.7)(0.4,0.1)
	\psbezier(0.0,1.4)(0.0,0.8)(0.4,0.8)(0.4,1.4)
	\psline(0.8,0.1)(0.8,1.4)\psline(1.2,0.1)(1.2,1.4)
\end{centerpict}}

\def\pictRelTuB{\begin{centerpict}(-0.1,0)(1.3,1.5)
	\pictRelTuCircles
	\psbezier(0.0,0.1)(0.0,1.0)(1.2,1.0)(1.2,0.1)
	\psbezier(0.4,0.1)(0.35,0.5)(0.85,0.5)(0.8,0.1)
	\psbezier(0.0,1.4)(0.0,1.0)(0.4,1.0)(0.4,1.4)
	\psbezier(0.8,1.4)(0.8,0.8)(1.2,0.8)(1.2,1.4)
\end{centerpict}}

\def\pictRelTuT{\begin{centerpict}(-0.1,0)(1.3,1.5)
	\pictRelTuCircles
	\psbezier(0.0,1.4)(0.0,0.5)(1.2,0.5)(1.2,1.4)
	\psbezier(0.4,1.4)(0.35,1.0)(0.85,1.0)(0.8,1.4)
	\psbezier(0.0,0.1)(0.0,0.5)(0.4,0.5)(0.4,0.1)
	\psbezier(0.8,0.1)(0.8,0.7)(1.2,0.7)(1.2,0.1)
\end{centerpict}}

\def\pictRelD{\begin{centerpict}(-0.6,-0.5)(0.6,0.5)
	\psset{linewidth=0.5pt,dash=1pt 1.5pt,dimen=middle,dotsize=4pt}%
	\pscircle(0,0){0.5}
	\psellipticarc(0,0)(0.5,0.15){-180}{0}
	\psdash\psellipticarc(0,0)(0.5,0.15){0}{180}
	\psdot(0,0.15)
\end{centerpict}}

\def\pictRelNeckI{\begin{centerpict}(-0.6,0)(0.6,2)
	\psset{linewidth=0.5pt,dash=1pt 1.5pt,dimen=middle,dotsize=4pt}%
	\psellipticarc(0,0.1)(0.4,0.1){-180}{0}
	\psdash\psellipticarc(0,0.1)(0.4,0.1){0}{180}
	\psellipse(0,1.9)(0.4,0.1)
	\psline(-0.4,0.1)(-0.4,1.9)
	\psline( 0.4,0.1)( 0.4,1.9)
\end{centerpict}}

\def\pictRelNeckB{\begin{centerpict}(-0.6,0)(0.6,2)
	\psset{linewidth=0.5pt,dash=1pt 1.5pt,dimen=middle,dotsize=4pt}%
	\psellipticarc(0,0.1)(0.4,0.1){-180}{0}
	\psdash\psellipticarc(0,0.1)(0.4,0.1){0}{180}
	\psellipse(0,1.9)(0.4,0.1)
	\psbezier(-0.4,0.1)(-0.4,0.9)(0.4,0.9)(0.4,0.1)
	\psbezier(-0.4,1.9)(-0.4,1.1)(0.4,1.1)(0.4,1.9)
	\psdot(0,0.4)
\end{centerpict}}

\def\pictRelNeckT{\begin{centerpict}(-0.6,0)(0.6,2)
	\psset{linewidth=0.5pt,dash=1pt 1.5pt,dimen=middle,dotsize=4pt}%
	\psellipticarc(0,0.1)(0.4,0.1){-180}{0}
	\psdash\psellipticarc(0,0.1)(0.4,0.1){0}{180}
	\psellipse(0,1.9)(0.4,0.1)
	\psbezier(-0.4,0.1)(-0.4,0.9)(0.4,0.9)(0.4,0.1)
	\psbezier(-0.4,1.9)(-0.4,1.1)(0.4,1.1)(0.4,1.9)
	\psdot(0,1.6)
\end{centerpict}}

\def\pictRelNeckM{\begin{centerpict}(-0.6,0)(0.6,2)
	\psset{linewidth=0.5pt,dash=1pt 1.5pt,dimen=middle,dotsize=4pt}%
	\psellipticarc(0,0.1)(0.4,0.1){-180}{0}
	\psdash\psellipticarc(0,0.1)(0.4,0.1){0}{180}
	\psellipticarc(0,1.0)(0.4,0.1){-180}{0}
	\psdash\psellipticarc(0,1.0)(0.4,0.1){0}{180}
	\psellipse(0,1.9)(0.4,0.1)
	\psbezier(-0.4,0.1)(-0.4,0.6)(0.4,0.6)(0.4,0.1)
	\psbezier(-0.4,1.9)(-0.4,1.4)(0.4,1.4)(0.4,1.9)
	\pscircle(0,1){0.4}
	\psdot(0,1.15)
	\psdot(0,0.75)
\end{centerpict}}

\def\pictRelTwoDots{\begin{centerpict}(0,0)(3,1)
	\psset{linewidth=0.5pt,dash=1pt 1.5pt,dimen=middle,dotsize=4pt}%
	\psline(0,0)(2,0)(3,1)(1,1)(0,0)
	\psdot(1.25,0.5)
	\psdot(1.75,0.5)
\end{centerpict}}

\def\pictRelB{\begin{centerpict}(-1.2,0)(1,2)
	\psset{linewidth=0.5pt,dash=1pt 1.5pt,dimen=middle,dotsize=4pt}%
	\psellipticarc(0,0.1)(0.6,0.1){-180}{0}
	\psdash\psellipticarc(0,0.1)(0.6,0.1){0}{180}
	\psellipse(0,1.9)(0.6,0.1)
	\psline(-0.6,0.1)(-0.6,1.9)
	\psline( 0.6,0.1)( 0.6,1.9)
	\psset{linewidth=1pt}%
	\psellipticarc(0,1)(0.6,0.1){-180}{0}
	\psellipticarc[linestyle=dashed,dash=2pt 3pt](0,1)(0.6,0.1){0}{180}
	\uput[l](-0.6,1){$\scriptstyle\gamma$}
	\psdot(0.2,1.4)
\end{centerpict}}

\newrgbcolor{membraneColor}{0.75 0.75 1}%
\def\membrane{%
	\begin{psset}{linewidth=0.7pt,linecolor=blue}%
		\psline{<-}(-0.0, 1.9)( 0.6, 1.9)
		\psline{<-}(-0.3, 1.7)( 0.3, 1.7)
		\psline{<-}(-0.6, 1.5)( 0.0, 1.5)
	\end{psset}%
	\pspolygon[linestyle=none,fillstyle=solid,fillcolor=membraneColor](-0.7,-2)(-0.7,1)(0.7,2)(0.7,-1)
}%

\def\cutpicture(#1,#2)(#3,#4)@(#5--#6)#7{%
	\begin{psclip}{\psframe[linestyle=none](#1,#2)(#5,#4)}\pscustom{#7}\end{psclip}%
	\ifx\relax#6\relax
		\begin{psclip}{\psframe[linestyle=none](#5,#2)(#3,#4)}\pscustom[linestyle=dashed]{#7}\end{psclip}%
	\else
		\begin{psclip}{\psframe[linestyle=none](#5,#2)(#6,#4)}\pscustom[linestyle=dashed]{#7}\end{psclip}%
		\begin{psclip}{\psframe[linestyle=none](#6,#2)(#3,#4)}\pscustom[fillstyle=solid,fillcolor=white]{#7}\end{psclip}%
	\fi
}

\def\frontellipse(#1){\psellipse(#1)(0.2,0.5)}
\def\backellipse(#1){%
	\psellipse[linestyle=none,fillstyle=solid,fillcolor=white](#1)(0.2,0.5)
	\fronthalfellipse(#1)
	\backhalfellipse(#1)
}
\def\fronthalfellipse(#1){\psellipticarc(#1)(0.2,0.5){90}{270}}
\def\backhalfellipse(#1){\psellipticarc[linestyle=dashed](#1)(0.2,0.5){-90}{90}}

\def\relRCapAtMembrane{\begin{centerpict}(-1.6,-2)(1,2)
	\psset{dimen=middle,dash=2pt 3pt}%
	\membrane
	\cutpicture(-2,-1)(1,1)@(-0.7--0){%
		\moveto(-1.2,-0.5)
		\lineto(0,-0.5)
		\curveto(1.2,-0.5)(1.2,0.5)(0,0.5)
		\lineto(-1.2, 0.5)
	}%
	\backellipse(0,0)
	\backellipse(-1.2,0)
\end{centerpict}}

\def\relRCapOffMembrane{\begin{centerpict}(-1.8,-2)(0.7,2)
	\psset{dimen=middle,dash=2pt 3pt}%
	\membrane
	\backhalfellipse(-1.5,0)
	\cutpicture(-2,-1)(0,1)@(-0.7--){%
		\moveto(-1.5,-0.5)
		\lineto(-1.2,-0.5)
		\curveto(0,-0.5)(0,0.5)(-1.2,0.5)
		\lineto(-1.5, 0.5)
		\fronthalfellipse(-1.5,0)
	}%
\end{centerpict}}

\def\relLCapAtMembrane{\begin{centerpict}(-1.2,-2)(1.3,2)
	\psset{dimen=middle,dash=2pt 3pt}%
	\membrane
	\cutpicture(-1,-1)(1,1)@(-0.7--0){%
		\moveto(1,0.5)
		\lineto(0,0.5)
		\curveto(-1.2,0.5)(-1.2,-0.5)(0,-0.5)
		\lineto(1,-0.5)
	}%
	\backellipse(0,0)
	\frontellipse(1,0)
\end{centerpict}}

\def\relLCapOffMembrane{\begin{centerpict}(-0.9,-2)(1.7,2)
	\psset{dimen=middle,dash=2pt 3pt}%
	\membrane
	\psbezier[fillstyle=solid,fillcolor=white]{c-c}(0.8,0.5)(-0.2,0.5)(-0.2,-0.5)(0.8,-0.5)
	\psline{c-c}(1.2, 0.5)(0.8, 0.5)
	\psline{c-c}(1.2,-0.5)(0.8,-0.5)
	\frontellipse(1.2,0)
\end{centerpict}}

\def\relDotLeftToMembrane{\begin{centerpict}(-1.8,-2)(1.4,2)
	\psset{dimen=middle,dash=2pt 3pt}%
	\membrane
	\pscustom[fillstyle=solid,fillcolor=white]{%
		\moveto(0.5,-0.8)\curveto(0.7,-0.2)(0.9,0.5)(1.4,0.9)
		\lineto(0.6,0.9)\curveto(0.1,0.5)(-0.1,-0.2)(-0.3,-0.8)
		\closepath
	}
	\psbezier(-1.8,-0.8)(-1.6,-0.2)(-1.4,0.5)(-0.9,0.9)
	\psline[linestyle=dashed](-0.3,-0.8)(-0.7,-0.8)\psline(-0.7,-0.8)(-1.8,-0.8)
	\psline[linestyle=dashed]( 0.7, 0.9)(-0.7, 0.9)\psline(-0.7, 0.9)(-0.9, 0.9)
	\psdot(-1,0)
\end{centerpict}}

\def\relDotRightToMembrane{\begin{centerpict}(-1.4,-2)(1.8,2)
	\psset{dimen=middle,dash=2pt 3pt}%
	\membrane
	\pscustom[fillstyle=solid,fillcolor=white]{%
		\moveto(0.9,-0.8)\curveto(1.1,-0.2)(1.3,0.5)(1.8,0.9)
		\lineto(0.6,0.9)\curveto(0.1,0.5)(-0.1,-0.2)(-0.3,-0.8)
		\closepath
	}%
	\cutpicture(-2,-2)(0.7,2)@(-0.7--){%
		\moveto(-0.3,-0.8)
		\lineto(-1.4,-0.8)
		\curveto(-1.2,-0.2)(-1.0,0.5)(-0.5,0.9)
		\lineto( 0.7, 0.9)
	}%
	\psdot(0.6,0)
\end{centerpict}}

\def\relSaddleBackArcs{%
	\psclip{\psframe[linestyle=none](-0.5,1.5)(0.5,0.7)}
		\psbezier{c-c}( 0.35, 1.2)( 0.05, 0.1)(-0.05, 0.1)( 0.35,-0.7)
	\endpsclip
	\psclip{\psframe[linestyle=none](-0.5,-1.5)(0.5,0.7)}
		\psbezier[linestyle=dashed]( 0.35, 1.2)( 0.05, 0.1)(-0.05, 0.1)( 0.35,-0.7)
	\endpsclip
	\psbezier{c-c}(-0.35, 0.7)(-0.05,-0.1)(-0.15,-0.1)(-0.35,-1.2)
}

\def\relSaddleLeftToMembrane{\begin{centerpict}(-1.9,-2)(1.3,2)
	\psset{dimen=middle,dash=2pt 3pt}%
	\membrane
	\pscustom[linestyle=none,fillstyle=solid,fillcolor=white]{%
		\moveto(-0.03,0.37)\lineto(0.82,0.37)
		\lineto(0.82,0.7)\lineto(-0.03,0.7)
		\closepath
		\moveto(0.03,-0.375)\lineto(0.88,-0.375)
		\lineto(0.88,-0.7)\lineto(0.03,-0.7)
		\closepath
	}%
	\pscustom[linestyle=none,fillstyle=solid,fillcolor=white]{%
		\moveto(1.2,1.2)\curveto(0.85, 0.2)(0.85, 0.2)(0.5, 0.7)
		\lineto(-0.35,0.7)\curveto(0,0.2)(0,0.2)(0.35,1.2)
		\closepath
		\moveto(1.2,-0.7)\curveto(0.85,-0.2)(0.85,-0.2)(0.5,-1.2)
		\lineto(-0.35,-1.2)\curveto(0,-0.2)(0,-0.2)(0.35,-0.7)
		\closepath
	}%
	\psbezier{c-c}(1.2, 1.2)(0.85, 0.2)(0.85, 0.2)(0.5, 0.7)
	\psbezier{c-c}(1.2,-0.7)(0.85,-0.2)(0.85,-0.2)(0.5,-1.2)
	\rput(-1.45,0){\relSaddleBackArcs}
	\psclip{\pspolygon[linestyle=none](-0.06,0.7)(0.06,-0.75)(0.6,-0.75)(0.6,0.7)}
		\psbezier[linestyle=dashed]{c-c}(0.35, 1.2)(0, 0.2)(0, 0.2)(-0.35, 0.7)
		\psbezier[linestyle=dashed]{c-c}(0.35,-0.7)(0,-0.2)(0,-0.2)(-0.35,-1.2)
	\endpsclip
	\psclip{\pspolygon[linestyle=none](-0.06,0.7)(0.06,-0.75)(0.06,-1.5)(-1,-1.5)(-1,1.5)(0.6,1.5)(0.6,0.7)}
		\psbezier{c-c}(0.35, 1.2)(0, 0.2)(0, 0.2)(-0.35, 0.7)
		\psbezier{c-c}(0.35,-0.7)(0,-0.2)(0,-0.2)(-0.35,-1.2)
	\endpsclip
	\psclip{\pspolygon[linestyle=none](-0.7,-2)(-0.7,1)(0.7,2)(0.35,1.2)(-0.35,0.7)(-0.35,-2)}%
		\psline[linestyle=dashed]{c-c}(-1.1, 1.2)(1.2, 1.2)
		\psline[linestyle=dashed]{c-c}(-1.8, 0.7)(0.5, 0.7)
		\psline[linestyle=dashed]{c-c}(-1.8,-1.2)(0.5,-1.2)
	\endpsclip
	\psclip{\pspolygon[linestyle=none](-1.9,-3)(-1.9,2)(0.7,2)(-0.7,1)(-0.7,-2)(-0.35,-2)(-0.35,0.7)(0.7,1.45)(1.3,1.45)(1.3,-3)}%
		\psline{c-c}(-1.1, 1.2)(1.2, 1.2)
		\psline{c-c}(-1.8, 0.7)(0.5, 0.7)
		\psline{c-c}(-1.8,-1.2)(0.5,-1.2)
	\endpsclip
	\psclip{\pspolygon[linestyle=none](0.5,-1.2)(0.8,-0.375)(-1.8,-0.375)(-1.8,-1.2)}
		\psline[linestyle=dashed]{c-c}(-1.1,-0.7)(1.2,-0.7)
	\endpsclip
	\psclip{\pspolygon[linestyle=none](0.5,-1.2)(0.8,-0.375)(2,-0.375)(2,-1.2)}
		\psline{c-c}(-1.1,-0.7)(1.2,-0.7)
	\endpsclip
	\psline{c-c}(-0.03,0.37)(0.82,0.37)
	\psline{c-c}( 0.03,-0.375)(0.88,-0.375)
	\psclip{\pspolygon[linestyle=none](-1.5,2)(-1.5,-1)(-0.7,-1)(-0.7,0.75)(-0.06,0.75)(0.06,-0.75)(1,-0.75)(1,2)}
		\psbezier(-0.03, 0.375)(-1.25, 0.37)(-1.25,-0.375)(0.03,-0.375)
	\endpsclip
	\psclip{\pspolygon[linestyle=none](-0.7,0.75)(-0.06,0.75)(0.06,-0.75)(-0.7,-0.75)}
		\psbezier[linestyle=dashed]{c-c}(-0.03, 0.375)(-1.25, 0.37)(-1.25,-0.375)(0.03,-0.375)
	\endpsclip
\end{centerpict}}

\def\relSaddleRightToMembrane{\begin{centerpict}(-1.6,-2)(1.6,2)
	\psset{dimen=middle,dash=2pt 3pt}%
	\membrane
	\pscustom[linestyle=none,fillstyle=solid,fillcolor=white]{%
		\moveto(0.8,-1.2)
		\lineto(1.18,-0.375)\curveto(0.2,-0.375)(0.2, 0.37)(1.12, 0.37)
		\lineto(0.8,0.7)\lineto(-0.35,0.7)
		\curveto(-0.05,-0.1)(-0.15,-0.1)(-0.35,-1.2)
		\closepath
		\moveto(0.35,1.2)
		\curveto( 0.05, 0.1)(-0.05, 0.1)( 0.35,-0.7)
		\lineto(0.35,0.5)\lineto(0.8,0.7)\lineto(0.8,1.2)
		\closepath
	}%
	\psbezier(1.5, 1.2)(1.15, 0.2)(1.15, 0.2)(0.8, 0.7)
	\psbezier(1.5,-0.7)(1.15,-0.2)(1.15,-0.2)(0.8,-1.2)
	\rput(-1.15,0){\relSaddleBackArcs}
	\relSaddleBackArcs
	\psclip{\pspolygon[linestyle=none](-0.7,-2)(-0.7,1)(0.7,2)(0.35,1.2)(-0.35,0.7)(-0.35,-2)}%
		\psline[linestyle=dashed]{c-c}(-0.8, 1.2)(1.5, 1.2)
		\psline[linestyle=dashed]{c-c}(-1.5, 0.7)(0.8, 0.7)
		\psline[linestyle=dashed]{c-c}(-1.5,-1.2)(0.8,-1.2)
	\endpsclip
	\psclip{\pspolygon[linestyle=none](-1.6,-3)(-1.6,2)(0.7,2)(-0.7,1)(-0.7,-2)(-0.35,-2)(-0.35,0.7)(0.7,1.45)(1.6,1.45)(1.6,-3)}%
		\psline{c-c}(-0.8, 1.2)(1.5, 1.2)
		\psline{c-c}(-1.5, 0.7)(0.8, 0.7)
		\psline{c-c}(-1.5,-1.2)(0.8,-1.2)
	\endpsclip
	\psclip{\pspolygon[linestyle=none](0.8,-1.2)(1.1,-0.375)(-1.5,-0.375)(-1.5,-1.2)}
		\psline[linestyle=dashed]{c-c}(-0.8,-0.7)(1.5,-0.7)
	\endpsclip
	\psclip{\pspolygon[linestyle=none](0.8,-1.2)(1.1,-0.375)(2,-0.375)(2,-1.2)}
		\psline{c-c}(-0.8,-0.7)(1.5,-0.7)
	\endpsclip
	\psbezier{c-c}(1.12, 0.37)(0.2, 0.37)(0.2,-0.375)(1.18,-0.375)
\end{centerpict}}

\makeatletter
\def\drawTorusMembrane{\@ifstar
	{\drawTorusMembraneX{\psellipticarc(0,0)(3,2){-45}{225}}}%
	{\drawTorusMembraneX{\psellipse(0,0)(3,2)}}%
}
\makeatother

\def\drawTorusMembraneXX{%
	\psset{dimen=middle,unit=0.75,dash=2pt 3pt}%
	\begin{psset}{linewidth=0.7pt,linecolor=blue}%
		\psline{<-}(-0.2, 3.50)( 0.9, 3.50)
		\psline{<-}(-0.5, 3.25)( 0.6, 3.25)
		\psline{<-}(-0.8, 3.00)( 0.3, 3.00)
	\end{psset}%
	\pspolygon[linestyle=none,fillstyle=solid,fillcolor=membraneColor](-0.9,-0.7)(-0.9,2.1)(0.9,3.6)(0.9,0.8)  % 18/15 = 6/5
}

\def\drawTorusMembraneX#1{%
	\drawTorusMembraneXX
	\psellipse[linestyle=none,fillstyle=solid,fillcolor=white](0,1.3)(0.28,0.7)
	\begin{psclip}{\pscustom[linestyle=none]{%
		\moveto(0,3)
		\psellipticarc(0,1)(2,1.5){270}{340}
		\lineto(3,3)
		\closepath
	}}%
		\pscustom[linestyle=none,fillstyle=solid,fillcolor=white]{%
			\moveto(3,2)
			\psellipticarc(0,0)(3,2){0}{90}
			\psellipticarcn(0,-0.9)(2,1.5){90}{0}
			\closepath
		}%
	\end{psclip}%
	\begin{psclip}{\psframe[linestyle=none](0,0)(-0.9,3)}%
		\psset{linestyle=dashed}%
		\psellipticarc(0,0)(3,2){90}{180}
		\psellipticarc(0,-0.9)(2,1.5){90}{180}
	\end{psclip}%
	\begin{psclip}{\pscustom[linestyle=none]{%
		\moveto(-3,3)
		\psellipticarc(0,1)(2,1.5){200}{340}
		\lineto(3,3)
		\lineto(0,3)
		\lineto(0,0.3)
		\lineto(-0.9,0.3)
		\lineto(-0.9,3)
		\closepath
	}}%
		\psellipticarc(0,-0.9)(2,1.5){0}{180}
	\end{psclip}%
	\begin{psclip}{\pspolygon[linestyle=none](-4,3)(-0.9,3)(-0.9,0)(0,0)(0,3)(4,3)(4,-3)(-4,-3)}#1\end{psclip}%
	\psellipticarc(0,1.3)(0.28,0.7){90}{270}
	\psellipticarc[linestyle=dashed](0,1.3)(0.28,0.7){-90}{90}
}

\def\pictCutTorusMembrane[#1#2]{\begin{centerpict}[midline=0](-2.5,-1.5)(2.5,2.7)
	\drawTorusMembrane*
	\pscustom[linestyle=none,fillstyle=solid,fillcolor=white]{%
		\psbezier{c-c}(-1.2,-0.2)(-0.65,-0.55)(-0.45,-0.8)(-0.5,-1.2)
		\lineto(-1.2,-1.2)
	}
	\psellipticarc(0,1)(2,1.5){200}{225}
	\psbezier{c-c}(-1.664,-1.664)(-0.8,-2)(-0.55,-1.5)(-0.5,-1.2)
	\psbezier{c-c}(-1.2,-0.2)(-0.65,-0.55)(-0.45,-0.8)(-0.5,-1.2)
	\ifx*#1\relax\psdot(-1.25,-1.05)\fi
	\psellipticarc(0,1)(2,1.5){315}{340}
	\psbezier{c-c}(1.664,-1.664)(0.8,-2)(0.55,-1.5)(0.5,-1.2)
	\psbezier{c-c}(1.2,-0.2)(0.65,-0.55)(0.45,-0.8)(0.5,-1.2)
	\ifx*#2\relax\psdot(1.25,-1.05)\fi
\end{centerpict}}

\def\pictTorusMembrane{\begin{centerpict}[midline=0](-2.5,-1.5)(2.5,2.7)
	\drawTorusMembrane
	\pscustom[linestyle=none,fillstyle=solid,fillcolor=white]{%
		\psellipticarc(0,1)(2,1.5){200}{340}
		\psellipticarcn(0,0)(2,1.5){340}{200}
	}
	\psellipticarc(0,1)(2,1.5){200}{340}
\end{centerpict}}

\def\pictTorusLeftSphere{\begin{centerpict}[midline=0](-2.5,-1.5)(0.8,2.7)
	\drawTorusMembraneXX
	\pscustom[linestyle=none,fillstyle=solid,fillcolor=white]{%
		\psbezier{c-c}(-1.2,-0.2)(-0.65,-0.55)(-0.45,-0.8)(-0.5,-1.2)
		\lineto(-1.2,-1.2)
	}
	\psellipticarc(0,1)(2,1.5){200}{225}
	\psbezier{c-c}(-1.664,-1.664)(-0.8,-2)(-0.55,-1.5)(-0.5,-1.2)
	\psbezier{c-c}(-1.2,-0.2)(-0.65,-0.55)(-0.45,-0.8)(-0.5,-1.2)
	\psdot(-1.25,-1.05)
	\psellipticarc{c-c}(0,0)(3,2){135}{225}
	\begin{psclip}{\psellipse[linestyle=none](0,1)(2,1.5)}
		\psellipticarc{c-}(0,-0.9)(2,1.5){135}{180}
	\end{psclip}%
	\cutpicture(-2,-1)(0,2)@(-0.9--){%
		\moveto(-1.664,1.664)
		\curveto(-0.8,2)(-0.55,1.5)(-0.5,1.2)
		\curveto(-0.45,0.8)(-0.65,0.65)(-1.2,0.3)
	}%
\end{centerpict}}

\def\pictTorusRightSphere{\begin{centerpict}[midline=0](-0.9,-1.5)(2.5,2.7)
	\drawTorusMembraneXX
	\psellipticarc(0,1)(2,1.5){-45}{-20}
	\psbezier{c-c}(1.664,-1.664)(0.8,-2)(0.55,-1.5)(0.5,-1.2)
	\psbezier{c-c}(1.2,-0.2)(0.65,-0.55)(0.45,-0.8)(0.5,-1.2)
	\psdot(1.25,-1.05)
	\psellipticarc{c-c}(0,0)(3,2){-45}{45}
	\begin{psclip}{\psellipse[linestyle=none](0,1)(2,1.5)}
		\psellipticarc{c-}(0,-0.9)(2,1.5){0}{45}
	\end{psclip}%
	\pscustom[fillstyle=solid,fillcolor=white]{%
		\moveto(1.664,1.664)
		\curveto(0.8,2)(0.55,1.5)(0.5,1.2)
		\curveto(0.45,0.8)(0.65,0.65)(1.2,0.3)
	}%
\end{centerpict}}

\newrgbcolor{sfceFillColor}{1 0.96 0.95}
\newrgbcolor{sfceBdryColor}{0.33 0.31 0.30}
\newpsobject{surgeryline}{psline}{linecolor=blue,linewidth=5\pslinewidth,linestyle=solid}%

\newpsstyle{mysurface}{%
	fillstyle=solid, fillcolor=sfceFillColor,
	linestyle=solid, linecolor=sfceBdryColor, linewidth=0.5pt
}

\newif\ifquantumsurface
\quantumsurfacefalse

\makeatletter

\newif\ifDrawSurgeryLine
\newif\ifDrawBefore
\newif\ifDrawRotated
\newif\ifDrawRotatedSurgeryLine

\def\@rotatecenterpict(#1,#2)(#3,#4){\centerpict(#2,#1)(#4,#3)\psset{swapaxes=true}}

\def\createSurgeryDiagrams{\@ifstar
	{\DrawRotatedSurgeryLinetrue\createSurgeryDiagrams@}%
	{\DrawRotatedSurgeryLinefalse\createSurgeryDiagrams@}%
}

\def\createSurgeryDiagrams@#1[#2][#3](#4)(#5)#6#7#8{%
	\enamedef{pict#1@}##1{%
		\noexpand\DrawRotatedfalse
		\noexpand\DrawSurgeryLinefalse
		\noexpand\DrawBeforetrue
		\expandonce\csname pict#1@options\endcsname##1%
	}
	\enamedef{pict#1@options}{%
		\noexpand\@ifstar{%
			\noexpand\DrawRotatedtrue
			\expandonce\csname pict#1@options\endcsname
		}{\noexpand\@ifnextchar+{%
			\noexpand\DrawBeforefalse
			\noexpand\expandafter\expandonce\csname pict#1@options\endcsname\noexpand\@gobble
		}{\noexpand\@ifnextchar-{%
			\noexpand\DrawBeforetrue
			\noexpand\expandafter\expandonce\csname pict#1@options\endcsname\noexpand\@gobble
		}{\noexpand\@ifnextchar|{%
			\noexpand\DrawSurgeryLinetrue
			\noexpand\expandafter\expandonce\csname pict#1@options\endcsname\noexpand\@gobble
		}{\noexpand\@ifnextchar[{%
			\expandonce\csname pict#1@draw\endcsname
		}{%
			\expandonce\csname pict#1@draw\endcsname[#3]%
		}}}}}%
	}%
	\ifDrawRotatedSurgeryLine
		\namedef{pict#1@draw}[#2]{%
			\begingroup
				\ifDrawRotated
					\@rotatecenterpict(#4)(#5)%
				\else
					\centerpict(#4)(#5)%
				\fi
					\psset{style=diagarc,linestyle=dashed,dash=2pt 3pt}#6%
					\psset{linestyle=solid}%
					\ifDrawSurgeryLine\begingroup
						\surgeryline(-0.5,0)(0.5,0)%
					\endgroup\fi
					\ifDrawBefore#7\else#8\fi
				\endcenterpict
			\endgroup
		}%
	\else
		\namedef{pict#1@draw}[#2]{%
			\begingroup
				\ifDrawRotated
					\@rotatecenterpict(#4)(#5)%
				\else
					\centerpict(#4)(#5)%
				\fi
					\psset{style=diagarc,linestyle=dashed,dash=2pt 3pt}#6%
					\psset{linestyle=solid}%
					\ifDrawSurgeryLine\begingroup
						\surgeryline(0,-0.5)(0,0.5)%
					\endgroup\fi
					\ifDrawBefore#7\else#8\fi
				\endcenterpict
			\endgroup
		}%
	\fi
	\namedef{pict#1}{\csname pict#1@\endcsname-|}%
	\namedef{pict#1Before}{\csname pict#1@\endcsname-}%
	\namedef{pict#1After}{\csname pict#1@\endcsname+}%
}%

\makeatother

\createSurgeryDiagrams{SplitCircle}[#1#2][00](-2.5,-1.5)(2.5,1.5){%
	\psellipticarc{c-c}(-1.5,0)(1,1){90}{270}
	\psellipticarc{c-c}( 1.5,0)(1,1){-90}{90}
	\ifnum1=#1\relax \psdot(-0.5,0) \fi
	\ifnum1=#2\relax \psdot( 0.5,0) \fi
}{%
	\psbezier{c-c}(-1.5, 1)(-0.5, 1)(-0.75, 0.5)(0, 0.5)
	\psbezier{c-c}(-1.5,-1)(-0.5,-1)(-0.75,-0.5)(0,-0.5)
	\psbezier{c-c}( 1.5, 1)( 0.5, 1)( 0.75, 0.5)(0, 0.5)
	\psbezier{c-c}( 1.5,-1)( 0.5,-1)( 0.75,-0.5)(0,-0.5)
}{%
	\psellipticarc{c-c}(-1.5,0)(1,1){-90}{90}
	\psellipticarc{c-c}( 1.5,0)(1,1){90}{270}
}%

\createSurgeryDiagrams{MergeCircles}[][](-1.5,-2.5)(1.5,2.5){%
	\psellipticarc{c-c}(0,-1.5)(1,1){180}{360}
	\psellipticarc{c-c}(0, 1.5)(1,1){0}{180}
}{%
	\psellipticarc{c-c}(0,-1.5)(1,1){0}{180}
	\psellipticarc{c-c}(0, 1.5)(1,1){180}{360}
}{%
	\psbezier{c-c}( 1,-1.5)( 1,-0.5)( 0.5,-0.75)( 0.5, 0)
	\psbezier{c-c}(-1,-1.5)(-1,-0.5)(-0.5,-0.75)(-0.5, 0)
	\psbezier{c-c}( 1, 1.5)( 1, 0.5)( 0.5, 0.75)( 0.5, 0)
	\psbezier{c-c}(-1, 1.5)(-1, 0.5)(-0.5, 0.75)(-0.5, 0)
}%

\createSurgeryDiagrams{MergeToArc}[][](-1.5,-3)(1.5,3){%
	\psellipticarc{c-c}( 0, 1.5)(1,1){0}{180}
}{%
	\psellipticarc{c-c}(0,-1.5)(1,1){0}{180}
	\psellipticarc{c-c}(0, 1.5)(1,1){180}{360}
}{%
	\psbezier{c-c}( 1,-1.5)( 1,-0.5)( 0.5,-0.75)( 0.5, 0)
	\psbezier{c-c}(-1,-1.5)(-1,-0.5)(-0.5,-0.75)(-0.5, 0)
	\psbezier{c-c}( 1, 1.5)( 1, 0.5)( 0.5, 0.75)( 0.5, 0)
	\psbezier{c-c}(-1, 1.5)(-1, 0.5)(-0.5, 0.75)(-0.5, 0)
}%

\createSurgeryDiagrams{SplitFromArc}[#1#2][00](-1.5,-1.5)(3,1.5){%
	\psellipticarc{c-c}( 1.5,0)(1,1){-90}{90}
	\ifnum1=#1\relax \psdot(-0.5,0) \fi
	\ifnum1=#2\relax \psdot( 0.5,0) \fi
}{%
	\psbezier{c-c}(-1.5, 1)(-0.5, 1)(-0.75, 0.5)(0, 0.5)
	\psbezier{c-c}(-1.5,-1)(-0.5,-1)(-0.75,-0.5)(0,-0.5)
	\psbezier{c-c}( 1.5, 1)( 0.5, 1)( 0.75, 0.5)(0, 0.5)
	\psbezier{c-c}( 1.5,-1)( 0.5,-1)( 0.75,-0.5)(0,-0.5)
}{%
	\psellipticarc{c-c}(-1.5,0)(1,1){-90}{90}
	\psellipticarc{c-c}( 1.5,0)(1,1){90}{270}
}%

\createSurgeryDiagrams{MergeAtPlatform}[][](-0.5,-3)(3,2){%
	\pscustom{%
		\moveto(-0.5,-2.5)
		\lineto(0,-2.5)
		\psellipticarc(0,-1.5)(1,1){-90}{0}
		\stroke
	}%
	\psline[style=diagplatform](-0.5,-3)(-0.5,2)
	\psline[style=diagplatform]( 3  ,-3)( 3  ,2)
}{%
	\psellipticarc{c-c}(0,-1.5)(1,1){0}{90} \psline{c-c}(-0.5,-0.5)(0,-0.5)
	\psellipticarc{c-c}(0, 1.5)(1,1){-90}{0}\psline{c-c}(-0.5, 0.5)(0, 0.5)
}{%
	\psbezier{c-c}( 1,-1.5)( 1,-0.5)( 0.5,-0.75)( 0.5, 0)
	\psbezier{c-c}( 1, 1.5)( 1, 0.5)( 0.5, 0.75)( 0.5, 0)
}%

\createSurgeryDiagrams{MergeArcsAtPlatform}[][](-0.5,-2.5)(3,2.5){%
	\psline(1.5, 1)(3, 1)
	\psline(1.5,-1)(3,-1)
	\psline[style=diagplatform](-0.5,-2.5)(-0.5,2.5)
	\psline[style=diagplatform]( 3  ,-2.5)( 3  ,2.5)
}{%
	\psbezier{c-c}( 1.5, 1)( 0.5, 1)( 0.75, 0.5)(0, 0.5)
	\psbezier{c-c}( 1.5,-1)( 0.5,-1)( 0.75,-0.5)(0,-0.5)
	\psline{c-c}(-0.5,-0.5)(0,-0.5)
	\psline{c-c}(-0.5, 0.5)(0, 0.5)
}{%
	\psellipticarc{c-c}( 1.5,0)(1,1){90}{270}
}%

\createSurgeryDiagrams{SplitAtPlatform}[#1][0](-0.5,-2.5)(3.5,2.5){%
	\psellipticarc(1.5, 0)(1, 1){-90}{90}
	\psline[style=diagplatform](-0.5,-2.5)(-0.5,2.5)
	\psline[style=diagplatform]( 3.5,-2.5)( 3.5,2.5)
	\ifnum1=#1\relax\rput{150}(1.5,0){\psdot(1,0)}\fi
}{%
	\psbezier{c-c}( 1.5, 1)( 0.5, 1)( 0.75, 0.5)(0, 0.5)
	\psbezier{c-c}( 1.5,-1)( 0.5,-1)( 0.75,-0.5)(0,-0.5)
	\psline{c-c}(-0.5,-0.5)(0,-0.5)
	\psline{c-c}(-0.5, 0.5)(0, 0.5)
}{%
	\psellipticarc{c-c}( 1.5,0)(1,1){90}{270}
}%

\createSurgeryDiagrams{MergeArcs}[][](-3,-2.5)(3,2.5){%
	\psline{c-c}(-3,-1)(-1.5,-1)
	\psline{c-c}( 3,-1)( 1.5,-1)
	\psline{c-c}(-3, 1)(-1.5, 1)
	\psline{c-c}( 3, 1)( 1.5, 1)
	\psline[style=diagplatform](-3,-2.5)(-3, 2.5)
	\psline[style=diagplatform]( 3,-2.5)( 3, 2.5)
}{%
	\psbezier{c-c}(-1.5, 1)(-0.5, 1)(-0.75, 0.5)(0, 0.5)
	\psbezier{c-c}(-1.5,-1)(-0.5,-1)(-0.75,-0.5)(0,-0.5)
	\psbezier{c-c}( 1.5, 1)( 0.5, 1)( 0.75, 0.5)(0, 0.5)
	\psbezier{c-c}( 1.5,-1)( 0.5,-1)( 0.75,-0.5)(0,-0.5)
}{%
	\psellipticarc{c-c}(-1.5,0)(1,1){-90}{90}
	\psellipticarc{c-c}( 1.5,0)(1,1){90}{270}
}%

\createSurgeryDiagrams{MergeHorArcs}[][](-3,-3)(3,3){%
	\psbezier(-3,-2.25)(-1.5,-2.25)(-1,-1.75)(-1,-1)
	\psbezier(-3, 2.25)(-1.5, 2.25)(-1, 1.75)(-1, 1)
	\psbezier( 3,-2.25)( 1.5,-2.25)( 1,-1.75)( 1,-1)
	\psbezier( 3, 2.25)( 1.5, 2.25)( 1, 1.75)( 1, 1)
	\begin{psset}{style=diagplatform}
		\psline(-3,-3)(-3,-1.5) \psline(-3, 1.5)(-3, 3)
		\psline( 3,-3)( 3,-1.5) \psline( 3, 1.5)( 3, 3)
	\end{psset}
	\psline[style=diaghline](-3,-1.5)(3,-1.5)
	\psline[style=diaghline](-3, 1.5)(3, 1.5)
}{%
	\psellipticarc(0, 1)(1,0.5){180}{360}
	\psellipticarc(0,-1)(1,0.5){0}{180}
}{%
	\psbezier(-1,-1)(-1,-0.5)(-0.5,-0.5)(-0.5,0)
	\psbezier(-1, 1)(-1, 0.5)(-0.5, 0.5)(-0.5,0)
	\psbezier( 1,-1)( 1,-0.5)( 0.5,-0.5)( 0.5,0)
	\psbezier( 1, 1)( 1, 0.5)( 0.5, 0.5)( 0.5,0)
}%

\createSurgeryDiagrams*{MergeVertArcs}[][](-3,-3)(3,3){%
	\psbezier(-3,-2.25)(-1.5,-2.25)(-1,-1.75)(-1,-1)
	\psbezier(-3, 2.25)(-1.5, 2.25)(-1, 1.75)(-1, 1)
	\psbezier( 3,-2.25)( 1.5,-2.25)( 1,-1.75)( 1,-1)
	\psbezier( 3, 2.25)( 1.5, 2.25)( 1, 1.75)( 1, 1)
	\begin{psset}{style=diagplatform}
		\psline(-3,-3)(-3,-1.5) \psline(-3, 1.5)(-3, 3)
		\psline( 3,-3)( 3,-1.5) \psline( 3, 1.5)( 3, 3)
	\end{psset}
	\psline[style=diaghline](-3,-1.5)(3,-1.5)
	\psline[style=diaghline](-3, 1.5)(3, 1.5)
}{%
	\psbezier(-1,-1)(-1,-0.5)(-0.5,-0.5)(-0.5,0)
	\psbezier(-1, 1)(-1, 0.5)(-0.5, 0.5)(-0.5,0)
	\psbezier( 1,-1)( 1,-0.5)( 0.5,-0.5)( 0.5,0)
	\psbezier( 1, 1)( 1, 0.5)( 0.5, 0.5)( 0.5,0)
}{%
	\psellipticarc(0, 1)(1,0.5){180}{360}
	\psellipticarc(0,-1)(1,0.5){0}{180}
}%

\createSurgeryDiagrams*{MergeArcsToTurnback}[][](-2.5,-2)(2.5,4){%
	\psbezier(-1,1)(-1,2.2)(-1.3,2.5)(-2.5,2.5)
	\psbezier( 1,1)( 1,3.0)( 0.3,3.5)(-2.5,3.5)
	\begin{psset}{style=diagplatform}
		\psline(-2.5, 4)(-2.5, 1.75)
		\psline( 2.5, 4)( 2.5, 1.75)
	\end{psset}
	\psline[style=diaghline](-2.5, 1.75)(2.5, 1.75)
}{%
	\psbezier(-1,-1)(-1,-0.5)(-0.5,-0.5)(-0.5,0)
	\psbezier(-1, 1)(-1, 0.5)(-0.5, 0.5)(-0.5,0)
	\psbezier( 1,-1)( 1,-0.5)( 0.5,-0.5)( 0.5,0)
	\psbezier( 1, 1)( 1, 0.5)( 0.5, 0.5)( 0.5,0)
}{%
	\psellipticarc(0, 1)(1,0.5){180}{360}
	\psellipticarc(0,-1)(1,0.5){0}{180}
}%

\def\pictDrawAnnulus(#1)#2#3{%
	\pscircle[style=mysurface](#1){#2}
	\pscircle[style=mysurface,fillcolor=white](#1){#3}
	\ifquantumsurface
		\rput(#1){\psline[linecolor=red,linewidth=0.5pt,linestyle=solid](#2,0)(#3,0)}%
	\fi
}

\createSurgeryDiagrams*{AnnMergeTrivial}[#1][-](-1.5,-3)(4.5,3){%
	\pictDrawAnnulus(1.5,0){3}{0.3}
	\psarc{c-c}(1.5,0){1.75}{-120}{120}
	\psellipticarc(-0.5,0)(0.5,1){90}{270}
}{\pscustom{%
		\moveto(0.625,1.516)
		\curveto(0.192,1.266)(0.5,0.75)(0.5,0)
		\curveto(0.5,-0.75)(0.192,-1.266)(0.625,-1.516)
		\moveto(-0.5,-1)
		\curveto(0,-1)(-0.5,-0.5)(-0.5,0)
		\curveto(-0.5,0.5)(0,1)(-0.5,1)
	}%
	\ifx d#1\relax
		\psline[arrowsize=6pt]{->}(0.445,0.617)(0.459,0.518)
	\else\ifx u#1\relax
		\psline[arrowsize=6pt]{->}(0.404,0.990)(0.392,1.090)
	\fi\fi
}{\pscustom{%
		\moveto(0.625,1.516)
		\curveto(0.192,1.266)(0.7,0.3)(0,0.3)
		\curveto(-0.4,0.3)(-0.1,1)(-0.5,1)
		\moveto(0.625,-1.516)
		\curveto(0.192,-1.266)(0.7,-0.3)(0,-0.3)
		\curveto(-0.4,-0.3)(-0.1,-1)(-0.5,-1)
	}%
	\ifx d#1\relax
		\psline[arrowsize=6pt]{->}(0.405,0.755)(0.397,0.656)
	\else\ifx u#1\relax
		\psline[arrowsize=6pt]{->}(0.448,1.198)(0.457,1.290)
	\fi\fi
}

\createSurgeryDiagrams{AnnSplitOffTrivial}[#1#2][0-](-1.55,-3)(4.45,3){%
	\pictDrawAnnulus(1.45,0){3}{0.3}
	\psarc{c-c}(1.45,0){1.75}{-120}{120}
	\psellipticarc(-0.55,0)(0.5,1){90}{270}
}{\pscustom{%
		\moveto(0.575,1.516)
		\curveto(0.142,1.266)(0.65,0.5)(-0.05,0.5)
		\curveto(-0.45,0.5)(-0.15,1)(-0.55,1)
		\moveto(0.575,-1.516)
		\curveto(0.142,-1.266)(0.65,-0.5)(-0.05,-0.5)
		\curveto(-0.45,-0.5)(-0.15,-1)(-0.55,-1)
	}%
	\ifx d#2\relax
		\psline[arrowsize=6pt]{->}(0.355,0.755)(0.347,0.656)
	\else\ifx u#2\relax
		\psline[arrowsize=6pt]{->}(0.395,1.268)(0.405,1.360)
	\fi\fi
}{\pscustom{%
		\moveto(0.575,1.516)
		\curveto(0.142,1.266)(0.45,0.75)(0.45,0)
		\curveto(0.45,-0.75)(0.142,-1.266)(0.575,-1.516)
		\moveto(-0.55,-1)
		\curveto(-0.05,-1)(-0.25,-0.5)(-0.25,0)
		\curveto(-0.25,0.5)(-0.05,1)(-0.55,1)
	}%
	\ifx d#2\relax
		\psline[arrowsize=6pt]{->}(0.395,0.617)(0.409,0.518)
	\else\ifx u#2\relax
		\psline[arrowsize=6pt]{->}(0.354,0.990)(0.342,1.090)
	\fi\fi
	\ifnum1=#1\relax\psdot(-0.25,0)\fi
}

\createSurgeryDiagrams*{AnnMergeEssentials}[#1#2][--](-1.25,-3)(4.75,3){%
	\pictDrawAnnulus(1.75,0){3}{0.3}
	\psarc(1.75,0){1.3}{-135}{135}
	\psarc(1.75,0){2.25}{-155}{155}
}{%
	\ifx d#1\relax
		\psarc[arrowsize=6pt]{->}(1.75,0){2.25}{155}{205}
	\else\ifx u#1\relax
		\psarcn[arrowsize=6pt]{->}(1.75,0){2.25}{205}{155}
	\else
		\psarc(1.75,0){2.25}{155}{205}
	\fi\fi
	\ifx d#2\relax
		\psarc[arrowsize=6pt]{->}(1.75,0){1.3}{135}{225}
	\else\ifx u#2\relax
		\psarcn[arrowsize=6pt]{->}(1.75,0){1.3}{225}{135}
	\else
		\psarc(1.75,0){1.3}{135}{225}
	\fi\fi
}{\pscustom{%
		\moveto(0.831,0.919)
		\curveto(0.615,0.707)(0.5,0.3)(0.1,0.3)
		\curveto(-0.4,0.3)(-0.416,0.679)(-0.289,0.951)
		\moveto(0.831,-0.919)
		\curveto(0.615,-0.707)(0.5,-0.3)(0.1,-0.3)
		\curveto(-0.4,-0.3)(-0.416,-0.679)(-0.289,-0.951)
	}%
	\ifx1#1\relax\psdot(0.1, 0.3)\fi
	\ifx1#2\relax\psdot(0.1,-0.3)\fi
}

\createSurgeryDiagrams*{AnnMergeEssentialsFull}[#1#2][--](-1.25,-3)(4.75,3){%
	\pictDrawAnnulus(1.75,0){3}{0.3}
}{%
	\pscircle(1.75,0){1.3}
	\pscircle(1.75,0){2.25}
	\ifx d#1\relax
		\psarc[arrowsize=6pt,linestyle=none]{->}(1.75,0){2.25}{145}{215}
	\else\ifx u#1\relax
		\psarcn[arrowsize=6pt,linestyle=none]{->}(1.75,0){2.25}{215}{145}
	\else
		\psarc(1.75,0){2.25}{155}{205}
	\fi\fi
	\ifx d#2\relax
		\psarc[arrowsize=6pt,linestyle=none]{->}(1.75,0){1.3}{125}{235}
	\else\ifx u#2\relax
		\psarcn[arrowsize=6pt,linestyle=none]{->}(1.75,0){1.3}{235}{125}
	\else
		\psarc(1.75,0){1.3}{135}{225}
	\fi\fi
}{\pscustom{%
		\moveto(3.05,0)
		\psarc(1.75,0){1.3}{-10}{135}
		\curveto(0.615,0.707)(0.5,0.3)(0.1,0.3)
		\moveto(4,0)
		\psarc(1.75,0){2.25}{-10}{155}
		\curveto(-0.416,0.679)(-0.4,0.3)(0.1,0.3)
		\moveto(3.05,0)
		\psarcn(1.75,0){1.3}{10}{-135}
		\curveto(0.615,-0.707)(0.5,-0.3)(0.1,-0.3)
		\moveto(4,0)
		\psarcn(1.75,0){2.25}{10}{-155}
		\curveto(-0.416,-0.679)(-0.4,-0.3)(0.1,-0.3)
	}%
	\ifnum1=#1\relax\psdot(0.1, 0.3)\fi
	\ifnum1=#2\relax\psdot(0.1,-0.3)\fi
}

\createSurgeryDiagrams{AnnSplitTrivial}[#1#2][--](-1.35,-3)(4.65,3){%
	\pictDrawAnnulus(1.65,0){3}{0.3}
	\psarc(1.65,0){1.3}{-135}{135}
	\psarc(1.65,0){2.25}{-155}{155}
}{\pscustom{%
		\moveto(0.731,0.919)
		\curveto(0.515,0.707)(0.4,0.5)(0.0,0.5)
		\curveto(-0.3,0.5)(-0.516,0.679)(-0.389,0.951)
		\moveto(0.731,-0.919)
		\curveto(0.515,-0.707)(0.4,-0.5)(0.0,-0.5)
		\curveto(-0.3,-0.5)(-0.516,-0.679)(-0.389,-0.951)
	}%
}{%
	\ifx d#1\relax
		\psarc[arrowsize=6pt]{->}(1.65,0){2.25}{155}{205}
	\else\ifx u#1\relax
		\psarcn[arrowsize=6pt]{->}(1.65,0){2.25}{205}{155}
	\else
		\psarc(1.65,0){2.25}{155}{205}
	\fi\fi
	\ifx d#2\relax
		\psarc[arrowsize=6pt]{->}(1.65,0){1.3}{135}{225}
	\else\ifx u#2\relax
		\psarcn[arrowsize=6pt]{->}(1.65,0){1.3}{225}{135}
	\else
		\psarc(1.65,0){1.3}{135}{225}
	\fi\fi
}

\createSurgeryDiagrams{AnnSplitTrivialFull}[#1#2][--](-1.35,-3)(4.65,3){%
	\pictDrawAnnulus(1.65,0){3}{0.3}
}{\pscustom{%
		\moveto(0.731,0.919)
		\curveto(0.515,0.707)(0.4,0.5)(0.0,0.5)
		\curveto(-0.3,0.5)(-0.516,0.679)(-0.389,0.951)
		\moveto(0.731,-0.919)
		\curveto(0.515,-0.707)(0.4,-0.5)(0.0,-0.5)
		\curveto(-0.3,-0.5)(-0.516,-0.679)(-0.389,-0.951)
	}%
	\psarc(1.65,0){1.3}{-135}{135}
	\psarc(1.65,0){2.25}{-155}{155}
}{%
	\pscircle(1.65,0){1.3}
	\pscircle(1.65,0){2.25}
	\ifx d#1\relax
		\psarc[arrowsize=6pt,linestyle=none]{->}(1.65,0){2.25}{155}{205}
	\else\ifx u#1\relax
		\psarcn[arrowsize=6pt,linestyle=none]{->}(1.65,0){2.25}{205}{155}
	\else
		\psarc(1.65,0){2.25}{155}{205}
	\fi\fi
	\ifx d#2\relax
		\psarc[arrowsize=6pt,linestyle=none]{->}(1.65,0){1.3}{135}{225}
	\else\ifx u#2\relax
		\psarcn[arrowsize=6pt,linestyle=none]{->}(1.65,0){1.3}{225}{135}
	\else
		\psarc(1.65,0){1.3}{135}{225}
	\fi\fi
}

\def\pictDrawBandCurve#1(#2,#3)(#4,#5)(#6,#7){
	\curveto(#4,#5)(#6,#7)(0,#1)
	\psarc(0,0){#1}{90}{270}%
	\curveto(#6,-#7)(#4,-#5)(#2,#3)
}

\def\pictDrawMoebiusBand(#1)#2#3#4(#5)(#6)(#7)(#8){%
	\rput(#1){%
		\pscustom[style=mysurface]{%
			\moveto(#4, 0)
			\pictDrawBandCurve{#2}(#4,0)(#5)(#6)
		}\pscustom[style=mysurface,fillcolor=white]{%
			\moveto(#4, 0)
			\pictDrawBandCurve{#3}(#4,0)(#7)(#8)
		}%
	}%
	\ifquantumsurface
		\rput(#1){\psline[linecolor=red,linewidth=0.5pt,linestyle=solid](#2,0)(#3,0)}%
	\fi
}

\createSurgeryDiagrams{MBBoundingToTrivial}[#1#2][--](-1.2,-3)(4.2,3){%
	\pictDrawMoebiusBand(1.8,0){3}{1}{1.8}(2.8,1.2)(2,3)(1.3,0.6)(1,1)
	\psbezier(3.6,0.2)(3.3,1.0)(1.731,1.919)(0.731,0.919)
	\psbezier(3.6,-0.2)(3.3,-1.0)(1.731,-1.919)(0.731,-0.919)
	\pscustom{%
		\moveto(3.6,0.2)
		\curveto(3.9,1.0)(3.5,2.25)(1.65,2.25)
		\psarc(1.65,0){2.25}{90}{155}
		\moveto(3.6,-0.2)
		\curveto(3.9,-1.0)(3.5,-2.25)(1.65,-2.25)
		\psarcn(1.65,0){2.25}{-90}{-155}
	}%
}{\pscustom{%
		\moveto(0.731,0.919)
		\curveto(0.515,0.707)(0.4,0.5)(0.0,0.5)
		\curveto(-0.3,0.5)(-0.516,0.679)(-0.389,0.951)
		\moveto(0.731,-0.919)
		\curveto(0.515,-0.707)(0.4,-0.5)(0.0,-0.5)
		\curveto(-0.3,-0.5)(-0.516,-0.679)(-0.389,-0.951)
	}%
}{%
	\ifx d#1\relax
		\psarc[arrowsize=6pt]{->}(1.65,0){2.25}{155}{205}
	\else\ifx u#1\relax
		\psarcn[arrowsize=6pt]{->}(1.65,0){2.25}{205}{155}
	\else
		\psarc(1.65,0){2.25}{155}{205}
	\fi\fi
	\ifx d#2\relax
		\psarc[arrowsize=6pt]{->}(1.65,0){1.3}{135}{225}
	\else\ifx u#2\relax
		\psarcn[arrowsize=6pt]{->}(1.65,0){1.3}{225}{135}
	\else
		\psarc(1.65,0){1.3}{135}{225}
	\fi\fi
}

\createSurgeryDiagrams*{MBTrivialToBounding}[#1#2][--](-1.1,-3)(4.3,3){%
	\pictDrawMoebiusBand(1.9,0){3}{1}{1.9}(2.9,1.2)(2.1,3)(1.4,0.6)(1.1,1)
	\psbezier(3.7, 0.2)(3.4, 1.0)(1.831, 1.919)(0.831, 0.919)
	\psbezier(3.7,-0.2)(3.4,-1.0)(1.831,-1.919)(0.831,-0.919)
	\pscustom{%
		\moveto(3.7,0.2)
		\curveto(4.0,1.0)(3.6,2.25)(1.75,2.25)
		\psarc(1.75,0){2.25}{90}{155}
		\moveto(3.7,-0.2)
		\curveto(4.0,-1.0)(3.6,-2.25)(1.75,-2.25)
		\psarcn(1.75,0){2.25}{-90}{-155}
	}%
}{%
	\ifx d#1\relax
		\psarc[arrowsize=6pt]{->}(1.75,0){2.25}{155}{205}
	\else\ifx u#1\relax
		\psarcn[arrowsize=6pt]{->}(1.75,0){2.25}{205}{155}
	\else
		\psarc(1.75,0){2.25}{155}{205}
	\fi\fi
	\ifx d#2\relax
		\psarc[arrowsize=6pt]{->}(1.75,0){1.3}{135}{225}
	\else\ifx u#2\relax
		\psarcn[arrowsize=6pt]{->}(1.75,0){1.3}{225}{135}
	\else
		\psarc(1.75,0){1.3}{135}{225}
	\fi\fi
}{%
	\pscustom{%
		\moveto(0.831,0.919)
		\curveto(0.615,0.707)(0.5,0.3)(0.1,0.3)
		\curveto(-0.4,0.3)(-0.416,0.679)(-0.289,0.951)
		\moveto(0.831,-0.919)
		\curveto(0.615,-0.707)(0.5,-0.3)(0.1,-0.3)
		\curveto(-0.4,-0.3)(-0.416,-0.679)(-0.289,-0.951)
	}%
	\ifx1#1\relax\psdot(0.1, 0.3)\fi
	\ifx1#2\relax\psdot(0.1,-0.3)\fi
}

\makeatletter

\def\pictCirclesOnAPlane#1{%
\begin{centerpict}(0.2,0.2)(3.5,2.6)
	\psccurve(1,1)(2,1.2)(2,1.8)(1.8,2.2)(2.2,2.6)(3.2,1)(2.2,0.8)(1,0.4)
	\pscircle(1,1.8){0.4}%
	\psellipse(2.6,1.5)(0.2,0.4)
	\ifnum#1>0\relax\psdot(1.8,2.2)\fi
	\ifnum#1>1\relax\psdot(1.4,1.8)\fi
	\ifnum#1>2\relax\psdot(3.2,1)\fi
\end{centerpict}}

\def\cupCob{\@ifstar{\drawCupCob{1}}{\drawCupCob{0}}}

\def\drawCupCob#1{%
	\begin{psset}{unit=1.5ex}
	\begin{centerpict}(-0.1,0,1)(2.1,3.1)
		\pscustom{%
			\moveto(0,2.6)
			\lineto(0,2.0)
			\psellipticarc(1,2.0)(1,1.5){180}{360}
			\lineto(2,2.6)
		}%
		\psellipse(1,2.6)(1,0.4)
		\ifcase#1\or\psdot(1,1.5)\fi
	\end{centerpict}%
	\end{psset}%
}

\newTangle{Tangle}{1.5ex}(0,-1)(3,1)<2>[][]{%
	\psbezier(1.5, 1)(1.5,0)(0.5,0)(0.5,1)
	\psbezier(1.5,-1)(1.5,0)(2.5,0)(2.5,1)
}

\def\pictTangleDisk{\@ifstar{\pictTangleDisk@0}{\pictTangleDisk@1}}

\def\pictTangleDisk@#1{%
\begin{centerpict}(-1,-3)(5,3.2)
	\psellipse[style=mysurface](1.5,0.1)(1.5,1.1)
	\psset{style=diagplatform}
	\psline(0,0.1)(-1,0.1)
	\psline(3,0.1)( 5,0.1)
	\psset{style=diagarc}%
	\pscustom{%
		\moveto( 4.5, 0.1)
		\psarc ( 3.0, 1.0){1.5}{0}{180}
		\curveto(1.5, 0.2)(0.5, 0.0)(0.5, 1.0)
		\psarc ( 0.0, 1.0){0.5}{0}{180}
		\ifx1#1\relax
			\psarc(2.0,-0.5){2.5}{180}{360}
			\closepath
		\else
			\lineto(-0.5, 0.1)
		\fi
		\stroke
		\newpath
		\moveto(3.5, 1.0)
		\psarc (3.0, 1.0){0.5}{0}{180}
		\curveto(2.5,0.0)(1.5, 0.0)(1.5,-1.0)
		\psarc (2.5,-1.0){1.0}{180}{360}
		\closepath
	}
\end{centerpict}}

\newTangle{VertTan}{2ex}(0,-0.75)(2,0.75)<1.5>[#3#4][--]{%
	\psline(0.5,-0.75)(0.5,0.75)
	\psline(1.5,-0.75)(1.5,0.75)%
	\ifx*##3\relax\psdot(0.5,0)\fi
	\ifx*##4\relax\psdot(1.5,0)\fi
}

\newTangle{VTCob}{2ex}(0,-0.75)(2,0.75)<1.5>[][]{%
	\surgeryline(0.5,0)(1.5,0)
	\psline(0.5,-0.75)(0.5,0.75)
	\psline(1.5,-0.75)(1.5,0.75)%
}

\newTangle{CupCap}{2ex}(0,-0.75)(2,0.75)<1.5>[#3#4][--]{%
	\psarc(1, 0.75){0.5}{180}{360}
	\psarc(1,-0.75){0.5}{0}{180}
}

\newTangle{CCCob}{2ex}(0,-0.75)(2,0.75)<1.5>[#3#4][--]{%
	\surgeryline(1,-0.25)(1,0.25)
	\psarc(1, 0.75){0.5}{180}{360}
	\psarc(1,-0.75){0.5}{0}{180}
}

\def\pictVertTanClosed{%
	\@ifnextchar[%
		{\pictVertTanClosed@}%
		{\pictVertTanClosed@[--]}%
}

\def\pictVertTanClosed@[#1#2]{%
\begin{centerpict}(-2,-1.5)(2,1.5)
	\pictVertTanClosedDraw{#1}{#2}%
	\ifx*#1\psdot(-0.5,0)\fi
	\ifx*#2\psdot( 0.5,0)\fi
\end{centerpict}}

\def\pictVTCobClosed{%
\begin{centerpict}(-2,-1.5)(2,1.5)
	\pictVertTanClosedDraw--%
	\surgeryline(-0.5,0)(0.5,0)
\end{centerpict}}

\def\pictCupCapClosed{%
\begin{centerpict}(-2,-1.5)(2,1.5)
	\pictCupCapClosedDraw
\end{centerpict}}

\def\pictCCCobClosed{%
\begin{centerpict}(-2,-1.5)(2,1.5)
	\pictCupCapClosedDraw
	\surgeryline(0,-0.5)(0,0.5)
\end{centerpict}}

\def\pictDrawDisk{%
	\psline[style=diagplatform](-1.1,0)(-2,0)
	\psline[style=diagplatform]( 1.1,0)( 2,0)
	\pscustom[style=mysurface]{%
		\moveto(1.1,0)
		\curveto( 1.1, 0.8)( 0.8, 1.1)( 0.0, 1.1)
		\curveto(-0.8, 1.1)(-1.1, 0.8)(-1.1, 0.0)
		\curveto(-1.1,-0.8)(-0.8,-1.1)( 0.0,-1.1)
		\curveto( 0.8,-1.1)( 1.1,-0.8)( 1.1, 0.0)
		\closepath
	}
	\psset{style=diagarc}%
}

\def\pictCupCapClosedDraw{%
	\pictDrawDisk
	\pscustom{%
		\moveto(-0.5,1)
		\psarc(-1, 1){0.5}{0}{180}
		\psarc(-1,-1){0.5}{180}{360}
		\psarcn(0,-1){0.5}{180}{0}
		\psarc( 1,-1){0.5}{180}{360}
		\psarc( 1, 1){0.5}{0}{180}
		\psarcn(0, 1){0.5}{0}{-180}
		\closepath
	}%
}

\def\pictVertTanClosedDraw#1#2{%
	\pictDrawDisk
	\pscustom{%
		\moveto(-0.5,1)
		\psarc(-1, 1){0.5}{0}{180}
		\psarc(-1,-1){0.5}{180}{360}
		\closepath
	}%
	\pscustom{%
		\moveto(1.5,1)
		\psarc(1, 1){0.5}{0}{180}
		\psarc(1,-1){0.5}{180}{360}
		\closepath
	}%
}

\makeatother

\newtheorem{theorem}{Theorem}[section]
\newtheorem{lemma}[theorem]{Lemma}

\newtheorem{proposition}[theorem]{Proposition}
\newtheorem{corollary}[theorem]{Corollary}

\theoremstyle{definition}

\newtheorem{definition}[theorem]{Definition}
\newtheorem{example}[theorem]{Example}

\theoremstyle{remark}

\newtheorem{remark}[theorem]{Remark}
\newtheorem*{remark*}{Remark}
\newtheorem*{notation*}{Notation}

\numberwithin{equation}{section}

\definecolor{internalLink}{rgb}{0.5,0,0}
\definecolor{citeLink}{rgb}{0,0.5,0}
\definecolor{urlLink}{rgb}{0,0,0.5}

\hypersetup{linkcolor=internalLink}
\hypersetup{citecolor=citeLink}
\hypersetup{urlcolor=urlLink}

\def\noic{\sb{}\kern-\scriptspace }

\diagarcset{%
	annular=false,%
	hlinecolor=sfceBdryColor,%
	hlinewidth=0.5pt,%
	hlinestyle=solid,%
	arccolor=black,%
	arcwidth=1pt,%
	arcstyle=solid,%
	arctype=flat,%
	dotsize=5pt,%
	platformcolor=red,%
	platformwidth=0.5pt,%
	platformstyle=dashed,%
	pt-dist=1,%
}%

\newrgbcolor{changecolor}{0.5 0.0 0.2}

\makeatletter

\def\blank{-}

\def\cocolon{%
	\nobreak\mskip6mu plus1mu\mathpunct{}%
	\nonscript\mkern-\thinmuskip{:}%
	\mskip2mu\relax
}

\begingroup
	\def\definering#1#2{\protected\gdef#1{\@ifstar{#2^*}{#2}}}%
	\definering\scalars\Bbbk
	\definering\qscalars{\Bbbk[q^{\pm1}]}
	\definering\escalars{\Bbbk'}
	\definering\Fld{\mathbb F}
	\definering\qFld{\mathbb F[q^{\pm1}]}
	\definering\Z{\mathbb Z}
	\definering\Zq{\mathbb Z[q^{\pm1}]}
	\definering\Q{\mathbb Q}
	\definering\R{\mathbb R}
	\definering\C{\mathbb C}
\endgroup

\newcommand*\LieSL{\mathfrak{sl}}

\newcommand*\qU{\mathcal{U}_q}

\newcommand*{\Uqsl}[1]{\qU(\LieSL_{#1})}

\def\S{\mathbb S^1}
\def\Ann{\mathbb A}
\def\sfce{F}
\def\MB{\mathbb M}

\def\RxI{\R\times I}

\def\AxI{\Ann\times I}
\def\MBxI{\MB\mkern\thinmuskip\lowtilde[1.5]\times\mkern\thinmuskip I}

\def\setcomponent@#1{%
	\@ifnextchar[{\setcomponent@weight#1}{\setcomponent@weight#1[]}
}
\def\setcomponent@weight#1[#2]#3{%
	\setcomponent@draw#1[#2]#3,,\relax
}
\def\setcomponent@parameter[#1]#2{\ifx\relax#1\relax#2\else#2,#1\fi}
\def\setcomponent@draw#1[#2]#3,#4,#5\relax{%
	\ifx\relax#4\relax
		\mathcal{#1}^{#3}%
	\else
		\mathcal{#1}^{#3}_{#4}%
	\fi
	\ifx\relax#2\relax\else(#2)\fi
}

\def\tangles@#1,#2,#3,#4\relax{%
	\ifx\relax#3\relax		%	points on one boundary only
		\mathcal{#1}^{#2}%
	\else									%	points on both boundaries
		\mathcal{#1}^{#3}_{#2}%
	\fi
}

\def\tangles#1{\tangles@ T,#1,,\relax}
\def\ftangles#1{\tangles@ FT,#1,,\relax}
\def\dtangles#1{\tangles@\widetilde T,#1,,\relax}
\def\matchings#1{\tangles@ M,#1,,\relax}
\def\gmatchings{\setcomponent@{GM}}

\def\revmatching#1{#1^!}

\DeclareMathOperator{\id}{id}		%	identity morphism
\DeclareMathOperator{\im}{im}		%	image
\DeclareMathOperator{\tr}{tr}		% trace of a function
\DeclareMathOperator{\dom}{dom}
\DeclareMathOperator{\cod}{cod}
\DeclareMathOperator{\iid}{\mathbf{1}}%

\DeclareMathOperator\adeg{adeg}
\DeclareMathOperator\qdeg{qdeg}
\DeclareMathOperator\weight{wt}

\def\tensor{\@ifstar\tensor@{\tensor@\otimes}}
\def\tensor@#1#2{\let\tensor@symb#1\tensor@@#2,\tensor@end,}
\def\tensor@@#1,{\tensor@expand#1\relax\tensor@@@}
\def\tensor@@@#1,{%
	\ifx\tensor@end#1\relax\else
		\tensor@symb\tensor@expand#1\relax
	\expandafter\tensor@@@\fi
}
\def\tensor@expand#1#2\relax{\ifx*#1#2\else#1_{#2}\fi}

\DeclareMathOperator\dcbtensor\heartsuit
\DeclareMathOperator\bdcbtensor\spadesuit

\def\utimes{\mytimes@\otimes}
\def\udtimes{\mytimes@{\hat\otimes}}

\def\mytimes@#1#2{%
	\mathchoice
		{\mytimes@@{#1}{#2}}%
		{\mytimes@@@{#1}{#2}\scriptstyle}%
		{\mytimes@@@{#1}{#2}\scriptscriptstyle}%
		{\mytimes@@@{#1}{#2}\scriptscriptstyle}%
}
\def\mytimes@clap#1#2{%
	\mathchoice
		{\mytimes@@clap{#1}{#2}}%
		{\mytimes@@@{#1}{#2}\scriptstyle}%
		{\mytimes@@@{#1}{#2}\scriptscriptstyle}%
		{\mytimes@@@{#1}{#2}\scriptscriptstyle}%
}

\def\mytimes@@#1#2{\mkern\thinmuskip\underset{\mkern-1mu#2}{#1}\mkern\thinmuskip}
\def\mytimes@@clap#1#2{\mkern\thinmuskip\underset{\mathclap{#2}}{#1}\mkern\thinmuskip}
\def\mytimes@@@#1#2#3{\mkern\thinmuskip{#1}_{\mkern-2mu\raisebox{-1pt}{$\m@th#3#2$}}\mkern\thinmuskip}

\newcommand\quotient[2]{%
	\raisebox{0.3\baselineskip}{$#1$}\!\Big/\!\raisebox{-0.3\baselineskip}{$#2$}%
}

\def\shadow@getindices{%
	\let\shadow@sup\relax
	\let\shadow@sub\relax
	\shadow@checkindices
}
\def\shadow@checkindices{%
	\@ifnextchar^{\shadow@readsup}{%
	\@ifnextchar_{\shadow@readsub}{%
	\mathchoice
		{\shadow@typeindices\scriptstyle{1.25ex}{-0.75ex}}%
		{\shadow@typeindices\scriptstyle{1.25ex}{-0.75ex}}%
		{\shadow@typeindices\scriptscriptstyle{0.75ex}{-0.65ex}}%
		{\shadow@typeindices\scriptscriptstyle{0.65ex}{-0.55ex}}%
	}}%
}
\def\shadow@readsup^#1{\def\shadow@sup{#1}\shadow@checkindices}%
\def\shadow@readsub_#1{\def\shadow@sub{#1}\shadow@checkindices}%
\def\shadow@typeindices#1#2#3{\begingroup
	\ifx\shadow@sup\relax
		\setbox128=\hbox{}%
	\else
		\setbox128=\hbox{$\m@th#1\mkern-2mu\shadow@sup$}%
		\raisebox{#2}{\copy128}%
	\fi
	\ifx\shadow@sub\relax\else
		\raisebox{#3}{\hskip-\wd128$\m@th#1\mkern-3mu\shadow@sub$}%
	\fi
	\setbox128=\hbox{}%
\endgroup}%

\def\shdw{\@ifnextchar^{\shdw@\blank}{\@ifnextchar_{\shdw@\blank}{\@ifnextchar[\shdw@@\shdw@}}}
\def\sshdw{\@ifnextchar^{\sshdw@\blank}{\@ifnextchar_{\sshdw@\blank}{\@ifnextchar[\sshdw@@\sshdw@}}}

\def\shdw@#1{%
	\langle\!\langle#1\rangle\!\rangle
	\shadow@getindices}
\def\sshdw@#1{%
	\langle\mkern-4mu\langle\mkern-4mu\langle#1\rangle\mkern-4mu\rangle\mkern-4mu\rangle
	\shadow@getindices}
\def\shdw@@[#1]#2{\shdw@{#2}_{#1}}
\def\sshdw@@[#1]#2{\sshdw@{#2}_{#1}}

\def\mcircm{\mskip\thinmuskip\circ\mskip\thinmuskip}

\def\circm{\circ\mskip\thinmuskip}

\def\coinv{\mathrm{coInv}}

\newtoks\cat@name
\newtoks\cat@upper

\def\cat@generate#1{\begingroup
	\countdef\parser@state=0
	\countdef\parser@newstate=1
	\parser@state=0		% 0: active or macro; 1: lowercase or symbol; 2: uppercase
	\cat@name={}%        % a token list with the generated name
	\cat@upper={}%		% a token list with uppercase characters
	\cat@parse#1\relax
	\expandafter\xdef\csname cat@declared:\detokenize{#1}\endcsname{\the\cat@name}%
\endgroup}

\def\cat@parse#1{%
	\ifx\relax#1\relax
		\parser@newstate=4
	\else\ifcat\noexpand#1\relax\relax
		\parser@newstate=1
	\else\ifnum`#1>`Z\relax
		\parser@newstate=3
	\else\ifnum`#1<`A\relax
		\parser@newstate=3
	\else
		\parser@newstate=2
	\fi\fi\fi\fi
	\ifnum\parser@state=\parser@newstate\else
		\ifcase\parser@state
		\or
			\ifnum\parser@newstate<4
				\edef\next{\cat@name={\the\cat@name\/}}\next
			\fi
		\or
			\edef\next{\cat@name={%
				\the\cat@name\noexpand\mathcal{\the\cat@upper}%
			}}\next
			\cat@upper{}%
		\fi
		\parser@state\parser@newstate
	\fi
	\ifnum\parser@state=2
		\edef\next{\cat@upper={\the\cat@upper#1}}\next
	\else
		\edef\next{\cat@name={\the\cat@name\noexpand#1}}\next
	\fi
	\ifnum\parser@state<4	\expandafter\cat@parse\fi
}

\def\DeclareCategory#1{\expandafter\def\csname cat@declared:\detokenize{#1}\endcsname}%
\def\DeclareBicategory#1{\expandafter\def\csname bicat@declared:\detokenize{#1}\endcsname}%

\newcommand*{\cat}[1]{%
	\ifcsname cat@declared:\detokenize{#1}\endcsname\else
		\cat@generate{#1}%
	\fi
	\csname cat@declared:\detokenize{#1}\endcsname
}
\newcommand*{\ccat}[1]{%
	\ifcsname bicat@declared:\detokenize{#1}\endcsname
		\csname bicat@declared:\detokenize{#1}\expandafter\endcsname
	\else
		\mathbf{#1}%
	\fi}

\newcommand*{\catAdd}[1]{{#1}^{\oplus}}

\newcommand*{\Kar}{\mathrm{Kar}}

\DeclareCategory{Tan}{\mathcal T\mkern-5mu\mathit{an}}
\DeclareCategory{OTan}{\mathcal{OT}\mkern-5mu\mathit{an}}

\DeclareBicategory{linTan}{\scalars\ccat{Tan}}

\newcommand*{\Tan}{\cat{Tan}}

\newcommand*{\TTan}{\ccat{Tan}}

\newcommand*{\sLinks}{\mathit{Links}}

\newcommand*{\Links}{\cat L\!\mathit{inks}}
\newcommand*{\aLinks}{\Links(\Ann)}

\newcommand*{\Cob}{\cat C\mkern-2mu\textit{ob}}

\newcommand*{\CCob}{\ccat{Cob}}

\newcommand*{\TL}{\mathcal T\mkern-3mu\mathcal L}						% category
\newcommand*{\KBSM}{\mathcal{S}}

\newcommand*{\BN}{\cat{BN}}

\newcommand*{\aBN}{\BN(\Ann)}

\newcommand*{\BNu}{\widetilde{\cat{BN}}}

\newcommand*{\BBN}{\ccat{BN}}

\newcommand*{\BNq}{\BN\!_q}
\newcommand*{\aBNq}{\BNq(\Ann)}

\newcommand*{\BoBN}{\cat{BBN}}
\newcommand*{\BoBNq}{\BoBN\!_q}

\newcommand*{\evBNq}{\BNq^{\textit{ev}}}

\DeclareCategory{SMod}{\mathcal{SM}\mkern-3mu\mathit{od}}
\DeclareCategory{HoSMod}{\cat{SMod}_{\!/_{\mkern-2mu h}}}

\DeclareCategory{HoCom}{\cat{Com}_{\!/_{\mkern-2mu h}}\@ifnextchar({\!}{}}

\DeclareCategory{ComB}{\cat{Com}^b}
\DeclareCategory{HoComB}{\cat{Com}^b_{\!/_{\mkern-2mu h}}\@ifnextchar({\!}{}}

\begingroup
	\def\defname#1{\expandafter\protected\expandafter\gdef\csname#1\endcsname}%
	\def\defmodulecats#1#2{%
		\defname{#1}{#2}%
		\defname{l#1}##1{%{##1\text{--}#2}%
			\raisebox{-0.35ex}{$\scriptstyle##1$}#2%
		}
		\defname{r#1}##1{%{#2\text{--}##1}%
			#2\raisebox{-0.35ex}{$\scriptstyle##1$}%
		}
		\defname{b#1}##1##2{%{##1\text{--}#2\text{--}##2}%
			\raisebox{-0.35ex}{$\scriptstyle##1$}%
			#2%
			\raisebox{-0.35ex}{$\scriptstyle##2$}%
		}
	}%
	\defmodulecats{Mod}{\cat M\mkern-3mu\mathit{od}}
	\defmodulecats{PMod}{\cat{Proj}}
	\defmodulecats{ModG}{\cat{grMod}}
	\defmodulecats{PModG}{\cat{grProj}}
	\defname{Bimod}{\cat{Bimod}}
	\defname{gBimod}{\mathrm g\cat{Bimod}}
	\def\defbicats#1#2{%
		\defname{#1}{\ccat{#2}}%
		\defname{D#1}{\mathcal D(\ccat{#2})}%
		\defname{DD#1}{\mathcal D^-(\ccat{#2})}%
		\defname{g#1}{\@ifstar{\ccat{g#2}_0}{\ccat{g#2}}}%
		\defname{gD#1}{\mathcal D(\@ifstar{\ccat{g#2}_0)}{\ccat{g#2})}}%
		\defname{gDD#1}{\mathcal D^-(\@ifstar{\ccat{g#2}_0)}{\ccat{g#2})}}%

		\defname{e#1}{%
			\mathrlap{\hskip-0.15em\widetilde{\phantom{\rule{1.25em}{1.5ex}}}}%
			\ccat{#2}}%
		\defname{eD#1}{\mathcal D(%
			\mathrlap{\hskip-0.15em\widetilde{\phantom{\rule{1.25em}{1.5ex}}}}%
			\ccat{#2})}%
		\defname{eDD#1}{\mathcal D^-(%
			\mathrlap{\hskip-0.15em\widetilde{\phantom{\rule{1.25em}{1.5ex}}}}%
			\ccat{#2})}%
		\defname{eg#1}{%
			\mathrlap{\phantom{\ccat g}\hskip-0.15em\widetilde{\phantom{\rule{1.25em}{1.5ex}}}}%
			\@ifstar{\ccat{g#2}_0}{\ccat{g#2}}}%
		\defname{egD#1}{\mathcal D(%
			\mathrlap{\phantom{\ccat g}\hskip-0.15em\widetilde{\phantom{\rule{1.25em}{1.5ex}}}}%
			\@ifstar{\ccat{g#2}_0)}{\ccat{g#2})}}%
		\defname{egDD#1}{\mathcal D^-(%
			\mathrlap{\phantom{\ccat g}\hskip-0.15em\widetilde{\phantom{\rule{1.25em}{1.5ex}}}}%
			\@ifstar{\ccat{g#2}_0)}{\ccat{g#2})}}%
	}%
	\defbicats{BBimod}{Bimod}
	\defbicats{CCom}{Com}
	\defbicats{RRep}{Rep}
	\defbicats{DiagBimod}{DB}
\endgroup

\DeclareCategory{Mod}{\mathcal M\!\mathit{od}}
\DeclareCategory{gMod}{\mathrm g\mathcal M\!\mathit{od}}

\DeclareCategory{Rep}{\mathcal R\!\mathit{ep}}
\DeclareCategory{gRep}{\mathrm g\mathcal R\!\mathit{ep}}

\protected\def\Rep{\cat R\!\mathit{ep}}
\protected\def\gRep{\mathrm g\Rep}

\protected\def\Com{\mathrm{Com}}

\def\Com{\@ifnextchar^{\Com@}{\Com@@}}
\def\Com@^#1(#2){\mathit{Com}^{#1}(#2)}
\def\Com@@(#1){\mathit{Com}(#1)}

\def\HCom{\@ifnextchar^{\HCom@}{\HCom@@}}
\def\HCom@^#1(#2){\mathit{Com}^{#1}_{\!\raisebox{0.3ex}{$\scriptstyle/$}\mkern-2mu h}(#2)}
\def\HCom@@(#1){\mathit{Com}_{\!\raisebox{0.3ex}{$\scriptstyle/$}\mkern-2mu h}(#1)}

\def\SimplMod(#1){\cat{SMod}_{#1}}
\def\HoSimplMod(#1){\cat{HoSMod}_{#1}}

\newcommand*{\Ob}{\mathrm{Ob}}

\newcommand*{\Aut}{\mathrm{Aut}}
\newcommand*{\AlgEnd}{\mathrm{End}_{\mathrm{Alg}}}
\newcommand*{\End}{\mathrm{End}}
\newcommand*{\Hom}{\mathrm{Hom}}

\newcommand*{\F}{\mathcal{F}}

\newcommand*{\Id}{\mathrm{Id}}

\newcommand*{\cone}{\mathrm{cone}}

\newcommand*\Tor{\mathrm{Tor}}
\newcommand*\Tot{\mathrm{Tot}}

\DeclareMathOperator\Tr{Tr}
\DeclareMathOperator\hTr{hTr}
\DeclareMathOperator\vTr{vTr}

\newcommand*\triTr{\Tr^{\triangle}}

\newcommand*{\qTr}{\Tr_q}
\newcommand*{\qhTr}{\hTr_q}
\newcommand*{\qvTr}{\vTr_q}

\def\endofun{\Sigma}

\def\HHSymb{\mathit{HH}}
\newcommand*{\HoHom}{\HHSymb}
\newcommand*{\CHoHom}{\mathit{CH}}

\def\HMSymb{\mathcal{HH}}
\newcommand*{\HMHom}{\HMSymb}
\newcommand*{\CHMHom}{\mathcal{CH}}

\def\qHHSymb{\mathit{qHH}}
\newcommand*{\qHoHom}{\qHHSymb}
\newcommand*{\qCHoHom}{\mathit{qCH}}

\def\qHMSymb{\mathcal{q\mkern-2muHH}}
\newcommand*{\qHMHom}{\qHMSymb}
\newcommand*{\qCHMHom}{\mathcal{q\mkern-2muCH}}

\def\Ksp{K_0^{\mathit{\vphantom psp}}}		% '\vphantom p' fixes the vertical position of 'sp'
\def\Kss{K_0^{\mathit{\vphantom pss}}}		% '\vphantom p' fixes the vertical position of 'ss'

\newcommand*{\FTL}{\mathcal{F}_{\mathit{TL}}}
\newcommand*{\FKh}{\mathcal{F}_{\!\textit{Kh}}}

\newcommand*{\FFCK}{\mathbf{F}_{\!\textit{CK}}}

\newcommand*{\extFFCK}{\lowtilde[1.2]{\mathbf F}_{\!\textit{CK}}}

\newcommand*{\FAKh}{\mathcal{F}_{\!\Ann}}
\newcommand*{\FMKh}{\mathcal{F}_{\MB}}
\newcommand*{\FqAKh}{\mathcal{F}_{\!\Ann_q}}
\newcommand*{\FqMKh}{\mathcal{F}_{\MB_q}}

\def\trivmod{W}
\def\essmod{V}
\def\sepmod{V}
\def\nsepmod{U}

\def\KhBracket{\@ifstar\KhBracketScaled\KhBracketSimple}
\def\KhCube{\@ifstar\KhCubeScaled\KhCubeSimple}
\def\KhCubeSigned{\@ifstar\KhCubeSignedScaled\KhCubeSignedSimple}

\newcommand*{\KhCubeScaled}[1]{\mathcal{I}\left(#1\right)}
\newcommand*{\KhCubeSimple}[1]{\mathcal{I}(#1)}
\newcommand*{\KhCubeSignedScaled}[2]{\mathcal{I}\left(#1,#2\right)}
\newcommand*{\KhCubeSignedSimple}[2]{\mathcal{I}(#1,#2)}
\newcommand*{\KhBracketScaled}[1]{\left\llbracket#1\right\rrbracket}
\newcommand*{\KhBracketSimple}[1]{\llbracket#1\rrbracket}

\newcommand*{\KhCom}{\mathit{CKh}}

\newcommand*{\MKhCom}{\KhCom_{\MB}}
\newcommand*{\qAKhCom}{\KhCom_{\mkern-2mu\Ann_q}}

\newcommand*{\Kh}{\mathit{Kh}}

\newcommand*{\AKh}{\Kh_{\mkern-2mu\Ann}}
\newcommand*{\MKh}{\Kh_{\MB}}
\newcommand*{\qAKh}{\Kh_{\mkern-2mu\Ann_q}}

\newcommand*{\ChKh}{\mathit{Kh}}

\def\twosubs#1#2#3{%
	\raisebox{-0.85ex}{$\scriptstyle#1$}%
	#3%
	\raisebox{-0.85ex}{$\scriptstyle#2$}%
}

\def\arcalg{\@ifnextchar|{\arcalg@comp}{\arcalg@whole}}
\def\arcalg@comp|#1#2|#3{\twosubs{#1}{#2}{\!\left(\arcalg@whole{#3}\right)\!}}
\def\arcalg@whole#1{H^{#1}}

\def\arcmod{\@ifnextchar|{\arcmod@comp}{\arcmod@whole}}
\def\arcmod@comp|#1#2|#3{\twosubs{#1}{#2}{\arcmod@whole{#3}}}
\def\arcmod@whole#1{C_{Kh}(#1)}

\def\arcmap{\@ifnextchar|{\arcmap@comp}{\arcmap@whole}}
\def\arcmap@comp|#1#2|#3{\twosubs{#1}{#2}{\arcmap@whole{#3}}}
\def\arcmap@whole#1{\@ifnextchar'{\varphi_{#1}}{\varphi_{#1}^{\mathstrut}}}

\def\CKalg{\@ifnextchar|{\CKalg@comp}{\CKalg@whole}}
\def\CKalg@comp|#1#2|#3{\twosubs{#1}{#2}{\!\left(\CKalg@whole{#3}\right)\!}}
\def\CKalg@whole#1{A^{#1}}

\def\CKmod{\@ifnextchar|{\CKmod@comp}{\CKmod@whole}}
\def\CKmod@comp|#1#2|#3{\twosubs{#1}{#2}{\CKmod@whole{#3}}}
\def\CKmod@whole{\@ifnextchar^\CKmod@whole@\CKmod@whole@@}
\def\CKmod@whole@^#1#2{C^{#1}_{CK}(#2)}
\def\CKmod@whole@@#1{C_{CK}(#1)}

\def\CKmap{\@ifnextchar|{\CKmap@comp}{\CKmap@whole}}
\def\CKmap@comp|#1#2|#3{\twosubs{#1}{#2}{\CKmap@whole{#3}}}
\def\CKmap@whole#1{\@ifnextchar'{\varphi_{#1}}{\varphi_{#1}^{\mathstrut}}}

\def\arctimes#1{\mytimes@clap\otimes{\arcalg{#1}}}
\def\CKtimes#1{\mytimes@clap\otimes{\CKalg{#1}}}

\def\arcdtimes#1{\mytimes@clap{\hat\otimes}{\arcalg{#1}}}
\def\CKdtimes#1{\mytimes@clap{\hat\otimes}{\CKalg{#1}}}

\def\carcalg{\@ifnextchar|{\carcalg@comp}{\carcalg@whole}}
\def\carcalg@comp|#1#2|#3{\twosubs{#1}{#2}{\!\left(\carcalg@whole{#3}\right)\!}}
\def\carcalg@whole#1{A\!H^{#1}}

\def\carcmod{\@ifnextchar|{\carcmod@comp}{\carcmod@whole}}
\def\carcmod@comp|#1#2|#3{\twosubs{#1}{#2}{\carcmod@whole{#3}}}
\def\carcmod@whole#1{A\!C_{Kh}(#1)}

\def\carcmap{\@ifnextchar|{\carcmap@comp}{\carcmap@whole}}
\def\carcmap@comp|#1#2|#3{\twosubs{#1}{#2}{\carcmap@whole{#3}}}
\def\carcmap@whole#1{\@ifnextchar'{\varphi_{#1}}{\varphi_{#1}^{\mathstrut}}}

\let\Lef\Lambda

\def\ev{\mathit{ev}}
\def\coev{\mathit{coev}}

\def\HS{t^{\mathit H\!S}}%

\def\lsupscript@#1#2{\vphantom{#2}^{#1}\!#2}
\def\rsupscript@#1#2{#2^{#1}}
\def\plsupscript@#1#2{\mathllap{\vphantom{#2}^{#1}\!}#2}
\def\prsupscript@#1#2{#2^{\mathrlap{#1}}}

\def\ldual {\@ifstar{\plsupscript@*}{\lsupscript@*}}%
\def\rdual {\@ifstar{\prsupscript@*}{\rsupscript@*}}%
\def\lddual{\@ifstar{\plsupscript@{**}}{\lsupscript@{**}}}%
\def\rddual{\@ifstar{\prsupscript@{**}}{\rsupscript@{**}}}%

\let\seam\mu
\def\membr{\mathcal M}

\makeatother

\DeclareCategory{Vect}{\mathcal V\!\mathit{ect}}
\DeclareCategory{gVect}{\mathit g\mathcal V\!\mathit{ect}}

\DeclareBicategory{LinCat}{\scalars\text{-}\ccat{Cat}}
\DeclareBicategory{gLinCat}{\scalars\text{-}\ccat{gCat}}
\DeclareBicategory{gEndoAddCat}{\ccat{gEndoCat}^\oplus}
\DeclareCategory{qLinks}{\cat{Links}_q}
\input{macros.tex}

\allowdisplaybreaks
\psset{unit=2.5ex}

\theoremstyle{plain}
\newtheorem{itheorem}{Theorem}

\begin{document}

\title{Quantum Link Homology via Trace Functor I}

\author{%
	Anna Beliakova
	\and
	Krzysztof K.\ Putyra
	\and
	Stephan M.\ Wehrli
}

\date{}
\maketitle

\begin{abstract}
Motivated by topology, we develop a~general theory of traces and shadows  for an~endobicategory,
which is a~pair: bicategory $\ccat C$ and  endobifunctor $\endofun\colon \ccat C \to\ccat C$.
For a graded linear bicategory and a~fixed invertible parameter $q$,
we quantize this theory by
 using the~endofunctor $\endofun_q$ such that
$\endofun_q \alpha:=q^{-\deg \alpha}\endofun\alpha$ for any 2-morphism $\alpha$ and coincides with $\endofun$ otherwise.

Applying the quantized trace to the~bicategory of Chen--Khovanov bimodules 
we get a new triply graded link homology theory
called {\it quantum annular link homology}.
If $q=1$
we reproduce Asaeda--Przytycki--Sikora (APS) homology for links in a~thickened annulus.
%If $\Sigma$ is a~reflection, we obtain the~APS homology for links in a~thickened M\"obius band.
% We call the~resulting invariant the~\emph{quantum link homology}.
We prove  that our homology carries an action of $\Uqsl 2$, which intertwines the action of cobordisms.
In particular, the~quantum annular homology of an $n$--cable admits an action of the braid group,
which commutes with the~quantum group action and factors through the~Jones skein relation.
This produces a~nontrivial invariant for surfaces knotted in four dimensions.
Moreover, a~direct computation for torus links shows that the~rank of quantum annular homology groups 
does depend on the~quantum parameter $q$.
%Hence, our quantum link homology has a~richer structure.
\end{abstract}

\setcounter{tocdepth}{2}
\tableofcontents

\section{Introduction}

\subsection*{Background and Overview}
Trace are an~important source of topological invariants.
Given a~category $\cat C$,
a~trace is a~collection of functions $t_x$ defined
on the~endomorphism spaces $\cat C(x, x)$,
indexed by objects $x\in \cat C$, such that
for any pair of composable morphisms $x \from^g y \from^f x$ in $\cat C$
the~relation
$$
	t_x(g\circ f) = t_y( f\circ g)
$$
holds. When $\cat C$ is the~category of vector spaces,
then the~above relation determines the~trace uniquely (up to a~scalar)
as the~sum of the~diagonal entries of a~matrix representing the~endomorphism.
The~well-known Reshetikhin-Turaev package of link invariants
is obtained by applying traces to linear endomorphisms associated with tangles.
In this context, the~cyclicity relation is interpreted topologically as 
passing to the annular closure of a tangle.
%\vspace{0.5mm}

Every category admits the~\emph{universal trace}, defined as
\begin{equation}\label{def-trace}
	\Tr(\cat{C})
		:= \quotient{\displaystyle\coprod_{x\in\Ob(\cat C)} \mskip-10mu \cat C(x, x)}
		            {g\circ f \sim  f\circ g}
\end{equation}
where $f$ and $g$ run through all pairs of composable morphisms.
Any trace on $\cat C$ factorizes through $\Tr(\cat C)$.
%If $\cat C$ is linear, then the~universal trace coincides with
%the~$0$th Hochschild--Mitchell homology.
If $\cat C$ is an~additive category, then it receives the~Chern
character map from the~Grothendieck group $K_0(\cat C)$,
and hence can be considered as an alternative decategorification
functor \cite{TraceDecat, SL2asTrace}.
Another important feature of the trace is its functoriality.

Let us consider the category $\cat{Tan}$ of tangles, objects of which
are points on the $x$-axis in $\R^2$ and morphisms are tangles in
$\R^2 \times I$. Any functor $F\colon \cat{Tan} \to \cat A$
induces a~map
$$
	\Tr(F)\colon \Tr(\cat{Tan})\to\Tr (\cat{A}),
$$
which is a~universal annular link invariant associated with $\cat A$.
In the Reshetikhin--Turaev case, $\cat A$ is a~category of representations
of a~quantum group. 

For quantum $\mathfrak{sl}_2$ the~Reshetikhin--Turaev
construction was categorified by Chen and Khovanov in \cite{ChenKhov}.
For this purpose $\cat{Tan}$ was extended to a~bicategory $\TTan$
with 2--morphisms given by tangle cobordisms.
Chen and Khovanov defined a~projective bifunctor 
$$
	\extFFCK\colon \TTan \to \HCom^b(\ccat{gBirep})
$$
valued in the~homotopy category of graded \emph{sweet} bimodules%
\footnote{
	Following \cite{KhTangles} we say that a~bimodule is \emph{sweet} if
	it is finitely generated and projective when considered as a~left and
	as a~right module separately.
}
$\ccat{Birep}$. Here, `projective' means that the~bifunctor is defined
on 2--morphisms only up to a~sign.

In this paper, we develop a~theory of traces in endobicategories and 
apply it to the Chen--Khovanov construction to obtain a~new 
triply graded \emph{quantum} annular link homology theory.
The sign issue will be resolved in a~second paper \cite{gl2foams} by redefining our
and Chen--Khovanov's constructions using $\mathfrak {gl}_2$--foams.
We will use it in \cite{Matt2} to show that all known definitions
of the \emph{colored} Khovanov homology coincide in the~\emph{quantum} annular
setting when the~quantization parameter is \emph{not} a root of unity.

\subsection*{Traces in Endobicategories}

\begin{figure}[ht]%
	\begin{displaymath}
		\vcenter{\hbox{\includegraphics[scale=0.7]{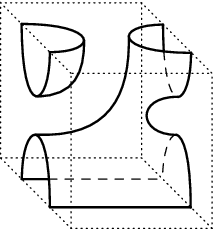}}}
		\quad\longmapsto\quad
		\vcenter{\hbox{\includegraphics[scale=0.7]{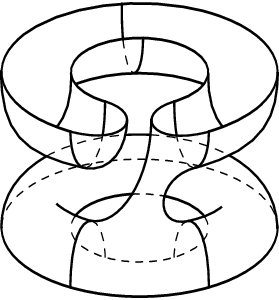}}}
		\quad\approx\quad
		\vcenter{\hbox{\includegraphics[scale=0.7]{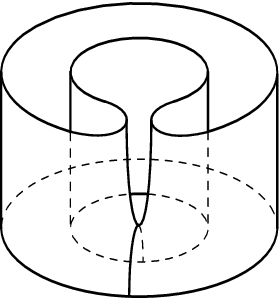}}}
	\end{displaymath}
	\caption{A horizontal trace of a~cobordism with corners.}%
	\label{fig:surface-closure}%
\end{figure}

Let us first present our categorical results.
There are two ways to define traces on a~bicategory:
%Note that in an $n$-category there are $n$ different compositions
%of morphisms on which a~cyclicity relation can be imposed.
%Hence there are $n$ different traces.
%These  traces 
%There is a natural filtration on the set of traces with the bigest one obtained by closing the 1-morphisms.
applying \eqref{def-trace} to morphism categories results in
the~\emph{vertical} trace, and a~more general \emph{horizontal} trace
was defined in \cite{SL2asTrace}.
Let us  visualize this construction for the bicategory $\CCob$, the~objects
of which are points on a~line, 1--morphisms are flat tangles in $\R\times I$,
and 2--morphisms are cobordisms in $(\RxI)\times I$.
The~\emph{horizontal trace} of a~cobordism between two different $(2,2)$--tangles 
is shown in Figure~\ref{fig:surface-closure},
where 1--morphisms are depicted horizontally.
Roughly speaking, the~horizontal trace is obtained
by imposing the~cyclicity relation on the~horizontal composition
of 1--morphisms in a way compatible with 2-morphisms
(see Section \ref{sec:categ-traces} for more details).

An~analogue of a~trace function for a~bicategory is called \emph{shadow},
and was defined by Ponto and Shulman \cite{Shadows}.
In this paper, inspired by topology, we further generalize the~notions
of the~universal trace of a~category,
as well as vertical and horizontal traces of a~bicategory.
%Furthemore, in the graded case, we 
% a quantization all the above mentioned 
%Then we apply our new theory to the Chen-Khovanov
%derived bicategory of bimodules to get a new 
%triply graded \emph{quantum} annular link homology. 
%
Recall that a~surface bundle $M$ over a~circle with monodromy $\phi$
is constructed from a~thickened surface $\sfce\times I$ by gluing
$\sfce\times \{0\}$ to $\sfce\times \{1\}$ along the~diffeomorphism $\phi$:
\begin{equation}\label{bundle}
	M := \ \quotient{\sfce\times I}{(p,1) \sim (\phi(p),0)} \; .
\end{equation}
To mimic this construction on a~more abstract level, we consider a~pair:
a~category $\cat C$ and an~endofunctor $\endofun\colon \cat C\to \cat C$. 
Then a~\emph{$\endofun$--twisted trace} of $\cat C$ is defined as
\begin{equation}\label{def-tw-trace}
	\Tr(\cat{C},\endofun) := \quotient
		{\displaystyle\coprod_{x\in\Ob(\cat C)} \mskip-10mu \cat C(x,\endofun x)}
		{\endofun f\circ g \sim g\circ f}
\end{equation}
where $f$ and $g$ run through all pairs of composable morphisms
$\endofun x \from^g y \from^f x$ in $\cat C$. This construction
is functorial with respect to $(\cat C,\endofun)$ and universal.
%in the~sense that any trace function factorizes through it.
We recover the usual trace of $\cat C$ when $\endofun$ is the~identity functor.

Let $\cat{Tan}$ be the~category of isotopy classes of tangles in $\R^2\times I$.
The~rotation of the~plane induces an~endofunctor $\endofun$ on $\cat{Tan}$,
for which $\Tr(\cat{Tan}, \endofun)$ recovers links in the thickened M\"obius
band.%, see Figure~\ref{fig:links-in-MB}.

We observe that any natural transformation of the endofunctor on an~endocategory
induces a new trace function.
For example, fixing a~framed $(1,1)$--tangle $\tau$ we can define a~natural
transformation of the~identity functor of $\cat{Tan}$ by sending a collection
of $n$ points to the~$n$--cabling $\tau^n$, viewed as a morphism in $\cat{Tan}$.
It is a~natural transformation because of the~equality $\tau^m\circ T =
T\circ \tau^n$, which holds for any $(m,n)$--tangle $T$. The same euqality also
implies that the~map sending an~$(n,n)$--tangle $T$ to the annular closure of the composition $T\circ \tau^n$ is well-defined trace function on $\cat{T\!an}$.
In particular, one can compose any $(n,n)$--tangle with a~full twist before closing.
Hence, traces on endocategories of tangles encode information not only about
a~simple annular closure $\widehat{T}$ of any $(n,n)$-tangle $T$, but also on all
satellites of $\widehat T$.

Triangulated categories are treated separately. They admit triangulated traces,
which are additive with respect to homomorphisms of distinguished triangles,
see \eqref{eq:trian-additivity}. We show that twisted triangulated traces on
the~homotopy endocategory of bounded complexes $\HCom^b(\cat C)$ are
Lefschetz traces, which are unique extensions of twisted traces on $\cat C$.
In particular, there is an~isomorphism between the universal triangulated trace
$\triTr(\HCom^b(\cat C),\endofun)$ and $\Tr(\cat C,\endofun)$.

Twisted traces admit a~quantization when $\cat C$ is pregraded.%
\footnote{
	A~category is \emph{pregraded} if its sets of morphisms are graded
	modules.}
For that fix an~invertible parameter $q$ and consider the~endofunctor
$\endofun_q$ that coincides with $\endofun$ on objects,
but $\endofun_q f:=q^{-|f|}\endofun f$ for a~homogeneous morphism
$f$ of degree $|f|$. This results in the~\emph{quantum universal trace}
$$
	\qTr(\cat C, \endofun) := \Tr(\cat C,\endofun_q) = \quotient
		{\displaystyle\coprod_{x\in\Ob(\cat C)} \mskip-10mu \cat C(x,\endofun x)}
		{\endofun f\circ g - q^{|f|} g\circ f} \; .
$$
%The~choice of $q$ is often limited, so that it is common to extend first
%the~coefficients of $\cat C$ to a~flat extension $\escalars$ of $\scalars$,
%such as $\qscalars$.
%
Further, we recall and generalize of quantize traces that appear commonly
in algebraic and topological contexts, including the~Hattori--Stallings trace,
the~Lefschetz number, traces in pivotal categories, and various closures
of tangles.

We also investigate the~quantum analogue of the~Hochschild--Mitchell
homology for endocategories and prove the~quantum version of the $K\noic$--theoretic
equivalence, generalizing that of Keller \cite{InvForCyclicHom}.
Roughly speaking, it states that the homotopy class of the~chain map
(between quantum Hochschild--Mitchell chain complexes)
induced by a~bimodule $M$ depends only on the~image of $M$
in the Grothendieck group, see Proposition \ref{prop:K-invariance-of-HH}.
This very powerful tool is used later to find the~quantum Hochschild homology
of Chen--Khovanov algebras.

\medskip

To categorify the~constructions described above,
we define the~notions of a~\emph{$\endofun$--twisted preshadow} and
a~\emph{$\endofun$--twisted horizontal trace} for a~bicategory $\ccat C$ and
its endobifunctor $\endofun$, generalizing those introduced in \cite{Shadows}
and \cite{SL2asTrace}. The~pair $(\ccat C,\endofun)$ is called
an~\emph{endobicategory}.
We prove that any twisted preshadow factors through the~twisted horizontal trace
under the~mild assumption that $\ccat C$ has left duals.
Hence, for bicategories with duals the~horizontal trace $\hTr(\ccat C, \endofun)$
is a~universal preshadow.

When $\ccat C$ is pregraded, a~deformation of $\endofun$, defined  by setting
$\endofun_q f:=q^{-|f|}\endofun f$ for any homogeneous 2-morphism $f$ of degree
$|f|$,  leads to a~theory of quantum preshadows and a~quantum horizontal trace.

The~horizontal trace can be viewed as a~categorification of the~trace of
an~endocategory in the~following sense.
Let $\Pi\colon n\ccat{Cat} \to (n{-}1)\ccat{Cat}$ be the~decategorification
functor, which forgets the $n$--morphisms and identifies isomorphic
$(n-1)$--morphisms. Then there is a~natural bijection of sets
$\Pi(\hTr(\ccat C,\endofun)) \approx \Tr(\Pi\ccat C,\Pi\endofun)$ for any small endobicategory $(\ccat C,\endofun)$ with both left and right duals
(see Section~\ref{sec:hTr-categ-Tr}).

Applying the~construction \eqref{def-tw-trace} to morphism categories in $\ccat C$
results in the~vertical trace $\vTr(\ccat C, \endofun)$.
There is a~functor from the~vertical to the~horizontal trace of an~endobicategory,
which is full and faithful, but not necessarily surjective on objects.
Quite often a~preshadow can be restricted to a~collection of traces on morphisms
sets of the~vertical trace, for instance the~shadow computing coinvariants in
bimodules restricts to the~Hattori--Stallings trace
(see Section~\ref{sec:shadow-coinv}).

Of particular interest to us are preshadows on bicategories of complexes.
We extend every preshadow on $(\ccat C,\endofun)$ to a~preshadow on
$\HCom^b(\ccat C, \endofun)$, which---by the~analogy to traces---we call
the~\emph{Lefschetz preshadow}. These satisfy a~higher analogue of additivity
with respect to distringuished triangles
(see Proposition~\ref{prop:additivity-of-trian-shadow} for a~precise statement).
%Althouth the~proposition can be stated for general triangulated preshadows,
%we were not able to prove it beyond Lefschetz traces. It seems that the~main
%problem is to find a~right set of axioms for a~triangulated preshadow.
%Although not stated in this way, the~problem was discussed by May in
%\cite{MonTriCat} with a~solution that includes quite large diagrams to
%be checked. Because of that we decided to avoid working with general
%triangulated preshadows and defining a~horizontal trace of a~locally
%triangulated bicategory.

Again we provide many examples of preshadows on endobicategories.
Among them are twisted spaces of coinvariants, quantum Hochschild homology,
and component-wise Hochschild homology of a~complex of bimodules (which is
an~example of a~Lefschetz preshadow mentioned above).
Furthermore, we compute explicitly twisted horizontal traces of $\ccat{Tan}(\sfce)$,
the~bicategory of oriented tangles in a~thickened surface $\sfce$ and oriented
tangle cobordisms (see Appendix~\ref{sec:def-tangles} for a~precise definition).

\begin{repeattheorem}{Theorem}{thm:hTr-of-FxI}
	Let $M$ be a~surface bundle with fiber $\sfce$ and monodromy
	$\phi\in\mathit{Diff}(\sfce)$. There is an~equivalence of categories
	\begin{equation}
		\hTr(\TTan(\sfce),\phi_*) \simeq \Links(M)
	\end{equation}
	where $\phi_*(S) := (\phi\times\id\times\id)(S)$
	for a~cobordism $S\subset\sfce\times I\times I$.
\end{repeattheorem}

\noindent
Here $\Links(M)$ stands for the~category of oriented links in $M$ and
oriented link cobordisms in $M\times I\noic$.
In this paper we are interested in the~case of when $M$ is a~solid torus,
seen as a~thickened annulus $\AxI$ or as a~thickened M\"obius band
$\MBxI$. We write shortly $\cat{Links}(\Ann)$ and $\cat{Links}(\MB)$
respectively for the~categories of links. Both arise by considering tangles
in a~thickened plane with $\phi$ being the~identity map or rotation
by 180 degrees.
Elements of $\cat{Links}(\Ann)$ are called \emph{annular links}.

An~analogous result holds for $\CCob$, the~bicategory of \emph{flat} tangles
embedded in a~stripe $\RxI$, and for its famous quotient by local relations
listed in Theorem \ref{thm:invariance-of-KhBracket}, known as the~Bar-Natan
bicategory $\BBN$.
Since a~real line admits only two diffeomorphisms up to isotopy,
the~identification space in this case is either an~annulus $\Ann$
or as a~M\"obius band $\MB$.
Hence, the~horizontal traces of $\BBN$ with the~corresponding endofunctors
recover the~Bar-Natan categories $\BN(\Ann)$ and $\BN(\MB)$.
The~identification $\BN(\Ann) \cong \hTr(\BBN)$ is also proven in
\cite{AnnKhRoz}.

%There is also an analogue of this Theorem for
Notice that 2-morphisms in $\TTan(F)$ are graded by the~Euler characteristic,
so that $\Links(M)$ can be quantized. Namely, we get the~category
$$
	\cat{qLinks}(M) := \qhTr(\scalars\TTan(F),\phi_*)
$$
where $\scalars$ is a~fixed ring containing the~value $q$, and
$\scalars\TTan(F)$ is the~linear extension of $\TTan(F)$, which 2--morphisms
are formal linear combinations of cobordisms with coefficients from $\scalars$.
This category it admits the~following graphical description. Represent the~image
of $F\times\{0\}$ in $M$ as a~cooriented membrane, with the~coorientation induced from the~orientation of $I$. The~objects of $\cat{qLinks}(M)$ are
oriented links in $M$ that intersect the~membrane transversely,
whereas morphisms are link cobordisms up to isotopies, where an~isotopy moving
a~cobordism of Euler characteristic $d$ through the~membrane scales
the~cobordism by $q^{\pm d}$:
\begin{equation*}
	\psset{unit=3ex}%
	\relRCapOffMembrane = q \relRCapAtMembrane
	\hskip 2cm
	\relSaddleLeftToMembrane = q^{-1} \relSaddleRightToMembrane 
\end{equation*}

Because the~local relations defining $\BBN$ are homogeneous,
we construct analogously the~\emph{quantum annular Bar-Natan category} 
as the~additive closure of the~quantum horizontal trace
\begin{equation}
	\BNq(\Ann) := \qhTr^\oplus(\BBN).
\end{equation}
In important feature of $\BNq(\Ann)$ is that a~torus wrapped once along
the~annulus evaluates to $q+q^{-1}$ instead of $2$, see \eqref{eq:torus-in-BNq}.

\subsection*{Link Homologies via Traces}

Our next goal is  to construct new functorial invariants of links in a~solid torus
and in a~M\"obius band.
The~first annular link homology theory was constructed
by Asaeda, Przytycki, and Sikora \cite{APSHom},
and it can be rephrased as applying the~APS TQFT functor
$$\FAKh\colon\aBN \to \Mod(\scalars)$$ to~the formal Khovanov bracket \cite{Cob}.
Grigsby, Licata, and Wehrli observed in \cite{KhHomSL2} that the~annular link homology admits an~action of $\LieSL_2$ commuting with the maps induced
by annular cobordisms. This motivates the~search of a~quantized annular homology,
on which the~quantum $\LieSL_2$ acts.

Our aim is to use $\BNq(\Ann)$ to deform the~APS construction.
Even though it is immediate that a~torus evaluates in $\BNq(\Ann)$ to $q+q^{-1}$
when it intersects the~membrane in a~circle, this does not extend
naively to a~TQFT functor.
Instead, we use the whole algebraic machinery developed in the first part of
the~paper to construct such a~TQFT, especially Theorem \ref{thm:hTr-of-FxI},
the~universality of the horizontal trace, and the~Lefschetz preshadow on
the~homotopy bicategory of Chen-Khovanov bimodules, induced by the~quantum Hochschild homology.

For this purpose, let us examine the~homology of tangles constructed by Chen
and Khovanov. In \cite{ChenKhov} they defined a~family of diagrammatic algebras
$\CKalg n$, commonly called \emph{arc algebras}, associated a~graded
$(\CKalg m, \CKalg n)$--bimodule with each flat $(m,n)$--tangle, and constructed
a~bimodule map for any cobordism of flat $(m,n)$--tangles.
This gives the~bifunctor
$$
	\FFCK\colon \BBN \to  \ccat{gBirep}
$$
valued in the~bicategory of sweet graded bimodules. A~precomposition with the~formal
Khovanov bracket $\KhBracket{\blank}$ results in
$$
	\extFFCK\colon \TTan \to  \HCom^b(\ccat{gBirep}),
$$
which assigns to an $(m,n)$-tangle $T$ a~chain complex $\CKmod T$ of~graded
$(\CKalg m, \CKalg n)$--bimodules, the~homotopy class of which is an~isotopy invariant of $T$.

The~algebra $\CKalg n$ categorifies the~$n$--th tensor power of the~fundamental
representation $V_1$ of $\Uqsl 2$. It was first introduced by Braden
\cite{SheavesOnGr} using generators and relations to describe the~category of
perverse sheaves on Grassmannians. The~arc algebras and their representations
were independently studied by Brundan and Stroppel \cite{ChenKhovIsCellular, ChenKhovKoszulity}.

We can now use our previous results to extend any given preshadow on
$\ccat{gBirep}$ to the~homotopy category of complexes, and then
pull it back along $\extFFCK$ to obtain a~preshadow on $\TTan$.
Since all involved bicategories have duals, such a~preshadow factors
through $\qhTr(\TTan)$ defining an~annular link homology theory.

An~immediate choice of a~preshadow on $\ccat{gBirep}$ is the~Hochschild homology
$\HoHom_{\!\bullet}$ of Chen--Khovanov bimodules.
However, we can also utilize the~grading and use the~\emph{quantum Hochschild
homology} $\qHoHom_{\!\bullet}$, a~one parameter deformation of $\HoHom_{\!\bullet}$
that factors through the~quantum horizontal trace of $\ccat{gBirep}$.
This results in the~following commuting diagram
\begin{equation}\label{diag-commut}
	\begin{diagps}(-11em,-0.5ex)(11em,14ex)
		\node tl(-11em,12ex)[\TTan]
		\node tm(  0em,12ex)[\HCom^b(\BBN)]
		\node tr( 11em,12ex)[\HCom^b(\ccat{gBirep})]
		\node bl(-11em, 0ex)[\qhTr(\TTan)]
		\node bm(  0em, 0ex)[\HCom^b(\qhTr(\BBN))]
		\node br( 11em, 0ex)[\HCom^b(\Uqsl 2)]
		\arrow|a{npos=0.45}|{->}[tl`tm;\KhBracket{\blank}]
		\arrow|a{npos=0.45}|{-->}[bl`bm;\KhBracket{\widehat\blank}_*]
		\arrow{->}[tm`tr;\extFFCK]
		\arrow{-->}[bm`br;\FqAKh]
		\arrow|b|{->}[tl`bl;\widehat{(\blank)}]
		\arrow|b|{->}[tm`bm;\widehat{(\blank)}]
		\arrow|a|{->}[tr`br;\mathit{qHHH}]
	\end{diagps}
\end{equation}
where the horizontal maps are functors and the vertical one are preshadows.
The~notation $\mathit{qHHH}$ means that we apply
$\qHoHom_{\!\bullet}$ to each bimodule in the~Chen--Khovanov complex separately.

To identify the~bottom right corner of this diagram
we apply the $K\noic$--theoretic invariance to arc algebras.
We deduce that the~Chern character map
$K_0(A^n)\to \qHoHom_{\!\bullet}(A^n)$ is actually an isomorphism.
Hence, the higher Hochschild homology of $A^n$ vanishes and
$\qHoHom_{\!0}(A^n) \cong V^{\otimes n}_1$. This relates
$\qhTr^\oplus(\ccat{Birep})=\BNq(\Ann)$
with the graded representation category $\gRep(\Uqsl 2)$.
Finally, to construct our TQFT functor
$$
	\FqAKh\colon \BNq(\Ann)\to \gRep(\Uqsl 2)
$$
we use the naturality of the Chern character map to identify 
$\BNq(\Ann)$ with (the~graded extension of) the~Temperley--Lieb category
$\TL$ and $\FqAKh$ with the~faithful functor from $\TL$ to $\gRep(\Uqsl 2)$.

This produces a~well-defined homology for annular links.
To any annular closure $L$ of an~$(n,n)$-tangle $T$ it assigns
a~chain complex
\begin{equation}\label{AKh-complex}
	C\noic\qAKh(L) := \mathit{qHHH}_{\!\bullet}(\CKalg n, \CKmod T)
		= \mathit{qHHH}_{\!0}(\CKalg n, \CKmod T) 
\end{equation} 
and to any annular cobordism an~induced chain map.
In addition to the~homological and quantum grading, $C\noic\qAKh(L)$
admits a~third grading, called the~annular grading, which comes
from the~weigh decomposition of the~Chen--Khovanov invariant.
 
The~chain maps induced by annular cobordisms are defined only up to multiplication
by $\pm q^{\pm1}$: the~sign comes from the~fact that $\extFFCK$ is merely
a~projective functor, and the~overall power of $q$ is not well-defined,
because it depends on the~presentation of an~annular cobordism as a~horizontal closure of a~surface with corners. We call this behavior `$q$--projective'.

\begin{repeattheorem}{Theorem}{thm:qAKh}
	The~quantum annular homology\/ $\qAKh(L)$ is a~triply graded invariant
	of an~annular link\/ $L$, which is $q$--projectively functorial with respect
	to annular link cobordisms.
	Moreover, it admits an~action of the~quantum group\/ $\Uqsl 2$ that
	commutes with the~differential and the~maps induced by annular link
	cobordisms intertwine this action.
\end{repeattheorem}

It was conjectured by Auroux, Grigsby, and the third author in
\cite{AnnKhAsTrace} that the~ASP homology of the~braid closure
coincides with the~Hochschild homology of the Chen--Khovanov complex 
associated to that braid. They checked this conjecture 
in the~next-to-top annular grading. 

Observe that the~quantum annular homology
arises as the~second page of the~spectral sequence associated
to $\qCHoHom_{\!\bullet}(\CKalg n, \CKmod T)$, the~bicomplex
computing the~quantum Hochschild complex for a~complex of bimodules.
Using the vanishing of the higher quantum Hochschild homology groups
for arc algebras, we can actually show that this spectral sequence 
collapses at the second page. Hence we have the~following.
 
\begin{repeattheorem}{Theorem}{thm:HH(CK)-vs-AKh}
	Let\/ $\widehat T\noic$ be the~annular closure of an~$(n,n)$--tangle\/ $T\noic$.
	Then there is an~isomorphism
	\begin{equation}
		\qAKh(\widehat T) \cong \qHoHom_{\!\bullet}(\CKalg n; \CKmod T),
	\end{equation}
	natural with respect to chain maps associated to tangle cobordisms.
	The~annular grading in $\AKh(\widehat T)$ corresponds to the~weight
	decomposition of\/ $\CKmod{T}$.
\end{repeattheorem}

When $q=1$, Theorem \ref{thm:HH(CK)-vs-AKh} proves Conjecture 1.1 from
\cite{AnnKhAsTrace} and motivates us to call our new link invariant
the~\emph{quantum annular link homology}.
 
We also show by explicit computation that the~rank
of the~quantum annular link homology for $(2,n)$ torus links does depend
on the~quantum parameter $q$. Hence, the quantized theory is richer than
APS annular link homology (which is the~case $q=1$).

There are many papers devoted to fixing the~sign issue to get
a~strictly functorial Khovanov homology \cite{VogelFun, BlanchetFoams, CMW,
sl2Webs}. In our next paper \cite{gl2foams} we will reconstruct the~Chen--%
Khovanov invariant using Blanchet foams to obtain a~strictly functorial
quantum annular link homology theory with an intertwining
$\cat U_q ({\mathfrak{gl}}_2)$-action.

\subsection*{Applications and Generalizations}
Let us first argue why the quantum annular link homology  
is actually more sensible to 4D topology than APS. 

A~link cobordism $W\colon L\to L$ induces a~(projective) map
$W_*\colon \Kh (L)\to\Kh(L)$ on Khovanov homology of $L$. 
The trace class of this endomorphism
is an invariant of the annular closure $\widehat W$ of $W$
and is characterised by its Lefschetz trace 
$$
	\Lambda (W_*)=\sum_{i,j} (-1)^i q^j\tr  W_*^{i,j} ,
$$
where $W_*^{i,j}$ is the component of $W_*$ in homological grading $i$ and
quantum grading $j$, and $\tr$ is the~Hattori--Stallings trace (which is
the~usual linear trace when we work over a~field). We prove the following result.

\begin{repeattheorem}{Theorem}{thm:inv-for-closed-sfce}
	Let\/ $\widehat W \subset \S\times\R^3$ be a~closed surface obtained
	as an~annular closure of a~link cobordism\/ $W\colon L\to L\noic$ with\/
	$L\subset\R^3$. Then $\qAKh(\widehat W) = \Lef(W_*)$ is the~graded Lefschetz
	trace of\/ $W_*\colon \Kh(L)\to \Kh(L)$, the~endomorphism of the~Khovanov
	homology of\/ $L$.
	In particular, $\qAKh(\S\times L)$ coincides with the~Jones polynomial\/ $J(L)$.
\end{repeattheorem}

Note that the APS invariant is trivial for closed surfaces.
Therefore, our invariant is a~nontrivial deformation.

Next, we establish a~nontrivial braid group action on
the~quantum annular homology of cablings of a~framed long knot. 
Consider a~framed annular knot $K\subset\AxI$. It defines an~embedding
$\nu_K\colon \AxI \to \AxI$ with the~tubular neighborhood 
of $K$ as its image, and hence induces a~functor $K_*\colon \Tan \to \aLinks$
mapping oriented points $B\subset \R^2$ to a collection of circles
$K^B\colon= \nu_K(\S\times B)$, and an~oriented tangle $T\subset \R^2\times I$ 
to the~oriented cobordism $K^T := \nu_K(\S\times T)$ between these circles.
Applying the~quantum annular homology produces a~map of homology
\[
	\qAKh(K^T)\colon \qAKh(K^B) \to \qAKh(K^{B'})
\]
for any oriented tangle $T \in \cat{Tan}(B,B')$, defined up to an~overall
power of $q$ (we work in characteristic 2 here to avoid the~sign issue).
This gives rise to a~$q$--projectively functorial action of $\Tan$, i.e.
\[
	\qAKh(K^{T'T}) = q^k \qAKh(K^{T'})\circ \qAKh(K^T),
\]
for any composable tangles $T$, $T'$, and some $k\in\Z$.
The~action was first observed in \cite{KhHomSL2} in the~non-quantized setting.
We compute this action and show that it factors through the~Jones skein relation.

\begin{repeattheorem}{Theorem}{thm:action-on-qAKh}
	Let\/ $K$ be a~framed annular link, considered as an~object in $\cat{qLinks}$.
	There is a~functorial action of\/ $\Tan\noic$ on the~quantum annular homology
	of oriented cablings of\/ $K\noic$, that takes a~tangle\/ $T\noic$ to the~chain
	map $\qAKh(K^T)$, and which intertwines the~action of\/ $\Uqsl 2$.
	It factors through the~Jones skein relation
	\begin{equation*}
	\tag{\ref{eq:skein-relation-for-K^T}}
		q^2 \qAKh(K^{\fntNegCr}) - q^{-2} \qAKh(K^{\fntPosCr})
		= (q-q^{-1}) \qAKh(K^{\fntNoCr})
	\end{equation*}
	if\/ $K$ intersects the~membrane in one point.
\end{repeattheorem}

The~condition on the~characteristic of the~ring of scalars can be dropped
once a~strictly functorial version of Chen--Khovanov is used. This will be
address in our next paper \cite{gl2foams}.

Let us now discuss few modifications and generalizations of Theorem
\ref{thm:qAKh}.
Recall that we have identified links in a~thickened M\"obius band with
twisted closures of tangles. Therefore, one may hope to get a~parallel
quantization of the~APS homology for links with diagrams on a~M\"obius band.
Indeed, such a~construction is described in Section~\ref{sec:quantum-MB}.
The~basic algebraic ingredient is the~\emph{twisted} quantum Hochschild homology,
defined for a~pair of a~graded algebra $A$ and its automorphism. Here, we consider
the~automorphism of arc algebras induced by the~horizontal flip of diagrams.
Again, setting $q=1$ recovers the~original APS homology.

We have already observed that any preshadow applied to $\CKmod T$ produces
an~annular link invariant.
In particular, we can precompose an~$(n,n$)--tangle $T$ with and $n$--cabling
of a~framed $(1,1$)--tangle, such as a~sequence of full twists. Geometrically,
these theories correspond to different embeddings of the~annulus into $\R^3$.
In the~other generalization we again fix a~$(1,1)$--tangle $T$ and we assign
the~quantum annular homology of either $\widehat T$ or its mirror image with
an~essential circle in an~annulus. We then use the~duality between the~homology
of a~link and its mirror image to define the~differential in the~chain complex.

In a~follow up \cite{Matt2} we will define a~quantum homology for colored
annular links. A very surprising feature of this construction is that
the~Cooper--Krushkal infinite complex categorifying the $n$-th Jones--Wenzl
idempotent becomes finite in the~annular closure for a~generic value of $q$,
and it does coincide with the~annular $n$-colored Khovanov complex.
%\comment{add references}.

\subsection*{Strategy for the Proof of Theorem \ref{thm:qAKh}}

To construct our new homology theory we follow the~general recipe for Khovanov
homology as described in \cite{Cob}. There, a~link $L$ in a~thickened surface
$\sfce$ is assigned a~formal complex $\KhBracket{L}$ in the~Bar-Natan skein
category $\BN(\sfce)$, whose objects are non-intersecting curves in $\sfce$
and morphisms are cobordisms in $\sfce\times I$. To construct actual
homology one then applies to $\KhBracket{L}$ a~certain TQFT functor
$\F_F\colon\BN(\sfce)\to \Mod(\scalars)$, where $\scalars$ is a~fixed ring
of scalars. The main body of this paper is devoted to the construction of such
a~functor when $F=\Ann$.
The~theory of preshadows and horizontal traces on endobicategories is needed to
guarantee the~existence of the~TQFT functor $\FqAKh$ in \eqref{diag-commut}.

Let us first explain why the two bottom arrows in \eqref{diag-commut} are dashed.
Strictly speaking, a~quantum preshadow on $\HCom^b(\BBN)$ factorizes through
$\qhTr(\HCom^b (\BBN))$ rather than $\HCom^b(\qhTr(\BBN))$.
Lemma~\ref{lem:lefshetz-shadow} resolves this problem by showing that we actually
can apply $\qHoHom_{\!\bullet}$ to each bimodule separately, rather than to
the complex of bimodules, and deal with flat tangles rather than formal complexes
of them. Hence, the~pullback shadow $(\FFCK)^*\qHoHom_{\!0}$ defined on
$\BBN$ factors through $\BNq(\Ann)$, inducing a~linear functor
$$
	\FqAKh\colon \BNq(\Ann)\to \Mod(\scalars).
$$
To identify the~target of this functor with $\Rep (\Uqsl 2)$ we need to find
an~isomorphism
\begin{equation}\label{Keller-iso}
	V_1^{\otimes n} \cong \qHoHom_{\!0}(\CKalg n)
\end{equation}
where $\CKalg n$ is the bimodule associated by $\FFCK$ with the trivial
$(n,n)$--tangle, and $\qHoHom_{\!0}(\CKalg n)$ is the~module assigned to
$n$ parallel essential circles in $\Ann$.

When $q=1$, the~existence of the~isomorphism \eqref{Keller-iso}
follows from the~invariance result for Hochschild homology due to Keller
\cite{InvForCyclicHom}: if $A$ is a~finite dimensional algebra,
$E\subset A$ is a~separable subalgebra such that $A=E\oplus\mathrm{rad}(A)$,
each simple $A$--module is one-dimensional, and $A$ has finite global dimension,
then $\HoHom_{\!\bullet}(E) \cong \HoHom_{\!\bullet}(A)$. It was proven by Brundan
and Stroppel that arc algebras have finite global dimension when $\scalars$ is
a~field \cite{ChenKhovIsCellular}, and other conditions hold after setting $E=A^n_0$,
the degree zero part of $A^n$. The~latter consists of idempotents, and hence is
isomorphic to $\scalars^{2^n}$. This provides the isomorphism \eqref{Keller-iso}
for $q=1$.

The~0th quantum Hochschild homology of $\CKalg n$ can be computed by hands,
which is all one needs to understand the~construction of our invariant. However,
computation of higher Hochschild homology is needed to identify the~invariant
with the~total Hochschild homology of Chen--Khovanov bimodules, as conjectured
in \cite{AnnKhAsTrace}. For that we reprove the~Keller's result for quantum
and---more generally---for twisted Hochschild homology of an~algebra $A$ by
identifying the~latter with quantum Hochschild--Mitchell homology of the~category
of finite dimensional representations of $A$, twisted by an~appropriate endofunctor.
The~advantage of replacing algebras and bimodules with representation categories
and functors is that the~latter provides a~more flexible framework, in which
the~action of bimodules on homology, and so the~$K\noic$--theoretic invariance,
is easier to understand, see Section~\ref{sec:K-invariance}.

Further we observe that the canonical embedding of $\qvTr(\BBN)$ into
$\qhTr(\BBN)$ is an equivalence of categories. This implies that every
cobordism in $\qhTr(\BBN)$ is a~linear combination of those of the form
$\S\times T$, where $T$ is a~Temperley-Lieb diagram. Therefore,
the~isomorphism \eqref{Keller-iso} determines $\FqAKh$ completely,
leading to a~commuting diagram of functors
\begin{equation}\label{comm-triangle}
	\begin{diagps}(0em,-0.5ex)(8em,11ex)
			\Vtriangle<8em,9ex>[%
				\TL`\BoBNq(\Ann)`\gRep(\Uqsl{2});
				\simeq`\F_{\mathit{TL}}`\FqAKh
			]
	\end{diagps}
\end{equation}
where the~horizontal functor sends a~flat tangle $T$ from $\TL$ to
the~annular cobordism $\S\times T$, and $\F_{\mathit{TL}}$ is
the~Reshetikhin--Turaev realization of the~Temperley--Lieb diagrams
as $\Uqsl 2$ intertwiners between tensor powers of $V_1$, the~fundamental representation of the~quantum group. The~category $\BoBNq(\Ann)$ is
a~quotient of $\BNq(\Ann)$ by a~certain local relation defined in \cite{Boerner}.

It follows from that our annular quantum TQFTs functor
$\FqAKh$, when restricted to $\BoBNq(\Ann)$, is faithful.
Once the~relation between $\FqAKh$ and $\F_{\mathit{TL}}$ is established,
we can directly see that setting $q=1$ recovers the~APS TQFT functor.

\subsection*{Outline}

The paper is organized as follows.
In Section \ref{sec:traces} we develop the~theory of traces in endocategories.
We discuss fundamental properties of the~universal twisted trace, such as
the~universal property, functoriality, and connection to the~additive
Grothendieck group. The~special cases of triangulated and graded endocategories
are treated separately. The~section ends with a~number of examples.

Section \ref{sec:categ-traces} is devoted to the~construction of categorical traces
in endobicategories. We first introduce the~notion a~(pre)shadow and construct
the~twisted horizontal trace. Its universal property and functoriality are
shown in Sections~\ref{sec:univeral-shadow} and \ref{sec:functoriality-of-hTr}.
Then we discuss the~connection with the~vertical trace and prove that the~horizontal
trace is a~categorification of the~universtal trace of an~endocategory.
Triangulated and quantum preshadows are discussed further. Again, all examples
are gathered at the~end of the~section. Theorem~\ref{thm:hTr-of-FxI} is proven
in Section~\ref{sec:links-in-sfce-bundle}.

After the~theory of categorified traces is established we take a~closer look
at Hochschild homology and its quantization in Section~\ref{sec:qHH}.
We show that it is a~preshadow on the~category of small linear categories
and prove its $K\noic$--theoretic invariance.

With Section \ref{sec:khov-hom} we move to the~topological part of the~paper.
Here we review link homology theories used in this paper: the~formal bracket
of Bar-Natan, the~Asaeda--Przytycki--Sikora homology for links in a~thickened
annulus and in a~thickened M\"obius band, and the~Chen--Khovanov homology for
tangles.

We quantize the~APS homology in Section \ref{sec:khov-hom}.
Theorems \ref{thm:qAKh} and \ref{thm:HH(CK)-vs-AKh}
are proven in Section \ref{sec:q-ann-hom}, whereas
\ref{sec:qAKh(torus-link)} presents a~computation of the~invariant
for torus links.

Section~\ref{sec:app-and-gen} contains further applications and generalizations:
an~extension to annular link cobordisms, the~action of oriented tangles on
cablings, quantization of the~APS homology for links in a~thickened M\"obius band,
and two generalizations of the~annular homology.

Finally, we gathered in Appendix basic definitions and results concerning
links and tangles, representations of $\Uqsl 2$, and elements of the~(bi)category theory.

\subsection*{Basic conventions and notation}

Throughout the~paper $\scalars$ is a~fixed unital commutative ring or field
such as $\Z$, $\Z_p$, or $\mathbb C$, and and linearity means linearity over
$\scalars$.
An~algebra means a~$\scalars$--algebra that is projective over $\scalars$
and likewise for modules and bimodules.

Graded means always $\Z$--graded. We denote by $\{d\}$ the~upwards degree shift,
i.e.\ $M\{d\}_i = M_{i-d}$ for a~graded module $M = \bigoplus_d M_d$.

A~differential in a~complex has homological degree $+1$ increases the~homological
degree. Thus we follow the~cohomological notation and put indices as superscripts.
The~only exception is the~Hochschild homology.
The~homological degree shift $[d]$ moves a~complex downwards,
i.e.\ $C[d]^i = C^{d+i}$.

Ordinary categories are typed with calligraphc letters ($\cat C$, $\cat{Vect}$, etc.),
whereas bold letters are reserved for bicategories ($\ccat C$, $\ccat{Rep}$, etc.).
We usually use small latin letters for objects ($x$, $y$, etc.) and for morphisms
($f$, $g$, etc.), whereas 2--morphisms are named by greek letters
($\alpha$, $\beta$, etc.), with the~exception of canonical isomorphisms
in bicategories, for which gothic letters are used ($\mathfrak a$, $\mathfrak l$,
$\mathfrak m$, etc.).
Capital letters are mostly reserved for functors.

\subsection*{Acknowledgements}
The~authors are grateful to Adrien Brochier, Matthew Hogancamp, Mikhail Khovanov, Slava Krushkal, Aaron Lauda, David Rose, and Paul Wedrich for stimulating discussions. During an~early stage of the~research Robert Lipshitz suggested to look on higher Hochschild homology of the~arc algebras and Ben Webster pointed a~connection between Hochschild homology and the~global dimension.
The~first two authors are supported by the~NCCR SwissMAP founded by the~Swiss National Science Foundation. The~third author was supported by the~NSF grant DMS-1111680.

\section{Generalized traces}
\label{sec:traces}

We start with a~discussion on twisted traces in general categories.
Section~\ref{sec:tr-trian} deals with traces on homotopy categories of complexes, whereas \ref{sec:tw-deformed} describes how twisted traces can be deformed
in case of graded categories. The~section ends with a~list of examples,
the~most important of which are the~twisted Hattori--Stallings trace
(\ref{sec:example-twisted-trace}) and the~annular closure of a~tangle
(\ref{sec:tr-annular-closure}). For a~brief list of constructions on categories
see Appendix~\ref{sec:categories}.

\subsection{Twisted traces}

Choose a~category $\cat C$ with an~endofunctor
$\endofun\colon \cat C\to \cat C$.
The~following definition extends the~usual notion of a~symmetric
trace.

%
% Definition: a trace function on an endocategory
%
\begin{definition}\label{def:trace-function}
	A~collection $t=\{t_x\}_{x\in\cat C}$ of functions
	$t_x\colon\cat C(x,\endofun x)\to S$ valued in a~set $S$ is
	a~\emph{$\endofun$--twisted trace} or a~\emph{trace on\/ $\cat C\noic$
	with a~monodromy\/ $\endofun$} if
	$t_y(\endofun f\circ g) = t_x(g\circ f)$
	for every pair of morphisms $\endofun x\from^g y\from^f x$.
\end{definition}

We shall sometimes refer to traces defined above as \emph{right traces},
whereas a~\emph{left trace} is defined dually as a~collection
of morphisms from $\cat C(\endofun x,x)$; the~trace condition
takes the~form $t_y(f\circ \endofun g) = t_x(g\circ f)$.
The~two definitions are clearly equivalent when $\endofun$ is invertible.
The~naming convention is motivated by traces in pivotal categories,
see Example~\ref{sec:pivotal-trace}.

\begin{lemma}\label{lem:sigma-inv-of-tr}
	A~trace $t$ with monodromy $\endofun$ is\/ $\endofun$--invariant, i.e.\ 
	$t_{\endofun x}(\endofun f) = t_x(f)$ for any $f\in\cat C(x, \endofun x)$.
\end{lemma}
\begin{proof}
	Take $y=\endofun x$ and $g=\id_{\endofun x}$ in the~trace relation.
\end{proof}

%
% The case of linear categories
%
We often write $t\colon(\cat C,\endofun)\to S$ for a~twisted trace
on $\cat C$, despite that it is not defined for all morphisms.
Further natural conditions on the~components of
twisted traces are imposed if $(\cat C,\endofun)$ has an~additional structure.
For instance, when both $\cat C$ and $\endofun$ are linear,
then $S$ is assumed to be a~module over the~ring of coefficients $\scalars$
and each component $t_x\colon \cat C(x,\endofun x) \to S$ to be
a~linear homomorphism.

%
% Universal trace is the coend
%
From the~point of view of Category Theory,
a~twisted trace is a~dinatural transformation from
$\cat C(\blank, \endofun(\blank))$ to the~constant bifunctor $\Delta_S$
that  assigns $S$ to any pair of objects and $\id_S$ to any pair of morphisms.
It follows that there exists a~\emph{universal\/ $\endofun$-twisted trace}
$\tr^\endofun\colon(\cat C,\endofun) \to \Tr(\cat{C},\endofun)$ when
$\cat C$ is small: the~coend of $\cat C(\blank, \endofun(\blank))$
\cite[Chapter IX.6]{Cats}. Explicitly,
\begin{gather}
	\Tr(\cat{C},\endofun)
		:= \quotient{\displaystyle\coprod_{x\in\Ob(\cat C)} \mskip-10mu \cat C(x,\endofun x)}
		            {\endofun f\circ g \sim g\circ f}
\intertext{or, when $\cat C$ is $\scalars$--linear,}
	\Tr(\cat{C},\endofun)
		:= \quotient{\displaystyle\bigoplus_{x\in\Ob(\cat C)} \mskip-10mu \cat C(x,\endofun x)}
		            {\Span_{\scalars}\{\endofun f\circ g - g\circ f\}}
\end{gather}
where $f$ and $g$ run through all pairs of morphisms
$\endofun x \from^g y \from^f x$ in $\cat C$.
Each component of $\tr^\endofun$ takes a~morphism $f\in\cat C(x,\endofun x)$
to its equivalence class $\tr^\endofun(f)\in\Tr(\cat C,\endofun)$, called
the~\emph{trace class of $f$}.
%
% The universal property of Tr
%
The~universality means that every twisted trace $t\colon (\cat C,\endofun)\to S$
factorizes uniquely through $\tr^\endofun$, i.e.\ there is a~unique
(linear) function $u\colon \Tr(\cat C,\endofun)\to S$
such that the~following triangle
\begin{equation}
\begin{gathered}
	\Rnode{C}{\cat C(x,\endofun x)}\hskip 6em \Rnode{S}{S}
	\\[1.5\baselineskip]
%	\hskip 1.5em
	\Rnode{Tr}{\Tr(\cat C,\endofun)}\rule{2em}{0pt}
	\begin{diagps}(0,0)(0,0)
		\arrow{->}[C`S;t_x]
		\arrow|b|{->}[C`Tr;\tr^\endofun_x]
		\arrow|b|{-->}[Tr`S;\exists ! u]
	\end{diagps}
\end{gathered}
\end{equation}
commutes for each object $x\in\cat C$.
Explicitly, $u(\tr^\endofun(f)) = t_x(f)$ for $f\in\cat C(x,\endofun x)$.

We shall often refer to $\Tr(\cat C,\endofun)$ as \emph{the} $\endofun$-twisted
trace of $\cat C$, keeping in mind the~components $\tr^\endofun_x$.
Moreover, following the~standard convention, we omit $\endofun$
from the~notation if it is the~identity functor, writing simply $\Tr(\cat C)$.
	
%
% Functoriality of Tr
%
The~trace $\Tr(\cat C,\endofun)$ can be seen as a~coend parametrized
by $\endofun$. Therefore, by the Parameter Theorem for coends
\cite[cp.\ Theorem IX.7.2]{Cats},
a~natural transformation of endofunctors induces a~map between traces.
In fact, such a~map exists in a~more general situation when $\cat C$
can vary as well. First let us create the~appropriate framework.

\begin{definition}
	An~\emph{endocategory} is a~pair $(\cat C,\endofun)$ consisting
	of a~category $\cat C$ and a~functor $\endofun\colon \cat C\to \cat C$.
	\emph{Functors} between endocategories
	$(\cat C,\endofun)$ and $(\cat C',\endofun')$ are pairs
	$(F,\omega)$ consisting of a~functor $F\colon \cat C\to \cat C'$
	and a~natural transformation $\omega\colon F\endofun \to \endofun'F$.
	It is finally understood that a~\emph{natural transformation}
	$\eta\colon (F,\omega) \to (F',\omega')$ between two such functors
	intertwines $\omega$ with $\omega'$.
\end{definition}

The~above describes $\cat{EndoCat}$, the~2-category of endocategories.
The~condition for a~natural transformation $\eta$ can be rephrased by
saying that the~square
\begin{equation}\label{diag:nat-trans-of-endo}
	\begin{diagps}(6em,12ex)
		\square<6em,10ex>[
			F\endofun x`F'\endofun x`\endofun'Fx`\endofun'F'x;
			\eta_{\endofun x}`\omega_x`\omega'_x`\endofun'\eta_x]
	\end{diagps}
\end{equation}
commutes for every $x\in\cat C$.
In the~following we argue that isomorphism classes of functors between
endocategories induce maps between traces.

Choose a~functor $(F,\omega)\colon (\cat C,\endofun)\to (\cat C',\endofun')$
and a~trace $t\colon (\cat C',\endofun')\to S$.
We check directly that $F^*t_x(f) := t_{Fx}(\omega_x\circ Ff)$ defines
a~$\endofun$--twisted trace on $\cat C$; we call it the~\emph{pullback
of\/ $t\noic$ along $F\noic$}.

\begin{proposition}\label{prop:pullbacks-of-traces}
	Choose a~trace\/ $t\colon (\cat C',\endofun')\to S$.
	Then $(GF)^*t = F^*(G^*t)$ for a~pair of functors\/
	$(F,\omega)\colon (\cat C,\endofun) \to (\cat C',\endofun')$ and\/
	$(G,\omega)\colon (\cat C',\endofun')\to (\cat C'',\endofun'')$.
	Morevoer, $F^*t = F'^*t$ if\/ $(F,\omega)$ and\/ $(F',\omega')$ are
	isomorphic functors.
\end{proposition}
\begin{proof}
	The~first equality follows from a~direct computation
	\begin{equation}
		F^*(G^*t)_x(f) = (G^*t)_{Fx}(\omega_x\circ Ff)
			= t_{GFx}(\omega'_{Fx}\circ G\omega_x \circ GFf)
			= (GF)^*t_x(f),
	\end{equation}
	while for the~second statement we use \eqref{diag:nat-trans-of-endo}
	to compute
	\begin{equation}
		\endofun'\eta_x \circ \omega_x \circ Ff \circ \eta_x^{-1}
		=
		\omega'_x\circ \eta_{\endofun x} \circ Ff \circ \eta_x^{-1}
		=
		\omega'_x \circ F'f
	\end{equation}
	where $\eta\colon(F,\omega)\to(F',\omega')$ is a~natural isomorphism.
\end{proof}

\begin{corollary}
	A~trace\/ $t\noic$ on $(\cat C,\endofun)$ can be twisted by
	a~natural transformation\/ $\eta\colon\endofun\to\endofun$
	to $(\eta^*t)_x := t_x(\eta_x\circ f)$.
\end{corollary}
\begin{proof}
	The~trace $\eta^*t$ is a~pullback of $t$ along
	the~functor of endocategories $(\Id,\eta)$.
\end{proof}

Let $\cat{EndoCat}^{\leqslant1}$ be the~truncation of $\cat{EndoCat}$
to a~category obtained by forgetting natural transformation and identifying
isomorphic functors. Proposition~\ref{prop:pullbacks-of-traces}, when
combined with the~universal property of $\Tr$, implies functoriality
of the~universal trace.

\begin{corollary}\label{cor:functoriality-of-Tr}
	There is a~functor\/ $\Tr\colon\cat{EndoCat}^{\leqslant1} \to \cat{Set}\noic$
	that sends an~endocategory\/ $(\cat C,\endofun)$ to its universal trace\/
	$\Tr(\cat C,\endofun)$ and a~functor\/
	$(F,\omega)\colon (\cat C,\endofun)\to (\cat C',\endofun')$ to the~map\/
	$\Tr(F,\omega)\colon \Tr(\cat C,\endofun) \to \Tr(\cat C', \endofun')$
	that takes\/ $\tr^\endofun(f)$ into\/ $\tr^{\endofun'}(\omega\circ Ff)$.
\end{corollary}

%
% Grothendieck group and Chern character.
%
Recall that the~split Grothendieck group $K_0(\cat C)$ of an~additive
category is generated by isomorphism classes $[x]$ of objects of $\cat C$
modulo the~relation $[x\oplus y] = [x] + [y]$.
There is a~natural additive map $h\colon K_0(\cat C)\to \Tr(\cat C)$
called the~\emph{Chern character} \cite{SL2asTrace} that takes $[x]$
to the~trace class of $\id_x$. This map is no longer well-defined for
twisted traces---identity morphisms does not represent trace classes
in $\Tr(\cat C, \endofun)$.
However, the map can be `twisted' in certain cases.

\begin{theorem}\label{thm:twisted-Chern-character}
	Let\/ $(\cat C,\endofun)$ be an~additive endocategory and
	$\eta\colon \Id\to \endofun$ a~natural transformation of functors.
	Then there is a~natural homomorphism of abelian groups
	\begin{equation}\label{eq:Chern-character}
		h^\eta\colon \Ksp(\cat C) \to \Tr(\cat C,\endofun),
			\qquad
		[x]_{\cong} \mapsto \tr^\endofun(\eta_x),
	\end{equation}
	called the~\emph{twisted Chern character}.
\end{theorem}
\begin{proof}
	By Corollary~\ref{cor:functoriality-of-Tr}
	there is a~map of traces
	$\eta_*\colon \Tr(\cat C)\to \Tr(\cat C,\endofun)$,
	which takes the~trace class $\tr(f)$ of $f\in\cat C(x,x)$
	into $\tr^\endofun(\eta_x\circ f)$.
	Compose it with the~usual Chern character to obtain $h^\eta$.
\end{proof}

%
% Additive closure and idempotent completion
%
Being defined for not necessarily additive categories is a~big
advantage of traces when compared to the~Grothendieck group.
In fact, $\Tr(\cat C,\endofun)$ can be often computed by
considering only a~certain subset of objects, which may
not be preserved by $\endofun$.
For a~precise statement, define the~\emph{trace of\/ $\cat C\noic$
restricted to $B\subset\Ob(\cat C)$} as the~quotient
\begin{equation}
	\Tr(\cat C|B,\endofun)
		:= \quotient{\displaystyle\bigoplus_{x\in B} \cat C(x,\endofun x)}
		            {\Span_{\scalars}\{\endofun f\circ g - g\circ f\}}
\end{equation}
where $f$ and $g$ run through all pairs of morphisms
$\endofun x \from^g y \from^f x$ in $\cat C$ with $x,y\in B$
(but $\endofun x$ may not belong to $B$).
There is a~natural homomorphism
\begin{equation}\label{eq:restr-to-entire}
	\Tr(\cat C|B,\endofun) \to \Tr(\cat C,\endofun).
\end{equation}
We say that $B$ \emph{generates additively} $\cat C$
if every object of $\cat C$ is a~direct summand of an~object
of the~form $b_1\oplus\dots\oplus b_r$ with $b_i\in B$.
For example, $\Ob(\cat C)$ generates additively
both the~additive closure $\catAdd{\cat C}$ and
the~idempotent completion $\Kar(\cat C)$ of $\cat C$
(see Section~\ref{sec:categories}).
Both categories come with natural endofunctors
induced by $\endofun$, for which we use the~same symbol.

\begin{proposition}\label{prop:tr-vs-additive-base}
	Suppose that\/ $\cat C\noic$ is generated additively by $B$.
	Then \eqref{eq:restr-to-entire} is an~isomorphism.
	In particular, the~inclusions of endocategories
	\begin{align}
		\label{fun:C-into-Cadd}
		(\cat C,\endofun)&\to \hbox to 8em{$\left(\cat C^{\oplus},\endofun\right)$,}
			x \mapsto (x) \\
		\label{fun:C-into-Kar(C)}
		(\cat C,\endofun)&\to \hbox to 8em{$\left(\Kar(\cat C),\endofun\right)$,}
			x \mapsto \id_x.
	\end{align}
	induce natural isomorphisms of the~universal traces.
\end{proposition}

This is a~generalization of \cite[Sections 3.4, 3.5]{SL2asTrace}
and a~special case of Proposition~\ref{prop:HH-vs-add-and-Kar}.
We leave the~proof as an~easy exercise.

\begin{remark}
	Although we discuss here only right traces, all the~results can be
	easily translated to left traces. In particular, an~endocategory
	$(\cat C,\endofun)$ admits the~universal \emph{left} trace
	$\Tr^\ell(\cat C,\endofun)$. Unless stated otherwise,
	the~universal trace always means the~right trace.
\end{remark}

\subsection{Triangulated traces}
\label{sec:tr-trian}

Let $\cat T$ be a~triangulated category and $\endofun$
its triangulated endofunctor,
meaning that there is a~natural isomorphism
$\alpha_x\colon \endofun (X[1]) \to^\cong (\endofun X)[1]$
and for each distinguished triangle
\begin{equation}
	X \to^f
	Y \to^g
	Z \to^h
	X[1]
\end{equation}
the~triangle
\begin{equation}
	\endofun X \to^{\endofun f}
	\endofun Y \to^{\endofun g}
	\endofun Z \to^{\alpha_X\circ\endofun h}
	\endofun X[1]
\end{equation}
is also distinguished.
We say that a~trace $t$ on $(\cat T, \endofun)$ is \emph{triangulated} if it is
additive with respect to homomorphisms of distinguished triangles, i.e.
\begin{equation}\label{eq:trian-additivity}
	t(b) = t(a) + f(c)
\end{equation}
for any commutative diagram
\begin{equation}\label{eq:dtri-twisted-endo}
	\begin{diagps}(14em,10ex)
		\node b0( 0em,0ex)[\endofun X]      \node t0( 0em,8ex)[X]
		\node b1( 4em,0ex)[\endofun Y]      \node t1( 4em,8ex)[Y]
		\node b2( 8em,0ex)[\endofun Z]      \node t2( 8em,8ex)[Z]
		\node b3(13em,0ex)[\endofun X\mathrlap{[1]}]
		\node t3(13em,8ex)[X\mathrlap{[1]}]
		\arrow{->}[b0`b1;\endofun f]   \arrow{->}[t0`t1;f]
		\arrow{->}[b1`b2;\endofun g]   \arrow{->}[t1`t2;g]
		\arrow{->}[b2`b3;\alpha_X\circ\endofun h]   \arrow{->}[t2`t3;h]
		\arrow{->}[t0`b0;a]
		\arrow{->}[t1`b1;b]
		\arrow{->}[t2`b2;c]
		\arrow{->}[t3`b3;{a[1]}]
	\end{diagps}
\end{equation}
with distinguished triangles as rows.

\begin{lemma}
	Assume that\/ $t\noic$ is a~triangulated trace on $(\cat T,\endofun)$.
	Then\/ $t(f[1]) = -t(f)$ for any $f\in \cat T(x,\endofun x)$.
\end{lemma}
\begin{proof}
	Apply the~additivity of $t$ to the~diagram \eqref{eq:dtri-twisted-endo}
	with $X \to 0 \to X[1] \to^{-\id} X[1]$ as the~top triangle
	and $(a,b,c) = (f,0,f[1])$.
\end{proof}

Likewise for usual traces,
triangulated traces are stable under pullbacks along triangulated functors,
and there is the~universal triangulated trace $\triTr(\cat T, \endofun)$,
through which all triangulated $\endofun$-twisted traces on $\cat T$ factor.

A~particular example of a~triangulated category is $\HCom^b(\cat C)$,
the~homotopy category of bounded formal complexes over $\cat C$,
i.e.\ sequences of objects and morphisms
\[
	\dots \to C^i \to^{d^i} C^{i+1} \to^{d^{i+1}} C^{i+2} \to \dots
\]
with $d^{i+1}\circ d^i = 0$ and $C^i=0$ except finitely many indices.
The~degree shift functor shifts the~complex leftwards and negates its differential,
whereas the~distinguished triangles are of the~form
\[
	C^\bullet
		\to^{f^\bullet}
	D^\bullet
		\to^{\mathit{in}^\bullet}
	\cone^\bullet(f)
		\to^{\mathit{pr}^\bullet}
	C[1]^\bullet,
\]
where $\cone^\bullet(f)$ stands for the~\emph{mapping cone} of $f\noic$,
see Section~\ref{sec:categories} for the~details.
It is implicitly assumed that $\cat C$ is additive (if not, take its additive
closure).

\begin{definition}
	A~\emph{Lefschetz trace} $\Lef^t$ on $\HCom^b(\cat C)$ induced by a~trace $t$
	on $\cat C$ is defined by the~formula
	\begin{equation}
		\Lef^t_x(f^\bullet) := \sum_i (-1)^i t_x(f^i).
	\end{equation}
\end{definition}

It is immediate that a~Lefschetz trace satisfies
the~$\endofun$-twisted trace relation
if it is induced by a~$\endofun$-twisted trace.
A~little longer computation shows that it is triangulated.
In fact, as explained below, all triangulated traces on $\HCom^b(\cat C)$ are Lefschetz traces.
In what follows we consider $\cat C$ as a~subcategory of $\HCom^b(\cat C)$ by
understanding an~object $X$ as the~complex $(0\to X\to 0)$
with $X$ at homological degree 0.

\begin{proposition}\label{prop:trian-trace-is-Lefschetz}
	Every triangulated\/ $\endofun$-twisted trace on\/ $\HCom^b(\cat C)$
	is a~Lefschetz trace, induced by its restriction to\/ $\cat C$.
	In particular, $\Tr(\cat C,\endofun) \cong \triTr(\HCom^b(\cat C),\endofun)$.
\end{proposition}
\begin{proof}
	Observe that a~bounded complex $C^\bullet$ is a~mapping cone of a~map between
	smaller complexes: $(0 \to C^{\min} \to 0)$, supported 	in homological degree
	$\min = \min\{i:C^i\neq 0\}$, and $C^{> \min} =
		(0\to C^{\min+1} \to C^{\min+2} \to \dots)$.
	The~first statement follows thus from the~induction on the~length of a~complex
	together with additivity of triangulated traces \eqref{eq:trian-additivity}.
	In particular, the~universal triangulated $\endofun$-twisted trace on
	$\HCom^b(\cat C)$ is the~Lefschetz trace induced by the~universal
	$\endofun$-twisted trace on $\cat C$, which implies the~second statement.
\end{proof}

\subsection{Deformation of twisted traces}
\label{sec:tw-deformed}

A~linear category $\cat C$ is \emph{graded} if it admits
an~autoequivalence $\{1\}\colon\cat C\to \cat C$, called the~\emph{degree shift}.
The~Grothendieck group $\Ksp(\cat C)$ of a~graded category has a~natural
structure of a~$\Zq$--module by setting $q\cdot [x] := [x\{1\}]$.
Motivated by Theorem~\ref{thm:twisted-Chern-character} we want to introduce
a~similar relation on the~trace. There are two ways to do it: by defining
the~action of $q$ explicitly or by deforming the~trace relation.
In what follows we explore both approaches.

%
% The trace of a graded endocategory is a Z[q,1/q]-module
%
A~\emph{graded endocategory} is a~quadruple $(\cat C,\endofun,\{1\},\alpha)$
consisting of a~category $\cat C$, endofunctors $\endofun$ and $\{1\}$,
where the~latter is an~equivalence, and a~natural isomorphism
$\alpha\colon \endofun(\blank)\{1\} \to^\cong \endofun(\blank\{1\})$.
The~equivalence $\{1\}$ is called the~degree shift functor.
We use it to equip the~universal trace $\Tr(\cat C,\endofun)$ with an~action
of $\Zq$ by setting
\begin{equation}
	q\cdot \tr^\endofun(f) := \tr^\endofun(\alpha_x\circ f\{1\})
\end{equation}
for any $f\in\cat C(x,\endofun x)$.
The~following is an~immediate consequence of this definition.

\begin{corollary}
	Let\/ $(\cat C,\endofun,\{1\},\alpha)$ be a~graded additive endocategory
	and\/ $\eta\colon \Id \to \endofun$ a~natural transformation satisfying\/
	$\eta_{x\{1\}} = \alpha_x\circ\eta_x\{1\}$.	
	Then the~twisted Chern character\/
	$h^\eta\colon \Ksp(\cat C) \to \Tr(\cat C,\endofun)$
	is\/ $\Zq$--linear.
\end{corollary}

%
% Pregraded categories
%
We shall now introduce the~action differently by deforming the~trace relation.
For that recall the~notion of a~\emph{pregraded category}:
it is a~category $\cat C$,
morphism sets of which are graded $\scalars$--modules
\begin{equation}
	\cat C(x,y) = \bigoplus_{d\in\Z} \cat C(x,y)^d,
\end{equation}
and the~degree is additive with respect to the~composition.
Naturally, we say that $(\cat C,\endofun)$ is a~\emph{pregraded category}
when $\endofun$ preserves the~grading,
i.e.\ $\endofun(\cat C(x,y)^d) \subset \cat C(\endofun x,\endofun y)^d$
for all objects $x$ and $y$.

%
% Pregraded -> graded
%
A~pregraded endocategory $(\cat C,\endofun)$ can be transformed into a~graded
one $(\cat C^{gr},\endofun^{gr})$ by introducing a~formal degree shift and taking
degree zero morphisms only.
Thence, objects of $\cat C^{gr}$ are symbols $x\{d\}$ with $x\in\Ob(\cat C)$
and $d\in\Z$, morphisms sets are $\cat C^{gr}(x\{a\},y\{b\}) := \cat C(x,y)^{a-b}$,
and the~degree shift functor increases the~number in brackets by one.
The~endofunctor $\endofun^{gr}$ takes $x\{a\}$ to $(\endofun x)\{a\}$.

%
% Graded -> pregraded
%
Dually, a~graded endocategory $(\cat C,\endofun)$ can be extended to a~pregraded
$(\cat C^{pre},\endofun^{pre})$ by forgetting the~degree shift functor, while
introducing morphisms of nonzero degrees by setting
$\cat C^{pre}(x,y)^d := \cat C(x\{d\},y)$ and extending the~composition
with the~formula
\begin{equation}
	\big(z \from^g y\{d'\}\big) \circ \big(y \from^f x\{d\}\big) :=
		\big(z \from^g y\{d'\} \from^{f\{d'\}} x\{d+d'\}\big)
\end{equation}
To define $\endofun^{pre}$ we use the~natural isomorphism $\alpha$:
\begin{align*}
	\endofun^{pre}x &:= \endofun x
		&&\text{for } x\in\Ob(\cat C^{pre}) = \Ob(\cat C)
\\
	\endofun^{pre}f &:= \endofun f \circ \alpha_x^d
		&&\text{for } f\in\cat C^{pre}(x,y)^d = \cat C(x\{d\},y)
\end{align*}
The~two construction described above are clearly inverse to each other.

%
% Quantum trace on a pregraded category
%
\begin{definition}
	Choose a~$\qscalars$--module $S$.
	A~\emph{quantum trace} on a~pregraded endocategory
	$(\cat C,\endofun)$ is a~linear trace
	$t\colon (\cat C,\endofun)\to S$ satisfying the~deformed trace condition
	$t_y(\endofun f\circ g) = q^{|f|} t_x(g\circ f)$ for
	every pair of homogeneous morphisms $\endofun x\from^g y\from^f x$,
	where we write $|f|$ for the~degree of $f$.
\end{definition}

Write $\cat C_q$ for the~category $\cat C$ with coefficients extended
to $\qscalars$.
There exists a~universal $\endofun$--twisted quantum trace
\begin{equation}
	\qTr(\cat{C},\endofun)
		:= \quotient{\displaystyle\bigoplus_{x\in\Ob(\cat C)} \mskip-10mu
				\cat C_q(x,\endofun x)}
		            {\Span_{\qscalars}\{\endofun f\circ g - q^{|f|}g\circ f \}}
\end{equation}
with the~obvious components
$\tr_q^\endofun\colon \cat C(x,\endofun x) \to \qTr(\cat C,\endofun)$,
through which every quantum trace on $(\cat C,\endofun)$ factorizes.
It is a~graded $\qscalars$--module,
because the~defining relation is homogeneous.
One can see it as a~one parameter deformation of $\Tr(\cat C,\endofun)$.
As before, we shall write $\qTr(\cat C)$ when $\endofun$ is the~identity functor.

%
% Quantum trace as a deformation
%
Quantum traces arise from a~deformation of the~endofunctor $\endofun$.
Indeed, a~quantum trace on $(\cat C, \endofun)$ is precisely a~trace
on $(\cat C_q,\endofun_q)$ as defined in the~previous section,
where $\endofun_q f := q^{-|f|}\endofun f$ for a~homogeneous morphism $f$.
In particular,
\begin{equation}
	\qTr(\cat C,\endofun) \cong \Tr(\cat C_q,\endofun_q)
\end{equation}
so that all the~properties of traces can be easily translated to the~quantum setting.

\begin{proposition}
	A~functor of pregraded endocategories\/
	$(F,\omega)\colon (\cat C,\endofun)\to (\cat C',\endofun')$
	induces a~map of universal quantum traces\/
	$\qTr(F,\omega)\colon \qTr(\cat C,\endofun) \to \qTr(\cat C', \endofun')$,
	which takes $\tr_q^\endofun(f)$ into $\tr_q^{\endofun'}(\omega\circ Ff)$.
\end{proposition}

Write as before $\qTr(\cat C|B,\endofun)$ for the~quantum trace restricted to
$B\subset\Ob(\cat C)$. We say that $\cat C$ is generated
\emph{additively by\/ $B\noic$ in the~graded sense} if it is
generated additively by $\tilde B := \{x\{d\}\ |\ x\in B, d\in\Z\}$.
The~following is a~special case of Proposition~\ref{prop:HH-vs-add-and-Kar}.

\begin{proposition}
	The~map $\qTr(\cat C|B,\endofun) \to \qTr(\cat C,\endofun)$
	is an~isomorphism if\/ $\cat C\noic$ is generated additively by\/ $B\noic$
	in a~graded sense. In particular,
	$\qTr(\cat C,\endofun) \cong \qTr(\cat C^{\oplus}, \endofun)$ and\/
	$\qTr(\cat C,\endofun) \cong \qTr(\Kar(\cat C),\endofun)$.
\end{proposition}

Because $\Sigma$ preserves the~degree of morphisms,
the~universal quantum traces if a~graded module.
Its degree $d$ component
\begin{equation}
	\qTr^d(\cat{C},\endofun)
		:= \quotient{\displaystyle
				\bigoplus_{x\in\Ob(\cat C)} \mskip-10mu\cat C_q(x,\endofun x)^d}
		       {\Span_{\qscalars}\{\endofun f\circ g - q^{|f|} g\circ f\}}
\end{equation}
is annihilated by $(1-q^d)$, which follows immediately from
Lemma~\ref{lem:sigma-inv-of-tr}.

%
% Comparing qTr(pregraded) with Tr(graded)
%
Traces on graded endocategories and quantum traces on pregraded endocategories
are closely related to each other.

\begin{proposition}
	Let\/ $(\cat C,\endofun)$ be a~pregraded extension of
	a~graded endocategory\/ $(\cat C^0, \endofun^0)$.
	Then $\Tr(\cat C^0, \endofun^0)$ and\/
	$\Tr^0_{q}(\cat C,\endofun)$
	are isomorphic $\qscalars$--modules.
\end{proposition}
\begin{proof}
	As the~notation suggest, $\cat C^0(x,\endofun^0 x)$
	are precisely the~degree zero morphisms in $\cat C(x,\endofun x)$.
	Therefore, there is an~obvious surjective map $i\colon \Tr(\cat C^0,\endofun^0)
	\to \qTr^0(\cat C,\endofun)$.
%	because morphisms of non-zero degree in $\cat C_0$ are not endomorphisms.
	We check first that it intertwines the~action of $\qscalars$.
	For that let $\iota_x\in\cat C(x,x\{1\})^1$ correspond to
	$\id_{x\{1\}} \in \cat C^0(x\{1\},x\{1\})$ and notice the~equality
	$f\{1\} = \iota_{\endofun x}\circ f\circ \iota_x^{-1}$
	for any morphism $f\colon x\to \endofun x$. Hence,
	\begin{equation*}
		\tr_q(\alpha_x\circ f\{1\})
			= \tr_q(\alpha_x\circ\iota_{\endofun x}\circ f\circ \iota_x^{-1})
			= q\cdot \tr_q(\endofun\iota_x^{-1} \circ \alpha_x \circ \iota_{\endofun x}\circ f )
			= q\cdot \tr_q(f)
	\end{equation*}
	in $\qTr^0(\cat C, \endofun)$. To show that $i$ is injective,
	and hence an~isomorphism, take any morphisms $f\in\cat C(x,y)^d$
	and $g\in\cat C(y,\endofun x)^{-d}$, and recall how $\endofun$ and
	the~composition of morphisms of nonzero degree is defined in $\cat C$:
	\begin{align*}
		g\circ f &:= g\circ f\{-d\} \\
		\endofun f\circ g &:= \endofun^0 f \circ \alpha^d \circ g\{d\}
	\end{align*}
	where the~left and right sides use the~composition in
	$\cat C$ and $\cat C^0$ respectively. Therefore,
	\begin{equation*}
			\tr( \endofun^0 f \circ \alpha^d \circ g\{d\} ) =
			\tr\left( \alpha^d \circ (g \circ f\{-d\})\{d\} \right) =
			q^d \tr( g\circ f\{-d\} )
	\end{equation*}
	in $\Tr(\cat C^0,\endofun^0)$, revealing that
	the~quantum trace relation in $\qTr^0(\cat C,\endofun)$
	is a~consequence of the~module structure on $\Tr(\cat C^0,\endofun^0)$.
	In particular, the~kernel of $i$ is trivial.
\end{proof}

\begin{remark}
	Triangulated traces also admit deformations as long as the~triangulated
	category, on which they are defined, has an~additional grading compatible
	with the~triangulated structure. An~example is the~homotopy category
	$\HCom(\cat C)$ of a~pregraded category $\cat C$, objects of which are formal
	complexes with a~graded differential, and formal chain maps have graded
	components.
\end{remark}

\subsection{Examples}
\label{sec:tr-examples}

In what follows we provide a~number of (twisted) traces
that appear in algebraic or topological contexts.
Of particular interest to us are the~Hattori--Stallings trace
(together with its twisted version) and annular closures
of tangles.

\subsubsection{The trace of a~linear map}

Let $\cat{Vect}$ be the~category of finite dimensional vector
spaces over a~field $\Fld$. It is a~classical result from Linear
Algebra that there is a~unique collection of maps
$\tr_V\colon \End(V)\to\Fld$ satisfying the~trace relation
and such that $\tr_V(\id) = \dim V$. Explicitly,
\begin{equation}
	\tr(f) := \sum_i f_{ii} \in \Fld
\end{equation}
where $(f_{ij})$ is a~matrix representation of $f$ in some basis
of $V$. This is the~universal trace: $\cat{Vect}$ is
generated additively by $B = \{\Fld\}$, so that
$\Tr(\cat{Vect}) \cong \Tr(\cat{Vect}|\{\Fld\}) = \Fld$.

\subsubsection{The~graded trace}

In case of $\Z$--graded vector spaces the~above trace admits
a~deformation
\begin{equation}
	\tr_q(f) := \sum_n q^n \tr(f|_{V_n}) \in \qFld,
\end{equation}
where $V_n\subset V$ is the~subspace generated by homogeneous
vectors of degree $n$. It satisfies the~quantum trace relation:
\begin{equation}
	q^{|f|}\tr_q(g\circ f)
	= \sum_{n\in\Z} q^{n+|f|}\tr(g\circ f|_{V_n})
	= \sum_{n\in\Z} q^{n+|f|}\tr(f\circ g|_{W_{n+|f|}})
	= \tr_q (f\circ g)
\end{equation}
for any pair of homogeneous maps $V\from ^g W\from^f V$.
Again, $\qTr(\cat{gVect}) \cong \qFld$ as $B = \{\Fld\}$
generates the~category additively in a~graded sense.

\subsubsection{The~Lefschetz number}

Recall that the~\emph{Lefschetz number} $\Lef(f^\bullet)$ of an~endomorphism
$f^\bullet$ of a~chain complex of vector spaces is the~alternating sum
\begin{equation}\label{eq:Lef-number}
	\Lef(f^\bullet) := \sum_i (-1)^i \tr(f^i).
\end{equation}
It is straightforward to check that $\Lef$ satisfies the~trace
relation. Even more, Lefschetz numbers of homotopic maps are equal,
so that $\Lef$ descends to a~trace on $\HCom^b(\Fld)$, the~homotopy
category of complexes.
It is the~universal triangulated trace as stated in
Proposition~\ref{prop:trian-trace-is-Lefschetz}.
When complexes of graded vector spaces are considered,
one constructs the~\emph{graded Lefschetz number} by replacing
$\tr$ with $\tr_q$ in \eqref{eq:Lef-number}.

\subsubsection{The~spectrum of a~linear map}

The~spectrum of an~endomorphism $f\colon V\to V$ of a~vector space
is the~set of nonzero eigenvalues. It is an~example of a~non-linear
trace on $\cat{Vect}$. Indeed,
\[
	(f\circ g)(x) = \lambda x
\quad\Rightarrow\quad
	(g\circ f)(g(x)) = \lambda g(x)
\]
and $g(x)\neq 0$ if $\lambda\neq 0$. Scaling a~homomorphism scales elements
of its spectrum accordingly, but the~spectum of $f+g$ cannot be expressed
in general in terms of eigenvalues of $f$ and $g$.

\subsubsection{Hattori--Stallings trace}
\label{sec:ex-trace-on-Mod}

Let $\Rep(A)$ be the~category of finitely generated projective right
modules over a~ring $A$.
The~\emph{Hattori--Stallings trace} \cite{Hattori,Stallings}
is a~homomorphism $\HS_P\colon\End_A(P) \to A/[A,A]$ defined as the~composition
\begin{equation}\label{eq:twisted-HS}
	\End_A(P) \to^\cong P\utimes{A} P^* \to^{\epsilon_P} A/[A,A],
\end{equation}
where $P^* = \Hom_A(P,A)$ is the~left $A$--module dual to $P$
and $A/[A,A]$ is the~quotient of $\scalars$--modules with $[A,A]$ generated
by differences $ab-ba$ for all $a,b\in A$.
The~left isomorphism takes an~endomorphism $f$ to $(f\otimes\id)(\coev(1))$
and the~right map is given as $\epsilon_P(x\otimes\alpha) := \alpha(x)$.%
\footnote{
	Notice the~difference between $\epsilon$ and the~evaluation map
	$\ev\colon P^*\otimes P\to A$, the~latter being
	a~homomorphism of $(A,A)$--bimodules.
}
In particular,
$\epsilon_P\circ(g\otimes\id) = \epsilon_Q\circ(\id\otimes g^*)$
for any $g\in\Hom_A(Q, P)$, which guarantees that the~trace relation
is satisfied. When $P$ is free of finite rank,
$\HS_P(f) = \sum_i f_{ii}$ for $f\in\End_A(P)$ represented by the~matrix
$(f_{ij})$.

Projective modules are direct summands of free modules. Hence, $\Rep(A)$
is generated additively by $A$, seen as a~right $A$--module,
and the~universal trace
\begin{equation}
	\Tr(\Rep(A)) \cong \Tr(\Rep{A}|\{A\}) = A/[A,A]
\end{equation}
coincides with the~Hattori--Stallings trace.

\subsubsection{Twisted traces}
\label{sec:example-twisted-trace}

Choose an~$(A,A)$--bimodule $M$ that is finitely generated and projective
as a~right module, so that $(\blank)\otimes_A M$ is an~endofunctor of $\Rep(A)$.
We construct a~twisted Hattori--Stallings trace $\tr^M$ as the~composition
\begin{equation}
	\Hom_A(P,P\utimes{A} M) \to^{c_P} P\utimes A M\utimes A P^*
		\to^{\epsilon_P} \quotient{M}{[A,M]} =: \coinv(M),
\end{equation}
where the~$\scalars$--module $[A,M]$ is generated by differences $am-ma$
with $a\in A$ and $m\in M$,
$c_P(f) = (f\otimes\id)(\coev(1))$, and
$\epsilon_P(x\otimes m\otimes\alpha) := \alpha(x)m$.
The~quotient $\coinv(M)$ is called the~\emph{space of coinvariants in $M\noic$}.
Again, the~trace condition is satisfied because
$\epsilon_P\circ(g\otimes\id\otimes\id) = \epsilon_Q\circ(\id\otimes\id\otimes g^*)$
for any $g\in\Hom_A(Q, P)$.

The~twisted Hattori--Stallings trace is the~universal trace
on $(\Rep(A),M)$. Indeed,
\begin{equation}
	\Tr(\Rep(A),M) \cong \Tr(\Rep(A)|\{A\}, M) \cong \coinv(M)
\end{equation}
by Proposition~\ref{prop:tr-vs-additive-base},
where the~right isomorphism arises from the~identifications
of $\End_A(A)$ and $\Hom_A(A,M)$ with $A$ and $M$ respectively,
each evaluating a~function at $1\in A$.

Given an~algebra endomorphism $\varphi\in\AlgEnd(A)$ one constructs a~twisted
trace on $\Rep(A)$ as follows. For a~right module $V$ define its
\emph{$\varphi$--twisting} $V_{\overline\varphi} := V \utimes{A} {}_\varphi A$,
where ${}_\varphi A$ is obtained from $A$ by redefining the~action as
$a\cdot x := \varphi(a)x$.  In case $\varphi$ is invertible,
there is an~isomorphism $V_{\overline\varphi} \cong V_{\varphi^{-1}}$
that takes $v\otimes a$ to $v\varphi^{-1}(a)$,
where $V_{\varphi^{-1}}$ is obtained from $V$ simply by twisting the~action of $A$.
If so, a~homomorphism $f\colon W\to V_{\overline\varphi}$
is precisely a~$\scalars$--linear map $f\colon W\to V$ satisfying
$f(w\varphi(a)) = f(w)a$. This leads to a~twisted Hattori--Stallings trace
\begin{equation}
	\tr^\varphi \colon \Rep(A) \to \quotient{A}{[A,A]_\varphi},
\end{equation}
where $[A,A]_\varphi = \Span_{\scalars}\{\varphi(a)b-ba\ |\ a,b\in A\}$.
Notice that $\tr^\varphi(f) = \sum_i f_{ii}$ for any endomorphism
$f$ of a~free $A$--module of finite rank.
This is the~same formula as in the~untwisted case, but
computed in a~different quotient of $A$.

Choose an~invertible scalar $q\in\scalars$.
If $A = \bigoplus_{d\in\Z} A_d$ is a~graded algebra and $M$ a~graded bimodule,
then the~twisted Hattori--Stallings trace can be deformed to take values in
the~\emph{quantum space of coinvariants} of $M\noic$
\begin{equation}\label{def:coinv-q}
	\coinv_q(M) := M/[A,M]_q,
\end{equation}	
where $[A,M]_q = \Span_{\scalars}\{am- q^dma\ |\ a\in A_d, m\in M \}$.
The~details are left to the~reader.

\subsubsection{Traces in a~pivotal category}
\label{sec:pivotal-trace}

A~\emph{pivotal category} is a~monoidal category $\cat C$
with left duals and a~natural isomorphism $\delta_x\colon x\to \lddual x$
for any $x\in \cat C$. Particular examples are
\begin{itemize}

	\item the~category $\cat{Vect}$ of finite dimensional
	vector spaces together with $\delta = \delta^{\mathit{vect}}$
	the~standard isomorphism, and
	
	\item the~representation category $\Rep(H)$ of a~Hopf algebra $H$
	with a~\emph{pivot}, a~group-like element $u\in H$ satisfying $S^2(h) = uhu^{-1}$;
	the~isomorphism $\delta_V$ is given by the~formula
	$\delta_V(x) = u\cdot\delta^{\mathit{vect}}(x)$.
	
\end{itemize}
It is common to represent $\cat C$ diagrammatically: products of objects
as sequences of points labeled by the~factors, and morphisms as graphs
with coupons, read from bottom to top and understood as follows:
\begin{center}
	\def\point(#1)#2{\rput[B](#1){$\scriptstyle#2$}}%
	\psset{linewidth=1pt,unit=1ex}%
	\begin{pspicture}(-5ex,-2ex)(18ex,10ex)%
		\psarc(0,0){4}{0}{180}%
		\point(-4,-2){x}%
		\point( 4,-2){\ldual*x}%
		\rput(8,4){$\equiv$}%
		\begin{diagps}%
			\node dom(14, -2)[\mathllap x\otimes \mathrlap{\ldual x}]%
			\node cod(14,  8)[I]%
			\arrow|*c|{->}[dom`cod;\mathrm{ev}]%
		\end{diagps}%
	\end{pspicture}
	\qquad
	\begin{pspicture}(-5ex,-2ex)(18ex,10ex)
		\psarc(0,8){4}{180}{360}%
		\point(-4, 8.5){\ldual*x}%
		\point( 4, 8.5){x}%
		\rput(8,4){$\equiv$}%
		\begin{diagps}%
			\node dom(14, -2)[I]%
			\node cod(14,  8.5)[\mathllap{\ldual x}\otimes\mathrlap{x}]%
			\arrow|*c|{->}[dom`cod;\mathrm{coev}]%
		\end{diagps}%
	\end{pspicture}%
	\qquad
	\begin{pspicture}(-2ex,-6ex)(10ex,6ex)
		\psline(0,-4)(0,4)
		\point(0,-6){x}%
		\point(0,4.5){y}%
		\coupon(0,0){\ f\ }%
		\rput(4,0){$\equiv$}%
		\begin{diagps}
			\node dom(8,-6)[x]
			\node cod(8, 4.5)[y]
			\arrow|*c|{->}[dom`cod;f]
		\end{diagps}
	\end{pspicture}
\end{center}
The~relations between evaluation and coevaluation translates into straightening zig-zags
\begin{equation*}
	\psset{linewidth=1pt,unit=4ex}%
	\begin{centerpict}(0,0)(13ex,12ex)
		\pscustom{%
			\moveto(0.5,0)
			\lineto(0.5,1.5)
			\psarcn(1,1.5){0.5}{180}{0}
			\psarc(2,1.5){0.5}{180}{360}
			\lineto(2.5,3)
			\stroke}%
	\end{centerpict}
\quad=\quad
	\begin{centerpict}(-2ex,0)(2ex,12ex)
		\psline(0,0)(0,3)
	\end{centerpict}
\quad=\quad
	\begin{centerpict}(0,0)(13ex,12ex)
		\pscustom{%
			\moveto(2.5,0)
			\lineto(2.5,1.5)
			\psarc(2,1.5){0.5}{0}{180}
			\psarcn(1,1.5){0.5}{360}{180}
			\lineto(0.5,3)
			\stroke}%
	\end{centerpict}
\end{equation*}
so that isotopic diagrams represent the~same morphism.
The~\emph{left} and \emph{right pivotal traces} $t_\ell(f)$
and $t_r(f)$ of an~endomorphism $f\in\cat C(x,x)$
are defined as the~compositions
\begin{equation}
	\psset{linewidth=1pt,unit=1ex}
	t_\ell(f) := \begin{centerpict}(-3ex,8ex)(8ex,-8ex)
		\pscustom{%
			\moveto(0,-4)
			\psarc(3,-4){3}{180}{360}
			\lineto(6,4)
			\psarc(3, 4){3}{0}{180}
			\closepath
			\stroke
		}
		\coupon(6,-2.5){\,f\,}
		\coupon(6, 2.5){\,\delta_x\,}
		\rput[r](-0.25, 1.0){$\scriptstyle \ldual x$}
		\rput[l]( 6.00, 5.5){$\scriptstyle \lddual x$}
		\rput[l]( 6.50, 0.0){$\scriptstyle x$}
		\rput[l]( 6.30,-5.5){$\scriptstyle x$}
	\end{centerpict}
	\hskip 0.15\textwidth
	t_r(f) := \begin{centerpict}(-3ex,8ex)(9ex,-8ex)
		\pscustom{%
			\moveto(0,-4)
			\psarc(3,-4){3}{180}{360}
			\lineto(6,4)
			\psarc(3, 4){3}{0}{180}
			\closepath
			\stroke
		}
		\coupon(0, 2.5){\,f\,}
		\coupon(0,-2.5){\,\delta_x^{-1}\,}
		\rput[l]( 6.30, 1.0){$\scriptstyle \ldual x$}
		\rput[r](-0.10,-5.5){$\scriptstyle \lddual x$}
		\rput[r](-0.50, 0.0){$\scriptstyle x$}
		\rput[r](-0.10, 5.5){$\scriptstyle x$}
	\end{centerpict}
\end{equation}
They satisfy left and right version of the~trace relation
(see the~discussion below Definition~\ref{def:trace-function})
\begin{equation}
	t_\ell( f \circ g^\delta) = t_\ell(g\circ f)
\hskip 0.15\textwidth
	t_r(f^{\delta} \circ g) = t_r(g\circ f)
\end{equation}
with $(\blank)^{\delta} := \delta\circ \lddual (\blank) \circ\delta^{-1}$.
The~two traces coincide when $\cat C$ is \emph{spherical},
but in general they are different. See \cite{SymInt} for more details.
%Actually, the~left trace is dual to the~right in the~sense of
%Remark~\ref{rmk:dual-trace}.

When $\cat C = \Rep(H)$ is the~representation category
of a~pivotal Hopf algebra, one can compute the~universal right trace
\begin{equation}
	\Tr(\Rep(H), \delta) \cong \quotient{H}{\Span_{\scalars}\{S^2(y)x - xy\ |\ x,y\in H\}}
\end{equation}
using the~same methods as before. The~resemblence with the~twisted trace
from Section~\ref{sec:example-twisted-trace} is not a~coincidence: there
is a~natural isomorphism of functors on $\Rep(H)$
\begin{equation}
	\eta\colon (\blank)^{\delta} \to (\blank)_{\overline u},
\qquad
	\eta_V(x) = xu^{-1},
\end{equation}
where $u$ acts on $H$ by conjugation.
Indeed, $\eta_V(xh) = (xu^{-1})(uhu^{-1}) = \eta_V(x) \cdot h$. Analogous statments hold
for the~left trace.

\subsubsection{The~annular closure of a~tangle}
\label{sec:tr-annular-closure}

\begin{figure}[ht]
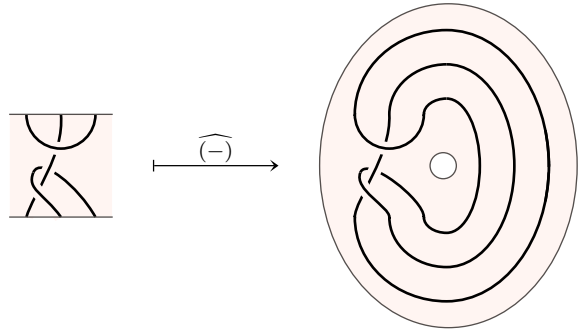

	\begin{displaymath}
	\begin{centerpict}(-0.5,0)(2.5,3)
		\psframe[style=mysurface,linestyle=none](-0.5,0)(2.5,3)
		\psset{style=diagarc}%
		\psbezier(0.2,1.3)(0.4,1.8)(1.5,0.9)(2,0)%
		\begin{psset}{border=2pt,bordercolor=sfceFillColor}
			\psbezier(1,3)(1.1,1.8)(0.2,0.7)(0,0)%
			\psbezier(1,0)(0.5,0.7)(0,0.8)(0.2,1.3)%
			\psarc(1,3){1}{180}{360}%
		\end{psset}
		\psline[style=diaghline](-0.5,0)(2.5,0)
		\psline[style=diaghline](-0.5,3)(2.5,3)
	\end{centerpict}
	\quad\mapsto<4em>^{\widehat{(\blank)}}\quad
	\begin{centerpict}(-1,-2.7)(6.5,5.7)
	\psellipse[style=mysurface](2.75,1.5)(3.75,4.75)
	\pscircle[style=mysurface,fillcolor=white](2.6,1.5){0.4}
	\psset{linewidth=1pt}%
	\pscustom{%
		\moveto(0.2,1.3)%
		\curveto(0.4,1.8)(2.04,0.5)(2.04,0)%
		\psellipticarc(2.7,0.0)(0.7,0.5){180}{270}%
		\psellipticarc(2.7,1.5)(1.0,2.0){-90}{ 90}}%
%		\stroke[linecolor=black,linewidth=1pt]}%
	\pscustom{%
		\moveto(2.7,4.5)%
		\psellipticarc(2.7,3.0)(1.7,1.5){ 90}{180}%
		\curveto(1.1,1.8)(0.04,0.5)(0.04,0)%
		\psellipticarc(2.7,0.0)(2.7,2.5){180}{270}%
		\psellipticarc(2.7,1.5)(3.0,4.0){-90}{ 90}%
		\stroke[linecolor=sfceFillColor,linewidth=5pt]}%
%		\stroke[linecolor=black,linewidth=1pt]}%
	\pscustom{%
		\moveto(0.2,1.3)
		\curveto(0,0.8)(1.04,0.3)(1.04,0)%
		\psellipticarc(2.7,0.0)(1.7,1.5){180}{270}%
		\psellipticarc(2.7,1.5)(2.0,3.0){-90}{ 90}%
		\stroke[linecolor=sfceFillColor,linewidth=5pt]}%
%		\stroke[linecolor=black,linewidth=1pt]}%
	\pscustom{%
		\moveto(5.7,1.5)
		\psellipticarc(2.7,1.5)(3.0,4.0){0}{90}%
		\psellipticarc(2.7,3.0)(2.7,2.5){ 90}{180}%
		\psarc(1.04,3){1}{180}{360}%
		\psellipticarcn(2.7,3.0)(0.7,0.5){180}{ 90}
		\stroke[linecolor=sfceFillColor,linewidth=5pt]}%
%		\stroke[linecolor=black,linewidth=1pt]}%
	\psellipticarc[linewidth=1pt](2.7,1.5)(3.0,4.0){-20}{20}%
	\end{centerpict}
	\end{displaymath}
	\caption{The~annular closure of a~tangle.}
	\label{fig:annular-closure}
\end{figure}

Let $q\colon \R^2\times I \to \R^2\times\S$ be the~quotient
map that identifies $(x,y,0)$ with $(x,y,1)$ for any $(x,y)\in\R^2$.
If the~tangle $T\subset\R^2\times I$ has the~same endpoints on
both boundary planes, then $q(T)$ is a~link in $\R^2\times\S$.
It is called the~\emph{annular closure of\/ $T\noic$}, because it can
visualized at the~level of tangle diagrams by connecting
the~endpoints inside an~annulus, see Figure~\ref{fig:annular-closure}.
We consider it as a~link $\widehat T$ in $\R^3$ by using
the~standard embedding of $\R^2\times\S \cong D^2\times\S$ into
the~3-space.

The~annular closure described above is a~trace function
on $\Tan(\R^2)$ valued in $\sLinks(\R^3)$. Indeed, the~links
$\widehat{T'T}$ and $\widehat{TT'}$ are isotopic by rotation.

\subsubsection{Satellite links as traces of tangles}

Embeddings of $\R^2\times\S$ into $\R^3$ are parametrized by
framed knots---the~tubular neighborhood of a~knot
is homeomorphic to $\R^2\times\S$, but it is the~framing
that fixes the~isotopy class of the~homeomorphism.
Let $h\colon\R^2\times\S \to \R^3$ be the~embedding
corresponding to a~framed knot $K$. The~image $h(L)$
of a~link $L\subset \R^2\times\S$ is called a~\emph{satellite link}
with \emph{companion} $K$.

\begin{figure}[ht]
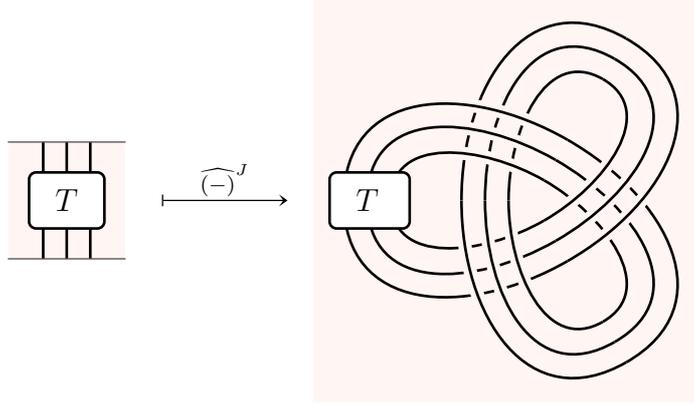

	\begin{displaymath}
	\begin{centerpict}(-1em,-2em)(2.5em,2em)
		\psset{linewidth=1pt,unit=1.5em}
		\psframe[style=mysurface,linestyle=none](-0.75,-1.25)(1.75,1.25)
		\psline(0.0,-1.25)(0.0,1.25)
		\psline(0.5,-1.25)(0.5,1.25)
		\psline(1.0,-1.25)(1.0,1.25)
		\pspolygon[fillstyle=solid,linearc=0.15](-0.3,-0.6)(-0.3,0.6)(1.3,0.6)(1.3,-0.6)
		\rput(0.5,0){$T$}
		\psline[style=diaghline](-0.75,-1.25)(1.75,-1.25)
		\psline[style=diaghline](-0.75, 1.25)(1.75, 1.25)
	\end{centerpict}	
	\quad\mapsto<4em>^{\widehat{(\blank)}^J}\quad
	\begin{centerpict}(-6.5em,-6em)(5em,6em)
		\psframe[style=mysurface,linestyle=none](-7em,-6.5em)(5.5em,6.5em)
		\psset{linewidth=1pt,unit=1.5em}
		\def\underarc{%
			\psbezier(0.25,-0.433)(-1.049,-1.183)(-3.0,-1.4)(-3.0,0)
			\psbezier(0.50,-0.866)(-1.232,-1.866)(-3.5,-2.0)(-3.5,0)
			\psbezier(0.75,-1.299)(-1.415,-2.549)(-4.0,-2.4)(-4.0,0)
		}%
		\def\overarc{%
			\psbezier(-3.0,0)(-3.0, 1.4)(-1.049, 1.183)(0.25, 0.433)
			\psbezier(-3.5,0)(-3.5, 2.0)(-1.232, 1.866)(0.50, 0.866)
			\psbezier(-4.0,0)(-4.0, 2.4)(-1.415, 2.549)(0.75, 1.299)
		}%

		\underarc
		\rput{120}(0,0)\underarc
		\rput{240}(0,0)\underarc

		\begin{psset}{border=2pt,bordercolor=sfceFillColor}
			\overarc
			\rput{120}(0,0)\overarc
			\rput{240}(0,0)\overarc
		\end{psset}
	
		\pspolygon[fillstyle=solid,linearc=0.15](-4.3,-0.6)(-4.3,0.6)(-2.6,0.6)(-2.6,-0.6)
		\rput(-3.5,0){$T$}
	\end{centerpict}
	\end{displaymath}
	\caption{The~closure of a~$(3,3)$--tangle along a~framed trefoil knot.}
	\label{fig:satellite}
\end{figure}

Write $K$ as an~annular closure of a~framed $(1,1)$--tangle $J$.
The~tangle $J$ describes a~tangled embedding of $\R^2\times I$
into itself.
Let $J_B$ be the~image under this embedding of the~trivial tangle
$B\times I$, where $B$ is any set of points on $\R^2$.
Then $J_{B'}T$ and $TJ_B$ are isotopic for any tangle $T$ with
$B$ and $B'$ as its bottom and top boundary, where the~isotopy
slides $T$ along $J$. In other words,
$J_\bullet$ is a~natural transformation of the~identity
functor on $\Tan(\R^2)$.
Using Proposition~\ref{prop:pullbacks-of-traces} we construct
a~new trace function $\widehat{(\blank)}^J$ on $\Tan(\R^2)$.
It takes a~tangle $T$ with same top and bottom boundary
to $\widehat{T}^J := h(\widehat T)$, a~satellite knot 
with compagnion $K$ (see Figure~\ref{fig:satellite}).

\subsubsection{Links in surface bundles}

Choose a~surface $\sfce$ with a~diffeomorphism $\phi\in\mathit{Diff}(\sfce)$.
The~mapping torus
\begin{equation}
	M := \ \quotient{\sfce\times I}{(p,1) \sim (\phi(p),0)}
\end{equation}
is called the~\emph{surface bundle} with \emph{fiber} $\sfce$
and \emph{monodromy} $\phi$.
The~annular closure of tangles can be generalized to a~procedure
producing links in $M$.

Given a~tangle $T\subset \sfce\times I$ with input $\phi(B)$ and output
$B$ we construct its closure $\widehat T := \pi(T)$
in $M$, where $\pi\colon \sfce\times I\to M$ is the~quotient map.
This operation satisfies the~left trace relation with respect to
the~endofunctor $\phi_*$ on $\Tan(\sfce)$ that takes a~tangle
$T\subset \sfce\times I$ into $(\phi\times\id)(T)$.
In another words, closing tangles describes a~twisted trace
$\widehat{(\blank)}\colon (\Tan(\sfce), \phi_*) \to \sLinks(M)$.

\begin{theorem}\label{thm:tr-of-tan}
	The~closure $\widehat{(\blank)}$ is the~universal left trace
	on $(\Tan(F),\phi_*)$.
\end{theorem}
\begin{proof}
	$\sLinks(M)$ is a~quotient set of $\Tr(\Tan(\sfce),\phi_*)$,
	because each link in $M$ is a~closure of a~certain tangle.
	Suppose that closures $\widehat T$ and $\widehat T'$ are isotopic
	in $M$. The~isotopy can be expressed as a~sequence of
	isotopies supported in small 3--balls, so that we may assume that
	it fixes some fiber $\sfce'\subset M$.
	Then the~cuts $\widehat T$ and $\widehat T'$ along $\sfce'$ are
	isotopic in $\sfce\times I$ and the~trace relation implies that
	the~images of $T$ and $T'$ in $\Tr^\ell(\Tan(\sfce),\phi_*)$ coincide.
\end{proof}

\subsubsection{Kauffman Bracket Skein Module}
\label{skein_module}

The~Kauffman Bracket Skein Module $\KBSM(M)$ of an~oriented 3--manifold $M$
is a~$\Z[A^{\pm1}]$--module generated by isotopy classes of framed tangles
in $M$ (we require the~isotopies to fix $\partial M$ if nonempty) modulo
the~local relations
\begin{align}
	\vcenter{\hbox{\bigNWSECr\,}}
		&= A \vcenter{\hbox{\,\bigHorRes\,}}
		+ A^{-1} \vcenter{\hbox{\,\bigVertRes\,}}\\
	\vcenter{\hbox{\bigCircle\,}}
		&= -A^2 - A^{-2}
\end{align}
where the~diagrams represent pieces of tangles picked by some
small ball, outside of which the~tangles coincide.
In particular, the~circle in the~second relation bounds a~disk.
When $\sfce$ is a~surface, then $\KBSM(\sfce\times I)$ is a~category
and any $\phi\in\mathit{Diff}(F)$ induces an~endofunctor $\phi_*$
as in the~previous example.
The~proof of Theorem~\ref{thm:tr-of-tan} can be easily modified to
show that the~universal left trace $\Tr^\ell(\KBSM(\sfce\times I),\phi_*)$
computes $\KBSM(M)$, the~skein module of the~surface bundle $M$ with fiber
$\sfce$ and monodromy $\phi$.

\section{Traces of bicategories}
\label{sec:categ-traces}

In this section we categorify twisted traces. We begin with the~definition
of a~preshadow and the~construction of the~(twisted) horizontal trace of
an~endobicategory. The~universal property and functoriality of the~latter
are shown in Sections \ref{sec:hor-trace} and \ref{sec:functoriality-of-hTr}
respectively, whereas Sections~\ref{sec:vTr-vs-hTr} and \ref{sec:hTr-categ-Tr}
discuss two ways to obtain traces from shadows: by restricting to the~vertical
trace or by decategorifying. These are followed by a~material on Lefschetz
preshadows and a~quantization of this framework. 
Examples are listed at the~end of the~section.

\subsection{Categorified traces}

Choose a~bicategory $\ccat C$ together with a~strong bifunctor
$\endofun\colon \ccat C\to \ccat C$.
These come together with natural 2-isomorphisms, the~\emph{coherence isomorphisms}
\begin{align*}
	\mathfrak a&\colon h\circ(g\circ f) \to/=>/^\cong (h\circ g)\circ f, &
	\mathfrak m&\colon \endofun g\circ\endofun f \to/=>/^\cong \endofun(g\circ f),
\\
	\mathfrak l&\colon \id_y\circ f \to/=>/^\cong f, &
	\mathfrak i&\colon \id_{\endofun x} \to/=>/^\cong \endofun(\id_x),
\\
	\mathfrak r&\colon f\circ\id_x \to/=>/^\cong f,
\end{align*}
which are often omitted for clarity, see Appendix~\ref{sec:bicats}.

%where $x,y$ are any objects and $f,g,h$ composable 1--morphisms in $\ccat C$.
%These isomorphisms satisfy certain conditions, which allows us to safely forget
%about them. For instance, whenever two compositions of 1-morphisms differ only
%by the~position of brackets or existence of identity 1-morphisms, there is
%a~unique 2-isomorphism between them obtained from $\mathfrak a$, $\mathfrak l$,
%and $\mathfrak r$ by composing them horizontally or vertically.
%Because of that we shall often omit them---especially in the~next sections---%
%as well as the~brackets to simplify the~exposition.

%The~following generalizes the~notion of a~(symmetric) shadow on a~bicategory.

%
%  Definition of a preshadow
%
\begin{definition}
	A~\emph{$\endofun$--twisted preshadow} or a~\emph{preshadow with monodromy\/
	$\endofun$} on $\ccat C$ valued in a~category $\cat T$ is a~collection of
	functors $\shdw_x\colon \ccat C(\endofun x,x) \to \cat T$,
	one per object $x\in \ccat C$, and natural morphisms
	$\theta_{g,f}\colon
		\shdw[y]{f\circ\endofun g} \to
		\shdw[x]{g\circ f}$ in $\cat T$,
	one for each pair of 1--morphisms $x\from^g y \from^f\endofun x$,
	such that the~following diagrams commute
	\begin{gather}
		\label{diag:preshadow-hex}
		\begin{diagps}(-14em,-1ex)(14em,21ex)
			\node f-Shg(-7.5em, 19ex)[{\shdw[y]{f\circ\endofun(h\circ g)}}]
			\node  hg-f( 7.5em, 19ex)[{\shdw[x]{(h\circ g) \circ f}}]
			\arrow{->}[f-Shg`hg-f;\theta]
			\node fSh-Sg(-10em, 12ex)[{\shdw[y]{(f\circ\endofun h)\circ \endofun g}}]
			\node  g-fSh( -5em,  0ex)[{\shdw[z]{g\circ (f\circ \endofun h)}}]
			\node  gf-Sh(  5em,  0ex)[{\shdw[z]{(g\circ f)\circ \endofun h}}]
			\node   h-gf( 10em, 12ex)[{\shdw[x]{h\circ (g\circ f)}}]
			\arrow|b{npos=0.6}|{->}[f-Shg`fSh-Sg;
				\shdw{\mathfrak a \ast (\iid\circ\mathfrak m^{-1})}]
			\arrow|b|{->}[fSh-Sg`g-fSh;\theta]
			\arrow{->}[g-fSh`gf-Sh;\shdw{\mathfrak a}]
			\arrow|b|{->}[gf-Sh`h-gf;\theta]
			\arrow|a{npos=0.6}|{<-}[hg-f`h-gf;\shdw{\mathfrak a}]
		\end{diagps}
	\\[2ex]
		\label{diag:preshadow-tri}
		\begin{diagps}(0em,-0.5em)(16em,12.5ex)
			\node tl( 0em,10ex)[{\shdw[x]{k\circ\id_{\endofun x}}}]
			\node tm( 8em,10ex)[{\shdw[x]{k\circ\endofun (\id_x)}}]
			\node tr(16em,10ex)[{\shdw[x]{\id_x\circm k}}]
			\node bm( 8em, 0ex)[{\shdw[x]{k}}]
			\arrow|a{npos=0.45}|{->}[tl`tm;\shdw{\iid\circ\mathfrak i}]
			\arrow{->}[tm`tr;\theta]
			\carrow*[arcangle=-20]|b|{->}[tl`bm;\shdw{\mathfrak r}]
			\carrow*[arcangle= 20]{->}[tr`bm;\shdw{\mathfrak l}]
		\end{diagps}
	\end{gather}
	for all $x \from^h z \from^g y \from^f\endofun x $ and $x\from^k\endofun x$.
	A~preshadow is a~\emph{shadow} if each $\theta_{g,f}$ is an~isomorphism.
	It is \emph{symmetric} when $\theta_{f,g}\circ\theta_{g,f} = \id$.
\end{definition}

\begin{remark}
	The~naturality of $\theta$ means that
		$\theta_{g',f'}\circ \shdw[y]{\alpha\circ\endofun\beta}
		= \shdw[x]{\beta\circ \alpha} \circ \theta_{g,f}$
	for any 2-morphisms $\alpha\colon f\dblto f'$ and $\beta\colon g\dblto g'$.
\end{remark}

We often write $\shdw\blank\colon (\ccat C,\endofun) \to \cat T$
for a~twisted preshadow on $\ccat C$, despite that it is not a~bifunctor.
It is understood that when morphism categories $\ccat C(x,y)$
are linear, then so are both $\cat T$ and the~functors $\shdw[x]\blank$.

Twisted shadows are categorified versions of traces in the~sense,
that the~equality $t(f\circ\endofun g) = t(g\circ f)$ defining a~(left)
twisted trace is replaced by a~choice of a~natural isomorphism
$\shdw{f\circ\endofun g} \cong \shdw{g\circ f}$.
When $\ccat C$ are $\endofun$ are strict, then the~diagrams
\eqref{diag:preshadow-hex} and \eqref{diag:preshadow-tri} simplify
to $\theta_{h\circ g,f} = \theta_{h,g\circ f} \theta_{g,f\circ\endofun h}$ and
$\theta_{\id,k} = \id_k$.

The~definition of a~preshadow coincides with the~one of a~categorified
trace from \cite{CategTrace} when $\endofun=\Id$ and $\ccat C$ has only
one object, i.e.\ it is a~monoidal category.
In comparison to \cite{Shadows} we do not require
$\theta_{f,g}\circ\theta_{g,f}=\id$.

%
% Definition: horizontal trace
%
In what follows we construct a~twisted preshadow for every bicategory%
---a~generalization of the~horizontal trace \cite{SL2asTrace}.
The~next section is devoted to show its functoriality
as well as universality.

\begin{definition}
	The~\emph{$\endofun$-twisted horizontal trace} of a~bicategory $\ccat C$
	is the~category $\hTr(\ccat C,\endofun)$, objects of which
	are 1--morphisms $f\in \ccat C(\endofun x, x)$,
	morphisms from $f\in\ccat C(\endofun x,x)$ to $g\in\ccat C(\endofun y,y)$
	are equivalence classes $[p,\alpha]$ of squares
	\begin{equation}\label{diag:morphism-in-shTr}
		\begin{diagps}(0em,-0.5ex)(5em,12.5ex)
			\square<5em,10ex>[\endofun x`x`\endofun y`y;f`\endofun p`p`g]
			\arrow[nodesep=1.65em]|b|{=>}[v1`v2;\alpha]
		\end{diagps}
	\end{equation}
	modulo the~relation
	\begin{equation}\label{rel:shTr-for-morphisms}
		\begin{diagps}(0em,-0.5ex)(6em,12.5ex)
			\square<5em,10ex>{->``->`->}[\endofun x`x`\endofun y`y;f``p`g]
			\carrow*[arcangle= 45]{->}[v0`v2;\endofun p\smash']
			\carrow*[arcangle=-45]|b|{->}[v0`v2;\endofun p]
			\morphism[nodesep=1.6em]{=>}|a{npos=0.53}
				|(5.5em,10ex)<-4em,-9.5ex>[`;\alpha]
			\morphism{<=}(-1em,5ex)<2em,0ex>[`;\endofun \tau]
		\end{diagps}
		\quad\sim\quad
		\begin{diagps}(-1em,-0.5ex)(5em,12.5ex)
			\square<5em,10ex>{->`->``->}[\endofun x`x`\endofun y`y;f`\endofun p'``g]
			\carrow*[arcangle= 45]{->}[v1`v3;p\smash']
			\carrow*[arcangle=-45]|b|{->}[v1`v3;p]
			\morphism[nodesep=1.6em]{=>}|b{npos=0.53}|
				(4em,10ex)<-4em,-9.5ex>[`;\alpha]
			\morphism{<=}(4em,5ex)<2em,0ex>[`;\tau]
		\end{diagps}
	\end{equation}
	with $[\id_x, \iid_f]$ the~identity on $f$, and the~composition
	\begin{equation}
		[q, \beta]\circ [p, \alpha] := [
			q\circ p,
			(\beta \circ \iid_{\endofun p}) \ast
			(\iid_q \circ \alpha)
		]
	\end{equation}
	can be visualized as stacking squares one on top of the~other:
	\begin{equation}
		\begin{diagps}(-1em,-0.5ex)(6em,12.5ex)
			\square<5em,10ex>[\endofun y`y`\endofun z`z;g`\endofun q`q`h]
			\arrow[nodesep=1.65em]|b|{=>}[v1`v2;\beta]
		\end{diagps}
	\quad\circ\quad
		\begin{diagps}(-1em,-0.5ex)(6em,12.5ex)
			\square<5em,10ex>[\endofun x`x`\endofun y`y;f`\endofun p`p`g]
			\arrow[nodesep=1.65em]|b|{=>}[v1`v2;\alpha]
		\end{diagps}
	\qquad:=\qquad
		\begin{diagps}(-1em,-0.5ex)(6em,14.5ex)
			\vsquares<5em,6ex>[
				\endofun x`x`\endofun y`y`\endofun z`z;
				f`\endofun p`p`g`\endofun q`q`h]
			\arrow[nodesep=1.5em]|b|{=>}[ur`ml;\alpha]
			\arrow[nodesep=1.5em]|b|{=>}[mr`dl;\beta]
		\end{diagps}
	\end{equation}
	Unitarity and associativity follows from \eqref{rel:shTr-for-morphisms}
	with an~appropriate composition of associators and unitors as $\tau$.
\end{definition}

%
% Construction: horizontal trace
%

\begin{proposition}\label{prop:hTr-is-preshadow}
	The~horizontal trace is the~target of a~preshadow
	$\shdw\blank^{\mathrm h}$ on $(\ccat C,\endofun)$ with components
	$\shdw[x]\blank^{\mathrm h}\colon
		\ccat C(\endofun x,x)\to\hTr(\ccat C,\endofun)$
	the~obvious functors and $\theta^{\;\!\mathrm h}_{g,f} := [g,\mathfrak a]$.
	It is a~shadow when $\ccat C$ has right duals.
\end{proposition}
\begin{proof}
	Commutation of both \eqref{diag:preshadow-hex}
	and \eqref{diag:preshadow-tri} is proven by a~direct
	computation, which is left to the~reader.
	When $\ccat C$ has right duals, consider a~morphism
	represented by
	\begin{equation}
		\begin{diagps}(0ex,-0.5ex)(10em,12.5ex)
			\hsquares(0em, 0ex)<5em,10ex>%
				|a`a`b`c`a`b`b|%
				{->`->`->``->`->`->}%
				[%
					\endofun x`y`x`\endofun y`\endofun x`y;%
					f`g`\endofun g^*``g^*`\endofun g`f%
				]
			\arrow|*c|{->}[ul`dm;\id_y]
			\arrow|*c|{->}[um`dr;\id_{\endofun x}]
			\node idY(2.5em,5ex)[\scriptstyle\phantom{\id_{\endofun x}}]
			\node idX(7.5em,5ex)[\scriptstyle\phantom{\id_y}]
			\arrow[nodesepA=3pt,nodesepB=4pt]|a{npos=0.45}|{=>}[idY`dl;\coev]
			\arrow[nodesepA=6pt,nodesepB=5pt]|b{npos=0.55}|{=>}[ur`idX;\ev]
		\end{diagps}
	\end{equation}
	where the~triangles are filled with the~coevaluation and evaluation 2-morphisms,
	and the~middle parallelogram with a~suitable composition of unitors.
	It is a~two-sided inverse of $\theta_{g,f}^{\;\!\mathrm h}$ due to the~relations
	\eqref{rel:left-duals} between $\mathit{ev}$ and $\mathit{coev}$.	
\end{proof}

%represented by the~square
%\begin{equation}\label{diag:def-of-theta}
%	\begin{diagps}(0ex,-0.5ex)(5em,12.5ex)
%		\square<5em,10ex>[%
%			x`\endofun x`y`\endofun y;
%			g\circ f`f`\endofun f`\endofun f\circ g]
%		\arrow[nodesep=1.65em]|b|{=>}[v1`v2;\mathfrak a]
%	\end{diagps}
%\end{equation}

\begin{figure}[ht]
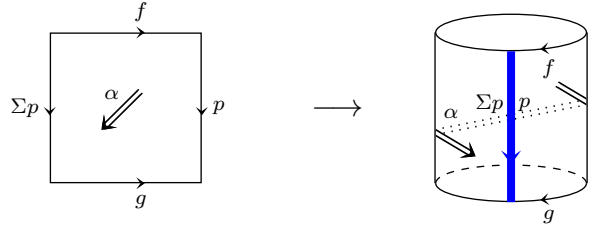
%
	\begin{displaymath}
	\psset{unit=1cm}%
		\begin{centerpict}(-1.5,0)(1.5,2)
			\usediagstyle
			\psframe(-1,0)(1,2)
			\rput[r](-1.15,1){$\scriptstyle \endofun p$}
			\rput[l]( 1.15,1){$\scriptstyle p$}
			\rput[b]( 0.2,2.15){$\scriptstyle f$}
			\rput[t]( 0.2,-0.15){$\scriptstyle g$}
			\morphism{=>}|b|(0.3,1.2)(-0.3,0.6)[`;\alpha]
			\psset{linestyle=none}
			\psline{->}(0,2)(0.25,2)
			\psline{->}(0,0)(0.25,0)
			\psline{->}(-1, 2)(-1, 0.9)
			\psline{->}( 1, 2)( 1, 0.9)
		\end{centerpict}
	\qquad\longrightarrow\qquad
		\begin{centerpict}(-1,0)(1,2.5)
			\begin{usediagstyle}
				\begin{psclip}{\psframe[linestyle=none](-1,0)(1,2)}
					\morphism{:}|b{npos=0.87}|(1.0,1.2)(-1.0,0.8)[`;\alpha]
					\morphism{=}(0.5,1.5)(1.5,0.9)%[`;]
					\morphism{=>}(-1.5,1.1)(-0.5,0.5)%[`;\alpha]
				\end{psclip}
				\psellipse(0,2.25)(1,0.25)
				\psellipticarc[linestyle=dashed,dash=3pt 3pt](0,0.25)(1,0.25){0}{180}
				\psellipticarc(0,0.25)(1,0.25){180}{360}
				\psline( 1, 0.25)( 1, 2.25)
				\psline(-1, 0.25)(-1, 2.25)
				\psline[linewidth=3pt,linecolor=blue](0,0)(0,2)
				\psline[linewidth=3pt,linecolor=blue,linestyle=none]{->}(0,2)(0,0.5)
				\psellipticarcn[linestyle=none]{->}(0,2.25)(1,0.25){0}{-30}
				\psellipticarcn[linestyle=none]{->}(0,0.25)(1,0.25){0}{-30}
				\rput[t](0.5, 1.9){$\scriptstyle f$}
				\rput[t](0.5,-0.1){$\scriptstyle g$}
				\rput[l](0.1, 1.3){$\scriptstyle p$}
				\rput[r](-0.1,1.3){$\scriptstyle \endofun p$}
			\end{usediagstyle}
		\end{centerpict}
	\end{displaymath}
	\caption{Visualization of the~horizontal trace.}%
	\label{fig:hTr-as-cylinders}%
\end{figure}

\begin{figure}[ht]
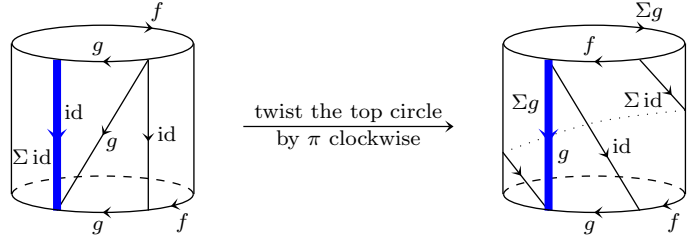
%
	\begin{displaymath}
	\psset{xunit=12mm,yunit=10mm}%
		\begin{centerpict}(-1.1,-1.5)(1.1,1.3)
			\begin{usediagstyle}
			% cylinder and lines
				\psellipticarc[linestyle=dashed,dash=3pt 3pt](0,-1)(1,0.25){0}{180}
				\psline( 0.5, 0.78349)(-0.5, -1.21651)
				\psline( 0.5, 0.78349)( 0.5, -1.21651)
				\psline[linewidth=3pt,linecolor=blue](-0.5, 0.78349)(-0.5, -1.21651)
				\psellipse[fillstyle=solid](0,1)(1,0.25)
				\psline(-1,1)(-1,-1)
				\psline( 1,1)( 1,-1)
				\psellipticarc(0,-1)(1,0.25){180}{360}
			% arrows and labels
				\psset{linestyle=none}
				\psellipticarcn{->}(0, 1)(1,0.25){ 0}{-90}
					\rput[b](-0.05, 0.85){$\scriptstyle g$}
				\psellipticarcn{->}(0, 1)(1,0.25){90}{20}
					\rput[b](0.6, 1.3){$\scriptstyle f$}
				\psellipticarcn{->}(0,-1)(1,0.25){ 0}{-90}
					\rput[t](-0.05,-1.35){$\scriptstyle g$}
				\psellipticarcn{->}(0,-1)(1,0.25){90}{-12}
					\rput[tl](0.8,-1.25){$\scriptstyle f$}
				\psline{->}( 0.5,0.5)( 0.5, -0.3)
					\rput[l]( 0.6,-0.2){$\scriptstyle\id$}
				\psline{->}( 0.5, 0.78349)(0, -0.21651)
					\rput[t]( 0.1,-0.25){$\scriptstyle g$}
				\psline[linewidth=3pt,linecolor=blue]{->}(-0.5, 0.5)(-0.5, -0.31651)
					\rput[l](-0.42, 0.1){$\scriptstyle\id$}
					\rput[r](-0.56,-0.5){$\scriptstyle\endofun\id$}
			\end{usediagstyle}
		\end{centerpict}
	\quad\to^{\text{twist the~top circle}}
	\diagbput{\scriptstyle\text{by $\pi$ clockwise}}\quad
		\begin{centerpict}(-1.1,-1.5)(1.1,1.3)
			\begin{usediagstyle}
			% cylinder and lines
				\psellipticarc[linestyle=dashed,dash=3pt 3pt](0,-1)(1,0.25){0}{180}
				\psline(-0.5, 0.78349)( 0.5, -1.21651)
				\psline( 0.5, 0.78349)( 1.0, 0.11349)
				\pscurve[linestyle=dotted](1, 0.11349)(0,-0.05)(-1, -0.44651)
				\psline(-1, -0.44651)(-0.5, -1.21651)
				\psline[linewidth=3pt,linecolor=blue](-0.5, 0.78349)(-0.5, -1.21651)
				\psellipse[fillstyle=solid](0,1)(1,0.25)
				\psline(-1,1)(-1,-1)
				\psline( 1,1)( 1,-1)
				\psellipticarc(0,-1)(1,0.25){180}{360}
			% arrows and labels
				\psset{linestyle=none}
				\psellipticarcn{->}(0, 1)(1,0.25){ 0}{-90}
					\rput[b](-0.05, 0.85){$\scriptstyle f$}
				\psellipticarcn{->}(0, 1)(1,0.25){90}{20}
					\rput[b](0.6, 1.3){$\scriptstyle\endofun g$}
				\psellipticarcn{->}(0,-1)(1,0.25){ 0}{-90}
					\rput[t](-0.05,-1.35){$\scriptstyle g$}
				\psellipticarcn{->}(0,-1)(1,0.25){90}{-12}
					\rput[tl](0.8,-1.25){$\scriptstyle f$}
				\psline{->}(-0.5, 0.78349)(0.15, -0.51651)
					\rput[l](0.2,-0.35){$\scriptstyle\id$}
				\psline[linewidth=3pt,linecolor=blue]{->}(-0.5, 0.5)(-0.5, -0.31651)
					\rput[l](-0.4,-0.5){$\scriptstyle g$}
					\rput[r](-0.6, 0.2){$\scriptstyle\endofun g$}
				\psline{->}( 0.5, 0.78349)(0.85, 0.31449)
				\psline{->}(-1,-0.44651)(-0.8,-0.75451)
					\rput[tr](0.75,0.35){$\scriptstyle\endofun\id$}
			\end{usediagstyle}
		\end{centerpict}
	\end{displaymath}
	\caption{Visualization of the~cyclicity morphism in $\hTr$.}%
	\label{fig:theta-as-twist}%
\end{figure}

The~idea underlying the~definition of the~horizontal trace
is to represent 1-endo\-mor\-phisms by boundary circles of a~cylinder,
the~interior of which carries a~2-morphism. Indeed, such
a~cylinder can be obtained by gluing the~vertical boundaries
of the~square \eqref{diag:morphism-in-shTr} together, but
keeping the~seam line, see Figure~\ref{fig:hTr-as-cylinders}.
The~seam is cooriented and a~2-morphism can ``slide'' through it,
in which case it is acted upon by $\endofun$
as described by \eqref{rel:shTr-for-morphisms}.
The~morphism $\theta^{\;\!\mathrm h}$ is represented in this picture
by the~cylinder obtained from the~identity on $g\circ f$ by
applying a~half-twist to its top boundary circle,
see Figure~\ref{fig:theta-as-twist}.

%
% Other names for hTr
%
The~symmetric horizontal trace (i.e.\ for $\endofun=\Id$) is also called
the~\emph{annularization functor} \cite{BlobComplex}. Its generalization
to any surface is known as \emph{factorization homology} \cite{FactHom}.

%% ============================================================================
%%  Section: horizontal trace
%% ============================================================================

\subsection{The~universal property of the~horizontal trace}\label{sec:hor-trace}
\label{sec:univeral-shadow}

Let $\shdw\blank\colon (\ccat C,\endofun) \to \cat T$ and
$\shdw\blank'\colon (\ccat C,\endofun) \to \cat T'$ be two preshadows.
We say that $\shdw\blank'$ \emph{factorizes} through $\shdw\blank$
if there exist a~functor $T\colon \cat T\to \cat T'$ and a~collection of natural
isomorphisms $\eta_x\colon T\circ \shdw[x]\blank \to \shdw[x]\blank\!\!'$
such that $\eta_x \circ T(\theta_{g,f}) = \theta_{g,f}'\circ \eta_y$
for any $x \from^g y \from^f\endofun x$.
Factorizations $(T,\eta)$ and $(T',\eta')$ are \emph{equivalent}
if there exists a~natural isomorphism $\epsilon\colon T\to T'$ such
that $\eta_x = \eta_x'\ast(\epsilon\circ \mathbf 1_{\shdw\blank})$, i.e.
\begin{equation}
	\begin{diagps}(0em,-0.5ex)(12em,14ex)
		\node  C(0em, 6ex)[\ccat C(\endofun x,x)]
		\node T1(8em,12ex)[\cat T]
		\node T2(8em, 0ex)[\cat T']
		\carrow*[arcangle= 15]|a|{->}[C`T1;{\shdw[x]\blank}]
		\carrow*[arcangle=-15]|b|{->}[C`T2;{\shdw[x]\blank\!\!'}]
		\carrow*[arcangle= 40]|a|{->}[T1`T2;T\phantom']
		\carrow*[arcangle=-40]|b|{->}[T1`T2;T']
		\morphism{=>}|b|(9em,6ex)<-2em,0ex>[`;\epsilon]
		\morphism{=>}|b{labelsep=-2pt}|(5.5em,8ex)<-1.7em,-3.2ex>[`;\eta']
	\end{diagps}
	=
	\begin{diagps}(-4em,-0.5ex)(8em,14ex)
		\node  C(0em, 6ex)[\ccat C(\endofun x,x)]
		\node T1(8em,12ex)[\cat T]
		\node T2(8em, 0ex)[\cat T']
		\carrow*[arcangle= 15]|a|{->}[C`T1;{\shdw[x]\blank}]
		\carrow*[arcangle=-15]|b|{->}[C`T2;{\shdw[x]\blank\!\!'}]
		\carrow*[arcangle= 40]|a|{->}[T1`T2;T\phantom']
		\morphism{=>}|b{labelsep=-1pt}|(6.5em,8.5ex)<-1.5em,-4ex>[`;\eta]
	\end{diagps}
\end{equation}
The~following result states that the~horizontal trace twisted by $\endofun$
is the~universal $\endofun$--twisted preshadow.

%
% Universality of hTr
%
\begin{theorem}\label{thm:hTr-is-universal}
	If\/ $\ccat C$ has left duals, then every preshadow on $(\ccat C,\endofun)$
	factorizes through\/ $\hTr(\ccat C,\endofun)$ uniquely up to an~equivalence.
%	Moreover, $\shdw\blank^{\!\mathrm h}$ is a~shadow when $\ccat C$ has right duals.
\end{theorem}
\begin{proof}
	Given a~preshadow $\shdw\blank \colon (\ccat C,\endofun) \to \cat T$
	we construct a~functor $T\colon \hTr(\ccat C,\endofun) \to \cat T$
	by taking $f\colon \endofun x\to x$ to $\shdw[x]f$ and a~morphism $[p,\alpha]\colon f\to g$ to the~composition
	\begin{equation}\label{eq:def-of-shdw-h}
	\begin{aligned}
			\shdw[x]f
				\to^{\shdw[x]{\coev\mcircm\iid}}
			\shdw[x]{\ldual p\circ p\circ f}
				\to^{\shdw[y]{\iid\circ\alpha}}
			&\shdw[x]{\ldual p\circ g\circ \endofun p}
		\\[1ex]
				\to^{\theta_{p,\ldual p\circ g}}
			&\shdw[y]{p\circ \ldual p \circ g}
				\to^{\shdw[y]{\ev\circ\iid}}
			\shdw[y]g,
	\end{aligned}
	\end{equation}
	where associators and unitors are omitted for clarity.
	Notice that $\shdw\blank = T\circ\shdw\blank^{\mathrm h}$ and
	$\theta_{g,f} = T(\theta^{\;\!\mathrm h}_{g,f})$.
	For uniqueness, suppose that
	$(T',\eta)\colon \shdw\blank^{\mathrm h} \to \shdw\blank$ is another
	factorization of a~morphism of shadows with each $\eta_x$ an~isomorphism.
	Then $\eta$ is a~natural isomorphism of functors $\eta\colon T'\to T$,
	because $\shdw f^{\mathrm h} = f$ and $Tf = \shdw f$;
	naturality follows from \eqref{eq:def-of-shdw-h} and
	the~observation that a~similar sequence determines $T'([p,\alpha])$.
\end{proof}

\begin{corollary}
	Suppose $\ccat C$ has both left and right duals.
	Then any preshadow on $\ccat C$ is a~shadow.
\end{corollary}
\begin{proof}
	It follows from Theorem~\ref{thm:hTr-is-universal} that
	a~preshadow $(\shdw\blank,\theta)$ factorizes through
	the~horizontal trace. In particular, $\theta$ is an~image
	of $\theta^{\;\!\mathrm h}\noic$, which is an~isomorphism by
	Proposition~\ref{prop:hTr-is-preshadow}.
\end{proof}

\subsection{Functoriality}
\label{sec:functoriality-of-hTr}

To discuss functoriality of $\hTr$ we begin with constructing
a~suitable 3-category of endobicategories.

%
%  Morphism of endobicategories
%
A~\emph{morphism of endobicategories}
$(\ccat F,\omega)\colon (\ccat C,\endofun)\to (\ccat C',\endofun')$
consists of a~strong bifunctor $\ccat F\colon \ccat C\to \ccat C'$,
so that the~coherence 2-morphisms
$\mathfrak m_{\ccat F}\colon \ccat Fg\circ \ccat Ff \dblto \ccat F(g\circ f)$ and
$\mathfrak i_{\ccat F}\colon \id_{\ccat Fx} \dblto \ccat F(\id_x)$ are invertible,
and a~weak natural transformation $\omega\colon \endofun'\ccat F \to
\ccat F\endofun$. The~latter comes with a~collection of 2-morphisms
\begin{equation}
	\begin{diagps}(-1em,-1ex)(7em,12ex)
		\square<6em,10ex>[%
			\endofun'\ccat F x`\ccat F\endofun x`
			\endofun'\ccat F y`\ccat F\endofun y;
			\omega_x`\endofun'\ccat F f`\ccat F\endofun f`\omega_y
		]
		\arrow[nodesep=1.4em]|b|{=>}[v1`v2;\omega_f]
	\end{diagps}
\end{equation}
parametrized with 1-morphisms $f\in\ccat C(x,y)$.
The 2-morphisms $\omega_f$ are not required to be invertible.
If they are, we say that $(\ccat F,\omega)$ is \emph{strong}.

%
%  Natural transformations
%
Suppose that $(\ccat F,\omega)$ and $(\ccat F',\omega')$ are morphisms
between the~same endobicategories. A~\emph{natural transformation}
$(\eta, \mathfrak n)\colon (\ccat F,\omega)\to (\ccat F',\omega')$ consists
of a~strong natural transformation of bifunctors $\eta\colon \ccat F\to \ccat F'$
and a~natural 2-morphism
$\mathfrak n_x\colon
	\eta_{\endofun x} \circ \omega_x \to/=>/
	\omega_x' \circ \endofun'\eta_x$
that---in addition to the~usual coherence conditions---fits into
the~commutative hexagon
\begin{equation}
	\begin{diagps}(-10em,-9ex)(10em,10.5ex)
		\node ml(-10em, 0ex)[\ccat F'\endofun f \circ\eta_{\endofun x} \circ\omega_x]
		\node mr( 10em, 0ex)[\omega'_y \circ \endofun'(\eta_y \circ \ccat Ff)]
		\node tl( -6em, 8ex)[\ccat F'\endofun f \circ\omega'_x \circ\endofun\eta_x]
		\node tr(  6em, 8ex)[\omega'_y \circ\endofun'\ccat F'f \circ\endofun'\eta_x]
		\node bl( -6em,-8ex)[\eta_{\endofun y} \circ\ccat F\endofun f \circ\omega_x]
		\node br(  6em,-8ex)[\eta_{\endofun y} \circ\omega_y \circ \endofun'\ccat Ff]
		\arrow{=>}[ml`tl;\iid\circ\mathfrak n_x]
		\arrow{=>}[tl`tr;\omega'_f\circ\iid]
%		\arrow|a{npos=0.6}|{=>}[tr`mr;\iid\circ(\endofun'\eta_f\ast\mathfrak m)]
		\arrow|a{npos=0.6}|{=>}[tr`mr;\iid\circ\endofun'\eta_f]
		\arrow|b{npos=0.5}|{=>}[ml`bl;\eta_{\endofun f}\circ\iid]
		\arrow{=>}[bl`br;\iid\circ\mskip\thinmuskip\omega_f]
%		\arrow|b|{=>}[br`mr;(\iid\circ\mathfrak m) \ast (\mathfrak n_y\circ\iid)]
		\arrow|b|{=>}[br`mr;\mathfrak n_y\circ\iid]
	\end{diagps}
\end{equation}
for any 1-morphism $f\in \ccat C(x,y)$, where canonical isomorphisms are omitted.

%
%  Modifications
%
Finally, we require that a~\emph{modification} $\Gamma\colon \eta \to \eta'$
of natural transformations $(\eta,\mathfrak n)$ and $(\eta', \mathfrak n')$
satisfies $(\iid\circ\endofun\Gamma_{\! x})\ast \mathfrak n_x =
\mathfrak n'_x \ast (\Gamma_{\!\endofun'x}\circ\iid)$ in addition to the~usual coherence condition.

\begin{definition}
	We write $\ccat{EndoBicat}$ for the~3-category of endobicategories,
	morphisms, natural transformations and modifications, as described above.
	Its restriction to endobicategories with left duals is denoted
	by $\ccat{EndoBicat}^*$.
\end{definition}

%
% Pullback of a preshadow
%
Preshadows can be pulled back along morphisms of endobicategories.
Indeed, given a~preshadow $\shdw\blank$ on $(\ccat C',\endofun')$
and a~morphism $(\ccat F,\omega)\colon (\ccat C,\endofun) \to (\ccat C',\endofun')$
we define a~its pullback $(\ccat F^*\!\shdw\blank, \ccat F^*\theta)$
on $(\ccat C,\endofun)$ as follows:
\begin{align*}
	\ccat F^*\!\shdw[x]{f} &:= \shdw[Fx]{\ccat Ff \circ \omega_x} \\
	\ccat F^*\!\shdw[x]{\alpha} &:= \shdw[Fx]{\ccat F\alpha \circ \iid}
\end{align*}
with the~cyclicity morphism $(\ccat F^*\theta)_{g,f}$ the~composition
\begin{equation*}
	\begin{split}
		\ccat F^*\!\shdw[y]{f\circ\endofun g} =
					\shdw[Fy]{\ccat F(f\circ \endofun g) \circ \omega_y}
			&\to<5em>^{\shdw{\mathfrak m_{\ccat F}^{-1}\circ\iid}}
					\shdw[Fy]{\ccat Ff\circ \ccat F\endofun g \circ \omega_y} \\
			&\to<5em>^{\shdw{\iid\circ\omega_g}}
					\shdw[Fy]{\ccat Ff\circ \omega_x \circ \endofun'\ccat Fg}  \\
			&\to<5em>^{\theta}
					\shdw[Fx]{\ccat Fg\circ \ccat Ff \circ \omega_x}	\\
			&\to<5em>^{\shdw{\mathfrak m_{\ccat F}^{\vphantom1}\circ\iid}}
					\shdw[Fx]{\ccat F(g\circ f)\circ\omega_x}
			= \ccat F^*\!\shdw[x]{g\circ f}.
	\end{split}
\end{equation*}
Notice that the~second arrow is invertible when $(\ccat F,\omega)$
is a~strong morphism.

\begin{proposition}\label{prop:shadow-from-nat-trans}
	The~datum\/ $(\ccat F^*\!\shdw\blank,\ccat F^*\theta)$ is a~preshadow
	on $(\ccat C,\endofun)$. It is a~shadow when $(\ccat F,\omega)$
	is a~strong morphism and\/ $(\shdw\blank,\theta)$ is a~shadow.
\end{proposition}
\begin{proof}
	Left as an~exercise.
\end{proof}

Let $\shdw\blank, \shdw\blank' \colon (\ccat C,\endofun) \to \cat T$
be two preshadows. A~collection of natural transformations
$\eta_x\colon \shdw\blank_x \to \shdw\blank_x'$ is a~\emph{morphism
of preshadows} if it is coherent with the~cyclicity morphisms:
$\theta'_{g,f}\circ\eta_y = \eta_x\circ \theta_{g,f}$.
There is a~morphism of pullback shadows associated
to any natural transformation of morphisms of endobicategories
\begin{equation}
	\begin{diagps}(-4em,-5ex)(4em, 7ex)
		\node[href= 0.8] dom(-4.5em,0ex)[(\ccat C,\endofun)]
		\node[href=-0.8] cod( 4.5em,0ex)[(\ccat C',\endofun')]
		\carrow*[arcangle=50]{->}[dom`cod;(\ccat F,\omega)]
		\carrow*[arcangle=-50]|b|{->}[dom`cod;(\ccat F',\omega')]
		\morphism{=>}|*c{npos=0.55}|(0,4ex)(0,-3ex)[`;(\eta,\mathfrak n)]
	\end{diagps}
\end{equation}
Indeed, there is a~morphism in $\hTr(\ccat C',\endofun')$
\begin{equation}\label{diag:nat-trans-on-hTr}
	\begin{diagps}(-5em,0ex)(5em,12ex)
		\hsquares(-5em,0ex)<5em,10ex>[%
			\endofun'\ccat F x`\ccat F \endofun x`\ccat F x`
			\endofun'\ccat F'x`\ccat F'\endofun x`\ccat F'x;
			\omega_x`\ccat F f`
			\endofun'\eta_x`\eta_{\endofun x}`\eta_x`
			\omega_x'`\ccat F'f
		]
		\arrow[nodesep=1.4em]|b|{=>}[um`dl;\mathfrak n_x]
		\arrow[nodesep=1.4em]|b|{=>}[ur`dm;\eta_f^{-1}]
	\end{diagps}
\end{equation}
which, assuming that $\ccat C'$ has left duals, descend to a~morphism
$(\eta,\mathfrak n)_f$ in $\cat T$ between $\ccat F^*\shdw f$
and $\ccat F'^*\shdw f$ for any shadow $\shdw\blank$ on $(\ccat C',\endofun')$.

\begin{lemma}\label{lem:hTr-for-trans}
	The~collection $(\eta,\mathfrak n)_*$ is a~morphism of shadows
	from $\ccat F^*\shdw\blank$ to\/ $\ccat F'^*\shdw\blank$.
	Morphisms induced by isomorphic transformations are equal.
\end{lemma}
\begin{proof}
	The~naturality and compatibility with cyclicity morphisms can be
	checked in $\hTr(\ccat C',\endofun')$. They follow directly from
	the~coherence conditions of $\eta_f$ and $\mathfrak n_x$.
	The~details are left as an~exercise.
\end{proof}

%
% Functoriality
%
Everything together implies functoriality of the~horizontal trace.
Let $\ccat{EndoBicat}^{*,\leqslant 2}$ be the~restriction of
$\ccat{EndoBicat}^*$ to a~2-category, obtained by forgetting
modifications and identifying isomorphic natural transformations.

\begin{theorem}[Functoriality of the~horizontal trace]\label{thm:hTr-functoriality}
	The~horizontal trace $\hTr$ extends to a~strict 2-functor\/
	$\hTr\colon\ccat{EndoBicat}^{*,\leqslant 2} \to \ccat{Cat}$,
	which assigns to a~morphism of endobicategories
	$(\ccat F,\omega)\colon (\ccat C,\endofun)\to(\ccat C',\endofun')$
	the~functor\/
	$\hTr(\ccat F,\omega)\colon \hTr(\ccat C,\endofun) \to \hTr(\ccat C',\endofun')$
	that takes an~object $f\colon \endofun x\to x$
	into\/ $\ccat Ff\circ\omega_x$ and a~morphism $[p,\alpha]$ to the~class
	of the~diagram
	\begin{equation}
		\begin{diagps}(0em,-0.5ex)(12em,12ex)
			\hsquares<6em,10ex>[
				\endofun'\ccat Fx`\ccat F\endofun x`\ccat Fx`
				\endofun'\ccat Fy`\ccat F\endofun y`\ccat Fy;
				\omega_x`\ccat Ff`
				\endofun'\ccat Fp`\ccat F\endofun p`\ccat Fp`
				\omega_y`\ccat Fg]
			\arrow[nodesep=1.4em]|b|{=>}[um`dl;\omega_p]
			\arrow[nodesep=1.5em]|b|{=>}[ur`dm;\ccat F\alpha]
		\end{diagps}
	\end{equation}
	and a~natural transformation $(\eta,\mathfrak n)$ is sent to
	the~natural transformation\/ $\hTr(\eta,\mathfrak n)$ with components
	\eqref{diag:nat-trans-on-hTr}.
\end{theorem}
\begin{proof}
	Apply Theorem~\ref{thm:hTr-is-universal} to the~pullback
	along $(\ccat F,\omega)$ of the~universal preshadow on
	$(\ccat C',\endofun')$ to construct $\hTr(\ccat F,\omega)$.
	The~coherence condition for a~natural transformation guarantees
	that $\hTr(\eta,\mathfrak n)$ is a~morphism of shadows.
\end{proof}

\begin{corollary}
	Let\/ $(\ccat F,\omega)\colon (\ccat C,\endofun)\to(\ccat C',\endofun')$
	be a~local equivalence, i.e.\ $\omega$ is invertible
	and\/ $\ccat F\noic$ restricts to equivalences of morphism categories
	$\ccat C(x,y) \approx \ccat C'(\ccat Fx,\ccat Fy)$ for all objects
	$x,y\in\ccat C$.
	Then\/ $\hTr(\ccat F, \omega)$ is full and faithful.
\end{corollary}
\begin{proof}
	By replacing $\ccat C'$ with the~image of $\ccat F$ we may assume
	that the~bifunctor is a~biequivalence.
	Then $\hTr(\ccat F,\omega)$ is an~equivalence of categories by
	Theorem~\ref{thm:hTr-functoriality} and Lemma~\ref{lem:hTr-for-trans},
	hence, full and faithful.
\end{proof}

%% ============================================================================
%%  Section: vertical trace
%% ============================================================================

\subsection{A~connection with the~vertical trace}
\label{sec:vTr-vs-hTr}

%
% Definition: vertical trace
%
Assume now that $\endofun$ fixes objects, i.e.\ $\endofun x = x$
for any object $x\in\ccat C$. Then $\endofun$ can be restricted
to each morphism category $\ccat C(x,y)$. Replacing each of them
by its trace $\Tr(\ccat C(x,y),\endofun)$ results in the~\emph{%
vertical trace} $\vTr(\ccat C,\endofun)$: a~category with the~same
objects as $\ccat C$, morphisms twisted trace classes of 2-morphisms
$\alpha\colon f\to/=>/ \endofun f$, and composition induced by
the~horizontal composition of 2-morphisms.
The~composition in $\vTr(\ccat C,\endofun)$ is both unital and
associative, because the~trace class of a~morphism is invariant
under conjugation.

%
% Natural transformation vTr --> hTr
%
There is a~natural functor $\vTr(\ccat C,\endofun)\to \hTr(\ccat C,\endofun)$, defined by expanding objects to identity morphisms. Explicitly,
\begin{equation*}
	x \mapsto (x\to^{\id_x} x)
		\qquad\text{and}\qquad
	\left(\begin{diagps}(-2em,1ex)(2em,9.5ex)
		\node t(0em,9ex)[x]
		\node b(0em,0ex)[y]
		\carrow*[arcangle= 45]{->}[t`b;f]
		\carrow*[arcangle=-45]|b|{->}[t`b;\endofun f]
		\morphism{<=}(-1em,4.5ex)<2em,0ex>[`;\sigma]
	\end{diagps}\right) \mapsto
	\begin{diagps}(-1.5em,-0.5ex)(4em,11ex)
		\square<4em,9ex>[x`x`y`y;\id_x`\endofun f`f`\id_y]
		\arrow[nodesep=1.25em]|b|{=>}[v1`v2;\sigma]
	\end{diagps}
\end{equation*}
It is clearly full and faithful, but not necessarily surjective on objects.

When $\ccat C$ has left duals, then the~above functor can be used
to restrict a~preshadow $\shdw\blank \colon (\ccat C,\endofun) \to \cat T$
to a~trace $t^{\shdw\blank} \colon (\ccat C(x,y), \endofun) \to
	\cat T(\shdw{\id_x},\shdw{\id_y})$
using the~formula
\begin{equation}\label{def:tr-from-shdw}
	t^{\shdw\blank}(p\to^\alpha \endofun p) := \shdw{[p,\alpha]},
\end{equation}
where we identify $\shdw\blank$ with the~corresponding functor on
$\hTr(\ccat C,\endofun)$. The~following is an~immediate consequence
of this construction.

\begin{proposition}
	The~restriction to $\ccat C(x,y)$ of the~universal\/ $\endofun$--twisted
	preshadow is the~universal\/ $\endofun$--twisted trace.
\end{proposition}

%% ============================================================================
%%  Section: hTr as a~categorification of Tr
%% ============================================================================

\subsection{Decategorification}
\label{sec:hTr-categ-Tr}

An~$n$--category $\ccat C$ can be truncated to an~$(n-1)$--category by forgetting its highest level morphisms and identifying isomorphic $(n-1)$--morphisms. Thence, we have a~functor $\Pi\colon n\ccat{Cat} \to (n{-}1)\ccat{Cat}$. Clearly, 	this construction applies to endo-$n$-categories as well.
We are mostly interested in the~case $n=1,2$.

Applying $\Pi$ to a~shadow $\shdw\blank \colon \ccat C \to \cat T$ results in a~(left) trace $\Pi\shdw\blank \colon \Pi\ccat C \to \Pi\cat T$. Conversely, any left trace on $\Pi\ccat C$ can be lifted to a~shadow if $\ccat C$ is small.\footnote{
	A~bicategory is \emph{small} if its 2-morphisms (and so objects as well as 1-morphisms) form a~set.
}

\begin{lemma}\label{lem:categ-trace}
	Choose a~small endobicategory\/ $(\ccat C,\endofun)$ and a~(left) trace
	$\tr \colon (\Pi\ccat C, \Pi\endofun) \to S$ for some set\/ $S\noic$.
	There exists a~shadow\/ $T\colon (\ccat C,\endofun) \to \cat S\noic$,
	such that\/ $\Pi\cat S = S\noic$ and\/ $\Pi T = \tr$.
\end{lemma}
\begin{proof}
	In what follows we write $\dom(\alpha) := f$ and $\cod(\alpha) := f'$
	for the~domain and codomain of a~2--morphism $\alpha\colon f\dblto f'$.
	Let $\tilde{\cat S}$ be the~category with $\Ob(\tilde{\cat S}) = S$
	and morphisms $\tilde{\cat S}(s,t)$ the~finite sequences
	$(\alpha_n,\dots,\alpha_1)$
	of 2--morphisms $\alpha_i\in\ccat C(\endofun x_i,x_i)$ satisfying
	$\tr(\dom(\alpha_{i+1})) = \tr(\cod(\alpha_{i}))$ for $i=1,\dots,n-1$,
	such that $s=\tr(\dom(\alpha_1))$ and $t=\tr(\cod(\alpha_n))$.
	Composition of morphisms is defined as concatenation of sequences.
	The~category $\cat S$ is a~quotient of $\tilde{\cat S}$ by the~relations
	\begin{align}
		(\dots, \beta, \alpha, \dots) &\sim (\dots, \beta \ast \alpha, \dots)
	\quad\text{whenever }\beta\ast\alpha\text{ exists, and}\\
		(\dots, \mathbf 1, \dots) &\sim (\dots,\dots).
	\end{align}
	They allow us to reduce a~given sequence of 2--morphism to a~sequence of noncomposable morphisms, none of which is the~identity. An~easy application of the~Bergman Diamond Lemma \cite{DiamondLemma} shows that this reduced sequence is unique. In particular,
	\begin{equation}
		(\beta_m,\dots,\beta_1)\circ (\alpha_1,\dots,\alpha_n) \sim ()
	\end{equation}		
	implies $m=n$ and each $\beta_i\circ\alpha_i = \mathbf 1$ if both sequences are reduced. Hence, each $\alpha_i$ is a~2--isomorphism in $\ccat C$ if $(\alpha_1,\dots,\alpha_n)$ is an~isomorphism in $\cat S$, which implies that $\tr(\dom(\alpha_i)) = \tr(\cod(\alpha_i))$ and 
	\begin{equation}
		\dom(\alpha_n,\dots,\alpha_1) = \tr(\dom(\alpha_1)) = \tr(\cod(\alpha_n)) = \cod(\alpha_n,\dots,\alpha_1).
	\end{equation}
	Therefore, $\cat S$ does not have isomorphisms between different objects, and $\Pi\cat S = S$. The~desired shadow $T\colon (\ccat C,\endofun)\to \cat S$ is defined as $T(f) := \tr(f)$ and $T(\alpha) := (\alpha)$.
\end{proof}

The~above result allows us to formally prove that the~horizontal trace
is a~categorification of the~universal trace of a~category.
Despite proving it only for small bicategories, we believe that
with a~slight modification of our argument the~result can
be generalized to locally small bicategories.

\begin{theorem}
	There is a~natural bijection\/
	$\Pi(\hTr(\ccat C,\endofun)) \approx \Tr^\ell(\Pi\ccat C,\Pi\endofun)$
	for each small endobicategory $(\ccat C,\endofun)$ with left and right duals.
%	In other words, the~diagram of functors
%	\begin{equation}
%		\begin{diagps}(0,0)(9em,12ex)
%			\square<9em,10ex>|a{npos=0.6}`b`a`a{npos=0.6}|[
%				\ccat{EndoBicat}^{*,\leqslant 2}`\ccat{Cat}`
%				\cat{EndoCat}^{\leqslant 1}`\cat{S\!et};
%				\hTr`\Pi`\Pi`\Tr^\ell
%			]
%		\end{diagps}
%	\end{equation}
%	commutes up to a~natural isomorphism.
\end{theorem}
\begin{proof}
	Because $\ccat C$ has right duals, $\hTr(\ccat C,\endofun)$ is a~shadow
	and we can apply $\Pi$ to get a~trace $\Pi(\hTr(\ccat C,\endofun))$.
	We must show that it is the~universal left trace,
	i.e.\ that every left trace $t\colon (\Pi\ccat C,\Pi\endofun) \to S$
	factorizes through it.
	For that use Lemma~\ref{lem:categ-trace} to find a~shadow
	$T\colon (\ccat C,\endofun)\to \cat S$ that lifts $t$.
	According to Theorem~\ref{thm:hTr-is-universal},
	$T$ factorizes through $\hTr(\ccat C,\endofun)$, and applying $\Pi$
	results in a~desired factorization of $t$.
\end{proof}

Thence, the~horizontal trace is a~categorification of the~universal left trace.
Right traces are categorified by a~dual version of shadows, defined on morphism
categories $\ccat C(x,\endofun x)$. The~convention for the~twisting of $\hTr$ used
in this paper is justified by the~relation of $\hTr$ with the~vertical trace
discussed in Section~\ref{sec:vTr-vs-hTr}.

\subsection{Categorified Lefschetz traces}
\label{sec:lef-shadow}

Let $\HCom^b(\ccat C)$ be a~bicategory obtained from a~locally additive
bicategory $\ccat C$ by replacing each morphism category $\ccat C(x,y)$
with the~corresponding homotopy category of bounded complexes
$\HCom^b(\ccat C(x,y))$. It is locally triangulated, i.e.\ the~functors
$f\circ(\blank)$ and $(\blank)\circ f$ preserve distinguished triangles
and commute with homological degree shifts. Having chosen an~endobifunctor
$\endofun$ on $\ccat C$ we extend it naturally over $\HCom^b(\ccat C)$.

\begin{definition}
	A~preshadow $\shdw\blank \colon (\HCom^b(\ccat C),\endofun) \to \cat T$ is
	\emph{triangulated} if $\cat T$ is a~triangulated category and the~components
	$\shdw[x]\blank$ are triangulated functors.
\end{definition}

A~twisted preshadow $\shdw\blank\colon (\ccat C,\endofun) \to \cat T$
can be extended to a~twisted triangulated preshadow
$\sshdw\blank \colon \HCom^b(\ccat C)(\endofun x,x) \to \HCom^b(\cat T)$,
which we call the~\emph{Lefschetz preshadow induced by $\shdw\blank$}.
For each object $x\in\ccat C$ the~functor $\sshdw\blank_x$ is constructed
by applying $\shdw[x]\blank$ component-wise:
\begin{equation}
	\sshdw{f^\bullet}_x = \big(
		\dots \to \shdw[x]{f^i} \to^{\shdw[x]d}
		    \shdw[x]{f^{i+1}} \to^{\shdw[x]d}
		    \shdw[x]{f^{i+2}} \to \dots
	\big)
\end{equation}
and likewise for 2-morphisms. Notice that
$\sshdw{g^\bullet\circ f^\bullet}^i_x = \bigoplus_{p+q=i} \shdw[x]{g^p\circ f^q}$
for a~pair of 1-morphisms $x\from^{g^\bullet} y\from^{f^\bullet} \endofun x$.
We check below that the~sums
\begin{equation}\label{eq:theta-for-complexes}
	\Theta_{g^\bullet,f^\bullet}^i := \sum_{p+q=i} (-1)^{pq} \theta_{g^p,f^q}
\end{equation}
are components of a~formal chain map
$\Theta_{g^\bullet,f^\bullet}\colon
	\sshdw{f^\bullet\circ \endofun g^\bullet}_x \to
	\sshdw{g^\bullet\circ f^\bullet}_y$.
It is invertible if $(\shdw\blank,\theta)$ is a~shadow.

\begin{lemma}\label{lem:lefshetz-shadow}
	The~datum\/ $(\sshdw\blank, \Theta)$ is a~triangulated preshadow
	on $\HCom^b(\ccat C)$. It is a~shadow if so is\/ $\shdw\blank$.
\end{lemma}
\begin{proof}
	The~functors $\sshdw\blank_x$ are additive and, as such,
	they preserve mapping cones. Hence, they are triangulated.
	To check that $\Theta_{g^\bullet,f^\bullet}$ is a~chain map, we compute
	\begin{align*}
		 d^i_{g^\bullet \circ f^\bullet} \circ \Theta^i_{g^\bullet,f^\bullet}
		&= \sum_{\mathclap{p+q=i}}\ (-1)^{pq}
			\shdw{
				d^p_{g^\bullet}\circ\id +
				(-1)^p\id\circ d^q_{f^\bullet}
			} \circ \theta_{g^p,f^q}
		\\
		&= \sum_{\mathclap{p+q=i}}\ (-1)^{pq} \left(
				\theta_{g^{p+1},f^q}\circ \shdw{\id\circ d^p_{\endofun g^\bullet}} +
				(-1)^p \theta_{g^p,f^{q+1}} \circ \shdw{d^q_{f^\bullet}\circ\id}
			\right)
		\\
		&= \sum_{\mathclap{p+q=i+1}}\ (-1)^{pq} \theta_{g^p,f^q} \circ
			\shdw{(-1)^q \id\circ d^p_{\endofun g^\bullet} + d^q_{f^\bullet}\circ\id}
		= \Theta^{i+1}_{g^\bullet,f^\bullet} \circ d^i_{f^\bullet \circ \endofun g^\bullet}
	\end{align*}
	and the~naturality follows from the~naturality of $\theta$:
	\begin{align*}
		\Theta^i_{g_1^\bullet,f_1^\bullet} \circ \sshdw{\alpha^\bullet\circ\endofun \beta^\bullet}^i
		&= \sum_{p+q=i} (-1)^{pq}
				\theta_{g_1^p,f_1^q} \circ \shdw{\alpha^p\circ\endofun\beta^q}
		\\
		&= \sum_{p+q=i} (-1)^{pq}
				\shdw{\beta^q\circ\alpha^p} \circ \theta_{g_0^p,f_0^q}
		= \sshdw{\beta^\bullet\circ\alpha^\bullet}^i \circ \Theta_{g_0^\bullet,f_0^\bullet}^i
	\end{align*}
	for any 2-morphisms $\alpha^\bullet\colon f_0^\bullet \to f_1^\bullet$
	and $\beta^\bullet \colon g_0^\bullet\to g_1^\bullet$. Finally,
	commutativity of \eqref{diag:preshadow-hex} and
	\eqref{diag:preshadow-tri} can be easily checked component-wise.
\end{proof}

%Below we list further properties of Lefschetz preshadows.
Recall that a~triangulated trace is additive with respect to (twisted) endomorphisms
of distinguished triangles. An~analoguous statement holds for Lefschetz preshadows.

\begin{proposition}\label{prop:additivity-of-trian-shadow}
	Let\/ $\ccat C$ have left duals and choose a~commuting diagram in
	$\HCom^b(\ccat C)$
	\begin{equation}\label{eq:dist-tri-in-hTr}
	\begin{diagps}(0em,-0.5ex)(21em,12ex)
		\psset{xunit=7em,yunit=10ex}%
		\node t0(0,1)[p^\bullet\circ f^\bullet]
		\node t1(1,1)[r^\bullet\circ f^\bullet]
		\node t2(2,1)[s^\bullet\circ f^\bullet]
		\node t3(3,1)[(p^\bullet\circ f^\bullet)\mathrlap{[1]}]
		\node b0(0,0)[g^\bullet\circ\endofun  p^\bullet]
		\node b1(1,0)[g^\bullet\circ\endofun  r^\bullet]
		\node b2(2,0)[g^\bullet\circ\endofun  s^\bullet]
		\node b3(3,0)[(g^\bullet\circ\endofun  p^\bullet)\mathrlap{[1]}]
		\arrow{->}[t0`t1;\alpha\circ\id]   \arrow{->}[b0`b1;\id\circ\endofun\alpha]
		\arrow{->}[t1`t2;\beta \circ\id]   \arrow{->}[b1`b2;\id\circ\endofun\beta]
		\arrow{->}[t2`t3;\gamma\circ\id]   \arrow{->}[b2`b3;\id\circ\endofun\gamma]
		\arrow{->}[t0`b0;\pi]
		\arrow{->}[t1`b1;\rho]
		\arrow{->}[t2`b2;\sigma]
		\arrow{->}[t3`b3;{\pi[1]}]
	\end{diagps}
	\end{equation}
	where $f^\bullet\in\HCom^b(\ccat C)(\endofun x, x)$,
	$g^\bullet\in\HCom^b(\ccat C)(\endofun y, y)$, and\/
	$p^\bullet\to^\alpha r^\bullet\to^\beta s^\bullet\to^\gamma p^\bullet[1]$
	is a~distinguished triangle in $\HCom^b(\ccat C)(x,y)$.
	Then $\sshdw{[r,\rho]} = \sshdw{[p,\pi]} + \sshdw{[s,\sigma]}$.
\end{proposition}
\begin{proof}
	The~1-morphism $s^\bullet$ is homotopy equivalent to
	the~mapping cone of $\alpha\colon p^\bullet\to r^\bullet$.
	In particular, $\sigma^n\colon (s^\bullet\circ f^\bullet)^n \to
	(g^\bullet\circ\endofun s^\bullet)^n$ consists of components
	\begin{equation}
	\begin{aligned}
		\sigma_{p,p}&\colon p^{j+1}\circ f^i \to g^{i'}\circ \endofun p^{j'+1} &
		\hskip 1cm
		\sigma_{r,p}&\colon r^j\circ f^i \to g^{i'}\circ \endofun p^{j'+1} \\
		\sigma_{p,r}&\colon p^{j+1}\circ f^i \to g^{i'}\circ \endofun r^{j'} &
		\sigma_{r,r}&\colon r^j\circ f^i \to g^{i'}\circ \endofun r^{j'}
	\end{aligned}
	\end{equation}
	for $i+j = i'+j' = n$.
	A~simple diagram chasing shows that $\sigma_{p,p}$ and $\sigma_{r,r}$
	are chain maps homotopic to $\pi[1]$ and $\rho$ respectively.
	We can thus take $\sigma_{p,p}^\bullet = \pi[1]^\bullet$
	and $\sigma_{r,r}^\bullet = \rho^\bullet$.
	
	Recall that the~left dual to $(C^\bullet, d)$ is given by the~complex
	\begin{equation}
		\dots \to
			\ldual C^{-i-1} \to^{\ldual d^{-i}}
			\ldual C^{-i  } \to^{\ldual d^{-i+1}}
			\ldual C^{-i+1} \to \dots
	\end{equation}	
	with $\ldual C^{-i}$ in homological degree $i$.
	The~evaluation and coevaluation maps do not vanish only in homological
	degree zero, where they are given by the~obvious maps
	\begin{align}
		A\to^{\sum\coev_{C^i}} \bigoplus_i \ldual C^i \otimes C^i
			= (\ldual C^\bullet \otimes C^\bullet)^0
	\intertext{and}
		(C^\bullet \otimes \ldual C^\bullet)^0
			= \bigoplus_i C^i\otimes\ldual C^i \to^{\sum\ev_{C^i}} A.
	\end{align}	
	In particular, $s^j = p^{j+1}\oplus r^j$ leads to the~decomposition
	\begin{equation}
		(\ldual s)^{-j} \circ s^j =
			(\ldual p^{j+1}\circ p^{i+j}) \oplus
			(\ldual r^j \circ r^j)  \oplus
			(\ldual p^{j+1}\circ r^j)  \oplus
			(\ldual r^j \circ p^{j+1}).
	\end{equation}
	The~coevaluation takes values only in the~first two terms,
	which are also the~only components on which the~evaluation does not vanish.
	Hence, the~only components of $\sshdw{[s,\sigma]}$ that contribute are
	\begin{equation}\label{eq:p-term}
		\begin{split}
			\shdw{f^i}
				\to^{\shdw{\coev\circ\iid}}
			\shdw{\ldual p^{j+1}\circ p^{j+1}\circ f^i}
				\to<4em>^{\shdw{\iid\circ\sigma_{p,p}}}&
			\shdw{\ldual p^{j+1} \circ g^i \endofun \circ p^{j+1}} \\
				\to^{(-1)^{(i+j)j} \theta}&
			\shdw{p^{j+1} \circ \ldual p^{j+1} \circ g^i}
				\to^{\shdw{\ev\circ\iid}}
			\shdw{g^i},
		\end{split}
	\end{equation}
	and those with $p^{j+1}$ replaced by $r^j$. Because
	$\sigma_{p,p}^\bullet = \pi[1]^\bullet$ and $\sigma_{r,r}^\bullet
	=\rho^\bullet$, the~compositions \eqref{eq:p-term} (resp.\ those
	with $r^j$ instead of $p^{j+1}$)  coincide with the~components of
	$\sshdw{[p[1],\pi[1]]}$ (resp.\ $\sshdw{[r,\rho]}$).
	Hence, $\sshdw{[s,\sigma]} = \sshdw{[r,\rho]} + \sshdw{[p[1],\pi[1]]}
	= \sshdw{[r,\rho]} - \sshdw{[p,\pi]}$.
\end{proof}

%A~direct consequence is that a~triangulated preshadow $\sshdw\blank$ on
%$\HCom^b(\ccat C)$ restricts to a~Lefschetz trace on $\HCom^b(\ccat C(x,y))$
%for any pair of objects $x,y\in\ccat C$.

\begin{corollary}\label{cor:restr-of-Lef-shadow}
	Assume that\/ $\ccat C$ has left duals and\/ $\endofun$ fixes objects of\/
	$\ccat C$. Then the~Lefschetz preshadow\/ $\sshdw\blank$ on\/ $\HCom^b(\ccat C)$
	induced by\/ $\shdw\blank$ restricts to the~Lefschetz trace on\/
	$\HCom^b(\ccat C(x,y))$ induced by\/ $t^{\shdw\blank}\noic$,
	the~trace defined as in \eqref{def:tr-from-shdw}.
\end{corollary}
\begin{proof}
	The~restricted trace is triangulated due to
	Proposition~\ref{prop:additivity-of-trian-shadow}
	applied to the~case $f=\id_x$ and $g=\id_y$.
	Hence, the~thesis follows from Proposition~\ref{prop:trian-trace-is-Lefschetz}.
\end{proof}

\subsection{Deformation of categorified traces}
\label{sec:quantum-shadows}

A~bicategory $\ccat C$ is \emph{locally pregraded} if each morphism category
$\ccat C(x,y)$ is pregraded, the~horizontal composition preserves the~degree
of 2-morphisms, and the~canonical 2-isomorphisms $\mathfrak a$, $\mathfrak l$,
and $\mathfrak r$ are homogeneous of degree 0.
As usual, we write $|\alpha|$ for the~degree of a~homogeneous 2-morphism $\alpha$.
We say that $\ccat C$ is \emph{locally graded} if each morphism category
$\ccat C(x,y)$ comes with a~translation functor $\{1\}$ that is an~equivalence,
such that for any pair of composable 1-morphism $x\to^f y\to^g z$
the~horizontal compositions $g\circ f\{1\}$, $g\{1\}\circ f$, and $(g\circ f)\{1\}$
are naturally isomorphic.
Likewise for categories, both frameworks are equivalent:
a~locally pregraded bicategory can be formally extended
to a~locally graded one, and vice versa.
In what follows we choose the~pregraded framework.

%Fix a~flat extension $\scalars'\supset\scalars$ with an~invertble element $q$.
Let $(\ccat C,\endofun)$ be a~locally pregraded endobicategory.
This means that $\ccat C$ is locally pregraded and $\endofun$ is
graded, i.e.\ it preserves the~degree of 2-morphisms and
the~structural 2-isomorphisms $\mathfrak m$ and $\mathfrak i$
are homogeneous in degree 0.

\begin{definition}
	A~\emph{$\endofun$--twisted quantum preshadow} on $\ccat C$ valued in
	a~pregraded $\qscalars$--linear category $\cat T$ is a~collection
	of graded functors
	$\shdw[x]\blank \colon \ccat C(\endofun x, x) \to \cat T$
	together with graded morphisms
	$\theta_{g,f}\colon
		\shdw[y]{f\circ \endofun g} \to
		\shdw[x]{g\circ f}$ in $\cat T$,
	one for each pair of 1--morphisms $x\from^g y \from^f \endofun x$,
	which are natural in the~graded sense, i.e.
	\begin{equation}
		\theta_{g',f'}\circ\shdw{\alpha\circ\endofun\beta} = q^{|\beta|}
			\shdw{\beta\circ\alpha} \circ \theta_{g,f}
	\end{equation}
	for homogeneous 2-morphisms $\alpha\colon f\to/=>/ f'$ and $\beta\colon g\to/=>/ g'$,
	and such that the~defining diagrams for preshadows \eqref{diag:preshadow-hex}
	and \eqref{diag:preshadow-tri} commute. We say that $(\shdw\blank, \theta)$
	is a~\emph{quantum shadow} if each $\theta_{f,g}$ is invertible.
\end{definition}

Likewise in the~case of traces, quantum preshadows arise by extending coefficients
of $\ccat C$ to $\qscalars$ and deforming $\endofun$ into $\endofun_q$ by redefining
it on homogeneous 2-morphisms as
\begin{equation}
	\endofun_q(\alpha) := q^{-|\alpha|} \endofun(\alpha).
\end{equation}
Quantum preshadows on $(\ccat C,\endofun)$ are then precisely
preshadows on $(\ccat C,\endofun_q)$.
It follows that there is a~deformed twisted horizontal trace
$\qhTr(\ccat C,\endofun)$, the~\emph{quantum\/ $\endofun$--twisted horizontal
trace} of $\ccat C$, obtained by deforming the~defining relation into
\begin{equation}\label{rel:shTrq-for-morphisms}
	\begin{diagps}(0em,-0.5ex)(6em,12.5ex)
		\square<5em,10ex>{->``->`->}[\endofun x`x`\endofun y`y;f``p`g]
		\carrow*[arcangle= 45]{->}[v0`v2;\endofun p\smash']
		\carrow*[arcangle=-45]|b|{->}[v0`v2;\endofun p]
		\morphism[nodesep=1.6em]{=>}|a{npos=0.53}
			|(5.5em,10ex)<-4em,-9.5ex>[`;\alpha]
		\morphism{<=}(-1em,5ex)<2em,0ex>[`;\endofun \tau]
	\end{diagps}
	\quad\sim\quad
	q^{|\tau|}
	\begin{diagps}(-1.5em,-0.5ex)(5em,12.5ex)
		\square<5em,10ex>{->`->``->}[\endofun x`x`\endofun y`y;f`\endofun p'``g]
		\carrow*[arcangle= 45]{->}[v1`v3;p\smash']
		\carrow*[arcangle=-45]|b|{->}[v1`v3;p]
		\morphism[nodesep=1.6em]{=>}|b{npos=0.53}|
			(4em,10ex)<-4em,-9.5ex>[`;\alpha]
		\morphism{<=}(4em,5ex)<2em,0ex>[`;\tau]
	\end{diagps}
\end{equation}
where $\alpha$ and $\tau$ are homogeneous.
When $\ccat C$ is locally graded, then, in an~analogy to traces of categories,
the~$\qscalars$--linear structure on $\qhTr(\cat C,\endofun)$ satisfies
$q\cdot [p,\alpha] := [p\{1\},\alpha\{1\}]$. Diagrammatically,
\begin{equation}
	q\cdot \begin{diagps}(-1em,-0.5ex)(6.5em,12.5ex)
		\square<5em,10ex>|a`*c`*c`a|[%
			\endofun x`x`\endofun y`y;
			f`\endofun p`p`g
		]
		\arrow[nodesep=1.5em]|b|{=>}[v1`v2;\alpha]
	\end{diagps}
	:=
	\begin{diagps}(-2em,-0.5ex)(5em,12.5ex)
		\square<5em,10ex>|a`*c`*c`a|[%
			\endofun x`x`\endofun y`y;
			f`\endofun p\{1\}`p\{1\}`g
		]
		\arrow[nodesep=1.5em]|b{npos=0.4}|{=>}[v1`v2;\alpha\{1\}]
	\end{diagps}
\end{equation}
The~following is a~direct generalization from the~undeformed framework.
In particular, every $\endofun$--twisted quantum preshadow factorizes
through $\qhTr(\ccat C,\endofun)$ uniquely up to an~equivalence.

\begin{proposition}
	Assume that\/ $\ccat C$ has left duals. Then\/ $\qhTr(\ccat C,\endofun)$
	with $\shdw[x]\blank^{\mathrm h}$ the~obvious functors 
	is the~universal\/ $\endofun$--twisted quantum preshadow.
	It is a~quantum shadow if\/ $\ccat C$ has right duals.
\end{proposition}

Locally pregraded endobicategories form a~3-category % $\cat{preEndoBicat}$
if we require all bifunctor to preserve degrees of 2-morphisms
and all structural 2-morphisms (such as $\omega_f$ in the~definition of
a~morphism or $\eta_f$ and $\mathfrak n_x$ in the~definition of a~natural
transformation) to be graded. A~pullback of a~twisted quantum preshadow
along a~morphism of locally pregraded endobicategories is again
a~twisted quantum preshadow, and natural transformations of such
morphisms induce morphisms of preshadows. The~following is
an~immediate consequence of that and the~universality of $\qhTr$.

\begin{theorem}[Functoriality of the~quantum horizontal trace]
	The~quantum horizontal trace $\qhTr$
	extends to a~strict 2-functor
%	$\qhTr\colon \cat{preEndoBicat}^{*,\leqslant 2} \to \cat{preCat}$
	from the~2-category of locally pregraded endobicategories
	with duals to pregraded categories.
\end{theorem}

We omit the~explicit construction of the~2-functor,
because it is a~verbatim copy of the~one from Theorem~\ref{thm:hTr-functoriality}.

%% ============================================================================
%%  Section: examples of traces
%% ============================================================================

\subsection{Examples}
\label{sec:examples}

\subsubsection{The~twisted trace}
\label{sec:shadow-Tr}

Choose a~pair of functors $\cat C\from^F \cat D\from^G \cat C$.
The~pair $(G,\id)$ is a~morphism of endocategories, from $(\cat C, FG)$
to $(\cat D, GF)$, and as such it induces a~map between their universal traces
\[
	\Tr(G,\id) \colon \Tr(\cat C,FG) \to \Tr(\cat D, GF),
\qquad
	\tr(f) \mapsto \tr(Gf).
\]
To show naturality, choose natural transformations
$\eta\colon F\to F'$ and $\mu\colon G\to G'$. They determine morphisms
of endocategories
\begin{align*}
	(\Id,\mu\circ F\eta)&\colon (\cat C, FG) \to (\cat C, F'G')
	\\
	(\Id,\eta\circ G\mu)&\colon (\cat D, GF)\to (\cat D, G'F')
\end{align*}
and hence maps on traces. Using $\mu\circ F\eta = F'\eta\circ\mu$
we check directly that
\begin{multline*}
	(\Tr(G',\Id)\circ \Tr(\Id,\mu\circ F\eta))(\tr(f)) =
		\tr(G'F'\eta\circ G'\mu\circ G'f) = \\
		\tr(G'(\mu\circ f)\circ \eta) =
		\tr(\eta\circ G(\mu\circ f)) =
	(\Tr(\Id,\eta\circ G\mu) \circ \Tr(G,\id))(\tr(f)).
\end{multline*}
Because $\tr(FGf) = \tr(f)$ by Lemma~\ref{lem:sigma-inv-of-tr},
the~above upgrades $\Tr(\blank)$ to a~symmetric shadow on $\ccat{Cat}$,
the~2-category of small categories. The~same argument shows that $\qTr(\blank)$
is a~quantum shadow on $\ccat{gCat}$, the~2-category of small graded categories,
with the~cyclicity map satisfying $\tr(FGf) = \tr(FGf\circ\id) = q^{|f|}\tr(f)$
for a~homogeneous $f$.

Although not every functor has a~left dual (left dualizable functors are
exactly those with a~left adjoint), $\Tr(\blank)$ factors through
$\hTr(\ccat{Cat})$ (resp.\ $\qTr(\blank)$ factors through $\qhTr(\ccat{gCat})$).
Indeed, a~representant $(F,\omega)$ of a~morphism in the~horizontal trace
is precisely a~morphism of endocategories, and it can be checked directly
that $\Tr(F,\omega)$ (resp.\ $\qTr(F,\omega)$) depends only on the~image of
$(F,\omega)$ in the~horizontal trace; see also Theorem~\ref{thm:HH-on-hTr(Cat)}.

\subsubsection{The~space of coinvariants}
\label{sec:shadow-coinv}

Write $\ccat{Bimod}$ for the~bicategory of $\scalars$--algebras,
bimodules, and bimodule maps; the~horizontal composition of bimodules
$N\in\ccat{Bimod}(A,B)$ and $N'\in\ccat{Bimod}(B,C)$ is defined as
their tensor product: $N'\circ N := N\utimes{B} N'$. It admits a~symmetric shadow
that assigns to an~$(A,A)$--bimodule $M$ the~space of coinvariants
$\coinv(M) := M/[A,M]$ and the~cyclicity morphism of which is the~standard
twist: $\theta_{N',N}([n'\otimes n]) = [n\otimes n']$,
where $N\in\ccat{Bimod}(A,B)$ and $N'\in\ccat{Bimod}(B,A)$.

Likewise, the~quantum space of coinvariants $\coinv_q(\blank)$ defined
in \eqref{def:coinv-q} determines a~shadow on $\gBBimod$, the~locally graded
bicategory of graded algebras, graded bimodules, and graded bimodule maps.
The~cyclicity morphism is a~noninvolutive deformation of
the~standard twist: $\theta_{N',N}([n'\otimes n]) = q^{|n'|}[n\otimes n']$
for homogeneous $n' \in N'$ and $n\in N$.

Let $\ccat{gRep}$ be the~restriction of $\ccat{gBimod}$ to those bimodules
that are finitely generated and projective as right modules. It can be realized
as a~sub-2-category of $\ccat{gCat}$ by interpreting an~algebra $A$ as its
representation category $\gRep(A)$, and an~$(A,B)$--bimodule $N$ as the~functor
$(\blank)\utimes{A} N$. Therefore, it admits a~preshadow as described in
the~previous example, which factorizes through $\qhTr(\ccat{gRep})$.
Explicitly, it assigns to
\begin{equation}
	\begin{diagps}(0em,-0.5ex)(5em,12.5ex)
		\square<5em,10ex>[A`A`B`B;M`P`P`M']
		\arrow[nodesep=1.65em]|b|{=>}[v1`v2;f]
	\end{diagps}
\end{equation}
the~linear map $[m]\mapsto \tr_q^P\!\left(p\mapsto f(m\otimes p)\right)$,
where $\tr_q^P$ is the~quantum deformation of the~twisted Hattori--Stallings
trace from Section \ref{sec:ex-trace-on-Mod}.
We claim that this shadow coincides with the~one restricted from $\ccat{Bimod}$.
Indeed, given $N\in\ccat{gRep}(A,B)$ and $N'\in\ccat{gRep}(B,A)$ with
homogeneous $n\in N$ and $n'\in N'$, the~cyclicity map takes the~class
of $[n'\otimes n]$ to
\begin{multline*}
	\tr_q^{N'}\!(n'\otimes n\otimes\id_{N'}) =
		\epsilon_{q,N'}(n'\otimes n\otimes\coev_{\!N'}(1)) = \\
		= q^{|n'|}[n\otimes(\id_{N'}\otimes\ev_{\!N'})(\coev_{\!N'}(1)\otimes n')]
		= q^{|n'|}[n\otimes n']
		= \theta_{N',N}([n'\otimes n]).
\end{multline*}

%Indeed, the~isomorphism $\Tr(\Rep(A),M) \cong \coinv(M)$ is shown in
%Section~\ref{sec:example-twisted-trace}, whereas the~cyclicity map $\theta_{N',N}$
%takes $[n'\otimes n]$ to the~twisted Hattori--Stallings trace of
%\[
%	N' \ni x \mapsto n'\otimes n\otimes x \in N'\utimes{A} (N\utimes{B} N'),
%\]
%which is exactly $[n\otimes n']$.

\subsubsection{Twisted coinvariants}
\label{sec:shadow-twisted-coinv}

Let $\eBBimod$ be an~extension of $\BBimod$, objects of which are pairs
$(A,\varphi)$ consisting of an~algebra $A$ and its automorphism $\varphi$,
and whose 1- and 2-morphisms are bimodules and bimodule maps respectively
with no restriction imposed. It can be upgraded to an~endobicategory via
an~endobifunctor $\endofun$ that affects only 1-morphisms:
\[
	\endofun N := {}_{\overline\varphi} N_{\overline\psi}
		\cong A_\varphi \utimes{A} N \utimes{B} {}_\psi B
\]
for $N\in\eBBimod(A,\varphi; B,\psi)$,
where $A_\varphi$ is the~$(A,A)$--bimodule $A$ with the~right action
$x\cdot a := x\varphi(a)$, and likewise for ${}_\psi B$.
For any object $(A,\varphi)$ there is an~isomorphism of bimodules
$A\cong {}_\varphi A_\varphi$ that takes $a$ to $\varphi(a)$. Hence,
$\endofun(N\utimes{B} N') \cong \endofun N\utimes{A} \endofun N'$ as desired.

Twisted spaces of coinvariants from Section~\ref{sec:example-twisted-trace}
constitute a~shadow on $\eBBimod$ with the~cyclicity map
\begin{multline}\label{eq:twisted-coinv-shadow}
	\coinv^\psi(\endofun N'\utimes{A} N) =
		\coinv({}_\psi B_\psi\utimes{B} N'\utimes{A} {}_\varphi A\utimes{A} N) \cong\\
		\coinv({}_\varphi A\utimes{A} N\utimes{B} N') =
		\coinv^\varphi(N\utimes{B} N')
\end{multline}
for $N\in\eBBimod(A,\varphi;B,\psi)$ and $N'\in\eBBimod(B,\psi;A,\varphi)$.
Recalling that $\endofun N'$ and $N'$ coincide as $\scalars$--modules,
we find out that \eqref{eq:twisted-coinv-shadow} is the~twist
$[n'\otimes n] \mapsto [n\otimes n']$.

Let $\eRRep\subset\eBBimod$ be the~restriction to bimodules that are finitely
generated and projective as right modules. Such an~$(A,B)$--bimodule $N$ has
a~left dual $\ldual N := \Hom_B(N,B)$, the~$(B,A)$--bimodule of right $B$--linear
functions on $N$.\footnote{
	Notice that $\ldual N$ is projective as a~left module, but not necessarily
	as a~right one. Hence, in general, it is not a~1-morphism in $\eRRep$.
}
Thus the~above shadow can be restricted to a~trace on
$\vTr(\eRRep(A,\varphi; B,\psi), \endofun)$ that takes
$f\in\Hom(N,\endofun N)$ to the~composition
\begin{equation}
\begin{aligned}
	\coinv^\varphi(A) \to^{\coev} \coinv^\varphi(N\utimes{B}\ldual N)
		\to<2.5em>^{f\otimes\id} &\coinv^\psi(\endofun N\utimes{B}\ldual N)\\
		\to<2.5em>^{\theta} &\coinv^\psi(\ldual N\utimes{A}N)
		\to^{\ev} \coinv^\psi(B),
\end{aligned}
\end{equation}
which, when evaluated at $[a]\in\coinv^\varphi(A)$, recovers the~$\psi$--twisted
Hattori--Stallings trace of $f_a(n) := f(an)$.

\subsubsection{Hochschild homology}

Let $M$ be an~$(A,A)$--bimodule and write $R_\bullet(A)$ for
the~\emph{bar resolution} of $A$ \cite[1.1.11]{Loday}.
The~\emph{Hochschild homology of\/ $A\noic$ with coefficients in $M\noic$} or,
shortly, the~\emph{Hochschild homology of\/ $M$} is the~homology of the~chain complex
$\CHoHom_{\!\bullet}(A,M) := \coinv(M\utimes{A} R_\bullet(A))$, to which we refer
as the~\emph{Hochschild complex}. Explicitly,
$\CHoHom_{\!n}(A, M) = M\otimes A^{\otimes n}$ and
the~Hochschild differential is the~alternating sum
\begin{equation}
	\begin{split}
	\partial(m\otimes a_1\otimes\dots\otimes a_n) = 
		ma_1\otimes a_2\otimes&\dots\otimes a_n \\
		+ \sum_{i=1}^{n-1} (-1)^i m\otimes a_1\otimes &\dots
					\otimes a_ia_{i+1} \otimes \dots \otimes a_n \\
		+ (-1)^n a_nm\otimes a_1\otimes&\dots\otimes a_{n-1}.
	\end{split}
\end{equation}
In particular, $\HoHom_{\!0}(A,M) \cong \coinv(M)$ is the~ordinary
space of coinvariants.

Hochschild homology is not directly a~shadow on $\BBimod$: for that
one has to extend it to derived categories of bimodules and use
the~derived tensor product as the~composition. It is not necessary
in the~case of $\RRep$, because an~$(A,A)$--bimodule of the~form $N\utimes{B} N'$
has two projective resolutions
\begin{gather}
	P_\bullet(N,N'):=R_\bullet(A) \utimes{A} (N\utimes B N')
\intertext{and}
	R_\bullet(N,N'):=\Tot(R_\bullet(A)\utimes{A}N\utimes{B}R_\bullet(B)\utimes{B}N'),
\end{gather}
which, by the~uniqueness of the~resolution, are homotopy equivalent.
Therefore, there is a~sequence of homotopy equivalences
\begin{equation}\label{eq:htpy-eq-for-HH}
	\CHoHom_{\!\bullet}(B, N'\utimes A N) \simeq 
		\coinv(R_\bullet(N',N)) \cong
		\coinv(R_\bullet(N,N')) \simeq
	\CHoHom_{\!\bullet}(A, N\utimes B N'),
\end{equation}
where the~middle map is the~standard twist
$[\underline b\otimes n'\otimes \underline a\otimes n] \mapsto
	[\underline a\otimes n\otimes \underline b\otimes n']$,
in which $\underline a$ stands for a~tensor product of elements of $A$,
and likewise for $\underline b$. Doing this twice results in a~map
chain homotopic to the~identity, so that $\CHoHom_{\!\bullet}$
is a~symmetric shadow on $\RRep$.

Let $\varphi\in\Aut(A)$. The~\emph{$\varphi$--twisted Hochschild complex}
$\CHoHom_{\!\bullet}^\varphi(A, M)$ of an~$(A,A)$--bimo\-dule $M$ is obtained
by replacing $M$ with ${}_\varphi M$
(or, equivalently, $\coinv$ with $\coinv^\varphi)$:
\[
	\CHoHom_{\!\bullet}^\varphi(A, M) :=
		\coinv({}_\varphi M \utimes{A} R_\bullet(A)).
\]
It coincides with $\CHoHom_{\!\bullet}(A, M)$ except that the~last
term in the~differential is replaced with
\begin{equation}
	(-1)^n \varphi(a_n)m\otimes a_1\otimes\dots\otimes a_{n-1}.
\end{equation}
Hochschild chains can be visualized by oriented circles divided into segments,
one labeled with $m\in M$ and the~others with $a_0,\dots,a_n$. Each of the~terms
of the~differential merges two segments multiplying their labels:
\[
		\HochschildChain{0,120,240}{m@180,a_0@-60,a_1@60}
	\longmapsto
		\HochschildChain{0,120}{ma_0@240,a_1@60}
		-
		\HochschildChain{120,240}{m@180,a_0a_1@0}
		+
		\HochschildChain{0,240}{a_1m@120,a_0@-60}
\]
In the~twisted case add a~mark on the~circle between segments labeled $m$ and
$a_n$. To merge these two segments, one has to move the~mark over $a_n$, acting
upon it with $\varphi$:
\[
		\HochschildChain[120]{0,120,240}{m@180,a_0@-60,a_1@60}
	\mapsto\quad
		\HochschildChain(-2.5em,-1.75em)(3.5em,1.75em)[0]{0,120,240}{m@180,a_0@-60,\varphi(a_1)@60}
	\mapsto\quad
		\HochschildChain(-2em,-1.75em)(3.5em,1.75em)[0]{0,240}{\varphi(a_1)m@120,a_0@-60}
\]
Defined \emph{a priori} for the~untwisted case, the~cyclicity map
\eqref{eq:htpy-eq-for-HH} can be easily extended to the~twisted homology
following \eqref{eq:twisted-coinv-shadow}:
\begin{multline*}
	\CHoHom_{\!\bullet}^\psi(B, \endofun N' \utimes{A} N) =
	\CHoHom_{\!\bullet}(B, {}_\psi B_\psi\utimes{B} N' \utimes{A} {}_\varphi N)
\simeq \\
	\CHoHom_{\!\bullet}(A, {}_\varphi N \utimes{B} {}_\psi B_\psi\utimes{B} N')
\cong
	\CHoHom_{\!\bullet}(A, {}_\varphi N \utimes{B} N') =
	\CHoHom_{\!\bullet}^\varphi(A, N \utimes{B} N'),
\end{multline*}
which upgrades the~twisted homology to a~symmetric shadow on $\eRRep$.

\subsubsection{Quantum Hochschild homology}
\label{sec:ex-qHH}

Fix an~invertible elements $q\in\scalars$.
We define the~\emph{quantum Hochschild complex} of a~graded $(A,A)$--bimodule $M$
by replacing $\coinv(M)$ with its deformation:
\begin{equation}
	\qCHoHom_{\!\bullet}(A, M) := \coinv_q(M\utimes{A} R_\bullet(A)).
\end{equation}
It has the~same chain groups as $\CHoHom_{\!\bullet}(A, M)$ with
coefficients extended to $\scalars'$, but the~last term of the~differential
is replaced with
\begin{equation}
	(-1)^n q^{-|a_n|}a_nm\otimes a_1\otimes\dots\otimes a_{n-1}.
\end{equation}
The~\emph{quantum Hochschild homology of $M\noic$},
denoted by $\qHoHom_{\!\bullet}(A, M)$, is the~homology of this complex.

The~quantum Hochschild complex is a~nonsymmetric shadow on $\gRRep$.
Indeed, replacing the middle isomorphism
in \eqref{eq:htpy-eq-for-HH} with 
$[\underline b\otimes n'\otimes \underline a\otimes n]
	\mapsto q^{|\underline b\otimes n'|}
		[\underline a\otimes n\otimes\underline b\otimes n']$
results in a~homotopy equivalence
\begin{equation}
	\theta_{N',N} \colon \qCHoHom_{\!\bullet}(B, N' \utimes A N)
		\to^\simeq \qCHoHom_{\!\bullet}(A, N\utimes B N')
\end{equation}
that is no longer involutive: $\theta_{N,N'}\circ\theta_{N',N}$
scales a~chain of degree $d$ by $q^d$.
The~twisted Hochschild homology can be deformed likewise,
leading to the~\emph{twisted quantum Hochschild complex}
$\qCHoHom_{\!\bullet}^\varphi(A, M)$ of a~graded $(A,A)$--bimodule M
and its homology $\qHoHom_{\!\bullet}^\varphi(A, M)$.
The~details are left to the~reader.

\subsubsection{The~functor \texorpdfstring{$\mathit{HHH}$}{HHH}}
\label{sec:shadow-HHH}

The~construction of Hochschild homology extends naturally to complexes
of bimodules, by letting
$\CHoHom_{\!\bullet}(A,C^\bullet) := \coinv(C^\bullet\utimes{A}R_\bullet(A))$
for a~complex of $(A,A)$--bimodules. It can be checked that it is a~triangulated
shadow on $\HCom^b(\RRep)$.

Instead of computing the~total Hochschild homology of the~complex $C^\bullet$
one can apply $\HoHom_{\!\bullet}$ component-wise. This leads to a~Lefschetz
shadow $\mathit{HHH}$ on $\HCom^b(\ccat{Rep})$ that replaces $C^\bullet$ with 
\begin{equation}\label{eq:HHH-complex}
	\dots\to \HoHom_{\!\bullet}(A, C^i)
	     \to \HoHom_{\!\bullet}(A, C^{i+1})
	     \to \HoHom_{\!\bullet}(A, C^{i+2})
	     \to \dots
\end{equation}
It is the~second page of the~spectral sequence associated to
$\coinv(C^\bullet\utimes{A}R_\bullet(A))$ seen as a~bicomplex. 
This functor appears quite often in categorification of link invariants,
when links are considered as closures of braids.% \comment{references}

Naturally, when dealing with complexes of graded bimodules, one can replace
$\HoHom_{\!\bullet}$ with $\qHoHom_{\!\bullet}$ to obtain $\mathit{qHHH}\!$,
which is a~quantum shadow.
According to Corollary~\ref{cor:restr-of-Lef-shadow}, it restricts to a~quantum
trace on $\HCom^b(\ccat{gRep})(A,B)$ for any graded algebras $A$ and $B$.
It takes a~particularly simple form when $A=B=\scalars$,
because then $\qHoHom_{\!0}(\scalars) = \scalars$
and higher homology vanishes.

\begin{corollary}\label{cor:HH-induces-Lef}
	Let\/ $C^\bullet \in \HCom^b(\gRRep(\scalars,\scalars))$ be a~bounded
	complex and\/ $f^\bullet\in\End(C^\bullet)$ its graded endomorphism.
	The~endomorphism of\/ $\qHoHom(\scalars) \cong \scalars\noic$
	induced by $(C^\bullet, f^\bullet)$ is	the~multiplication by
	$\Lef_q(f^\bullet) = \sum_{i,j} (-1)^i q^j \tr(f^{i,j})$,
	where $i$ and $j$ are the~homological and internal gradings in\/ $C^\bullet$ respectively, and\/ $\tr$ is the~Hattori--Stallings trace on $\Rep(\scalars)$.
\end{corollary}

\subsubsection{Cobordisms of links in surface bundles}
\label{sec:links-in-sfce-bundle}

It is generally a~hard problem to identify the~horizontal trace
$\hTr(\ccat C,\endofun)$ for a~given bicategory $\ccat C$,
but the~answer is very natural for $\TTan(\sfce)$,
the~bicategory of points in a~surface $\sfce\!$,
tangles in $\sfce\times I$, and tangle cobordisms in $\sfce\times I\times I$.
Recall that the~quotient space
\begin{equation}
	M := \quotient{\sfce\times I}{(p,1) \sim (\varphi(p),0)},
\end{equation}
where $\varphi$ is a~diffeomorphism of $\sfce\!$, is a~manifold called
the~\emph{surface bundle with fiber $\sfce\noic$ and monodromy $\varphi$}.

\begin{figure}
	\centering
	\psset{unit=3ex}%
	\begin{pspicture}(-4,-1.5)(4,5)
		\newrgbcolor{cylinderLighter}{0.9 0.9 1}
		\newrgbcolor{cylinderLight}{0.85 0.85 1}
		\newrgbcolor{cylinderNormal}{0.75 0.75 1}
		\newrgbcolor{cylinderDark}{0.6 0.6 1}
		%
		% background colors
		%
		\pscustom[linestyle=none,fillstyle=solid,fillcolor=cylinderLight]{%
			\moveto(4,0)
			\psellipticarc(0,0)(4,1){0}{180}
			\lineto(-4,4)
			\psellipticarc(0,4)(4,1){-180}{0}
			\closepath
		}
		\psellipse[linestyle=none,fillstyle=solid,fillcolor=cylinderLighter](0,0)(4,1)
		%
		% the cylinder
		%
		\psset{linewidth=0.5pt,arrowsize=3pt,arrowinset=0.5,arrowlength=1.2}
			\psellipse[fillstyle=solid,fillcolor=cylinderNormal](0,4)(4,1)
			\psellipticarc(0,0)(4,1){-180}{0}
			\psellipticarc[linewidth=0.3pt,linestyle=dashed,dash=3pt 3pt](0,0)(4,1){0}{180}
			\psline(-4,0)(-4,4)
			\psline( 4,0)( 4,4)
			\psline[linestyle=dotted,linewidth=1pt](1.6905,-0.9063)(1.6905,3.0937)
		%
		% arrows visuzalizing the deformation
		%
			\psellipticarc{->}(0,3.4)(4,1){215}{280}
			\psellipticarc{->}(0,3.0)(4,1){215}{325}
			\psellipticarc{->}(0,2.6)(4,1){215}{325}
			\psellipticarc{->}(0,2.2)(4,1){215}{325}
			\psellipticarc{->}(0,1.8)(4,1){215}{325}
			\psellipticarc{->}(0,1.4)(4,1){215}{325}
			\psellipticarc{->}(0,1.0)(4,1){215}{325}
			\psellipticarc{->}(0,0.6)(4,1){215}{260}
		%
		% projections of the membranes
		%
		\psset{linewidth=1.5pt}
			\psline(-1.6905,-0.9063)(-1.6905,3.0937)
			\psbezier{-c}(-1.6905,-0.9063)(-1.6905,-0.4063)(1.6905,-0.9063)(1.6905,0.0937)
			\psline{c-c}(1.6905,0.0937)(1.6905,2.0937)
			\psbezier{c-}(1.6905,2.0937)(1.6905,2.9937)(-1.6905,2.5937)(-1.6905,3.0937)
		%
		% labels of the membranes
		%
		\rput[t](-1.7,-1.2){$\sfce_0$}
		\rput[t]( 1.7,-1.2){$\sfce\smash'$}
	\end{pspicture}
	\caption{A~visualization of the~isotopy pushing $\sfce_0\times I$ onto $\sfce'\times I$. Each point of the~cylinder represents a~fiber of the~$\sfce$--bundle $M\times I \protect\to \S\times I$. The~thick straight line is the~standard membrane $\sfce_0\times I$, whereas the~curve is its isotopic deformation. The~flat part of the~curve corresponds to the~piece of the~deformed membrane contained in $\sfce'\times I$.}
	\label{fig:isot-of-membrane}
\end{figure}
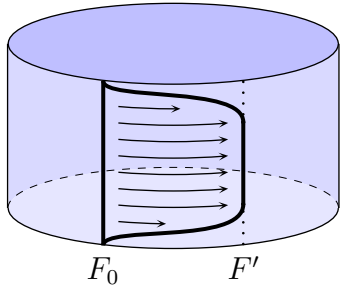

\begin{itheorem}\label{thm:hTr-of-FxI}
	Let\/ $M\noic$ be a~surface bundle with fiber\/ $\sfce\noic$ and monodromy
	$\phi\in\mathit{Diff}(\sfce)$. There is an~equivalence of categories
	\begin{equation}
		\hTr(\TTan(\sfce),\phi_*) \simeq \Links(M)
	\end{equation}
	where $\phi_*(S) := (\phi\times\id\times\id)(S)$
	for a~cobordism\/ $S\subset\sfce\times I\times I\noic$.
\end{itheorem}
\begin{proof}
	Let $\pi\colon\sfce\times I \to M$ be the~quotient map 	and consider
	the~fiber $\sfce_0 := \pi(\sfce\times\{0\})$ along which $M$ can be
	cut open to $\sfce\times I$. Objects of $\hTr(\TTan(\sfce),\phi_*)$
	can be identified with links in $M$ transverse to $\sfce_0$,
	while morphisms are represented by link cobordisms in $M\times I$
	transverse to the~3-dimensional \emph{membrane} $\sfce_0\times I$.
	The~cobordism can be deformed by an~ambient isotopy that fixes
	the~membrane, and the~trace relation allows us to isotope the~embedding
	of the~membrane (although $\sfce_0\times\{i\}$ is fixed for $i=0,1$).
	
	There is an~obvious functor $\hTr(\TTan(\sfce),\phi_*)\to \Links(M)$
	that forgets the~membrane. By the~transversality argument it is essentially
	surjective on objects (each link is isotopic to a~link transverse to $\sfce_0$)
	and full on morphisms (each surface between links transverse to $\sfce_0$
	can be isotoped to be transverse to the~standard membrane $\sfce_0\times I$).
	It remains to show that if two surfaces $S,S'\subset \sfce\times I\times I$
	represent isotopic cobordisms $\widehat S$ and $\widehat S'$ in $M\times I$,
	then their images in $\hTr(\TTan(\sfce),\phi_*)$ coincide.

	Assume there is an~isotopy $\varphi_t$ of $M\times I$ taking $\widehat S$
	to $\widehat S'$ with support\footnote{
		A~\emph{support} of an~isotopy is the~closure of the~set of points
		that are not stationary under the~isotopy.
	}
	$\mathrm{supp}(\varphi)$ disjoint from $M\times \partial I$ and a~membrane
	$\sfce'\times I$ for some fiber $\sfce'\subset M$. It is enough to consider
	only such isotopies, because every two isotopic surfaces in $M\times I$ are
	connected by a~sequence of them. If $\sfce'=\sfce_0$, then $S$ and $S'$ are
	already isotopic in $\sfce\times I\times I$ and we are done. Otherwise, let
	$p\in\S$ be the~point over which $\sfce'$ lives and consider a~bump function
	$\beta\colon I\to \S$ with $\beta(0) = \beta(1) = 1$ and $\beta(t) = p$ for
	$t\in[\epsilon,1-\epsilon]$ for some $\epsilon>0$ such that
	$\mathrm{supp}(\varphi) \subset M\times [\epsilon,1-\epsilon]$,
	see Fig.~\ref{fig:isot-of-membrane}. The~preimage in $M\times I$ of the~graph of $\beta$
	is a~membrane isotopic to $\sfce_0\times I$; it can be visualized as pushing
	the~interior of $\sfce_0\times I$ onto $\sfce'\times I$. Because the~new membrane
	is disjoint from the~support of $\varphi$, the~cuts of $\widehat S$ and
	$\widehat S'$ along it are isotopic. This proves the~faithfulness,
	because the~cuts represent in $\hTr(\TTan(\sfce),\phi_*)$ the~same morphisms
	as the~surfaces $S$ and $S'$.
\end{proof}

Each orientation preserving diffeomorphism $\phi$ of $\R^2$ is isotopic to identity, so that $\hTr(\TTan,\phi_*)$ is equivalent to the~category of links in a~solid torus. Because tangles in a~thickened plane $\R^2\times I$ can be represented by diagrams on the~stripe $\{0\}\times\R\times I$, is it worth to consider those diffeomorphisms that preserve the~line $\{0\}\times\R$. There are two of them:
\begin{itemize}
	\item the~identity, in which case the~stripe is closed to an~annulus $\Ann$, and
	\item the~rotation by $180$ degrees, for which the~image of the~stripe is a~M\"obius band $\MB$.
\end{itemize}
The~solid torus is a~trivial line bundle over $\Ann$ and a~twisted one over $\MB$ respectively. Hence, given an~invariant of tangles computed from their diagrams, there are two ways to get invariants of links in a~solid torus.

\begin{corollary}\label{cor:annular-invariants}
	There~are equivalences of categories
	$$
		\hTr(\TTan) \simeq \Links(\S\times\R^2)  \simeq \hTr(\TTan,\rho_*),
	$$
	where $\rho\in\mathit{Diff}(\R^2)$ is the~half--rotation.
	Thence, a~bifunctor\/ $\ccat I\colon\TTan\to \ccat C$ induces invariants of links in a~solid torus
	\begin{align*}
		\hTr(\ccat I)\colon& \Links(\R^2\times \S) \to \hTr(\ccat C), \text{ and}\\
		\hTr(\ccat I,\rho_*)\colon&\Links(\R^2\times \S) \to \hTr(\ccat C,\endofun),
	\end{align*}
	where $\endofun$ is an~endofunctor of\/ $\ccat C$ satisfying $\endofun\circ\ccat I \cong \ccat I\circ\rho_*$.
	In particular, $\shdw{\hTr(\ccat I)}$ is a~link invariant for any symmetric preshadow $\shdw\blank$ on\/ $\ccat C$.
\end{corollary}

\section{Quantum Hochschild homology}
\label{sec:qHH}

Here we develop a~machinery to show that higher quantum Hochschild homology
of Chen--Khovanov algebras vanish, which is used to proof the~conjecture
of Auroux, Grigsby, and Wehrli. In particular, we prove the~invariance
of quantum Hochschild homology under $K\noic$--theoretic equivalences.
This section can be skipped by a~reader interested only in the~construction
of the~quantum link homology.
Hereafter we fix an~invertible $q\in\scalars$.

%% ============================================================================
%%  Section: the cyclic set
%% ============================================================================

\subsection{Twisted Hochschild--Mitchell homology}

Choose a~small pregraded endocategory $(\cat C, \endofun)$.
The~\emph{(twisted) quantum Hochschild--Mitchell complex} of $(\cat C,\endofun)$
is the~chain complex $\qCHMHom_{\!\bullet}(\cat C,\endofun)$ with components
\begin{equation}\label{def:HM-chains}
	\qCHMHom_{\!n}(\cat C, \Sigma) :=
		\bigoplus_{\quad\mathclap{x_0,\dots,x_n\in\Ob(\cat C)}\quad}
			\cat C(x_0,\endofun x_n) \otimes
			\cat C(x_1,x_0) \otimes \dots \otimes
			\cat C(x_n,x_{n-1})
\end{equation}
and differential the~alternating sum $\partial_n = \sum_{i=0}^n (-1)^i d_n^i$, where
\begin{align}
	\label{def:cat-HH-face-map}
		d_n^i(\tensor{f0,*\dots,fn}) &:= \begin{cases}
			\tensor{f0,*\dots,*(f_i\circ f_{i+1}),*\dots,fn}
				& \textrm{ if }i<n,\\
			q^{-|f_n|}\tensor{*(\endofun f_n\circ f_0),f1,*\dots,f{n-1}}
				& \textrm{ if }i=n.
		\end{cases}
\end{align}
Its homology is called the~\emph{(twisted) quantum Hochschild--Mitchell homology
of $(\cat C,\endofun)$} and denoted by $\qHMHom_{\!\bullet}(\cat C, \endofun)$.
We also use the~symbols $\CHMHom_{\!\bullet}$ and $\HMHom_{\!\bullet}$ when
$q=1$, dropping the~adjective `quantum'. In particular,
$\HMHom_{\!\bullet}(\cat C,\Id)$ recovers the~usual Hochschild--Mitchell
homology of a~category \cite{HochMitchell}.
Furthermore, $\HMHom_{\!0}(\cat C,\endofun)$ and $\qHMHom_{\!0}(\cat C,\endofun)$
are precisely the~universal $\endofun$--twisted trace of $\cat C$ and its
quantum deformation.

%
% A note on homotopies
%
Expressing the~differential as an~alternating sum makes it evident
that the~Hochschild--Mitchell complex arises from a~presimplicial module
\cite[1.0]{Loday}. Despite not working in this framework, we will often
use the~following characterisation of chain homotopies.

\begin{lemma}\label{lem:simpl-htpy}
	Choose a~collection 
	$\{	h_n^k\colon \qCHMHom_{\!n}(\cat C,\endofun)
				\to \qCHMHom_{\!n+1}(\cat C', \endofun')
		\}_{0\leqslant k\leqslant n}$
	of linear maps, satisfying
	\begin{equation}\label{def:simpl-htpy}
		d_{n+1}^i \circ h_n^k =
		\begin{cases}
			h_{n-1}^{k-1}\circ d_n^i &\quad\textrm{for }i<k,\\
			d_{n+1}^k\circ h_n^k     &\quad\textrm{for }i=k+1,\\
			h_{n-1}^k\circ d_n^{i-1} &\quad\textrm{for }i>k+1.
		\end{cases}
	\end{equation}
	Then the~alternating sums\/ $h_n := \sum_{k=0}^n (-1)^k h_n^k$
	are components of a~chain homotopy from $f_n:=d_{n+1}^0\circ h_n^0$ to\/
	$g_n:=d_{n+1}^{n+1}\circ h_n^n$.
\end{lemma}

%
% Simplicial maps from functors
%
Given a~morphism of endocategories
$(F,\omega)\colon (\cat C,\endofun)\to (\cat C',\endofun')$
define a~chain map
$(F,\omega)_*\colon \qCHMHom_{\!\bullet}(\cat C,\endofun)
	\to \qCHMHom_{\!\bullet}(\cat C',\endofun')$
by the~formula
\begin{equation}\label{def:HH-for-functor}
	(F,\omega)_*(\tensor{f0,*\dots,fn}) :=
		\tensor{*(\omega\circ F\!f_0),{F\!f}1,*\dots,{F\!f}n},
\end{equation}
where $\omega$ is used to fix the~codomain of $F\!f_0$.
It is straightfoward to check that the~assignment $(F,\omega) \mapsto (F,\omega)_*$
is functorial and it can be shown that $(F,\omega)_*$ and $(F',\omega')_*$ are
chain homotopic if the~morphisms are naturally isomorphic
(this also follows from the~results of the~following section).

%
% Additive closure and idempotent completion
%
For a~set of objects $B\subset\Ob(\cat C)$ define a~subcomplex
$\qCHMHom_{\!\bullet}(\cat C|B,\endofun)
	\subset \qCHMHom_{\!\bullet}(\cat C,\endofun)$
by assuming that all $x_i$ in \eqref{def:HM-chains} are from $B\noic$.
Recall that $B$ generates additively $\cat C$ in a~graded sense if each object
of $\cat C$ admits a~homogeneous isomorphism, possibly of nonzero degree,
to a~direct summand of some $x_1\oplus\dots\oplus x_r$ with $x_i\in B\noic$.

\begin{proposition}\label{prop:HH-vs-add-and-Kar}
	The~inclusion\/ $I\colon\qCHMHom_{\!\bullet}(\cat C|B,\endofun)
		\to\qCHMHom_{\!\bullet}(\cat C,\endofun)$
	is a~homotopy equivalence if\/ $B\noic$ generates additively\/ $\cat C$
	in a~graded sense. In particular, there are homotopy equivalences
	\begin{align*}
		\qCHMHom_{\!\bullet}(\cat C,\endofun)
			&\to^\simeq \qCHMHom_{\!\bullet}(\catAdd{\cat C},\endofun),
		\quad\text{and}\\
		\qCHMHom_{\!\bullet}(\cat C,\endofun)
			&\to^\simeq \qCHMHom_{\!\bullet}(\Kar(\cat C),\endofun)
	\end{align*}
	induced by the~canonical inclusions.
\end{proposition}
\begin{proof}
	We proof the~statement in two steps. First assume each object $x\in\cat C$
	admits a~decomposition $x = x_1\oplus\dots\oplus x_r$ with $x_i\in B$, together
	with inclusions $\iota_i\colon x_i\to x$ and projections $\pi_i\colon x\to x_i$,
	both homogeneous of degree 0.
	Pick the~trivial decomposition if $x$ is already in $B$. This leads to
	a~chain map $P\colon \qCHMHom_{\!\bullet}(\cat C,\endofun)
		\to \qCHMHom_{\!\bullet}(\cat C|B,\endofun)$
	\begin{equation*}
		P_{\!n}(f_0\otimes\dots\otimes f_n) := \sum_{i_0,\dots,i_n}
			(\endofun\pi_{i_n}\circ f_0\circ\iota_{i_0}) \otimes
			(\pi_{i_0}\circ f_1\circ\iota_{i_1}) \otimes\dots\otimes
			(\pi_{i_{n-1}}\circ f_n \circ\iota_{i_n}),
	\end{equation*}
	which satisfies $P\circ I = \id$. To show that $I\circ P \simeq \id$,
	consider the~family of maps
	\begin{multline*}
		h_n^k(f_0 \otimes\dots\otimes f_n) := \sum_{i_0,\dots,i_k}
			(f_0\circ\iota_{i_0}) \otimes
			(\pi_{i_0}\circ f_1 \circ \iota_{i_1})
				\otimes\cdots \\ \otimes
			(\pi_{i_{k-1}}\circ f_k \circ \iota_{i_k})
			\otimes \pi_{i_k} \otimes f_{k+1}
			 	\otimes\cdots\otimes f_n
	\end{multline*}
	and apply Lemma~\ref{lem:simpl-htpy}.
	Indeed, one checks directly that $d_{n+1}^0\circ h_n^0 = \id$,
	$d_{n+1}^n\circ h_n^n = I_n\circ P_n$,
	and the~conditions \eqref{def:simpl-htpy} are satisfied.
	
	Hence, by extending $B$, we can assume that for each object $x\in\cat C$
	there is $b(x)\in B$ together with homogeneous morphisms
	$\iota_x \colon x\to b(x)$ and $\pi_x\colon b(x)\to x$, such that
	$\pi_x\circ \iota_x = \id_x$.
	Consider now the~chain map $P\colon \qCHMHom_{\!\bullet}(\cat C,\endofun)
		\to \qCHMHom_{\!\bullet}(\cat C|B,\endofun)$ 
	\[
		P_{\!n}(f_0\otimes\dots\otimes f_n) := q^{|\iota_{x_n}|}
			(\endofun \iota_{x_n}\circ f_0 \circ \pi_{x_0}) \otimes
			(\iota_{x_0}\circ f_1\circ\pi_{x_1})
				\otimes\dots\otimes
			(\iota_{x_{n-1}}\circ f_n\circ \pi_{x_n})
	\]
	where $x_i = \dom(f_i)$. Again, $P\circ I = \id$, whereas the~other composition
	is chain homotopic to the~identity by Lemma~\ref{lem:simpl-htpy} applied
	to the~collection of linear maps
	\begin{multline*}
		h_n^k(\tensor{f0,*\cdots,fn}) :=
			(f_0\circ\pi_{x_0}) \otimes (\iota_{x_0}\circ f_1\circ\pi_{x_1})
				\otimes\cdots \\ \otimes
			(\iota_{x_{i-1}}\circ f_i\circ\pi_{x_i}) \otimes \iota_{x_i} \otimes
			f_{i+1} \otimes\cdots\otimes f_n.
	\end{multline*}
	Indeed, $\pi_x\circ\iota_x = \id_x$ implies that $d_{n+1}^0 \circ h_n^0 = \id$,
	and $d_{n+1}^{n+1}\circ h_n^n = I_n\circ P_n$ is immediate.
\end{proof}

A~graded $(A,A)$--bimodule $M$ can be understood as an~endofunctor
$(\blank)\utimes{A} M$ on $\gRep(A)$ if it is finitely generated
and projective as a~right module. This leads to the~identification
of the~quantum Hochschild--Mitchell homology of the~representation category
of $A$ with the~quantum Hochschild homology of the~algebra.

\begin{corollary}\label{cor:HH(A)-is-HH(RepA)}
	Choose a~graded algebra $A$ and a~graded $(A,A)$--bimodule $M$
	that is finitely generated and projective as a~right module.
	Then the~chain complexes $\qCHoHom_{\!\bullet}(A, M)$ and
	$\qCHMHom_{\!\bullet}(\gRep(A),M)$ are homotopy equivalent.
\end{corollary}
\begin{proof}
	The~category $\gRep(A)$ is generated additively by $B=\{A\}$ in a~graded sense,
	where $A$ is seen as a~right $A$--module. Therefore,
	\[
		\qCHMHom_{\!\bullet}(\gRep(A), M)
			\simeq \qCHMHom_{\!\bullet}(\gRep(A)|\{A\}, M)
			\cong  \qCHoHom_{\!\bullet}(A, M),
	\]
	where the~second isomorphism identifies $f\in \End_A(A)$
	with $f(1)\in A$ and $g\in \Hom_A(A,M)$ with $g(1)\in M$.
\end{proof}

\begin{remark}
	Considering $\Rep(A)$ as a~trivially graded category, the~above result
	can be rephrased to say that $\CHMHom_{\!\bullet}(\Rep(A),M)$
	is homotopy equivalent to $\CHoHom_{\!\bullet}(A,M)$.
\end{remark}

\subsection{Twisted Hochschild--Mitchell complex is a~preshadow}

%
% Endocategories as object in hTr(Cat)
%
Let $\ccat{gLinCat}$ be the~2-category of pregraded linear categories.
It is locally pregraded with a~natural transformation $\eta$ being
homogeneous of degree $d$ if each component $\eta_x$ is homogeneous of degree $d$.
Objects of $\qhTr(\ccat{gLinCat})$ are precisely pregraded endocategories.
Recall that a~morphism in the~horizontal trace from $(\cat C,\endofun)$ to
$(\cat C',\endofun')$ is represented by a~pair $(F,\omega)$ consisting
of a~functor $F\colon \cat C\to \cat C'$ and a~natural transformation
$\omega\colon F\endofun \to \endofun'F$. Hence, $(F,\omega)$ is actually
a~functor of endocategories. It is then natural to ask, whether
the~chain map $(F,\omega)_*\colon \qCHMHom_{\!\bullet}(\cat C,\endofun)
\to \qCHMHom_{\!\bullet}(\cat C',\endofun')$ depends only on the~image of
the~functor in the~horizontal trace.

\begin{theorem}\label{thm:HH-on-hTr(Cat)}
	The~quantum Hochschild--Mitchell complex descends to a~functor
	$$
		\qCHMHom\colon\qhTr(\ccat{gLinCat})\to \HCom^b(\scalars).
	$$
	In particular, it is a~quantum preshadow on\/ $\ccat{gLinCat}$.
\end{theorem}
\begin{proof}
	Consider the~following two diagrams in $\ccat{gLinCat}$
	\[
		\begin{diagps}(0em,-0.5ex)(6em,12.5ex)
			\square<5em,10ex>{->``->`->}[
				\cat C`\cat C`\cat C'`\cat C';
				\endofun``G`\endofun'
			]
			\carrow*[arcangle= 45]{->}[v0`v2;F]
			\carrow*[arcangle=-45]|b|{->}[v0`v2;G]
			\morphism[nodesep=1.6em]{=>}|a{npos=0.53}
				|(5.5em,10ex)<-4em,-9.5ex>[`;\sigma]
			\morphism{<=}(-1em,5ex)<2em,0ex>[`;\nu]
		\end{diagps}
	\hskip\dimexpr -6em + 0.33\textwidth\relax
		\mathllap{q^{|\nu|}}
		\begin{diagps}(-1em,-0.5ex)(5em,12.5ex)
			\square<5em,10ex>{->`->``->}[
				\cat C`\cat C`\cat C'`\cat C';
				\endofun`F``\endofun'
			]
			\carrow*[arcangle= 45]{->}[v1`v3;F]
			\carrow*[arcangle=-45]|b|{->}[v1`v3;G]
			\morphism[nodesep=1.6em]{=>}|b{npos=0.53}|
				(4em,10ex)<-4em,-9.5ex>[`;\sigma]
			\morphism{<=}(4em,5ex)<2em,0ex>[`;\nu]
		\end{diagps}
	\]
	Composing the~2--morphisms results in morphisms
	$(G,(\nu\circ\iid)\ast\sigma)$ and $(F,q^{|\nu|}\sigma\ast(\iid\circ\nu))$
	of endocategories. We check directly that the~family of maps
	\begin{equation*}
		h_n^k(\tensor{f0,*\dots,fn}) := q^{|\nu|}\tensor{%
			*(\sigma_{x_n} \circ Gf_0), {Gf}1, *\dots, {Gf}k,
			*\nu_{x_k}, {Ff}{k+1}, *\dots, {Ff}n}
	\end{equation*}
	where $x_i = \dom(f_i)$,
	defines a~chain homotopy between the~induced chain maps.
	Indeed, $d_{n+1}^0\circ h_n^0 = q^{|\nu|}(F, \sigma\ast(\iid\circ\nu))_n$
	by the~naturality of $\nu$, whereas
	$d_{n+1}^{n+1}\circ h_n^n=(G, (\nu\circ\iid)\ast\sigma)_n$
	is straightforward.
\end{proof}

\begin{corollary}
	If $(F,\omega)$ and $(F',\omega')$ are isomorphic functors,
	then $(F,\omega)_* = (F',\omega')_*$.
\end{corollary}
\begin{proof}
	Suppose $\eta\colon (F,\omega)\to (F',\omega')$ is an~isomorphism.
	Then
	\[
		\begin{diagps}(-1em,-0.5ex)(6em,12.5ex)
			\square<5em,10ex>[
				\cat C`\cat C`\cat C'`\cat C';
				\endofun`F'`F'`\endofun'
			]
			\morphism[nodesep=1.6em]{=>}|b{npos=0.5}
				|(4.75em,10ex)<-4em,-9.5ex>[`;\omega']
		\end{diagps}
		\quad=\quad
		\begin{diagps}(-2em,-0.5ex)(8em,12.5ex)
			\node tl(0em,10ex)[\cat C ]	\node bl(0em, 0ex)[\cat C']
			\node tr(6em,10ex)[\cat C ]	\node br(6em, 0ex)[\cat C']
			\arrow{->}[tl`tr;\endofun ]	\arrow{->}[bl`br;\endofun']
			
			\carrow*[arcangle= 45]|a|{->}[tl`bl;F ]
			\carrow*[arcangle=-45]|b|{->}[tl`bl;F']
			\morphism{<=}(-1em,5ex)<2em,0ex>[`;\eta]

			\carrow*[arcangle= 45]|a|{->}[tr`br;F']
			\carrow*[arcangle=-45]|b|{->}[tr`br;F ]
			\morphism{<=}|a{npos=0.6}|( 5em,5ex)<2em,0ex>[`;\eta^{-1}]
			
			\morphism[nodesep=1.6em]{=>}|b{npos=0.53}
				|(5.3em,10ex)<-4em,-9.5ex>[`;\omega]
		\end{diagps}
		\quad=\quad
		\begin{diagps}(-1em,-0.5ex)(6em,12.5ex)
			\square<5em,10ex>[
				\cat C`\cat C`\cat C'`\cat C';
				\endofun`F`F`\endofun'
			]
			\morphism[nodesep=1.6em]{=>}|b{npos=0.5}
				|(4.75em,10ex)<-4em,-9.5ex>[`;\omega]
		\end{diagps}
	\]
	in $\qhTr(\ccat{gLinCat})$, where the~second equality follows
	from the~horizontal trace relation.
\end{proof}

It follows from Theorem~\ref{thm:HH-on-hTr(Cat)} that the~quantum
Hochschild--Mitchell homology is a~preshadow on $\ccat{gLinCat}$,
which generalizes Example~\ref{sec:shadow-Tr}. In the~view
of Corollary~\ref{cor:HH(A)-is-HH(RepA)} it induces a~preshadow on $\ccat{gRep}$ 
with \emph{a~priori} a~different cyclicity map.

\begin{proposition}\label{prop:qHM-is-qHH-as-shadows}
	The~homotopy equivalences
	$\qCHMHom_{\!\bullet}(\gRep(A),M) \simeq \qCHoHom_{\!\bullet}(A,M)$
	constitutes an~isomorphism of preshadows.
\end{proposition}
\begin{proof}
	We have to check that the~square
	\[
		\begin{diagps}(0,-1ex)(15em,13ex)
			\square<15em,10ex>{->`->`<-`->}[%
				\qCHoHom_{\!\bullet}(B, N'\utimes{A} N)`
				\qCHoHom_{\!\bullet}(A, N\utimes{B} N')`
				\qCHMHom_{\!\bullet}(\gRep(B), N'\utimes{A} N)`
				\qCHMHom_{\!\bullet}(\gRep(A), N\utimes{B} N');
				\theta_{N',N}`I`P`N'_*
			]
		\end{diagps}
	\]
	commutes up to a~chain homotopy for bimodules $N\in\gRRep(A,B)$ and
	$N'\in\gRRep(B,A)$, where $I$ and $P$ are mutually inverse homotopy equivalences.
	Due to the~naturality of both $\theta_{N',N}$ and $N_*'$ we can assume that $N'$
	is free as a~right $A$--module with a~graded basis $\{n'_1,\dots,n'_s\}$.
	Then the~left action of $B$ determines an~algebra map $B\to \mathit{Mat}_s(A)$
	that assigns to $b\in B$ the~matrix $(a_{ij}(b))_{1\leqslant i,j\leqslant s}$
	determined uniquely by the~formula
	\begin{equation}\label{eq:left-b-to-right-mat(a)}
		bn'_i = \sum_{j} n'_j a_{ij}(b).
	\end{equation}
	One checks that $a_{ij}(b'b) = \sum_k a_{kj}(b')a_{ik}(b)$ for any $b,b'\in B$.
	Likewise, $n'_k\otimes n\otimes\id_{N'}$ is represented, as a~right $A$--linear
	map, by the~matrix $(\delta_{jk} n\otimes n'_j)_{1\leqslant i,j\leqslant s}$.
	It follows from the~proof of Proposition~\ref{prop:HH-vs-add-and-Kar} that
	the~composition $P\circ N'_*\circ I$ is the~chain map
	\begin{equation}\label{eq:explicit-N*}
		(n'_{i_0}\otimes n)\otimes b_r\otimes\cdots\otimes b_1
	\mapsto
		\sum_{i_1,\dots,i_r = 1}^s q^{|n'_{i_r}|}
			(n\otimes n'_{i_r}) \otimes a_{i_{r-1},i_r}(b_r)
				\otimes \cdots \otimes
			a_{i_0,i_1}(b_1).
	\end{equation}
	To compute $\theta_{N',N}$, we need to lift the~identity
	on $N'\utimes{A} N$ to a~pair of graded chain maps
	\[
		p^{BA}\colon
			R_\bullet(B)\utimes{B} N' \utimes{A} R_\bullet(A) \utimes{A} N
		\twoways
			R_\bullet(B) \utimes{B} N'\utimes{A} N
		\cocolon s^{BA}\!.
	\]
	By the~uniqueness of a~lift, $p^{BA}$ and $s^{BA}$
	are mutually inverse homotopy equivalences.
	
	The~map $p^{BA}$ is defined as the~standard projection
	on $R_r(B)\utimes{B} N'\utimes{A} R_0(A) \utimes{A} N
	= R_r(B) \utimes{N} N'\otimes N$
	and vanishes otherwise. The~chain map $s^{BA}$ replaces each $b_\ell$,
	one-by-one, with entries of the~matrix $a_{ij}(b_\ell)$:
	\begin{multline}
		s^{BA}(\tensor{br,*\cdots,b1,*(n'_{i_0}\otimes n)}) := \\
			\sum_{\ell=0}^r \sum_{i_1,\dots,i_\ell=1}^s
				\tensor{br,*\cdots,b{\ell+1},{n'}{i_\ell},
					*a_{i_{\ell-1},i_\ell}(b_{\ell}),*\cdots,
					*a_{i_0,i_1}(b_1),*n}.
	\end{multline}
	The~equality \eqref{eq:left-b-to-right-mat(a)} makes $s^{BA}$ commute
	with the~differential.
	
	Up to the~isomorphisms
	$\coinv_q(M\utimes{A} R_\bullet(A)) \cong \coinv_q(R_\bullet(A)\utimes{A}M)$
	we can express $\theta_{N',N}$ as the~composition
	$p^{AB}\circ\tau\circ s^{BA}$,
	where $\tau$ it the~graded twist from \ref{sec:ex-qHH}.
	Because $|a_{ij}(b_k)| = |b_k|$ and $p^{AB}$ vanishes unless none of
	$b_k$ appears, the~composition coincides with \eqref{eq:explicit-N*}.
\end{proof}

\subsection{\texorpdfstring{$K$--}{K-}theoretic invariance}
\label{sec:K-invariance}

We say that a~sequence of linear functors
\begin{equation}\label{eq:sses-of-functors}
	(F',\omega') \to^{\iota} (F,\omega) \to^{\pi} (F'',\omega'')
\end{equation}
\emph{semisplits} if
$0 \to F'x \to^{\iota_x} Fx \to^{\pi_x} F''x \to 0$ splits
for every object $x$. Clearly, a~split exact sequence semisplits,
but not otherwise.
Moreover, the~sequence \ref{eq:sses-of-functors} remains semisplit when
each of the~three functors is pre- or post-composed with another
functor, because linear functors preserve direct sums.

Let $\ccat{gEndoAddCat}$ be the~bicategory of small graded additive
endocategories. We construct its \emph{semisplit Grothendieck category}
$\Kss(\ccat{gEndoAddCat})$ by replacing functor categories with
$\Z[q^{\pm1}]$--modules generated by isomorphism classes $[F,\omega]$ of functors
modulo the~relations $[F\{1\},\omega\{1\}] = q[F,\omega]$ and
$[F,\omega] = [F',\omega'] + [F'',\omega'']$ for every semisplit
exact sequence \eqref{eq:sses-of-functors}. The~composition is well-defined
by the~discussion above.
A~functor $(F,\omega)\colon (\cat C,\endofun) \to (\cat C',\endofun')$
is a~\emph{$K\noic$--theoretic equivalence} if it is an~isomorphism in
$\Kss(\ccat{gEndoAddCat})$. In other words, there is a~functor
$(F',\omega')\colon (\cat C',\endofun') \to (\cat C,\endofun)$,
called the~\emph{$K\noic$--theoretic inverse} of $(F,\omega)$,
such that the~compositions $(F',\omega')\circ(F, \omega)$ and
$(F,\omega)\circ(F',\omega')$ coincide in the~corresponding Grothendieck groups
with the~images of the~identity functors.

\begin{theorem}[$K\noic$--theoretic invariance]\label{thm:K-invariance-of-HM}
	The~homotopy class of\/ the~chain map $(F,\omega)_*$ assigned to
	a~functor of linear graded endocategories
	$(F, \omega)\colon (\cat C, \endofun)\to (\cat C',\endofun')$
	depends only on\/ $[F,\omega] \in \Kss(\ccat{gEndoAddCat})$.
	In particular, $(F,\omega)_*$ is a~homotopy equivalence if
	the~functor is a~$K$--theoretic equivalence.
\end{theorem}
\begin{proof}
	The~relation $(F\{1\},\omega\{1\})_* = q(F,\omega)_*$ follows from
	the~proof of Proposition~\ref{prop:HH-vs-add-and-Kar}. Thence, we have
	to show that $(F,\omega)_* = (F',\omega')_* + (F'',\omega'')_*$
	for a~semisplit sequence of functors \eqref{eq:sses-of-functors}.
	For that fix an~isomorphism $\varphi_x \colon Fx \cong F'x \oplus F''x$
	for each $x\in\cat C$ that results from a~splitting of the~sequence
	$0\to F'x\to Fx \to F''x\to 0$.
	Then every morphism $g\in\cat C(Fx,Fy)$
	can be identified with a~$2\times 2$ matrix $(g_{ij})$.
	In particular, the~proof of Proposition~\ref{prop:HH-vs-add-and-Kar}
	implies that $\eta_*$ is chain homotopic to
	\begin{equation}\label{eq:F-map-flattened}
		\tensor{f0,*\dots,fn} \mapsto \sum_{i_0,\dots,i_n}
			\tensor{
				{(\omega \circ Ff_0)}{i_0,i_n},
				{(Ff_1)}{i_1,i_0},
				*\dots,
				{(Ff_n)}{i_n,i_{n-1}}
			}.
	\end{equation}
	A~simple diagram chasing reveals that for any $f\in\cat C(x,y)$
	the~matrix of $F\noic f$ is upper triangular with $F'\!f$ and $F''\!f$
	on the~diagonal, and likewise for $\omega_x$.
	Therefore, the~only nontrivial summands in \eqref{eq:F-map-flattened}
	are the~two with $i_0 = \dots = i_n$, which are precisely
	$(F',\omega')_*(\tensor{f0,*\dots,fn})$ and
	$(F'',\omega'')_*(\tensor{f0,*\dots,fn})$.
\end{proof}

\subsection{An~application to Hochschild homology of algebras}

%We shall now translate the~$K$--theoretic invariance to Hochschild homology
%of algebras. This gives a~very powerful tool that reduces the~problem of
%computing the~homology of an~algebra to a~simpler one.

Let $A$ be a~graded algebra and $\varphi\in\Aut(A)$ a~graded automorphism.
Following \ref{sec:example-twisted-trace} construct the~bimodule ${}_\varphi A$
by redefining the~left action as $a\cdot x := \varphi(a)x$. The~bimodules
$A_\varphi$ and ${}_\varphi A_\varphi$ are constructed likewise.
The~automorphism $\varphi$ can be seen as an~isomorphism of bimodules
$A \to^\cong {}_\varphi A_\varphi \cong {}_\varphi A\utimes{A} A_\varphi$.
Tensoring with ${}_\varphi A$ is an~endofunctor on $\cat{gMod}(A)$ denoted
in \ref{sec:example-twisted-trace} by $(\blank)_{\overline\varphi}$.
A~module $V$ coincides with $V_{\overline\varphi}$,
except that the~action of $A$ is twisted: $v\cdot a = v\varphi^{-1}(a)$.

Let $B$ and $\psi\in\Aut(B)$ be another graded algebra together with a~graded
automorphism. Following \ref{sec:shadow-twisted-coinv} we put
$\endofun M := A_\varphi\utimes{A} M\utimes{B} {}_\psi B$ for a~graded
$(A,B)$--bimodule $M$. It coincides with $M$, except that $\varphi$ and $\psi$
are used to twist both actions. Suppose there is a~bimodule homomorphism
$\omega\colon M\to \endofun M$. Then there is a~$B$--linear map
\[
	V_{\overline\varphi} \utimes{A} M
		\to^{\id\otimes\omega}
	V_{\overline\varphi} \utimes{A} {}_{\overline\varphi} M_{\overline\psi}
		\to^\cong
	(V\utimes{A} M)_{\overline\psi}
\]
for any right $A$--module $V\noic$, where the~second map is the~isomorphism
${}_\varphi A_\varphi \cong A$. Hence, $(M,\omega)$ determines a~functor
of endocategories.

Let $\gBimod_0(A,\varphi; B,\psi)$ be the~category with objects
graded pairs $(M,\omega_M)$ as above and morphisms between
$(M,\omega_M)$ and $(N,\omega_N)$ graded bimodule maps $f\colon M\to N$
that intertwine the~structure, i.e.\ $\omega_N\circ f = f\circ\omega_M$.
One checks directly that $\omega_M(\ker f)\subset \ker f$ and
$\omega_N(\im f)\subset \im f$, so that $\gBimod_0(A,\varphi; B,\phi)$ has kernels
and cokernels. Thence, it is an~abelian category.
Restricting bimodules to those that are finitely generated and projective
as left modules picks the~subcategory $\gRep_0(A,\varphi; B,\psi)$.
It is additive, but not abelian. Note that the~categories are slightly
different from morphism categories of $\egBBimod$ and $\egRRep$
from Section~\ref{sec:shadow-twisted-coinv}.

Following the~usual convention we shall write $G_0(A,\varphi; B,\psi)$
and $K_0(A,\varphi; B,\psi)$ for the~(exact) Grothendieck groups%
\footnote{
	That is the~isomorphism class of a~bimodule $[M]$ is identified with
	$[M']+[M'']$ whenever there is an~exact sequence $0\to M'\to M\to M''\to 0$.
}
of $\gBimod_0(A,\varphi; B,\psi)$ and $\gRep_0(A,\varphi; B,\psi)$ respectively.
The~inclusion of categories induces a~linear map $K_0(A,\varphi; B,\psi)
\to G_0(A,\varphi; B,\psi)$, which is neither surjective nor injective in general.
Because not all exact sequence in $\gRep_0(A,\varphi; B,\psi)$ splits,
$K_0(A,\varphi; B,\psi)$ is usually a~proper quotient of the~split
Gro\-then\-dieck group.

The~tensor product $(M\utimes B M', \omega\utimes B\omega')$ 
of bimodules $(M,\omega)\in\gRep_0(A,\varphi; B,\psi)$
and $(M',\omega')\in\gRep_0(B,\psi; C,\rho)$,
belongs to $\gRep_0(A,\varphi; C,\rho)$. Hence,
the~categories $\gRep_0(A,\varphi; B,\psi)$ constitute a~bicategory.
Moreover, tensoring with $(M,\omega)$ is exact.
Indeed, exact sequences in $\gRep_0(B,\varphi; C,\rho)$
split when considered as sequence of left modules, which is enough to prove
the~exactness of $(M,\omega)\utimes{B} (\blank)$, whereas the~exactness of
$(\blank)\utimes{A}(M,\omega)$ follows from $M$ being projective as a~left module.
Hence, the~Grothendieck groups $K_0(A,\varphi; B,\psi)$ form a~category as well.
Following \cite{InvForCyclicHom} we say that $(M,\omega)$ is
a~\emph{$K\noic$--theoretic equivalence} if $[M,\omega]$ is invertible
in this category. The~following is an~immediate translation
of Theorem~\ref{thm:K-invariance-of-HM}.

\begin{proposition}[Algebraic $K\noic$--theoretic invariance]
\label{prop:K-invariance-of-HH}
	The~homotopy class of the~chain map\/
	$(M,\omega)_*\colon \qCHoHom_{\!\bullet}^\psi(B)
		\to\qCHoHom_{\!\bullet}^\varphi(A)$
	induced by\/ $(M,\omega)\in\gRep_0(A,\varphi;B,\psi)$
	depends only on\/ $[M,\omega]\in K_0(A,\varphi;B,\psi)$.
	In particular, $(M,\omega)_*$ is a~homotopy equivalence if\/
	$(M,\omega)$ is a~$K\noic$--theoretic equivalence.
\end{proposition}

Choose now a~positively graded algebra $A=\bigoplus_{d\in\mathbb N} A_d$.
The~degree zero subalgebra $A_0\subset A$ is preserved by
any graded automorphism $\varphi\in\Aut(A)$.
Hence, the~inclusion and projection induce chain maps
$\qCHoHom_{\!\bullet}^{\varphi_0}(A_0) \to \qCHoHom_{\!\bullet}^\varphi(A)$
and
$\qCHoHom_{\!\bullet}^\varphi(A) \to \qCHoHom_{\!\bullet}^{\varphi_0}(A_0)$
respectively, where $\varphi_0\in\Aut(A_0)$ is the~restriction of $\varphi$.
One of the~compositions is clearly the~identity map, but not the~other.

These chain maps may not be homotopy equivalences.
For instance, it is known that the~algebra of dual numbers $\scalars[x]/(x^2)$,
where $\deg x =2$, has unbounded Hochschild homology,
whereas the~homology of $\scalars$ is one-dimensional.
The~situation changes drastically when $A$ has finite global dimension:
the~inclusion $A_0\to A$ induces an~isomorphism on Hochschild homology
under some other technical conditions \cite{InvForCyclicHom}.
Here we reprove this result for the~twisted homology.

\begin{theorem}\label{thm:HH-for-fgdim}
	Let\/ $\scalars$ be a~field and\/ $A=\bigoplus_{d\in\mathbb N} A_d$
	a~positively graded\/ $\scalars$--algebra of finite dimension.
	If each simple\/ $A\noic$--module is one dimensional and\/ $A\noic$
	has finite global dimension, then the~inclusion\/ $A_0\subset A\noic$
	induces a~homotopy equivalence
	$\qCHoHom_{\!\bullet}^{\varphi_0}(A_0) \to \qCHoHom_{\!\bullet}^\varphi(A)$
	for any graded automorphism $\varphi\in\Aut(A)$.
\end{theorem}
\begin{proof}
	The~chain map
	$\qCHoHom_{\!\bullet}^{\varphi_0}(A_0) \to \qCHoHom_{\!\bullet}^\varphi(A)$
	is induced by $(A,\varphi)$, seen as an~object of
	$\gRep_0(A,\varphi; A_0,\varphi_0)$. Likewise, the~other chain map is induced by
	$(A_0,\varphi_0) \in \gRep_0(A_0,\varphi_0; A, \varphi)$,
	where the~right action of $A$ on $A_0$ is given by the~projection.
	A~quick computation shows that $(A_0,\varphi_0)\utimes{A} (A,\varphi) \cong
	(A_0,\varphi_0)$ in $\gRep_0(A_0,\varphi_0; A_0,\varphi_0)$. Hence,
	the~composition
	$\qCHoHom_{\!\bullet}^{\varphi_0}(A_0) \to
		\qCHoHom_{\!\bullet}^\varphi(A) \to
		\qCHoHom_{\!\bullet}^{\varphi_0}(A_0)$
	is homotopic to the~identity.
	In the~view of Proposition~\ref{prop:K-invariance-of-HH}
	it is now enough to show that $(A,\varphi)$ and
	$(A,\varphi)\utimes{A_0} (A_0,\varphi_0)$ coincide in
	$K_0(A,\varphi; A,\varphi)$.
	
	We check first that the~bimodules have the~same image in
	$G_0(A,\varphi; A,\varphi)$. Indeed,
	\begin{equation*}
		A_{\geqslant d}/A_{>d} \cong
			(A_{\geqslant d}\utimes{A_0} A_0)/A_{>d}\utimes{A_0} A_0),
	\end{equation*}
	which shows that the~graded associate of the~bimodules coincide,
	and so their composition series (which are finite due to $A$
	being of finite dimension). It remains to show that every $(A,A)$--bimodule
	of finite length can be represented in $K_0(A,\varphi; A,\varphi)$.
	
	We first notice that $A^e:=A\otimes A^{op}$ has finite global dimension:
	each simple $A^e$--module is of the~form $L'\otimes L^{\!*}$ for certain simple
	$A$--modules $L'$ and $L$, and as such it has a~uniformly bounded
	projective resolution.
	Thence, every $(A,A)$--bimodule $M$ of finite length
	has a~finite projective resolution $P^\bullet(M)$.
	Because $A_\varphi$ and ${}_\varphi A$ are projective as left and right
	modules, $\endofun$ is exact and preserves projective bimodules. 
	Hence, $P^i(\endofun M) := \endofun(P^i(M))$
	is a~projective resolution of $\endofun M$. Lifting $\omega$ to a~chain map
	$P^\bullet(\omega)$ between these complexes, we construct a~resolution of
	$(M,\omega)$ in $\cat{gRep}_0	(A,\varphi; A,\varphi)$. The~assignment
	\begin{equation}
		[M,\omega] \mapsto \sum_i (-1)^i [P^i(M),P^i(\omega)]
	\end{equation}
	is then a~two-sided inverse of
	$K_0(A,\varphi; A,\varphi) \to G_0(A,\varphi; A,\varphi)$.
	Hence, the~bimodules $(A,\varphi)$ and $(A,\varphi)\utimes{A_0} (A_0,\varphi_0)$
	coincide in $K_0(A,\varphi; A,\varphi)$ as desired.
\end{proof}

\section{Khovanov Homology}
\label{sec:khov-hom}

With this section we move from the~algebraic to the~topological part of the~paper.
It starts with a~brief description of the~formal bracket, a~very generic approach
to Khovanov-type homology due to Bar-Natan \cite{Cob}, which is followed by a~list
of TQFT functors producing link homologies used in this paper. The~main purpose
of this section is to fix notation and conventions. In particular, the~reader
should be aware that our convention for the~quantum grading makes Chen--Khovanov
algebras negatively graded.

\subsection{The~formal bracket}
\label{sec:bracket}

Let us start with a~brief reminder of the~formal Khovanov bracket following
\cite{Cob}. In this section $\sfce$ stands for a~surface, possibly with boundary.

%
% Cube of resolutions
%

Choose a~diagram $D\subset\sfce$ of an~oriented tangle $T \subset \sfce\times\R$ and let $n$ stand for the~number of crossings in $D$. To compute the~formal Khovanov bracket of $T$ one begins with creating the~$n$-dimensional \emph{cube of resolutions} $\KhCube{D}$ of $D$ defined as follows:
\begin{itemize}
	\item each vertex $\xi$ of $\KhCube{D}$ is decorated with the~resolution $D_\xi$ of $D$,
	i.e.\ the~collection of circles and proper intervals in $\sfce$ obtained from $D$ by forgetting
	the~orientation and replacing each $i$--th crossing \fntNWSECr\ by its horizontal \fntHorRes\ or
	vertical \fntVertRes\ smoothing for $\xi_i=0$ and $\xi_i=1$ respectively, and
	\item each edge is decorated by a~cobordism with a~unique saddle point over the~smoothing being changed
	and directed from the~resolution with less vertical smoothings to the~one with more of them.
\end{itemize}
We can view $\KhCube{D}$ as a~commutative cubical diagram in $\Cob(\sfce)$,
the~category with objects (non-oriented) flat tangles in $\sfce$ and
morphisms the~isotopy classes of surfaces in $\sfce\times I$. The~degree
of a~surface $S \subset \sfce\times I$ is given by the~formula
\begin{equation}\label{def:degree-in-Cob}
	\deg S := \chi(S) - \frac{\#B}{4}
\end{equation}
where $B$ is the~set of corners and $\chi(S)$ the~Euler characteristic of $S$.
Alternatively, $\deg S$ counts the~critical points of the~natural height
function $h\colon S\to I$ with signs: a~point of index $\mu$ contributes
$(-1)^\mu$ towards $\deg$. Thence, to ensure that each morphism in $\KhCube{D}$
has degree $0$, we introduce formal degree shifts%
\footnote{
	Formally, objects in $\Cob(\sfce)$ are now symbols $T\{i\}$
	formed by a~flat tangle $T\subset\sfce$ and $i\in\Z$.
}
and place at each vertex $\xi$ the~shifted resolution $D_\xi\{|\xi|\}$
for $|\xi| := \xi_1 + \dots + \xi_n$.

%
% Definition of the bracket
%

Let $\catAdd{\Cob(\sfce)}$ be the~additive closure of $\Cob(\sfce)$. The~\emph{formal Khovanov bracket} $\KhBracket{D}$ of the~tangle diagram $D$ is the~complex in $\catAdd{\Cob(\sfce)}$ obtained from the~cube by distributing signs over some edges to make every square anticommute, then taking direct sums along diagonal sections of the~cube, and finally applying suitable degree shifts:
\begin{equation}
	\KhBracket{D}^i := \bigoplus_{|\xi|=i+n_-}\mskip-12mu D_\xi\{i+n_+-n_-\}
\end{equation}
where $n_{\pm}$ stand for the~amount of positive or negative crossings in $D$.

%
% Invariance
%

\begin{theorem}[cf.\ \cite{Cob}]\label{thm:invariance-of-KhBracket}
	The~homotopy type of the~complex\/ $\KhBracket D$ is a~tangle invariant after imposing on\/ $\catAdd{\Cob(\sfce)}$ the~following local relations%
	\footnote{
		\emph{Locality} means that each picture represents a~part of a~cobordism inside a~ball in $\sfce\times I$.
	}
	\begin{align*}
		\textnormal{(\textit{S})}&\psset{unit=1cm}\quad\pictRelS = 0
		\hskip 2cm
		\textnormal{(\textit{T})}\quad\pictRelT = 2\\[1ex]
		\textnormal{(\textit{4Tu})}&\quad\psset{unit=0.75cm}%
			\pictRelTuL + \pictRelTuR = \pictRelTuB + \pictRelTuT
	\end{align*}
	introduced first by Bar-Natan \cite{Cob}.
\end{theorem}

%
% Viro's exact sequence
%

Let \fntPosCr, \fntNegCr, \fntNoCr, and \fntHorRes\ be tangle diagrams
that coincide everywhere except a~small disk, in which they look as indicated
by the~pictures. It follows from the~construction that the~formal brackets
of the~first two diagrams are mapping cones of chain maps between the~formal
complexes of the~other two diagrams.

\begin{proposition}[cp.\ \cite{KhHom,Cob}]\label{prop:Viro-ses}
	Let\/ \fntPosCr, \fntNegCr, \fntNoCr, \fntHorRes\ be four link diagrams as above,
	where we choose any orientation for\/ \fntHorRes, and write
	$e = n_-(\fntHorRes) - n_-(\fntNoCr)$.
	Then there are isomorphisms of formal complexes
	\begin{align}
	\label{ses:Viro-pos}
		\KhBracket\fntPosCr \cong&\, \cone\big(
			\Rnode{tl}{\KhBracket\fntNoCr \{1\}[-1]} & \hbox to 4em{} &
			\Rnode{tr}{\KhBracket\fntHorRes \{3e+2\}[-1-e]}
			\big), \text{and}
		\\
	\label{ses:Viro-neg}
			\KhBracket\fntNegCr \cong&\, \cone\big(
			\Rnode{bl}{\KhBracket\fntHorRes \{3e-2\}[-e]} &&
			\Rnode{br}{\KhBracket\fntNoCr \{-1\}}
			\big),
		\begin{diagps}(0,0)
			\arrow|a{npos=0.40}|{->}[tl`tr;\KhBracket\fntVHCob]
			\arrow|a{npos=0.55}|{->}[bl`br;\KhBracket\fntHVCob]
		\end{diagps}
	\end{align}
	where the~chain maps are induced by saddle cobordisms.
\end{proposition}

The~isomorphisms \eqref{ses:Viro-pos} and \eqref{ses:Viro-neg} were
first observed in \cite[Section 4.2]{KhHom} for knots in $\R^3$ and
then by Bar-Natan in the~framework of the~formal bracket \cite[Lemma 4.4]{Cob}.
They lead to distinguished triangles in the~homotopy category of
$\Cob(\sfce)^\oplus$, which can be seen as categorified versions
of the~Kauffman skein relation.

%
% TQFTs and homology groups
%

To retrieve homology groups from $\KhBracket D$ one has to replace
$\catAdd{\Cob(\sfce)}$ with an~abelian category. This is done by applying to
the~above construction a~graded TQFT functor that preserves the~relations
\textit{S}, \textit{T}, and \textit{4Tu}. Most known functors of this type
factorize through the~universal Bar-Natan skein category $\BNu(\sfce)$ \cite{Cob},
the~additive linear category generated by circles and proper intervals in $\sfce$
as objects, and formal linear combinations of cobordisms in $\sfce\times I$
decorated by dots (each dot decreases the~degree of a~cobordism by 2) as morphisms,
subject to the~local relations
\begingroup\psset{unit=7mm}%
\begin{align*}
		\textnormal{(\textit{S})}\quad&\pictRelS = 0, \hskip 1.5cm
		\textnormal{(\textit{D})}\quad\pictRelD = 1,\\[2ex]
		\textnormal{(\textit{N})}\quad&\pictRelNeckI = \pictRelNeckB + \pictRelNeckT - \pictRelNeckM.
\end{align*}
\endgroup
The~last relation is commonly referred as the~\emph{neck cutting relation}
and it allows to reduce any surface to a~linear combination of surfaces of genus
zero. Both \textit{T} and \textit{4Tu} follows, and a~contractible circle can be
replaced with a~pair of shifted copies of an~empty set. In particular, the~empty
set generates $\BNu(\sfce)$ when $\sfce$ is either $\R^2$ or $\mathbb S^2$.

%
% Delooping
%

\begin{proposition}[Delooping, cf.\ \cite{DotCob}]\label{prop:delooping}
	There is a~pair of mutually inverse isomorphisms in\/ $\BNu(\sfce)$
	\begin{equation}\label{eq:delooping}\psset{unit=5ex}%
		\begin{diagps}(-3,-1.1)(3,1.3)
			\node l(-2.8,0)[\fntCircle]
			\node r( 2.8,0)[\fntCircle]
			\node[vref=0.3] t(0, 0.8)[\emptyset\{-1\}]
			\node[vref=0.0] b(0,-0.8)[\emptyset\{+1\}]
			\rput[B](0,0){$\oplus$}
			\psset{dash=1pt 1.5pt,linewidth=0.5pt}
			\arrow|a|{->}[l`t;\pictDeloopTL]
			\arrow|b|{->}[l`b;-\,\pictDeloopBLb\,+\,\pictDeloopBLa]
			\arrow|a|{->}[t`r;\pictDeloopTR]
			\arrow|b|{->}[b`r;\pictDeloopBR]
		\end{diagps}
	\end{equation}
	for every circle $\fntCircle$ bounding a~disk in $\sfce$.
\end{proposition}

%If $\F\colon\Cob(\sfce)\to\cat A$ produces an~invariant finite dimensional chain complex, $\F(\fntCircle)$ must be a~rank two module over $\F(\emptyset)$ \cite{FrobExt}. In particular, $\F$ factors through $\BNu(\sfce)$ if $\F(\fntCircle)$ is a~free module, so that little information is lost when we decide to work with the~universal Bar-Natan category instead of $\Cob(F)$.

%
% Additional relations
%

The~TQFT functors we consider in this paper factor through certain quotients of
$\BNu(\sfce)$. Consider the~following two relations:
\begingroup\psset{unit=7mm}
\begin{align*}
	\textnormal{(\textit{TD})}\quad&\pictRelTwoDots = 0 \\[1ex]
	\textnormal{(\textit{B})} \quad&\hskip 0.4\psxunit\pictRelB\hskip 0.4\psxunit = 0 \quad\text{if}\quad [\gamma] \neq 0\text{ in }\pi_1(\sfce\times I).
\end{align*}
\endgroup
The first one asserts that two dots annihilate a~cobordism when placed on the~same
connected component. Together with the~neck cutting relation it implies further
that a~cobordism is annihilated when a~dot is placed on its component of positive
genus. The~second relation prohibits a~component of a~cobordism to carry a~dot if
the~component contains a~closed curve that is nontrivial in $\pi_1(\sfce\times I)$.
In particular, a~component of a~cobordism cannot carry a~dot if its boundary curves
are not contractible in $\sfce$.

We write $\BN(\sfce)$ when only \textit{TD} is imposed and $\BoBN(\sfce)$
when both. $\BN(\sfce)$ is commonly called the~\emph{Bar-Natan skein category}
and it was first defined in \cite{Cob}. The~case of the~annulus was extensively
studied by Russell \cite{AnnBN}. The~relation \textit{B} was introduced by Boerner
\cite{Boerner} for any surface $\sfce$ and we call $\BoBN(\sfce)$
the~\emph{Boerner--Bar-Natan skein category}.

% ==================================================================
%  Section: Functoriality
% ==================================================================

\subsubsection*{Functoriality}

The~construction of the~formal bracket $\KhBracket{T}$ is functorial up to signs: given a~cobordism $S$ between tangle diagrams $T$ and $T'$ there is a~chain map $\KhBracket{S}\colon \KhBracket{T} \to \KhBracket{T'}$ defined up to the~factor $\pm1$, such that $\KhBracket{SS'} = \pm\KhBracket{S}\circ\KhBracket{S'}$. Hence, the~same type of functoriality holds for any TQFT functor $\F\colon\Cob(\sfce)\to\cat{A}$ satisfying the~relations \textit{S}, \textit{T}, and \textit{4Tu}.

The~chain map $\KhBracket{S}$ is computed from a~\emph{movie presentation} of $S$, a~sequence of generic sections $S_t = S\cap(\sfce\times\{t\})$ called \emph{movie clips}, such that the~part $S|_{[t,t']}$ of $S$ between two consecutive clips $S_t$ and $S_{t'}$ is either one of the~Reidemeister moves, a~saddle cobordism, a~cap, or a~cup \cite{CarterSaito}. There is a~well-defined chain map for each of the~parts and $\KhBracket{S}$ is defined as the~composition of these pieces. A~cobordism $S$ admits many movie presentations and it is proven that up to sign $\KhBracket{S}$ does not depend on the~presentation chosen \cite{Cob}. Unfortunately, direct computation shows that the~sign of $\KhBracket{S}$ does depend on the~presentation \cite{JacobFun}.

There are a~few approaches to attack the~sign issue. In case $\sfce = \R^2$
one can use the~Lee deformation of the~Khovanov homology \cite{Lee} to define
\emph{canonical generators}, which are preserved by $\KhBracket{S}$ up to sign
\cite{CanGens}. We can then redefine $\KhBracket{S}$ so that the~generators are
actually preserved. This approach was used in \cite{KhHomSL2} to fix signs in
certain cases.

A~different idea is to replace $\Cob(\sfce)$ with another category. This was done
successfully for $\sfce=\R^2$ by Clark, Morrison and Walker using cobordisms with
seams and coefficients in the~ring of Gaussian integers $\Z[i]$ \cite{CMW},
then by Blanchet using nodal foams \cite{BlanchetFoams} and by Vogel using
\emph{mixed cobordisms}: locally oriented cobordisms with certain disorientation
curves \cite{VogelFun}. The~first strictly functorial construction for
$(2m,2n)$--tangles is due to Caprau \cite{sl2Webs}, and is defined over Gaussian
integers. It assigns to a~tangle a~direct summand of the~corresponding invariant
due to Chen and Khovanov, which is used in this paper. In a~forthcoming paper
\cite{gl2foams} we address the~functoriality of the~Chen--Khovanov invariant by
rephrasing it in terms of $\mathfrak{gl}_2$ foams and constructing an~explicit
isomorphism between the~new and the~original invariant.
Finally, the~functoriality of the~bracket for links in thickened surfaces has been
fixed in a~recent work of Queffelec and Wedrich \cite{KhForSfces}, who constructed
certain foam categories that extend the~nodal foams introduced by Blanchet.

\subsection{Khovanov homology for links in \texorpdfstring{$\R^3$}{R3}}

Let $R$ be a~commutative algebra and $A$ a~Frobenius algebra over $R$ of rank 2. This datum determines a~TQFT functor with $\F(\emptyset):=R$ and $\F(\fntCircle):=A$, and it produces an~invariant chain complex for link diagrams on a~plane \cite{FrobExt}. For instance, Khovanov's functor $\FKh\colon\Cob(\R^2)\to\Mod(\scalars)$ is defined this way by taking $R:=\scalars$ and equipping $A:=Rw_+\oplus Rw_-$ with the~structure maps
\begin{align}
	m&\colon A\otimes A\to A &&
	\left\{\begin{array}{l}
		\tensor{w+,w+} \mapsto w_+, \\
		\tensor{w\pm,w\mp} \mapsto w_-, \\
		\tensor{w-,w-} \mapsto 0,
	\end{array}\right.
	\\
	\Delta&\colon A\to A\otimes A &&
	\left\{\begin{array}{l}
		w_+ \mapsto \tensor{w+,w-} + \tensor{w-,w+}, \\
		w_- \mapsto \tensor{w-,w-},
	\end{array}\right.
	\\
	\eta&\colon \scalars\to A &&
	\left\{\begin{array}{l}
		1\mapsto w_+,
	\end{array}\right.
	\\
	\epsilon&\colon A\to\scalars &&
	\left\{\begin{array}{l}
		w_+ \mapsto 0, \\
		w_- \mapsto 1.
	\end{array}\right.
\end{align}
The~functor is graded if we set $\deg w_\pm := \pm1$. We shall write $\Kh(D)$ for the~homology of $\KhCom(D):=\FKh\KhBracket{D}$, where $D$ is a~diagram of a~link $L$; it is called the~\emph{Khovanov homology} of the~link $L$.

Khovanov's functor factorizes through $\BN(\R^2)$, where a~dot is understood as multiplication by $w_-$. In particular, both generators are images of $1\in\scalars$ under cup cobordisms
\begin{align}
	\FKh\left(\cupCob\right)\colon& 1\mapsto w_+ &
	\FKh\left(\cupCob*\right)\colon& 1\mapsto w_-
\end{align}
which motivates the~following graphical description of $\FKh$. Given a~\emph{collection of curves} $\Gamma\subset\R^2$ we identify $\FKh(\Gamma)$ with the~module generated freely by all diagrams obtained from $\Gamma$ by decorating some curves with dots and imposing the~relation that two dots on a~single curve annihilate the~diagram.%
\footnote{
	Clearly, those diagrams in which each circle carries at most one dot form a~free basis for $\FKh(\Gamma)$.
}
For example,
\begin{equation*}
		\FKh\left(\pictCirclesOnAPlane{0}\right) := \raisebox{2ex}{$%
			\Span_{\scalars}\left\{
				\pictCirclesOnAPlane{0},\ 
				\pictCirclesOnAPlane{1},\ 
				\pictCirclesOnAPlane{3},\ 
			\dots\right\}
		$}
		\!\!\Bigg/\!\!
		\raisebox{-2ex}{$%
			\begin{centerpict}(-2,1)(0,2)
				\psarc(0,0){2}{100}{150}
				\rput{115}(0,0){\psdot(2,0)}
				\rput{135}(0,0){\psdot(2,0)}
			\end{centerpict}
			= 0
		$}
\end{equation*}
The~generators $w_+$ and $w_-$ of the~algebra $A$ are represented by the~circle without and with a~dot respectively. To redefine $\FKh$ on a~\emph{cobordism} $S$ we use the~following rules:
\begin{itemize}
	\item if $S$ creates a~circle, then $\FKh(S)$ modifies a~diagram by inserting the~new circle with no dot on it,
	\item if $S$ contracts a~circle, then $\FKh(S)$ removes the~circle from a~diagram if it was decorated by a~dot, or takes the~diagram to $0$ otherwise, and
	\item we use the~following local surgery formulas to define $\FKh(S)$ if $S$ is a~merge or a~split
	\begingroup\psset{unit=2ex,dash=2pt 3pt}
	\begin{align}
		\label{surg:merge}
		\pictMergeCircles*&\quad\mapsto\quad \pictMergeCirclesAfter*
		\\
		\label{surg:split}
		\pictSplitCircle &\quad\mapsto\quad \pictSplitCircleAfter[10] \quad+\quad \pictSplitCircleAfter[01]
	\end{align}
	\endgroup
 	where the~blue thick arcs visualize the~saddle of $S$.
\end{itemize}
Notice that a~merge of two curves is zero, when each curve carries a~dot, as the~surgery \eqref{surg:merge} produces a~curve with two dots. Likewise, a~split of a~curve with a~dot results in one diagram, as one of the~two terms at the~right hand side of \eqref{surg:split} vanishes.

\subsection{Annular link homology}\label{sec:ann-kh}

There are two types of closed curves in the~annulus: \emph{trivial curves} bounding disks in $\Ann$, and \emph{essential curves}, parallel to the~core of $\Ann$. The~value of an~annular TQFT functor $\F\colon\aBN\to\Mod(\scalars)$ on trivial curves is determined by Bar-Natan's relations, but not the~value on essential curves.

The~first construction of an~annular TQFT functor is due to Asaeda, Przytycki, and Sikora \cite{APSHom}. The~APS functor $\FAKh\colon\aBN \to \Mod(\scalars)$ assigns to a~trivial and an~essential curve the~free modules
\begin{align}
	\label{def:ann-spaces}
	W &:= \Span_{\scalars}\{w_+,w_-\},
	&
	V &:= \Span_{\scalars}\{v_+,v_-\},
\\
\intertext{respectively, with the~degree defined on generators as}
	\label{def:ann-qdeg}
	&\deg w_\pm = \pm1,
	&
	&\deg v_\pm = 0.
\intertext{%
	This degree is denoted by $j'$ in \cite{KhHomSL2} and differs from the~one used in \cite{AnnHom}. In addition, the~modules admit the~\emph{annular grading}, denoted $\adeg$ and defined as%
}
	\label{def:ann-adeg}
	&\adeg w_\pm = 0,
	&
	&\adeg v_\pm = \pm1.
\end{align}
One can define $\FAKh$ by comparing it to $\FKh$. Indeed, $V$ and $W$ are
isomorphic as ungraded modules, but \eqref{def:ann-qdeg} induces a~filtration
on $\FKh$ and the~functor $\FAKh$ can be constructed as the~graded associate
\cite{AnnHom}. For completeness we write down the~maps corresponding to
the~elementary saddle moves. A~\emph{merge} is assigned one of the~maps
\begin{center}
	\begin{minipage}[t]{\widthof{$W\otimes W\to W$\quad}}
		\underline{$W\otimes W\to W$\quad}\newline
		$\tensor{w+,w+}\mapsto w_+$\newline
		$\tensor{w\pm,w\mp}\mapsto w_-$\newline
		$\tensor{w-,w-}\mapsto 0$
	\end{minipage}%
	\hskip 1cm
	\begin{minipage}[t]{\widthof{$V\otimes W\to V\quad$}}
		\underline{$V\otimes W\to V$\quad}\newline
		$\tensor{v\pm,w+}\mapsto v_\pm$\newline
		$\tensor{v\pm,w-}\mapsto 0$
	\end{minipage}
	\hskip 1cm
	\begin{minipage}[t]{\widthof{$V\otimes V\to W$\quad}}
		\underline{$V\otimes V\to W$\quad}\newline
		$\tensor{v\pm,v\pm}\mapsto 0$\newline
		$\tensor{v\pm,v\mp}\mapsto w_-$
	\end{minipage}%
\end{center}
depending on the~curves involved, whereas for \emph{splits} we choose one of
\begin{center}
	\begin{minipage}[t]{\widthof{$w_+\mapsto\tensor{w+,w-}+\tensor{w-,w+}$}}
		\underline{$W\to W\otimes W$\quad}\newline
		$w_-\mapsto\tensor{w-,w-}$\newline
		$w_+\mapsto\tensor{w+,w-}+\tensor{w-,w+}$
	\end{minipage}
	\hskip 5mm
	\begin{minipage}[t]{\widthof{$v_\pm\mapsto \tensor{v\pm,w-}$}}
		\underline{$V\to V\otimes W$\quad}\newline
		$v_\pm\mapsto \tensor{v\pm,w-}$
	\end{minipage}
	\hskip 1cm
	\begin{minipage}[t]{\widthof{$w_+\mapsto\tensor{v+,v-}+\tensor{v-,v+}$}}
		\underline{$W\to V\otimes V$\quad}\newline
		$w_-\mapsto 0$\newline
		$w_+\mapsto\tensor{v+,v-}+\tensor{v-,v+}$
	\end{minipage}
\end{center}
The~value of $\FAKh$ on caps and cups is unchanged.

%
% Graphical description
%
The~graphical description of $\FKh$ can be extended to the~annular case. Trivial curves can again carry dots, but the~essential ones cannot, because the~merge cobordism takes $\tensor{v\pm,w-}$ to zero. Therefore, we shall visualize the~two generators of $V$ by choosing an~orientation of the~essential curve, anticlockwise for $v_+$ and clockwise for $v_-$:
\begin{equation}
	\begin{centerpict}(-1.7,-1.7)(1.7,1.7)
		\pictDrawAnnulus(0,0){1.7}{0.3}
		\pscircle[dimen=middle](0,0){1}
		\rput{30}(0,0){\rput(1,0){\psline[arrowsize=6pt]{->}(-0.15pt,0pt)(0.15pt,3pt)}}
	\end{centerpict}
	\quad\to/<->/\quad v_+
	\hskip 1cm
	\begin{centerpict}(-1.7,-1.7)(1.7,1.7)
		\pictDrawAnnulus(0,0){1.7}{0.3}
		\pscircle[dimen=middle](0,0){1}
		\rput{30}(0,0){\rput(1,0){\psline[arrowsize=6pt]{->}(-0.15pt,0pt)(0.15pt,-3pt)}}
	\end{centerpict}
	\quad\to/<->/\quad v_-
\end{equation}
We use the~usual surgery formulas for merging a~trivial curve to an~essential one or splitting it off, keeping in mind that an~essential curve cannot carry dots:
\begin{equation}
	\psset{unit=0.75}%
	\pictAnnMergeTrivial[d] \:\mapsto\: \pictAnnMergeTrivialAfter[d]
		\hskip 3\psunit
	\pictAnnSplitOffTrivial \:\mapsto\: \pictAnnSplitOffTrivialAfter[1d]	
\end{equation}
Two essential curves can be merged together only if they have opposite orientations, in which case we decorate the~resulting trivial curve with a~dot, and otherwise we have zero:
\begin{equation}
	\psset{unit=0.75}%
	\pictAnnMergeEssentials[ud] \:\mapsto\: \pictAnnMergeEssentialsAfter[10] 
		\hskip 6\psunit
	\pictAnnMergeEssentials[uu] \:\mapsto\: \mathrlap 0 \hskip 3\psunit
\end{equation}
Finally, a~surgery from a~trivial curve to two essential ones is assigned the~map
\begin{equation}
	\psset{unit=0.75}%
	\pictAnnSplitTrivial \quad\mapsto\quad \pictAnnSplitTrivialAfter[ud] \quad+\quad \pictAnnSplitTrivialAfter[du]
\end{equation}
which can be viewed as performing the~surgery on the~trivial curve considered with both orientations at the~same time.

%
% Action of SL2
%
\subsubsection*{Action of $\LieSL_2$}

It has been recently observed in \cite{KhHomSL2} that the~annular link homology
admits an~action of $\LieSL_2$ if we consider $\trivmod$ as a~trivial
representation and $\essmod$ is identified with the~fundamental one
$V_1   = \Span_{\scalars}\{v_+, v_-\}$ or its dual
$V_1^* = \Span_{\scalars}\{v_+^*, v_-^*\}$,
depending on the~nestedness of the~associated essential curve,
i.e.\ $V_1$ and $V_1^*$ are assigned alternatively. The~following
tables describe the~action of $\LieSL_2$.
\begin{equation}
	\begin{tabular}{llp{1cm}ll}
		\multicolumn{2}{c}{The~action on $V_1$} && \multicolumn{2}{c}{The~action on $V_1^*$}\\
		\cline{1-2} \cline{4-5}
		  $Ev_+ = 0$   & $Ev_- = v_+$ && $Ev_-^* = 0$      & $Ev_+^* = -v_-^*$ \\
	    $Fv_+ = v_-$ & $Fv_- = 0$   && $Fv_-^* = -v_+^*$ & $Fv_+^* = 0$
	\end{tabular}
\end{equation}
There is an~obvious isomorphism of $\LieSL_2$--modules $V_1\cong V_1^*$, which identifies $v_\pm$ with $\pm v_\mp^*$. However, the~action on the~annular chain complex is defined using instead the~linear isomorphism $V \cong V_1^*$ that sends $v_\pm$ to $v_\mp^*$, so that the~action depends on the~position of $V$ in the~tensor product. For instance, two essential curves are assigned $V\otimes V \cong V_1^*\otimes V_1$ with $\LieSL_2$ acting in the~following way:
\begin{align*}
	E(\tensor{v+,v+}) &= \phantom-0              &
	F(\tensor{v+,v+}) &= \phantom-\tensor{v-,v+}-\tensor{v+,v-} \\
	E(\tensor{v+,v-}) &= -\tensor{v+,v+}         &
	F(\tensor{v+,v-}) &= \phantom-\tensor{v-,v-} \\
	E(\tensor{v-,v+}) &= \phantom-\tensor{v+,v+} &
	F(\tensor{v-,v+}) &= -\tensor{v-,v-}         \\
	E(\tensor{v-,v-}) &= \phantom-\tensor{v+,v-}-\tensor{v-,v+} &
	F(\tensor{v-,v-}) &= \phantom-0
\end{align*}
It follows that the~maps $V\otimes V\to W$ and $W\to V\otimes V$ intertwine the~action and the~annular TQFT functor is upgraded to $\FAKh\colon\aBN\to\gRep(\LieSL_2)$. In particular, $\LieSL_2$ acts on the~triply graded annular homology.

\begin{remark*}
The~action admits the~following graphical description: each clockwise oriented curve in a~diagram $w$ contributes to $Ew$ a~diagram obtained from $w$ by reversing the~curve, and scaling it by $(-1)$ if it is separated from the~outer boundary by an~odd number of curves. Likewise for $F$ we reverse orientations of anticlockwise oriented curves.
\end{remark*}

\subsection{Homology for links in a~thickened M\"obius band}
\label{sec:khov-hom-mobius}

Beyond links in a~thickened annulus, we also consider links in the~twisted line
bundle over the~M\"obius band (``twisted'' means that the~monodromy along
the~orientation reversing curve is $-id$, so that the~bundle is an~orientable
3-manifold). Let us recall the~construction of the~APS functor in this case.

A~M\"obius band $\MB$ admits three types of curves: \emph{trivial curves} bounding disks, \emph{separating curves} cutting an~annulus out of $\MB$, and \emph{nonseparating ones}. The~APS functor $\FMKh\colon \BN(\MB)\to \Mod(\scalars)$ assigns to them the~following free modules
\begin{equation}
	\trivmod := \Span_{\scalars}\{w_+,w_-\},
\quad
	\sepmod := \Span_{\scalars}\{v_+, v_-\},
\quad\text{and}\quad
	\nsepmod := \Span_{\scalars}\{u_+,u_-\},
\end{equation}
with the~degree function vanishing on both $V$ and $U$, and $\deg w_\pm = \pm1$ as usual. There are more types of saddle cobordisms in $\MB\times I$ than in $\Ann\times I$, and $\FMKh$ vanishes on those without trivial circles in the~boundary. Otherwise, it is defined as in the~annular case for merges and splits (where both $V$ and $U$ can play the~role of the~``annular'' $V$) and one of the~maps
\begin{center}
	\begin{minipage}[t]{\widthof{$w_+\mapsto v_+ + v_-$\quad}}
		\underline{\rlap{$W\to V$}\phantom{$w_+\mapsto v_+ + v_-$}}\newline
		$w_+\mapsto v_+ + v_-$\newline
		$w_-\mapsto 0$
	\end{minipage}
	\hskip 1cm
	\begin{minipage}[t]{\widthof{$V\to W$\quad}}
		\underline{$V\to W$\quad}\newline
		$v_\pm\mapsto w_-$
	\end{minipage}
\end{center}
for a~saddle cobordism between a~trivial circle and a~separating one. We can represent the~latter graphically as
\begin{gather}
	\psset{unit=0.75}
	\pictMBBoundingToTrivial
	\quad\mapsto\quad
	\pictMBBoundingToTrivialAfter[uu]
	\quad+\quad
	\pictMBBoundingToTrivialAfter[dd]
\intertext{and}
	\psset{unit=0.75}
	\pictMBTrivialToBounding[uu]
	\quad,\quad
	\pictMBTrivialToBounding[dd]
	\quad\mapsto\quad
	\pictMBTrivialToBoundingAfter[1-]
\end{gather}

We shall write $\MKhCom(D) := \FMKh\KhBracket{D}$ for the~chain complex for a~link diagram $D$ on $\MB$, and $\MKh(D)$ for its homology. A~resolution of $D$ can have at most one non-separating curve---such a~curve cuts the~band into an~annulus. In particular, $U$ only appears in $\MKhCom(D)$ when $D$ meets any cross section of $\MB$ in an~odd number points. If so, $\MKhCom(D) \cong \overline\KhCom_{\MB}(D) \otimes U$, where we write $\overline\KhCom_{\MB}(D)$ for the~chain complex of $D$ computed with $\scalars$ assigned to non-separating curves instead of $U$.

\subsection{Chen--Khovanov homology for tangles}\label{sec:chen-khov}

Tangle diagrams in a~thickened stripe $\RxI$ form a~category. Thence, the~formal Khovanov bracket of a~tangle is a~chain complex built over the~bicategory $\CCob = \Cob(\RxI)$ of points on a~line, flat tangles in a~stripe $\RxI$, and surfaces in $(\RxI)\times I$. This bicategory is graded with the~degree of a~surface $S$ defined in \eqref{def:degree-in-Cob}. To preserve this richer structure, the~homology for tangles is constructed by Chen and Khovanov \cite{ChenKhov} using a~2--functor $\FFCK\colon \CCob\to \ccat{gBirep}$ valued in the~bicategory of graded bimodules that are \emph{sweet}, i.e.\ finitely generated and projective
as left and as right modules, but not necessarily as bimodules.
In fact, $\FFCK$ factors through the~Bar-Natan skein bicategory $\BBN = \BN(\RxI)$.
We begin with describing the~modules assigned to tangles, then the~algebras assigned to points, and finally reconstructing the~bimodule structure.

%
%  Cup diagrams with platforms
%
\subsubsection*{Cup diagrams with platforms}

\begin{figure}[ht]
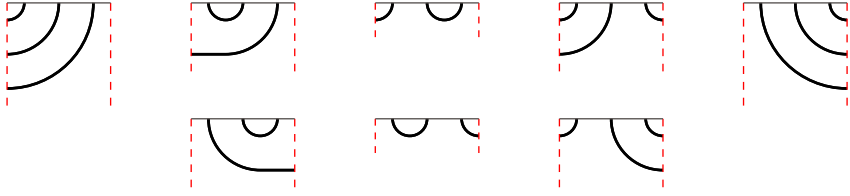
%
	\centering
	\begin{tabular}{*5{p{2cm}}}
		\planarmatching{|000|} &
		\planarmatching{|100|} &
		\planarmatching{|010|} &
		\planarmatching{|001|} &
		\planarmatching{|111|} \\ &
		\planarmatching{|110|} &
		\planarmatching{|101|} &
		\planarmatching{|011|}
	\end{tabular}\qquad
	\caption{The~generalized cup diagrams with three termini.}%
	\label{fig:cup-and-platforms-example}%
\end{figure}

A~\emph{crossingless matching} between $2n$ points in a~line is a~collection of $n$ disjoint arcs attached to the~points. We shall draw the~arcs in the~lower half-plane $\R\times\R_-$ and refer to them as a~\emph{cup diagram}. Following \cite{ChenKhov} we generalize cup diagrams to allow semi-infinite arcs, each attached to one point only and going left or right towards infinity. This can be visualized by drawing two vertical platforms going out of the~horizontal line, one to the~left and one to the~right of all the~points, and attaching semi-infinite arcs to them. In particular, odd number of points are allowed. We shall call the~points on the~line \emph{termini} to distinguish them from the~endpoints on the~platforms. Figure~\ref{fig:cup-and-platforms-example} presents all cup diagrams with three termini.

%
%  Figure: examples of cup diagrams
%

%
%  Weight of a diagram
%
Let $\gmatchings n$ be the~set of such diagrams with $n$ termini. We define the~\emph{weight} of a~diagram $a\in\gmatchings n$ as $\weight(a) := r-\ell$, where $r$ and $\ell$ count respectively the~arcs terminating on the~right and on the~left platform. In what follows we shall write $\gmatchings[\lambda]{n} \subset \gmatchings{n}$ for the~subset of diagrams of weight $\lambda$. Notice that $\gmatchings[\lambda]{n}$ is empty unless $\lambda$ has the~same parity as $n$.

%
%  Cap diagrams
%
Dually we define the~set $\gmatchings{,n}$ of \emph{cap diagrams with platforms} with arcs drawn in the~upper half-plane. The~reflection along the~horizontal line induces a~bijection of sets
\begin{equation}
	\gmatchings[\lambda]{n} \ni a\mapsto \revmatching a\in\gmatchings[\lambda]{,n}
\end{equation}
for every $n$ and $\lambda$.

%
% Extending Khovanov's TQFT to diagrams with platforms
%
\subsubsection*{An~extension of $\FKh$}

\begin{figure}[ht]
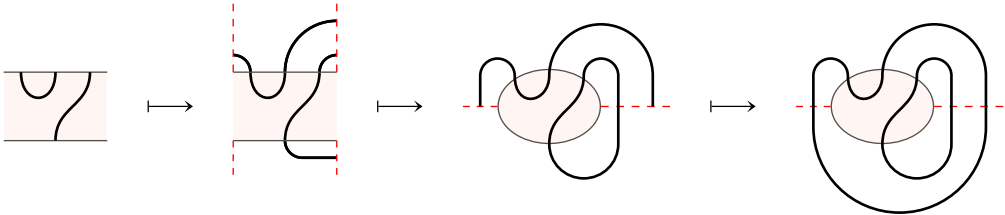
%
	\begin{equation*}
		\pictTangle*
		\quad\mapsto\quad	\pictTangle*<|.1.|><|011|>
		\quad\mapsto\quad	\pictTangleDisk*
		\quad\mapsto\quad	\pictTangleDisk
	\end{equation*}
	\caption{The~construction of a~planar closure of a~$(1,3)$--tangle.}%
	\label{fig:tangle-closure}%
\end{figure}

A~pair of cup diagrams $a\in\gmatchings[\lambda]{m}$ and $b\in\gmatchings[\lambda]{n}$ can be used to produce a~planar closure $\revmatching bTa$ of any flat tangle $T\in\CCob(m,n)$. The~closure is constructed by gluing $a$ to the~bottom of $T$ and $\revmatching b$ to the~top, then turning the~platforms towards themselves and identifying the~endpoints of the~arcs from inside out. In case $m\neq n$ there will be unmatched endpoints, the~same number at each side, because $a$ and $b$ have equal weights. We connect them with half-circles, see Figure~\ref{fig:tangle-closure}.

%
%  Figure: constructing a closure of a tangle
%

\begin{figure}[ht]
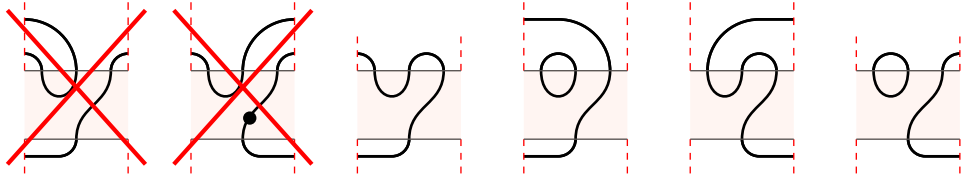
%
	\centering
	\pictTangle*<|.0.|><|001|>%
	\psline[linecolor=red,linewidth=2\pslinewidth](-3.5,-1.5)( 0.5,3)
	\psline[linecolor=red,linewidth=2\pslinewidth]( 0.5,-1.5)(-3.5,3)
	\qquad
	\pictTangle*<|.1.|><|011|>%
	\psdot[style=diagarc](-1.3,-0.15)
	\psline[linecolor=red,linewidth=2\pslinewidth](-3.5,-1.5)( 0.5,3)
	\psline[linecolor=red,linewidth=2\pslinewidth]( 0.5,-1.5)(-3.5,3)
	\qquad
	\pictTangle*<|.0.|><|010|>\qquad
	\pictTangle*<|.0.|><|100|>\qquad
	\pictTangle*<|.1.|><|110|>\qquad
	\pictTangle*<|.1.|><|101|>
	\caption{%
		Examples of generators of $\FFCK\left(\fntTangle\right)$. The~first two
		diagrams are nonadmissible, because they contain either a~turnback or
		a~dot on an~open arc.}%
	\label{fig:arc-tangle}%
\end{figure}

%
%  Conditions for diagrams generating F(curves)
%
The~functor $\FKh\colon\Cob(\R^2)\to \Mod(\scalars)$ is extended to collections of curves with platforms by assigning to such a~collection a~module generated by all possible decorations of the~curves with dots as before, but with more restrictions:
\begin{enumerate}
	\item a~diagram vanishes when it contains a~curve intersecting any of the~platforms twice,
	\item a~dot annihilates a~diagram when placed on a~curve that intersects a~platform, and
	\item as before, two dots on one curve annihilate the~diagram.
\end{enumerate}
A~diagram is \emph{nonadmissible} if one of the~above situations happens,
see Figure~\ref{fig:arc-tangle}. Surgeries \eqref{surg:merge} and \eqref{surg:split}
on nonadmissible diagrams produce nonadmissible ones, so that $\FKh(S)$ is
well-defined for any surface $S$, see also \cite{ChenKhov}. The~Chen--Khovanov
functor assigns to a~flat tangle $T\in\CCob(m,n)$ the~module
\begin{equation}
	\FFCK(T) := \bigoplus_\lambda \FFCK(T;\lambda),
\quad\text{with}\quad
	\FFCK(T;\lambda) := \bigoplus_{\substack{a\in\gmatchings[\lambda]{m}\\ b\in\gmatchings[\lambda]{n}}}\FKh(\revmatching bTa),
\end{equation}
and to a~cobordism $S$ between flat tangles $T_0$ and $T_1$ the~linear map
\begin{equation}\label{eq:FFCK-for-surfaces}
	\FFCK(S) := \bigoplus_\lambda \FFCK(S;\lambda),
\quad\text{with}\quad
	\FFCK(S;\lambda) := \bigoplus_{\substack{a\in\gmatchings[\lambda]{m}\\ b\in\gmatchings[\lambda]{n}}}\FKh(\revmatching bSa),
\end{equation}
where $\revmatching bSa$ stands for the~surface
$(\revmatching b\times I)\cup S\cup(a\times I)$.

\begin{remark*}
	The~closures of \fntTangle\ in Figure~\ref{fig:arc-tangle} are drawn without
	identifying the~platforms of cup diagrams. Not only makes this smaller diagrams,
	but also easier to describe the~module structure on $\FFCK(T)$ once
	the~Chen--Khovanov algebras are introduced.
\end{remark*}

% ==================================================================
%  Examples
% ==================================================================
\begin{example}\label{ex:ChKh-V-to-C}
	Consider the~saddle cobordism $S:=\fntVTCob\colon\fntVertTan\to\fntCupCap$ between the~identity $(2,2)$--tangle and the~tangle consisting of a~cap followed by a~cup. The~module $\FFCK(\fntVertTan)$ has seven generators, on which $\FFCK(S)$ takes the~following values:
	\begin{align*}
		\pictVertTan*<|11|><|11|> &\:\mapsto\: 0 &
		\pictVertTan*<|01|><|10|>      &\:\mapsto\: \pictCupCap*<|01|><|10*={60}|> &
		\pictVertTan*<|10|><|10*={60}|>&\:\mapsto\: \pictCupCap*<|10*={240}|><|10*={60}|> &
		\pictVertTan*<|01|><|01|>      &\:\mapsto\: \pictCupCap*<|01|><|01|>
		\\[2ex]
		\pictVertTan*<|00|><|00|> &\:\mapsto\: 0 &
		\pictVertTan*<|10|><|01|> &\:\mapsto\: \pictCupCap*<|10*={240}|><|01|> &
		\pictVertTan*<|10|><|10|> &\:\mapsto\: \pictCupCap*<|10|><|10*={60}|>\mathrlap{\ +\ \pictCupCap*<|10*={240}|><|10|>}
	\end{align*}
	For example, the~top right component of $\FFCK(S)$ is a~merge when the~platforms are identified:
	\begin{equation*}
		\pictVTCobClosed \colon \pictVertTanClosed \;\mapsto\; \pictCupCapClosed
	\end{equation*}
\end{example}

\begin{example}\label{ex:ChKh-C-to-V}
	Consider now the~cobordism $S:=\fntCCCob\colon\fntCupCap\to\fntVertTan$ going in the~other direction. The~module $\FFCK(\fntCupCap)$ has eight generators, on which $\FFCK(S)$ is defined as below:
	\psset{unit=2.5ex}%
	\begin{align*}
		\pictCupCap*<|10|><|10|> &\:\mapsto\: \pictVertTan*<|10|><|10|> &
		\pictCupCap*<|01|><|10|> &\:\mapsto\: \pictVertTan*<|01|><|10|> &
		\pictCupCap*<|10*={240}|><|10*={60}|>,\ \pictCupCap*<|01|><|01|> &\:\mapsto\: 0
		\\[1ex]
		\pictCupCap*<|10|><|10*={60}|>,\  \pictCupCap*<|10*={240}|><|10|> &\:\mapsto\: \pictVertTan*<|10|><|10*={60}|> &
		\pictCupCap*<|10|><|01|> &\:\mapsto\: \pictVertTan*<|10|><|01|> &
		\pictCupCap*<|01|><|10*={60}|>,\  \pictCupCap*<|10*={240}|><|01|> &\:\mapsto\: 0
	\end{align*}
	The~two arcs in the~top right corner are mapped to zero, because the~corresponding closure of $S$ in this case is a~split with each circle in its output touching a~platform:
	\begin{equation}\label{eq:merge-of-arcs-vanishes}
		\pictCCCobClosed \colon \pictCupCapClosed \;\mapsto\; \pictVertTanClosed[*-] \;+\; \pictVertTanClosed[-*] = 0.
	\end{equation}
\end{example}

\subsubsection*{Arc algebras and diagrammatic bimodules}

Let $c\in\gmatchings{n}$ be a~generalized cup diagram and write $S_c$ for
the~cobordism from $c\sqcup \revmatching c$ to $2n$ vertical lines obtained by
a~sequence of $n$ surgeries, one per arc in $c$, see Figure~\ref{fig:n-surgeries}.
The~collection of such cobordisms defines linear maps
\begin{equation}\label{eq:CK-gluing-map}
	\mu_{T',T}\colon \FFCK(T) \otimes \FFCK(T') \to \FFCK(T'T),
\end{equation}
one per a~pair of tangles $T\in\CCob(m,n)$ and $T'\in\CCob(n,k)$. Explicitly,
$\mu_{T',T}(x\otimes y)=0$ for $x\in\FFCK(\revmatching bTa)$ and
$y\in\FFCK(\revmatching dT'c)$ unless $b=c$, in which case $\mu_{T',T} =
\FKh(\revmatching d T^! S_c T a)$.

\begin{figure}[ht]%
	\centering
	\psset{unit=2ex,style=diagarc}%
	\begin{movie}[hc](8,6){4}
		\movieclip{%
			\surgeryline(-2,-1.5)(-2,1.5)
			\psarc(-2, 2){0.5}{180}{360}
			\psarc( 1, 2){0.5}{180}{360}
			\psarc( 1, 2){1.5}{180}{360}
			\psarc(-2,-2){0.5}{0}{180}
			\psarc( 1,-2){0.5}{0}{180}
			\psarc( 1,-2){1.5}{0}{180}
		}
		\movieclip{%
			\psline(-2.5,-2)(-2.5,2)
			\psline(-1.5,-2)(-1.5,2)
			\surgeryline(1,-0.5)(1,0.5)
			\psarc( 1, 2){0.5}{180}{360}
			\psarc( 1, 2){1.5}{180}{360}
			\psarc( 1,-2){0.5}{0}{180}
			\psarc( 1,-2){1.5}{0}{180}
		}
		\movieclip{%
			\psline(-2.5,-2)(-2.5,2)
			\psline(-1.5,-2)(-1.5,2)
			\psline(-0.5,-2)(-0.5,2)
			\psline( 2.5,-2)( 2.5,2)
			\surgeryline(1,-1.5)(1,1.5)
			\psarc( 1, 2){0.5}{180}{360}
			\psarc( 1,-2){0.5}{0}{180}
		}
		\movieclip{%
			\psline(-2.5,-2)(-2.5,2)
			\psline(-1.5,-2)(-1.5,2)
			\psline(-0.5,-2)(-0.5,2)
			\psline( 2.5,-2)( 2.5,2)
			\psline( 0.5,-2)( 0.5,2)
			\psline( 1.5,-2)( 1.5,2)
		}
	\end{movie}
	\caption{A~sequence of surgeries replacing a~disjoint union of a~cup diagram and its vertical flip with vertical lines.}%
	\label{fig:n-surgeries}%
\end{figure}

The~\emph{Chen--Khovanov algebra} $\CKalg n$ is the~module assigned to the~tangle
formed by $n$ vertical lines, with $x\cdot y := \mu(x\otimes y)$.
It admits a~\emph{weight decomposition}
\begin{equation}
	\CKalg n = \bigoplus_\lambda \CKalg{n}(\lambda),
\end{equation}
which is related to that from \cite{ChenKhov} by setting $\CKalg{n-k,k} = \CKalg{n}(n-2k)$. There is a~unique primitive idempotent $e_c\in \CKalg n$ for each closure $c\in\gmatchings{n}$ given by the~diagram $\revmatching cc$ with no dots. The~idempotents are mutually orthogonal, and their sum is a~unit in $\CKalg n$.

\begin{example}
	The~algebra $\CKalg 2$ has generators in weights $-2$, $0$, and $2$. Both $\CKalg 2(-2)$ and $\CKalg 2(2)$ are one dimensional, whereas $\CKalg 2(0)$ has five generators:
	\begin{equation}
		\pairofmatchings{|01|}{|01|}\qquad
		\pairofmatchings{|10|}{|10|}\qquad
		\pairofmatchings{|01|}{|10|}\qquad
		\pairofmatchings{|10|}{|01|}\qquad
		\pairofmatchings{|10|}{|10*|}
	\end{equation}
	of which the~first two are idempotents and the~other square to zero. Furthermore,
	\begin{equation}
		\pairofmatchings{|10|}{|01|}
		\ \cdot\ 
		\pairofmatchings{|01|}{|10|}
		\ = \
		\pairofmatchings{|10|}{|10*|}
	\qquad\text{and}\qquad
		\pairofmatchings{|01|}{|10|}
		\ \cdot\ 
		\pairofmatchings{|10|}{|01|}
		\ = \
		0.
	\end{equation}
\end{example}

The~product in $\CKalg n$ can be described explicitly using generalized surgeries as in \cite{ChenKhovIsCellular}. Because we do not identify platforms when drawing diagrams, each diagram has four platforms drawn vertically. The~product $x\cdot y$, when nonzero, can be then computed graphically by placing $y$ on $x$, connecting the~platforms in between, and following the~two steps below.

\smallskip\noindent\textit{Step I: surgeries at platforms.}
Replace two opposite arcs touching one of the~platforms with a~vertical line and decorate with a~dot each closed loop created that way:
\begin{equation}\label{srg:CK-at-platform}
	\psset{unit=1.8ex}%
		\pictMergeAtPlatform \mapsto \pictMergeAtPlatformAfter
	\hskip 15mm
		\pictMergeArcsAtPlatform \mapsto \pictMergeArcsAtPlatformAfter
	\hskip 15mm
		\pictSplitAtPlatform \mapsto \pictSplitAtPlatformAfter[1]
\end{equation}

\smallskip\noindent\textit{Step II: surgeries on half-circles.}
When no arc at inner platforms is left, perform surgeries on the~remaining arcs using the~usual surgery formulas \eqref{surg:merge} and \eqref{surg:split}, except that a~merge of two open arcs is zero and a~diagram with a~dot on an~open arc vanishes (in particular, the~first term in the~result of the~second surgery may vanish):
\begin{equation}\label{srg:CK-inside}
	\psset{unit=1.8ex}%
		\pictMergeToArc \mapsto \pictMergeToArcAfter
	\hskip 1cm
		\pictSplitFromArc \mapsto \pictSplitFromArcAfter[10] + \pictSplitFromArcAfter[01]
	\hskip 1cm
		\pictMergeArcs \mapsto 0.
\end{equation}
The~merge of two arcs vanishes in the~second step, because the~two arcs have endpoints on the~outer platforms and they belong to the~same circle when the~platforms are identified (compare with \eqref{eq:merge-of-arcs-vanishes}).

It follows from the~construction that the~maps $\mu_{T',T}$ are natural with respect to tangle cobordisms.
Furthermore, \eqref{eq:FFCK-for-surfaces} can be extended to dotted surfaces,
because $\FFCK$ is defined by $\FKh$, in which case $\FFCK(S) = 0$ if any
component of $S$ carries two dots. Thence, the~following result holds,
see also \cite{ChenKhov}.

\begin{proposition}
	$\FFCK(T)$ is a~sweet $(\CKalg n, \CKalg m)$--bimodule
	for any flat tangle $(m,n)$--tangle $T$, where the~actions of
	the~algebras are given by $\mu$. Moreover, \eqref{eq:CK-gluing-map}
	descend to natural isomorphisms of bimodules
	\begin{equation}\label{eq:CKmod-tensor-isom}
		\FFCK(T) \CKtimes{n} \FFCK(T') \to^\cong \FFCK(T'T),
	\end{equation}
	so that there is a~strong bifunctor $\FFCK\colon \BBN\to \ccat{Birep}$.%
	\footnote{%
		Recall that $M\circ N := N\otimes M$ in $\ccat{Birep}$.
	}
\end{proposition}

Throughout the~paper we call $\FFCK(T)$ a~\emph{diagrammatic bimodule}.
They are called \emph{geometric} in \cite{ChenKhov}.
Each $\FFCK(T)$ has a~two sided dual $\FFCK(T^!)$, where $T^!$ is the~vertical
flip of $T$. It is also known that each weight component $\FFCK(T; \lambda)$ is
indecomposable when $T$ contains no loops \cite[Theorem 4.14]{ChenKhovKoszulity},
and otherwise $\FFCK(T; \lambda) \cong \FFCK(\widetilde T; \lambda)\otimes
\scalars^{2\ell}$, where $\widetilde T$ is the~tangle $T$ with $\ell$ loops removed.
Therefore, the~category of diagrammatic bimodules has duals and is closed under
direct summands.

\subsubsection*{Grading}

The~grading on $\CKalg n$ is defined in \cite{ChenKhov} by shifting by $n$
the~grading induced by the~functor $\FKh$. This does not work well for bimodules
assigned to tangles, though. For instance, regarding a~cup diagram
$c\in\gmatchings{n}$ as a~flat $(0,n)$--tangle, there is an~isomorphism of graded
bimodules $\FFCK(c) \otimes \FFCK(\revmatching c)\cong e_c\CKalg n e_c$,
but $e_c\CKalg n e_c$ is graded differently from $\FFCK(\revmatching cc)$.

In \cite{ChenKhovIsCellular, ChenKhovKoszulity} a~grading is computed differently.
It agrees with that from \cite{ChenKhov} for arc algebras, and it is coherent with
tensor products of tangle bimodules. On the~other hand, it depends on a~Morse
decomposition of a~tangle, although it is well understood how the~degree changes
under planar isotopies, see \cite[Lemma 2.4]{ChenKhovKoszulity}.
This motivates the~following definition.

Let $T$ be a~flat tangle with $\ell$ loops and $c$ arcs connecting bottom endpoints. Given a~diagram $x\in\FFCK(T)$ orient all its curves counterclockwise. The~platforms and boundary lines of $T$ split some curves of $x$ into vertical lines, caps, and cups. Let $a$ be the~number of cups and caps with clockwise orientation. Then the~degree of $x$ is given by the~formula,
\begin{equation}\label{def:CK-degree}
	\deg x := \ell + c - a - 2d,
\end{equation}
where $d$ is the~number of dots. For example,
\begin{equation}
	\diagarcset{arccolor=gray,arcwidth=0.5pt}%
	\psset{style=diagarc}%
	\deg\left(\begin{centerpict}(-3.2,-2)(3.2,2)
		\drawpairofmatchings(0,0){|101001|}{|0110*10|}
		\psarc[linecolor=black,linewidth=1.5pt]( 2,0){0.5}{0}{180}
		\psarc[linecolor=black,linewidth=1.5pt](-2,0){0.5}{180}{360}
	\end{centerpict}\right) = -4,
\qquad\text{and}\qquad
	\deg\left(\begin{centerpict}(-0.4,-3.5)(4.4,2.5)
		\psframe[style=mysurface,linestyle=none](0,1.5)(4,-1.5)
		\drawtangleclosure(0,0)<3>{|1001|}{%
			\psbezier(0.5,-1.5)(0.5,0.5)(1.5,0.5)(1.5,1.5)
			\psarc(2,-1.5){0.5}{0}{180}
			\psbezier(2.5,1.5)(2.5,0.5)(3.5,0.5)(3.5,-1.5)
			\pscircle(2,0){0.4}
		}{|.01.|}%
		\psarc[linecolor=black,linewidth=1.5pt]( 1,-1.5){0.5}{180}{360}
	\end{centerpict}\right) = 1,
\end{equation}
where the~arcs with clockwise orientation are thickened. A~quick look on the~surgery formulas \eqref{surg:merge} and \eqref{surg:split} reveals that for any tangle cobordism $S$ the~map $\FFCK(S)$ is homogeneous of degree $\deg S$. Furthermore, $\mu_{T',T}\colon \FFCK(T)\otimes\FFCK(T')\to \FFCK(T'T)$ preserves the~degree, see \cite[Theorem 3.5 (iii)]{ChenKhovKoszulity}, so that \eqref{eq:CKmod-tensor-isom} is an~isomorphism of graded bimodules. However, $\FFCK(T^!)$ is dual to $\FFCK(T)$ only up to a~degree shift.

\subsubsection*{Chain complex and homology}

Assume now that $T$ is an~oriented tangle with $m$ points at the~bottom and $n$
at the~top. The~\emph{Chen--Khovanov complex} $\CKmod{T} := \FFCK\KhBracket{T}$ is
a~chain complex of graded sweet $(\CKalg m, \CKalg n)$--bimodules, obtained from
the~formal bracket by applying the~bifunctor $\FFCK$ component-wise. We refer to
the~homology $\ChKh^\bullet(T) := H_{\!\bullet}(\CKmod T)$ as the~\emph{Chen--Khovanov homology}. It is a~triply graded theory: beyond the~homological and quantum grading $\ChKh^\bullet(T)$ admits a~weight decomposition.

\subsubsection*{Decategorification}

Let $V$ be the~fundamental representation of $\Uqsl 2$. It is known
that flat tangles can be intrepreted as intertwiners between tensor
powers of $V$, see Appendix~\ref{sec:sl2-rep}. More precisely, there
is a~linear Temperly--Lieb category $\TL$, the~morphisms of which are
generated by flat tangles, and a~functor $\F_{TL}\colon \TL\to \Rep(\Uqsl 2)$
that takes a~collection of $n$ points to $V^{\otimes n}$, whereas a~cap and a~cup
to the~evaluation and coevalution map \eqref{eq:ev-coev-for-V} respectively.
The~Chen--Khovanov construction categorifies this functor.

\begin{theorem}[cf.\ \cite{ChenKhov}]\label{thm:ChenKhov-categ-TL}
	There are isomorphisms\/
	$\gamma_n\colon K_0(\CKalg n)\utimes\Zq \scalars \to^\cong V^{\otimes n}$\!,
	such that\/ $\gamma_n\circ [\FFCK(T)]\circ\gamma_m^{-1} = \F_{TL}(T)$
	for any flat $(m,n)$--tangle $T\noic$.
\end{theorem}

Let us briefly recall from \cite{ChenKhov} how $\gamma_n$ is constructed.
The~group $K_0(\CKalg n)$ is generated freely by indecomposable projectives
$P_a := e_a \CKalg n$, one projective for each generalized cup diagram
$a\in\gmatchings{n}$.
On the~other hand, cup diagrams are flat tangle.
The~isomorphism $\gamma_n$ is defined by
$\gamma_n([P_a]) := \F_{TL}(a)(v_-^{\otimes\ell} \otimes v_+^{\otimes r})$,
where $\ell$ and $r$ count points of $a$ on the~left and right platform
respectively.

\section{Quantization of the~annular link homology}
\label{sec:quantum}

This section discusses the~construction and properties of the~quantum annular
link homology. We start with redefining the~action of $\LieSL_2$ in a~way that
motivates the~search for the~quantization. Section~\ref{sec:qnt-annulus} contains
a~detailed construction of the~quantum annular TQFT, which is then applied to
(a~quantization of) the~formal bracket to obtain the~new invariant.

\subsection{The~action of \texorpdfstring{$\LieSL_2$}{sl2} revisited}

The~action of $\LieSL_2$ on the~annular link homology can be understood
already at the~level of the~skein category.
Consider now an~operation that takes a~flat tangle $T$ into the~surface
$\S\times T \subset \Ann\times I$. A~closed loop in $T$ corresponds to
a~toroidal component of the~surface, which is evaluated to 2 in $\BN(\Ann)$
due to the~neck cutting relation. Hence, there is a~well-defined
functor $\S\times(\blank)\colon\TL|_{q=1} \to \BN(\Ann)$, where
$\TL|_{q=1}$ is the~Temperly--Lieb category specialized at $q=1$,
see Appendix~\ref{sec:sl2-rep}.
In fact, it takes values in $\BoBN(\Ann)$, the~quotient of $\BN(\Ann)$
by the~Boerner's relation that forces a~dot to annihilate a~connected
surface with an~essential circle in its boundary.

\begin{remark*}
	Because $\BoBN(\Ann)$ is both graded and additive, from now on we make
	two modification to $\TL$: we introduce a~formal degree shift despite all
	morphisms in $\TL$ having degree 0, and formal direct sums, so that $\TL$
	becomes a~graded additive category. The~functor $\FTL$ extends
	naturally to a~faithful functor $\FTL\colon \TL\to\gRep(\Uqsl 2)$	
	and $\S\times(\blank)$ is still well-defined when $q=1$.
\end{remark*}

\begin{proposition}\label{prop:BBN-is-TL}
	The~functor\/ $\S\times(\blank) \colon \TL|_{q=1}\to\BoBN(\Ann)$ is an~equivalence of categories.
\end{proposition}
\begin{proof}
	By the~Delooping Lemma, each object in $\BoBN(\Ann)$ is isomorphic
	to a~collection of essential curves with shifted degree,
	whereas the~neck cutting relation implies that morphisms between such
	collections are generated by incompressible surfaces,
	i.e.\ annuli \cite{BN-for-sfces}. These are graded morphisms only when
	the~degree shifts of the~collection of curves at the~bottom and top agree.
	Hence, the~functor $\S\times(\blank)$ is full and essentially surjective.
	Faithfulness follows, because the~annuli are linearly independent in
	$\BoBN(\Ann)$, see \cite{AnnBN}.
\end{proof}

Specializing $q=1$ makes $\FTL$ valued in $\gRep(\LieSL_2)$.
We check directly that the~triangle of functors
\begin{equation}
	\begin{diagps}(0em,-0.5ex)(8em,11.5ex)
		\Vtriangle<8em,9ex>[%
			\TL|_{q=1}`\BoBN(\Ann)`\gRep(\LieSL_2);
			\S\times(\blank)`\FTL`\FAKh
		]
	\end{diagps}
\end{equation}
commutes. This equips the~annular TQFT with an~action of $\LieSL_2$
that coincides with the~one from~\cite{KhHomSL2}.

\subsection{Deformation of the~annular skein category}
\label{sec:qnt-annulus}

Our goal is to quantize the~annular skein category,
so that Proposition~\ref{prop:BBN-is-TL} holds for all values
of $q$. We achieve this by identifying $\aBN$ with
the~additive closure of the~horizontal trace of the~bicategory $\BBN$
(taking the~additive closure is necessary, because $\hTr(\BBN)$ is only
linear), then deforming the~trace relation.
Because the~horizontal traces are compared with Bar-Natan skein categories,
it is understood that they are made graded by introducing formal degree shifts.

\begin{proposition}
	There is an~equivalences of categories $\hTr^\oplus(\BBN) \simeq \BN(\Ann)$.
\end{proposition}
\begin{proof}
	The~proof follows the~same argument as the~one of Theorem~\ref{thm:hTr-of-FxI}.
\end{proof}

The~bicategory $\BBN$ is locally pregraded and each cobordism $S$ is
a~homogeneous morphism of degree
\begin{equation}
	\deg S = \chi(S) - \frac{\#B}{4} - 2d,
\end{equation}
where $B$ is the~set of corners of $S$ and $d$ the~number of dots.
Thence, we can deform $\BN(\Ann)$ by taking quantum horizontal traces.
Namely, we define the~\emph{quantum Bar--Natan skein category} of the~annulus
as the~additive closure of the~quantum horizontal trace%
%\footnote{
%	The~coefficients in $\BBN$ are extended to the~flat ring $\scalars$
%	before computing the~trace.
%}
\begin{equation}
	\BNq(\Ann) := \qhTr^\oplus(\BBN).
\end{equation}
This category admits the~following graphical description.
The~identified boundaries of $(\RxI)\times I$ form a~membrane
in the~resulting solid torus, and the~orientation of the~core
of $\Ann$ equips the~membrane with a~coorientation.
Isotopic cobordisms are identified whenever the~isotopy fixes
the~membrane. Otherwise, we scale the~target cobordism according
to the~following rules:
\begin{equation}
\label{rel:qBN}%
	\psset{unit=3ex}%
	\begin{gathered}
		\relRCapOffMembrane = q \relRCapAtMembrane
		\hskip 2cm
		\relSaddleLeftToMembrane = q^{-1} \relSaddleRightToMembrane
	\\[1ex]
		\relLCapAtMembrane = q \relLCapOffMembrane
		\hskip 2cm
		\relDotLeftToMembrane = q^{-2}\relDotRightToMembrane
	\end{gathered}%
\end{equation}
where, in each equality, the~tangles draw on the~membrane by the~cobordism
have the~same formal shift.
For instance, a~torus wrapped once around the~annulus evaluates
to $q+q^{-1}$:
\begin{equation}
\label{eq:torus-in-BNq}
	\pictTorusMembrane
	= \pictCutTorusMembrane[*-] + \pictCutTorusMembrane[-*]
	= q^{-1}\pictTorusLeftSphere +  q\pictTorusRightSphere.
\end{equation}
It follows also that $(1-q^2)$ annihilates a~surface $S$ that has
a~connected component with both an~essential boundary and a~dot.
Indeed, moving the~dot along the~essential boundary curve does
not change the~isotopy class of $S$, yet it requires to pass
the~dot through a~membrane.
Therefore, from now on we impose the~Boerner relation, writing
\begin{equation}\label{def:BBNq}
	\BoBNq(\Ann) := \quotient{\aBNq}{\!B} = \quotient{\qhTr^\oplus(\BBN)}{\!B}
\end{equation}
for the~quotient category.
According to \eqref{eq:torus-in-BNq}, the~Cartesian product with a~circle
is a~well-defined functor
\begin{equation}\label{eq:S1x(-):TL-->BNq}
	\S\times(\blank) \colon \TL \to \BoBNq(\Ann)
\end{equation}
for any value of $q$.
The~rest of this section is devoted to prove the~following statement.

\begin{theorem}\label{thm:TL-is-aBNq}
	There is a~commuting diagram
	\begin{equation}\label{diag:TL-aBNq-Rep}
		\begin{diagps}(0em,-0.5ex)(8em,11.5ex)
			\Vtriangle<8em,9ex>|a{npos=0.45}`b`a|[%
				\TL`\BoBNq(\Ann)`\gRep(\Uqsl 2);
				\S\times(\blank)`\F_{TL}`\FqAKh
			]
		\end{diagps}
	\end{equation}
	with the~horizontal functor an~equivalence of categories.
\end{theorem}

We first show the~surjectivity of $\S\times(\blank)$.
Hereafter we write $\seam = \{1\}\times\R\subset\S\times\R$ for the~\emph{seam},
the~arc formed by identifying the~boundaries of $\RxI$.

\begin{figure}[ht]
\centering
\begin{pspicture}(-5,0)(6,7)
	\psframe[style=mysurface,linestyle=none](-5,0)(6,7)
	\begin{psset}{linestyle=none,fillstyle=solid,fillcolor=membraneColor}
		\psarc(1,3){1}{0}{180}
		\psarc(3,3){1}{180}{360}
	\end{psset}
	\begin{psset}{linecolor=red,linestyle=solid,linewidth=0.5pt}
		\psline(-5,3)(6,3)
		\psline{->}(5.7,3)(5.7,4)
		\psline{->}(5.3,3)(5.3,4)
		\psline{->}(4.9,3)(4.9,4)
	\end{psset}
	\uput[r](6,3){\textcolor{red}{$\mu$}}
	\pscustom[style=diagarc]{%
		\moveto(1,0)
		\curveto(0,1)(0,2)(0,3)
		\psarcn(1,3){1}{180}{0}
		\psarc(3,3){1}{180}{360}
		\psarc(1,3){3}{0}{180}
		\curveto(-2,2)(-2,1)(-3,0)
	}%
	\uput[120](-3,0){$\Gamma$}
\end{pspicture}
\caption{An~example of a~positive (to the~left) and a~negative
	(to the~right) retractible arc, each with the~retracting disk.
	The~top arc is not retractible, but it will be after the~positive
	arc is retracted.}%
\label{fig:retractible-arcs}%
\end{figure}
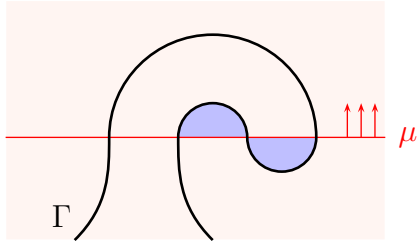

\begin{lemma}\label{lem:vTr(BN)=hTr(BN)}
	The~canonical embedding\/
	$\mathfrak i\colon \qvTr^\oplus(\BBN) \to \qhTr^\oplus(\BBN)$
	is an~equivalence of categories.
\end{lemma}
\begin{proof}
	The~functor $\mathfrak i$ is full and faithful,
	see Section~\ref{sec:vTr-vs-hTr}. Therefore, we need
	only to show that every object of $\qhTr^\oplus(\BBN) = \aBNq$ is
	isomorphic to one from the~image of $\mathfrak i$.
	
	Let $\Gamma$ be a~collection of curves in $\Ann$.
	If $\alpha\subset \Gamma$ is an embedded arc such that
	$\partial\alpha=\seam\cap\alpha$, then we shall say that
	$\Gamma$ is \emph{retractible} if there is an embedded
	disk $D\subset\Ann$ with interior disjoint from $\Gamma$,
	and which boundary is formed by $\alpha$ and the~subarc
	of $\seam$ that lies between the two endpoints of $\alpha$,
	see Figure~\ref{fig:retractible-arcs}.
	
	Assume first that $\Gamma$ has no retractible arcs,
	so that it is a~collection of essential and trivial circles.
	Essential circles are in the~image of $\mathfrak i$,
	and the~trivial ones can be removed using the~delooping
	isomorphism from Proposition~\ref{prop:delooping}.

	It remains to show that collections of curves with
	no retractible arcs generate $\aBNq$.
	For that pick a~retractible arc $\alpha$ in $\Gamma$
	with a~corresponding disk $D$ and apply an~isotopy that
	pushes $\alpha$ across $D$, taking it off $\seam$.
	This procedure reduces the geometric intersection number
	of $\Gamma$ with $\seam$ by $2$. Hence, applying this step
	several times, we end up with $\Gamma'$ that contains
	no retractible arcs, hence from the~image 	of $\mathfrak i$.
\end{proof}

\begin{corollary}\label{col:Sx(-) is full}
	The~functor\/ $\S\times(\blank)\colon \TL \to \BoBNq(\Ann)$ is full and
	essentially surjective.
\end{corollary}
\begin{proof}
	In the~view of \eqref{def:BBNq} and Lemma~\ref{lem:vTr(BN)=hTr(BN)}
	we only need to show that every cobordism in $\BoBNq(\Ann)$ is a~linear
	combination of those of the~form $\S\times T\noic$, where $T$ is
	a~Temperly--Lieb diagram. We achieve that by using the~Bar-Natan
	and the~trace relations.
	
	Using the~neck-cutting relation we reduce first a~surface $S$
	to a~linear sum of surfaces $S_i$ of genus 0. Because all closed
	components evaluate to scalars, we may assume each $S_i$ has
	a~boundary component, which is an~essential curve in $\Ann$
	intersecting the~seam once. Hence, $S_i$ is an~annulus that intersects
	the~membrane in an~arc and, perhaps, in a~collection of circles.
	The~latter can be removed at the~cost of some power
	of $q$ with the~left relations in \eqref{rel:qBN}.
	The~resulting surface is isotopic to $\S\times T\noic$,
	where $T$ is the~intersection of $S_i$ with the~membrane.
	Furthermore, the~Boerner relation prohibits the~surface from
	carrying dots.
\end{proof}

To show faithfulness of $\S\times(\blank)$, so that $\BoBNq(\Ann)$ does not collapse,
we construct the~functor $\FqAKh\colon\BoBNq(\Ann)\to \cat{Mod}(\scalars)$ using
the~diagram
\begin{equation}
	\begin{diagps}(0em,-0.5ex)(11em,12ex)
		\square<11em,10ex>{->`=`->`-->}[%
			\qhTr(\BBN)`\qhTr(\ccat{gBirep})`\aBNq`\cat{Mod}(\scalars);%
			\qhTr(\FFCK)``\qHoHom_{\!0}`\FqAKh
		]
	\end{diagps}
\end{equation}
where $\qhTr(\FFCK)$ is the~functor induced by the~Chen--Khovanov TQFT.
Here $\ccat{gBirep}$ is the~bimodule bicategory restricted to sweet bimodules,
i.e.\ those that are finitely generated and projective as both left and right
modules. This bicategory has duals, so that the~shadow $\qHoHom_{\!0}$ factorizes
through the~horizontal trace by Theorem~\ref{thm:hTr-is-universal},
providing the~right vertical map. Explicitly,
\begin{equation}
	\FqAKh(\widehat T) := \qHoHom_0(\CKalg n, \FFCK(T))
\end{equation}
for a~flat $(n,n)$--tangle $T$. Equivalently, $\FqAKh$ can be constructed
by pulling back $\qHoHom_{\!0}$ to $\BBN$ along $\FFCK$ and factorizing it
through the~horizontal trace.
The~following justifies taking only the~0th
Hochschild homology.

\begin{proposition}\label{prop:HH-of-CKalg}
	Suppose\/ $\scalars$ is flat over $\Zq$.
	Then the~inclusion\/ $\CKalg n_0\subset\CKalg n$ induces an~isomorphism
	of quantum Hochschild homology.
	In particular, the~Chern character\/
	$h\colon K_0(\CKalg n)\utimes\Zq\scalars \to \qHoHom_{\!0}(\CKalg n)$
	is an~isomorphism and\/ $\qHoHom_{\!i}(\CKalg n) = 0$ for $i>0$.
\end{proposition}
\begin{proof}
	It is enough to check the~case $\scalars=\Zq$.
	Consider quantum Hochschild homology as abelian groups and let $R$
	be any ring. The~Universal Coefficient Theorem provides a~commuting diagram
	\[
		\begin{diagps}(0,-0.5ex)(33.3em,12ex)
			\psset{xunit=9em,yunit=10ex}
			\node t0(0.0,1)[0]
			\node t1(0.8,1)[\Tor(\qHoHom_{\!i-1}(\CKalg n_0),R)]
			\node t2(2.0,1)[\qHoHom_{\!i}(\CKalg n_0\otimes R)]
			\node t3(3.1,1)[\qHoHom_{\!i}(\CKalg n_0)\otimes R]
			\node t4(3.7,1)[0]
			\node b0(0.0,0)[0]
			\node b1(0.8,0)[\Tor(\qHoHom_{\!i-1}(\CKalg n),R)]
			\node b2(2.0,0)[\qHoHom_{\!i}(\CKalg n\otimes R)]
			\node b3(3.1,0)[\qHoHom_{\!i}(\CKalg n)\otimes R]
			\node b4(3.7,0)[0]
			\arrow{->}[t0`t1;]\arrow{->}[t1`t2;]\arrow{->}[t2`t3;\cong]\arrow{->}[t3`t4;]
			\arrow{->}[b0`b1;]\arrow{->}[b1`b2;]\arrow{->}[b2`b3;]\arrow{->}[b3`b4;]
			\arrow{->}[t1`b1;\alpha_i]
			\arrow{->}[t2`b2;\beta_i]
			\arrow{->}[t3`b3;\gamma_i]
		\end{diagps}
	\]
	with exact rows, where the~vertical homomorphisms are induced by the~inclusion
	of algebras.
	Theorem~\ref{thm:HH-for-fgdim} implies that $\beta_i$ is an~isomorphism if
	$R = \Z_p$ for prime $p$, because the~global dimension of $\CKalg n\otimes\Z_p$
	is finite \cite{ChenKhovIsCellular} and simple modules are one-dimensional.
	Thus $\gamma_0$ is an~isomorphism, because the~left groups vanish when $i=0$. Hence, it is an~isomorphism for $R=\Z$, and so must be $\alpha_1$.
	Using 5-Lemma we can now prove by induction that $\beta_i$ is an~isomorphism
	for $R=\Z$ and any $i$, which shows the~first claim.
	
	The~second claim follows from a~direct computation.
	The~algebra $\CKalg n_0 \cong \scalars^{2^n}$ is generated by $2^n$ orthogonal idempotents, so that
	$\qHoHom_{\!0}(\CKalg n_0) \cong \scalars^{2^n} \cong
		K_0(\CKalg n)\utimes\Zq \scalars$
	and higher homology vanishes.
\end{proof}

\begin{corollary}\label{cor:qHH-of-CK(T)}
	We have $\qHoHom_{\!i}(\CKalg n, \FFCK(T)) = 0$ for a~flat\/ $(n,n)$--tangle\/
	$T\noic$ and\/ $i>0$.
\end{corollary}
\begin{proof}
	In the~view of Lemma~\ref{lem:vTr(BN)=hTr(BN)} every flat $(n,n)$--tangle $T$
	is isomorphic in $\aBNq$ to the~disjoint union of the~identity tangle on $m$
	points and $\ell$ trivial loops for some $m,\ell\geqslant 0$.
	Hence, $\qHoHom_{\!i}(\CKalg n, \FFCK(T)) \cong
		\qHoHom_{\!i}(\CKalg m) \otimes W^{\otimes\ell}$,
	where $W = \Span_{\scalars}\{w_+, w_-\}$ is the~module assigned by $\FFCK$
	to a~circle. The~thesis follows from Proposition~\ref{prop:HH-of-CKalg}.
\end{proof}

We are now ready to show the~main result of this section.

\begin{proof}[Proof of Theorem~\ref{thm:TL-is-aBNq}]
	The~naturality of the~Chern character provides a~commuting square
	\begin{equation}\label{diag:K0(CKalg)-vs-qHH(CKalg)}
		\begin{diagps}(0em,-0.75ex)(12em,12.5ex)
			\square<12em,10ex>[%
				K_0(\CKalg n)`%\utimes\Z\escalars`%
				K_0(\CKalg m)`%\utimes\Z\escalars`%
				\qHoHom_{\!0}(\CKalg n)`%\qTr(\gRep(\CKalg n))`%
				\qHoHom_{\!0}(\CKalg m);%\qTr(\gRep(\CKalg m));%
				{[\FFCK(T)]}`h`h`\FFCK(T)_*]%
		\end{diagps}
	\end{equation}
	in which both compositions take the~class $[P]$ of a~projective to
	the~quantum Hattoring--Stallings trace of the~identity morphism on
	$P\CKtimes{n}\FFCK(T)$, see \ref{sec:shadow-coinv}.
	The~vertical arrows are isomorphisms by Proposition~\ref{prop:HH-of-CKalg}
	and the~top horizontal map can be identified with $\F_{TL}(T)$
	due to Theorem~\ref{thm:ChenKhov-categ-TL}. Finally, the~bottom horizontal map
	coincides with $\FqAKh$, because it is the~result of applying the~pullback
	shadow $(\FFCK)^*\qHoHom_{\!0}$ to the~morphism in $\qvTr(\BBN) \subset
	\qhTr(\BBN)$ represented by the~square $[T, I\times T]$.
\end{proof}

% ==================================================================
%  Section: Hochschild homology of Chen-Khovanov's 2-functor
% ==================================================================

\subsection{Quantum annular homology}
\label{sec:q-ann-hom}

We use the~deformed skein category together with the~functor $\FqAKh$
to produce a~quantization of the~annular link homology. The~construction
follows the~usual pattern.
\begin{enumerate}

	\item Choose an~annular link diagram $D$ that is transverse to
	the~seam $\mu\subset\Ann$. In particular, all crossings are
	assumed to be away from $\mu$.

	\item Construct the~formal bracket $\KhBracket D_q$ in
	$\HCom^b(\BoBNq(\Ann))$ as explained in Section~\ref{sec:bracket},
	where we use the~subscript $q$ to emphasize that the~resolutions
	of $D$ are taken in the~quantized skein category

	\item Apply the~functor $\FqAKh$ component-wise to $\KhBracket D_q$
	to get the~\emph{quantum annular chain complex}
	$\qAKhCom(D) := \FqAKh\KhBracket D_q$.

\end{enumerate}
The~\emph{quantum annular homology} $\qAKh(L)$ of an~annular link $L$
is defined as the~homology of $\qAKhCom(D)$.
It is a~triply graded $\scalars$--module with a~homological grading,
a~quantum grading coming from the~grading on $\FFCK$,
and an~annular grading arising from the~$\LieSL_2$ weight decomposition.
In what follows we prove that it is well-defined and examine functoriality
with respect to annular link cobordisms.
We begin with a~detailed look on the~quantized formal bracket.
In what follows, we say that a~functor or a~shadow is \emph{projective}
or \emph{$q$-projective} if it is defined on 2-morphisms up to scaling
by $\pm1$ or $\pm q^{\pm k}$ for some $k\in\Z$ respectively.

\begin{proposition}\label{prop:qKhBracket-functorial}
	The~quantized formal bracket\/ $\KhBracket\blank_q\colon \aLinks
	\to \HCom^b(\aBNq)$ is a~$q$--projective functor. In other words,
	the~homotopy class of the~formal chain complex\/ $\KhBracket D_q$
	is an~invariant of annular links and an~annular link cobordism\/
	$W\colon L\to L'$ induces a~formal chain map\/
	$\KhBracket W_q\colon \KhBracket D_q \to \KhBracket{D'}_q$,
	defined up to an~overall sign and power of\/ $q\noic$, where\/
	$D\noic$ and\/ $D'$ are diagrams of\/ $L\noic$ and\/ $L'$ respectively.
\end{proposition}
\begin{proof}
	Let $\ccat{linTan}$ be the~linear extension of the~tangle
	bicategory, 2-morphisms of which are finite sums of tangle
	cobordisms with coefficients in $\scalars$. Motivated
	by Theorem~\ref{thm:hTr-of-FxI} we define
	$\cat{qLinks}(\Ann) := \qhTr(\ccat{linTan})$. It is a~deformed
	linear extension of the~category of annular links, in which
	an~annular link cobordism gets scaled by a~power of $q$ when
	isotoped through the~membrane. Clearly, the~formal bracket
	on $\TTan$ extends to $\ccat{linTan}$.
	
	The~universal shadow $\widehat{(\blank)} \colon \BBN \to \aBNq$,
	defined as the~annular closure of flat tangles,
	induces a~Lefschetz shadow on $\HCom^b(\BBN)$,
	see Lemma~\ref{lem:lefshetz-shadow}.
	Its pullback to $\ccat{linTan}$ along $\KhBracket\blank$ is thence a~projective
	shadow. The~universality of the~horizontal trace gives then
	a~projective functor
	$\KhBracket\blank'_q\colon \cat{qLinks}(\Ann) \to \HCom^b(\aBNq)$.
	One can easily check that $\KhBracket{W}'_q = \KhBracket{W}_q$
	for a~fixed presentation of a~link cobordism $W$\!,
	and it remains to check how dropping the~membrane affects the~bracket.
	For that choose isotopic link cobordisms
	$W$ and $W'$ that are related by a~trace move.
	Regarded as morphisms in $\cat{qLinks}$, they satisfy a~relation
	$W = q^k W'$ for some $k$. Hence,
	$\KhBracket{W}_q = \KhBracket{W}'_q =
		\pm q^k\KhBracket{W'}'_q = \pm q^k \KhBracket{W'}_q$,
	which ends the~proof.
\end{proof}

\begin{itheorem}\label{thm:qAKh}
	The~quantum annular homology $\qAKh(L)$ is a~triply graded invariant
	of an~annular link $L$, which is $q$--projectively functorial with respect
	to annular link cobordisms.
	Moreover, it admits an~action of the~quantum group $\Uqsl 2$ that
	commutes with the~differential and the~maps induced by annular link
	cobordisms intertwine this action.
\end{itheorem}
\begin{proof}
	It follows from Proposition~\ref{prop:qKhBracket-functorial}
	that the~functor $\qAKhCom(\blank) := \FqAKh\KhBracket\blank_q$
	is $q$--projective functor and the~action of $\Uqsl 2$ follows
	from Theorem~\ref{thm:TL-is-aBNq}.
\end{proof}

It was conjectured in \cite{AnnKhAsTrace} that the~Hochschild homology of
Chen--Khovanov complexes recovers the~annular chain complexes, with a~proof
for the~algebras $\CKalg{1,n-1}$, the~next--to--highest weight subalgebra
of $\CKalg n$. The~conjecture follows from Theorem~\ref{thm:TL-is-aBNq}
and the~computation of Hochschild homology in Corollary~\ref{cor:qHH-of-CK(T)}.

\begin{itheorem}\label{thm:HH(CK)-vs-AKh}
	Let\/ $\widehat T\noic$ be the~annular closure of an~$(n,n)$--tangle\/ $T\noic$.
	Then there is an~isomorphism
	\begin{equation}\label{eq:aAKh-is-qHH(cplx)}
		\qAKh(\widehat T) \cong \qHoHom_{\!\bullet}(\CKalg n, \CKmod T),
	\end{equation}
	natural with respect to chain maps associated to tangle cobordisms.
	The~annular grading in $\AKh(\widehat T)$ corresponds to the~weight
	decomposition of\/ $\CKmod{T}$.
\end{itheorem}
\begin{proof}
	The~right hand side of \eqref{eq:aAKh-is-qHH(cplx)} is the~total homology
	of the~bicomplex
	$$
		\coinv_q(\CKmod^\bullet T \CKtimes{n} R_\bullet(\CKalg n)),
	$$
	where $R_\bullet(\CKalg n)$ is the~bar resolution of $\CKalg n$.
	Consider the~associated spectral sequence $E^r\!$ that starts
	with the~Hochschild differential.
	The~first page $E^1_{ij} = \qHoHom_{\!j}(\CKalg n, \CKmod^i T)$
	has vanishing rows except the~one with $j=0$ due to
	Corollary~\ref{cor:qHH-of-CK(T)}.
	Hence, the~sequence collapses at the~second page with
	$E^2_{i0} = \qAKh^i(\widehat T)$ and $E^2_{ij} = 0$ for $r\neq 0$,
	which ends the~proof.
\end{proof}

\noindent\wrapfigure[l]<2>{%
		\psset{unit=1.75em}%
		\begin{pspicture}[shift=-3.5](-1.5,-1.5)(6.1,3.5)
			\psline{<->}(-0.5,3)(-0.5,-0.5)(5,-0.5)
			\uput[l](-0.5,2.85){$\scriptstyle j$}
			\uput[d](4.85,-0.5){$\scriptstyle i$}
			\uput[l](-0.5,2){$\scriptstyle 1$}
			\uput[l](-0.5,1){$\scriptstyle 0$}
			\uput[l](-0.5,0){$\scriptstyle-1$}
			\uput[d](0,-0.5){$\scriptstyle 0$}
			\uput[d](1,-0.5){$\scriptstyle 1$}
			\uput[d](2,-0.5){$\scriptstyle 2$}
			\uput[d](3,-0.5){$\scriptstyle 3$}
			\uput[d](4,-0.5){$\scriptstyle 4$}
			\rput(0,0){\Rnode{00}{$\Z^{\mathrlap 2}$}}
			\rput(0,1){\Rnode{01}{$\Z^{\mathrlap 3}$}} \psdot(0,2)
			\psdot(1,0) \psdot(1,1) \psdot(1,2)
			\psdot(2,0) \rput(2,1){\Rnode{21}{$\Z^{\mathrlap{\phantom 2}}$}} \psdot(2,2)
			\psdot(3,0) \psdot(3,1) \psdot(3,2)
			\psdot(4,0) \psdot(4,1) \psdot(4,2)
			\diagline[linecolor=red,nodesepB=5pt]{->}{21}{00}
			\rput[bl](0,2.5){${}^\backprime E^2_{ij} =
				\HoHom_{\!i}(\Kh^j_{CK}(\text{\fntNegCr}))$}%
		\end{pspicture}%
}%
\indent There is another spectral sequence $\{{}^\backprime E^r\}$ associated to
$\coinv_q(\CKmod^\bullet T\CKtimes{n} R_\bullet(\CKalg n))$, where the~Chen--Khovanov
differential is computed first. Its second page
\begin{equation}
	{}^\backprime E^r_{ij} = \qHoHom_{\!i}(\CKalg n, \Kh^j_{CK}(T))
\end{equation}
is already an~invariant of the~annular link $\widehat T$.
Contrary to the~previous case,
this sequence may not collapse immediately.
The~second page for $T=\text{\fntNegCr}$ has three non-trivial entries
and there is a~non-trivial differential that kills two generators, see
the~diagram to the~left for the~case $\scalars=\Z$ and $q=1$. It is the~third page
that agrees with the~annular homology of $\widehat T$,
which is the~homology of $(\trivmod\to\essmod^{\otimes 2})$ with $\trivmod$ in
homological degree $-1$.

% ==================================================================
%  Section: Hochschild homology of T(2,n)
% ==================================================================

\subsection{Homology for \texorpdfstring{$(2,n)$}{(2,n)} torus links}
\label{sec:qAKh(torus-link)}

Consider the~subcategory $\cat E_2\subset\BoBNq(\Ann)$ generated by objects
intersecting the~seam $\seam$ in exactly two points. Note that every such object
is of the form $\Gamma_I\cup \Gamma_N$, where $\Gamma_N$ is a~(possibly empty)
union of trivial circles not intersecting $\seam$, and $\Gamma_I$ is either
a~trivial circle intersecting $\seam$ in two points, or a~pair of essential
circles, each intersecting $\seam$ in a~single point.

In what follows we shall write $W$ or $W_{\!I}$ for the~module assigned to
a~trivial circle depending on whether it is disjoint from $\seam$ or not.
$W$ is freely generated by $w_+$ and $w_-$, the~images of $1\in\scalars$
under the~maps induced by a~cup cobordism disjoint from the~membrane
$\membr = \seam\times I$, without and with a~dot respectively.
To pick generators for $W_{\!I}$ consider a~cup cobordism that intersects
$\membr$ in a~single arc $\alpha$ and define $w_+$, $w_-^-$ and $w_-^+$
as the~images of $1\in\scalars$ under the~maps induced by the~cobordism
respectively without any dot, a~dot on the~negative side of $\membr$,
and a~dot on the~positive side of $\membr$. All three generate $W_{\!I}$,
but they are not linearly independent: $w_-^- = q^2 w_-^+$.

We can represent elements of the~modules graphically as in
Section~\ref{sec:ann-kh}: the~generators of $V$ are visualized by orienting
the~essential circles, and those of $W$ and $W_{\!I}$ are given as the~trivial
circle without or with a~dot. The~relation in $W_{\!I}$ between $w_-^-$ and
$w_-^+$ can be then written diagrammatically as
\begin{equation}
	\quantumsurfacetrue
	\begin{centerpict}(-1.8,-1.7)(1.8,1.7)
		\pictDrawAnnulus(0,0){1.7}{0.3}
		\psarc[linewidth=0.5pt,arrowinset=0.4,arrowlength=1,arrowsize=5pt,linecolor=red]{->}(0,0){2}{25}{65}
		\psarc{c-c}(0,0){1.3}{-30}{30}
		\psarc{c-c}(0,0){0.7}{-30}{30}
		\rput{ 30}(0,0){\psarc{c-c}(1,0){0.3}{0}{180}}
		\rput{-30}(0,0){\psarc{c-c}(1,0){0.3}{180}{360}}
		\rput{-20}(0,0){\psdot[dotsize=5pt](1.3,0)}
%		\psline[linecolor=red,linewidth=0.5pt](0.3,0)(1.7,0)
%		\psarc[linewidth=0.5pt,arrowsize=4pt]{->}(0,0){1}{130}{230}
	\end{centerpict}
	\quad=\quad q^2
	\begin{centerpict}(-1.8,-1.7)(1.8,1.7)
		\pictDrawAnnulus(0,0){1.7}{0.3}
		\psarc[linewidth=0.5pt,arrowinset=0.4,arrowlength=1,arrowsize=5pt,linecolor=red]{->}(0,0){2}{25}{65}
		\psarc{c-c}(0,0){1.3}{-30}{30}
		\psarc{c-c}(0,0){0.7}{-30}{30}
		\rput{ 30}(0,0){\psarc{c-c}(1,0){0.3}{0}{180}}
		\rput{-30}(0,0){\psarc{c-c}(1,0){0.3}{180}{360}}
		\rput{ 20}(0,0){\psdot[dotsize=5pt](1.3,0)}
%		\psline[linecolor=red,linewidth=0.5pt](0.3,0)(1.7,0)
%		\psarc[linewidth=0.5pt,arrowsize=4pt]{->}(0,0){1}{130}{230}
	\end{centerpict}
\end{equation}
Here we choose the~counter-clockwise orientation of the~core of the~annulus,
so that the~seam $\seam$ is cooriented upwards. Because the~intertwners
\eqref{eq:ev-coev-for-V} are not symmetric, the~essential circles must
be ordered. We choose the~left--to--right ordering read from the~seam.

Capping off a~trivial circle touching the~seam vanishes on $w_+$ and takes
$w_-^\pm$ to $q^{\mp1}$, as the~result is a~sphere without or with a~dot
respectively, but intersecting the~membrane, and isotoping it off the~membrane
using \eqref{rel:qBN} introduces a~power of $q$.
It is now straightforward to compute the~saddle cobordisms in $\cat E_2$
using the~comparison with $\F_{\mathit{TL}}$. A~merge of two essential circles
followed by capping off the~trivial circle is the~evaluation map, which implies
the~following surgery rules :
\begin{equation}
	\psset{unit=0.8}%
	\quantumsurfacetrue
	\pictAnnMergeEssentialsFull[ud] \,\mapsto\, \pictAnnMergeEssentialsFullAfter[01]
\hskip 2cm
	\pictAnnMergeEssentialsFull[du] \,\mapsto\, q\,\pictAnnMergeEssentialsFullAfter[10]
\end{equation}
The~other saddle cobordism is even easier to find out. For degree reasons it must vanish if the~trivial circle carries a~dot, and otherwise it is the~coevaluation map:
\begin{equation}
	\psset{unit=0.8}%
	\quantumsurfacetrue
	\pictAnnSplitTrivialFull \,\mapsto\,
	\pictAnnSplitTrivialFullAfter[ud] \,+\,q^{-1}\,
	\pictAnnSplitTrivialFullAfter[du] 
\end{equation}

Fix $n>0$ and let $T_{2,n}$ denote the~annular $(2,n)$ torus link: the~annular closure of the braid $\sigma^{-n}$, where $\sigma$ is the positive generator of the $2$--strand braid group.

\begin{proposition}\label{prop:qAKh-for-torus-link}
	The quantum annular Khovanov homology of the~annular\/ $(2,n)$
	torus link is given by
	\begin{equation}
		\qAKh^{i,j}(T_{2,n}) =
		\begin{cases}
			V_2                & \text{if }i=0\text{ and }j=-n,\\
			V_0/(q^2 + (-1)^i) & \text{if }-n+1\leqslant i\leqslant -1\text{ and }j=2i-n,\\
			K(q^2 - (-1)^i)    & \text{if }-n\leqslant i \leqslant -2\text{ and }j = 2i - n + 2,\\
			V_0                & \text{if }i=-n\text{ and }j=-3n,\\
			0                  & \text{else},
		\end{cases}
	\end{equation}
	where $V_2:=(V\otimes V) / \Span_{\scalars} \{ \tensor{v+,v-} + q^{-1}\tensor{v-,v+} \}$
	is the~simple representation of\/ $\Uqsl 2$ of dimension 3,
	and\/ $K(a) := \{v\in V_0\ |\ av = 0\}$ for any $a\in\scalars$.
\end{proposition}

Here $i$ and $j$ represent respectively the~homological and the~quantum degree
as defined in \eqref{def:ann-qdeg}. The~latter is denoted by $j'$ in \cite{KhHomSL2}.

\begin{proof}
	Let $D\subset\Ann$ be a~standard diagram for $T_{2,n}$ such that
	cutting $D$ along $\seam$ results in a~diagram for $\sigma^{-n}$.
	Then each resolution of $D$ belongs to $\cat E_2$. We introduce
	the~notations $u_q, l_q\colon W_{\!I}\to W_{\!I}$ for the~maps that put
	a~dot on the~positive and negative side of the~circle respectively,
	and $w_q\colon W_{\!I}\to V\otimes V$ for the~split map. Explicitly,
	\begin{align}
		u_q(w_+) &= w_-^+, &
		u_q(w_-) & =0, \\
		l_q(w_+) &= w_-^- = q^2 w_-^+, &
		l_q(w_-) &= 0, \\
		w_q(w_+) &= \tensor{v+,v-} + q^{-1}\tensor{v-,v+}, &
		w_q(w_-^\pm) &= 0.
	\end{align}
	Let $\{m\}$ denote the grading shift functor which raises the $j$--degree by $m$.
	Arguing as in \cite[Proposition~26]{KhHom}, one can show that $\qAKhCom(D)$
	is quasi-isomorphic to the~chain complex
	\begin{multline*}
		0 \to
		W_{\!I}\{-3n+1\} \to^{\partial^{-n}}
		W_{\!I}\{-3n+3\} \to^{\partial^{-n+1}} \ldots \\
		\to^{\partial^{-3}}
		W_{\!I}\{-n-3\} \to^{\partial^{-2}}
		W_{\!I}\{-n-1\} \to^{\partial^{-1}}
		V_1\otimes V_1^*\{-n\} \to 0
	\end{multline*}
	where $\partial^{-1}=w_q$ and $\partial^i = u_q - (-1)^i l_q$ for $-n\leq i\leq -2$. One can write the above complex more explicitly by writing each $W_{\!I}$ as a~direct sum
	\begin{equation*}
		W_{\!I}=\Span_{\scalars}\{w_+,w_-\}=V_0\{+1\}\oplus V_0\{-1\},
	\end{equation*}
	and by noting that the~map $u_q - (-1)^i l_q$ is given by
	\begin{equation*}
		u_q-(-1)^il_q=
		\begin{pmatrix}
			0 & 0 \\
			0 & 1-(-1)^iq^2
		\end{pmatrix}
	\end{equation*}
	with respect to this direct sum decomposition. It then follows that the above complex is isomorphic to a direct sum of complexes
	\begin{align*}
		&\Bigl(
			0
				\to
			V_0
				\to<3em>^{w_q}
			V_1\otimes V_1
				\to
			0
		\Bigr)\{-n\}
		\\
		&\Bigl(
			0
				\to
			V_0
				\to^{1+(-1)^i q^2}
			V_0
				\to
			0
		\Bigr)\{2i-n\}
			\qquad
		\text{for }-n+1\leq i\leq -1
		\\
		&\Bigl(
			0 \to V_{0} \to 0
		\Bigr)\{-3n\}
	\end{align*}
	where, in each of these complexes, the~bidegree of the~rightmost
	nonzero term is supported on the diagonal $j=2i-n$.
	The~proposition now follows by passing to homology.
\end{proof}

If $\scalars$ is a~field and $q^2 \neq \pm 1$, then $q^2 + (-1)^i$ is invertible. It follows from the~above proof that in such a~case
\begin{equation}
	\qAKh^{i,j}(D) = 
	\begin{cases}
		V_2 & \text{if }(i,j)=(0,-n),\\
		V_0 & \text{if }(i,j)=(-n,-3n),\\
		0   & \text{else}.
	\end{cases}
\end{equation}
On the~other hand, the~quantum homology contains additional copies of $V_0$ if $q^2=\pm 1$, see Table~\ref{tab:qAKh-for-T(2,5)}. This illustrates that quantum annular homology is in general richer than the~non-quantized theory, which corresponds to the case $q=1$.

\begin{table}
\centering
	\begin{tabular}{||c||c|c|c|c|c|c||}\hhline{|=======|}
		\backslashbox{$j$}{$i$}
	     & $-5$ & $-4$ & $-3$ & $-2$ & $-1$ & $0$\\\hhline{#=#=|=|=|=|=|=#}
		$-5$ & & & & & & $V_2$\\\hline
		$-7$ & & & & $V_0^\dagger$ & $V_0^\dagger$ &\\\hline
		$-9$ & & & $V_0^\ddagger$ & $V_0^\ddagger$ &  &\\\hline
		$-11$& & $V_0^\dagger$ & $V_0^\dagger$ & &  &\\\hline
		$-13$& $V_0^\ddagger$ & $V_0^\ddagger$ &  & & &\\\hline
		$-15$& $V_0$ & &  & & &\\\hhline{|=======|}
	\end{tabular}
	\caption{%
		The~quantum annular homology for the~torus knot $T_{2,5}$.
		The~representations marked by a dagger $(\dagger)$ only occur
		if $q^2=1$, and the ones marked by a double dagger ($\ddagger$)
		only occur if $q^2=-1$. The unmarked representations are always there.%
	}%
	\label{tab:qAKh-for-T(2,5)}
\end{table}

\section{Applications and generalizations}
\label{sec:app-and-gen}

Here we discuss functorial properties of the~quantum annular link homology,
as well as generalizations, such as the~quantized APS homology for links
in a~thickened M\"obius or the~twisted annular homology (which provides a~way
to describe homology of satellite knots).

\def\qHHH{\mathit{qHHH}}%
\subsection{Invariants for annular link cobordisms}

It follows from Theorem~\ref{thm:qAKh} that the~quantum
annular homology assigns a~chain map to annular link cobordisms.
Let us now recall how this map is computed.

Let $L$ and $L'$ by annular closures of an~$(m,m)$--tangle $T$ and
an~$(n,n)$--tangle $T'$ respectively.
A~link cobordism $W\colon L\to L'$ intersects the~membrane
in an $(n,m)$--tangle $P$, so that it can be represented
by a~tangle cobordism $\widetilde W\colon P T \to T'P$. Applying
the~Chen--Khovanov functor results in a~square
\begin{equation}\label{diag:map-for-cob}
	\begin{diagps}(0em,-0.5ex)(6em,12.5ex)
		\square<6em,10ex>|a`b`a`a|[%
			\CKalg m`\CKalg m`\CKalg n`\CKalg n;
			\CKmod{T}`\CKmod{P}`\CKmod{P}`\CKmod{T'}
		]%
		\arrow[nodesep=1.75em]|b|{=>}[v1`v2;\widetilde W_*]
	\end{diagps}
\end{equation}
The~quantum annular complex is computed by replacing each component
of the~Chen--Khovanov complex with its Hochschild homology. This
is an~example of a~Lefschetz shadow discussed in Section~\ref{sec:lef-shadow}.
Hence, the~component of the~chain map $W_*$ at homological degree $i$
and quantum degree $j$ is given by summing up the~compositions
\begin{equation}\label{eq:component-of-W*}
	\begin{split}
		\qHoHom_{\!0}(\CKalg m, \CKmod^{i,j} T)
		&\to<5em>^{\coev_*}
		\qHoHom_{\!0}(\CKalg m, \CKmod^{i,j} T \CKtimes{m}
						\CKmod^{i',j} P \CKtimes{n} \CKmod^{i',j} P^*) \\
		&\to<5em>^{\widetilde W_*}
		\qHoHom_{\!0}(\CKalg m, \CKmod^{i',j} P \CKtimes{n} \CKmod^{i,j}{T'}
						\CKtimes{n} \CKmod^{i',j} P^* ) \\
		&\to<5em>^{(-1)^{i'(i+i')}\theta}
		\qHoHom_{\!0}(\CKalg n, \CKmod^{i,j}{T'} \CKtimes{n} \CKmod^{i',j} P^*
						\CKtimes{m} {\CKmod^{i',j} P}) \\
		&\to<5em>^{\ev_*}
		\qHoHom_{\!0}(\CKalg n, \CKmod^{i,j}{T'})
	\end{split}
\end{equation}
over all indices $i'$. Notice that the~third map takes the~class
of $x\otimes y\otimes\alpha$ to the~class of
$(-1)^{i'(i+i')}q^{\qdeg(x)}y\otimes\alpha\otimes x$.

\begin{proposition}
	Let the~cobordism\/ $W\colon L\to L$ trace the~rotation of\/ $L\noic$
	along the~annulus.
	The~induced map $\qAKh(W)$ scales $x\in\qAKh(L)$ by\/ $q^{\qdeg(x)}\!$.
\end{proposition}
\begin{proof}
	Write $L$ as the~annular closure of $T$.
	Then $W$ is represented by the~identity cobordism from $TT$ to itself,
	and the~only components \eqref{eq:component-of-W*} that contribute
	are those with $i=i'$. In particular, the~second map scales the~argument
	by $q^j$, but does not change the~sign.
\end{proof}

It follows that quantum annular homology detects twists: the~chain
map $W_*$ in the~above proposition is not scaled identity. Even more,
if a~resolution $T_\zeta$ of $T$ is a~collection of $n$ vertical lines,
the~restriction of the~twist to $T_\zeta$ is $q^d$, where $d$ is the~degree
shift of $T_\zeta$ in the~complex. Indeed, the~component of the~twist
is represented by
\begin{equation}
	\begin{diagps}(0em,-0.5ex)(8em,12.5ex)
		\square<8em,10ex>|a`b`a`a|[%
			\CKalg n`\CKalg n`\CKalg n`\CKalg n;
			\FFCK(T_\zeta)\{d\}`\FFCK(T_\zeta)\{d\}`
			\FFCK(T_\zeta)\{d\}`\FFCK(T_\zeta)\{d\}
		]%
		\arrow[nodesep=2.75em]|b|{=>}[v1`v2;\id]
	\end{diagps}
\end{equation}
whereas the~identity map has $\FFCK(T_\zeta)$ with no degree shifts
as vertical maps.

Because $\qAKh$ is defined using a~Lefschetz shadow, we can use
Corollary~\ref{cor:HH-induces-Lef} to compute the~invariant for
closed surfaces.

\begin{itheorem}\label{thm:inv-for-closed-sfce}
	Let\/ $\widehat W \subset \S\times\R^3$ be a~closed surface obtained
	as an~annular closure of a~link cobordism\/ $W\colon L\to L\noic$ with\/
	$L\subset\R^3$. Then $\qAKh(\widehat W) = \Lef(W_*)$ is the~graded Lefschetz
	trace of\/ $W_*\colon \Kh(L)\to \Kh(L)$, the~endomorphism of the~Khovanov
	homology of\/ $L$.
	In particular, $\qAKh(\S\times L)$ coincides with the~Jones polynomial\/ $J(L)$.
\end{itheorem}
\begin{proof}
	The~first statement is an~immediate consequence
	of Corollary~\ref{cor:HH-induces-Lef}. Applying it to
	$W = L\times I$ results in computing the~graded Euler
	characteristic of $\Kh(L)$, the~Khovanov homology of $L$.
\end{proof}

The~APS homology does not distinguish closed surfaces.
Hence, $\qAKh$ is a~nontrivial deformation of the~annular $\LieSL_2$--homology.

\subsection{An~action of tangles on cablings}
\label{sec:action}

The~action of the~braid group on the~annular homology of a~cabling of a~framed
annular knot $K$ was studied in \cite{KhHomSL2}. In what follows we compute
the~action of the~entire category of oriented tangles on the~quantum annular
homology of cablings of $K$. We work in this section in characteristic 2, so that
the~Chen--Khovanov functor is strictly functorial.

%
%  Functoriality gives the braid group action
%
To a~framed knot $K\subset\AxI$ we can associate an~embedding
$\nu_K\colon \AxI \to \AxI$ with the~tubular neighborhood as its image.
It determines a~functor $K_*\colon \Tan \to \aLinks$
such that
\begin{itemize}
	\item a~collection of oriented points $B\subset \R^2$ is taken
	to $K^B := \nu_K(\S\times B)$, an~oriented cabling of $K$, and
	\item an~oriented tangle $T\subset \R^2\times I$ is mapped to the~oriented
	cobordism $K^T := \nu_K(\S\times T)$ between the~cablings.
\end{itemize}
Applying the~quantum annular homology produces a~map of homology
\[
	\qAKh(K^T)\colon \qAKh(K^B) \to \qAKh(K^{B'})
\]
for any oriented tangle $T \in \cat{Tan}(B,B')$, defined up to an~overall
power of $q$. These give a~projectively functorial action of $\cat{Tan}$, i.e.
\[
	\qAKh(K^{T'T}) = q^k\left( \qAKh(K^{T'})\circ \qAKh(K^T) \right)
\]
for composable tangles $T$, $T'$, and some $k\in\Z$.

This action has an~interpretation in the~framework of horizontal
traces when $K$ is the~annular closure of a~framed $(1,1)$--tangle $\tau\noic$.
Denote by $\nu_\tau\colon \R^2\times I \to \R^2\times I$ the~associated
parametrization of the~tubular neighborhood of $\tau$ and
choose a~tangle $T$ with collections of points $B$ and $B'$ as its bottom
and top boundary respectively. Then the~cobordism $K^T$ intersects the~membrane
in $T$ and is represented by the~square
\begin{equation}
	\begin{diagps}(0,-0.5ex)(6em,12ex)
		\square<6em,10ex>[%
			B`B`B'`B';
			\tau^B`T`T`\tau^{B'}
		]
		\arrow[nodesep=1.7em]|b|{=>}[v1`v2;\tau^T]
	\end{diagps}
\end{equation}
where $\tau^B := \nu_\tau(B\times I)$ and $\tau(B') := \nu_\tau(B'\times I)$
are cablings of $\tau$, and $\tau^T := \mu_\tau(T\times I)$ traces the~isotopy
that slides $T$ along $\tau$. We use this to prove that the~action satisfies
the~Jones relation.

\begin{itheorem}\label{thm:action-on-qAKh}
	Let\/ $K$ be a~framed annular link, considered as an~object in
	$\cat{qLinks}(\Ann)$.
	There is a~functorial action of\/ $\Tan\noic$ on the~quantum annular homology
	of oriented cablings of\/ $K\noic$, that takes a~tangle\/ $T\noic$ to the~chain
	map $\qAKh(K^T)$, and which intertwines the~action of\/ $\Uqsl 2$.
	It factors through the~Jones skein relation
	\begin{equation}\label{eq:skein-relation-for-K^T}
		q^2 \qAKh(K^{\fntNegCr}) - q^{-2} \qAKh(K^{\fntPosCr})
		= (q-q^{-1}) \qAKh(K^{\fntNoCr})
	\end{equation}
	if\/ $K$ intersects the~membrane in one point.
\end{itheorem}
\begin{proof}
	Functoriality has been discussed before and compatibility with
	the~action of $\Uqsl 2$ follows from the~construction---the~action
	of $\cat{Tan}$ is already defined at the~level of the~quantized formal
	bracket. It remains to show that the~Jones relation holds. We prove it
	for the~formal bracket.
	
	Let $K$ be the~annular closure of a~$(1,1)$--tangle, represented
	by a~tangle diagram $\tau$. 
	Given a~finite sequence $\underline\epsilon = (\epsilon_1,\dots,\epsilon_k)$
	with $\epsilon_i\in\{-1,+1\}$ we write $\tau^{\underline\epsilon}$ for
	the~oriented $k$--cabling of $\tau\noic$, in which each $i$-th cable is
	oriented parallel to $\tau$ when $\epsilon_i=+1$ and opposite to $\tau$ otherwise.
	We write $\tau^k$ when each $\epsilon_i=1$.
	Up to degree shifts, the~formal bracket depends only on $k$. Explicitly,
	\begin{equation}\label{eq:degree-shift-for-cablings}
		\KhBracket{\tau^{\underline\epsilon}} =
		\KhBracket{\tau^k}
			\left[
				\tfrac{w(\tau)}{2}(k^2 - \|\underline\epsilon\|^2)
			\right]
			\left\{
				\tfrac{3w(\tau)}{2}(k^2 - \|\underline\epsilon\|^2)
			\right\},
	\end{equation}
	where we write $w(T) = n_+(T) - n_-(T)$ for the~writhe of\/ $T$
	and $\|\underline\epsilon\| := \sum\epsilon_i$. Notice that
	$\|\underline\epsilon\| \equiv k\mod 2$, so that both numbers are integers.
	To derive this formula it is enough to compare degree shifts of a~single
	resolution, which is left as an~exercise.

	It follows from functoriality of the~action and
	\eqref{eq:degree-shift-for-cablings} that the~Jones relation
	has to be checked only for the~first three tangles shown below:
	\begin{center}
	\setlength{\tabcolsep}{1.5em}
	\begin{tabular}{cccc}
		\begin{pspicture}(-3em,-3.5ex)(3em,2ex)
			\psset{linewidth=1pt}%
			\psline{->}(-2.5em,-2ex)(-2.5em,2ex)
			\rput(-1.5em,0ex){$\cdots$}
			\psline{->}( 0.5em,-2ex)(-0.5em,2ex)
			\psline[border=3pt]{->}(-0.5em,-2ex)( 0.5em,2ex)
			\rput( 1.5em,0ex){$\cdots$}
			\psline{->}( 2.5em,-2ex)( 2.5em,2ex)
			\rput[B](-2.5em,-3.5ex){$\scriptstyle 0$}
			\rput[B](-0.5em,-3.5ex){$\scriptstyle i$}
			\rput[B]( 0.75em,-3.5ex){$\scriptstyle i+1$}
			\rput[B]( 2.5em,-3.5ex){$\scriptstyle k$}
		\end{pspicture}
		&
		\begin{pspicture}(-3em,-3.5ex)(3em,2ex)
			\psset{linewidth=1pt}%
			\psline{->}(-2.5em,-2ex)(-2.5em,2ex)
			\rput(-1.5em,0ex){$\cdots$}
			\psline{->}(-0.5em,-2ex)( 0.5em,2ex)
			\psline[border=3pt]{->}( 0.5em,-2ex)(-0.5em,2ex)
			\rput( 1.5em,0ex){$\cdots$}
			\psline{->}( 2.5em,-2ex)( 2.5em,2ex)
			\rput[B](-2.5em,-3.5ex){$\scriptstyle 0$}
			\rput[B](-0.5em,-3.5ex){$\scriptstyle i$}
			\rput[B]( 0.75em,-3.5ex){$\scriptstyle i+1$}
			\rput[B]( 2.5em,-3.5ex){$\scriptstyle k$}
		\end{pspicture}
		&
		\begin{pspicture}(-3em,-3.5ex)(3em,2ex)
			\psset{linewidth=1pt}%
			\psline{->}(-2.5em,-2ex)(-2.5em,2ex)
			\rput(-1.5em,0ex){$\cdots$}
			\psline{->}(-0.5em,-2ex)(-0.5em,2ex)
			\psline{->}( 0.5em,-2ex)( 0.5em,2ex)
			\rput( 1.5em,0ex){$\cdots$}
			\psline{->}( 2.5em,-2ex)( 2.5em,2ex)
			\rput[B](-2.5em,-3.5ex){$\scriptstyle 0$}
			\rput[B](-0.5em,-3.5ex){$\scriptstyle i$}
			\rput[B]( 0.75em,-3.5ex){$\scriptstyle i+1$}
			\rput[B]( 2.5em,-3.5ex){$\scriptstyle k$}
		\end{pspicture}
		&
		\begin{pspicture}(-3em,-3.5ex)(3em,2ex)
			\psset{linewidth=1pt}%
			\psline{->}(-2.5em,-2ex)(-2.5em,2ex)
			\rput(-1.5em,0ex){$\cdots$}
			\psbezier{->}(0.5em,-2ex)(0.5em,0ex)(-0.5em,0ex)(-0.5em,-2ex)
			\psbezier{->}(-0.5em,2ex)(-0.5em,0ex)(0.5em,0ex)( 0.5em, 2ex)
			\rput( 1.5em,0ex){$\cdots$}
			\psline{->}( 2.5em,-2ex)( 2.5em,2ex)
			\rput[B](-2.5em,-3.5ex){$\scriptstyle 0$}
			\rput[B](-0.5em,-3.5ex){$\scriptstyle i$}
			\rput[B]( 0.75em,-3.5ex){$\scriptstyle i+1$}
			\rput[B]( 2.5em,-3.5ex){$\scriptstyle k$}
		\end{pspicture}
		\\
		$T_+$ & $T_-$ & $T_0$ & $T_h$
	\end{tabular}\end{center}
	We use the~forth one to express the~formal brackets of $T_+$ and $T_-$
	as mapping cones following Proposition~\ref{prop:Viro-ses},
	obtaining distinguished triangles
	\begin{equation*}
		\setlength\arraycolsep{1.25em}
		\begin{array}{cccc}
	    		\Rnode{t0}{\KhBracket{T_0}[-1]\{1\}} &
		  	\Rnode{t1}{\KhBracket{T_h}[-1]\{2\}} &
			\Rnode{t2}{\KhBracket{T_+}} &
			\Rnode{t3}{\KhBracket{T_0}\{1\}}
		\\[1ex]
	    		\Rnode{b0}{\KhBracket{T_h}\{-2\}} &
	    		\Rnode{b1}{\KhBracket{T_0}\{-1\}} &
	    		\Rnode{b2}{\KhBracket{T_-}} &
	    		\Rnode{b3}{\KhBracket{T_h}[1]\{-2\}}
		\end{array}
		\begin{diagps}(0,0)
			\arrow{->}[t0`t1;\mathit{sd}]
			\arrow|a{npos=0.6}|{->}[t1`t2;\mathit{in}]
			\arrow|a{npos=0.45}|{->}[t2`t3;\mathit{pr}]
			\arrow{->}[b0`b1;\mathit{sd}]
			\arrow|a{npos=0.55}|{->}[b1`b2;\mathit{in}]
			\arrow|a{npos=0.35}|{->}[b2`b3;\mathit{pr}]
		\end{diagps}
	\end{equation*}
	in which the~left morphisms are saddle cobordisms in $\BBN$.
	Consider now the~diagram
	\begin{equation}\label{diag:nat-of-dist-tris}
	\begin{diagps}(0,-0.5ex)(27em,12ex)
		\psset{xunit=2em,yunit=10ex}%
		\node b0( 0  ,0)[{\KhBracket{T_0}[-1] \otimes \KhBracket{\tau^k}}]
		\node b1( 5  ,0)[{\KhBracket{T_h}[-1] \otimes \KhBracket{\tau^k}}]
		\node b2( 9.5,0)[{\KhBracket{T_+} \otimes \KhBracket{\tau^k}}]
		\node b3(13.5,0)[{\KhBracket{T_0} \otimes \KhBracket{\tau^k}}]
		\node t0( 0  ,1)[{\KhBracket{\tau^k}\otimes\KhBracket{T_0}[-1]}]
		\node t1( 5  ,1)[{\KhBracket{\tau^k}\otimes\KhBracket{T_h}[-1]}]
		\node t2( 9.5,1)[{\KhBracket{\tau^k}\otimes\KhBracket{T_+}}]
		\node t3(13.5,1)[{\KhBracket{\tau^k}\otimes\KhBracket{T_0}}]
		\arrow{->}[t0`t1;\mathit{sd}]
		\arrow{->}[t1`t2;\mathit{in}]
		\arrow{->}[t2`t3;\mathit{pr}]
		\arrow{->}[b0`b1;\mathit{sd}]
		\arrow{->}[b1`b2;\mathit{in}]
		\arrow{->}[b2`b3;\mathit{pr}]
		\arrow{->}[t0`b0;{\KhBracket{\tau^{T_0}}[-1]}]
		\arrow{->}[t1`b1;{\KhBracket{\tau^{T_h}}[-1]}]
		\arrow{->}[t2`b2;\KhBracket{\tau^{T_+}}]
		\arrow{->}[t3`b3;\KhBracket{\tau^{T_0}}]
	\end{diagps}
	\end{equation}
	where the~quantum degree shifts are dropped for clarity.
	The~left and right vertical maps are identities, whereas those
	in the~middle are induced by a~sequence of Reidemeister moves.
	The~left square commutes up to a~formal chain homotopy $H$ by
	the~functoriality of $\KhBracket\blank$. In particular,
	$\KhBracket{\tau^{T_0}}$ and $\KhBracket{\tau^{T_h}}$
	together with $H$ induce a~chain map
	\[
		\alpha = \begin{pmatrix}
			\KhBracket{\tau^{T_0}} & 0 \\
			-H & \KhBracket{\tau^{T_h}}[-1]
		\end{pmatrix}\colon
		\KhBracket{\tau^k}\otimes\cone(\mathit{sd})
		\to \cone(\mathit{sd}) \otimes \KhBracket{\tau^k}
	\]
	that would make the~other squares in \eqref{diag:nat-of-dist-tris}
	commute when placed as the~third vertical map. Hence, it is enough
	to show that $\alpha$ and $\KhBracket{\tau^{T_+}}$ are chain homotopic.
	
	It is shown in \cite[Section 8]{Cob} that every degree 0 automorphism
	of $\KhBracket{T}$ is homotopic to $\pm\id$ if $T$ is obtained from
	a~crossingless tangle by twisting its endpoints. Therefore, there exists
	at most one homotopy equivalence $\KhBracket{T} \simeq \KhBracket{T'}$
	when $T$ and $T'$ are two such tangles (as we work over $\Z_2$).
	In particular, there are unique homotopy equivalences
	\[
		\KhBracket*{\vcenter{\hbox{\begin{pspicture}(-2.5em,-5ex)(2.5em,5ex)
			\psset{linewidth=1pt}%
			\psline(-1em,-5ex)(-1em,5ex)
			\psline( 1em,-5ex)( 1em,5ex)
			\pscustom{%
				\moveto(2em,-5ex)
				\lineto(2em, 0ex)
				\curveto(2em, 4ex)(-2em, 1ex)(-2em,5ex)
				\stroke[linewidth=7pt,linecolor=white]
			}
			\psframe[fillstyle=solid](-1.5em,-2.5ex)(1.5em,0.75ex)
			\rput[c](0,-1ex){$T_*$}
			\rput(0,-4ex){$\cdots$}
			\rput(0, 4ex){$\cdots$}
		\end{pspicture}}}}
		\simeq
		\KhBracket*{\vcenter{\hbox{\begin{pspicture}(-2.5em,-5ex)(2.5em,5ex)
			\psset{linewidth=1pt}%
			\psline(-1em,-5ex)(-1em,5ex)
			\psline( 1em,-5ex)( 1em,5ex)
			\pscustom{%
				\moveto(-2em, 5ex)
				\lineto(-2em, 0ex)
				\curveto(-2em,-4ex)(2em,-1ex)(2em,-5ex)
				\stroke[linewidth=7pt,linecolor=white]
			}
			\psframe[fillstyle=solid](-1.5em,2.5ex)(1.5em,-0.75ex)
			\rput[c](0, 0.75ex){$T_*$}
			\rput(0,-4ex){$\cdots$}
			\rput(0, 4ex){$\cdots$}
		\end{pspicture}}}}
	\quad\text{and}\quad
		\KhBracket*{\vcenter{\hbox{\begin{pspicture}(-2.5em,-5ex)(2.5em,5ex)
			\psset{linewidth=1pt}%
			\pscustom{%
				\moveto(2em,-5ex)
				\lineto(2em, 0ex)
				\curveto(2em, 4ex)(-2em, 1ex)(-2em,5ex)
			}
			\psline[border=3pt](-1em,-5ex)(-1em,5ex)
			\psline[border=3pt]( 1em,-5ex)( 1em,5ex)
			\psframe[fillstyle=solid](-1.5em,-2.5ex)(1.5em,0.75ex)
			\rput[c](0,-1ex){$T_*$}
			\rput(0,-4ex){$\cdots$}
			\rput(0, 4ex){$\cdots$}
		\end{pspicture}}}}
		\simeq
		\KhBracket*{\vcenter{\hbox{\begin{pspicture}(-2.5em,-5ex)(2.5em,5ex)
			\psset{linewidth=1pt}%
			\pscustom{%
				\moveto(-2em, 5ex)
				\lineto(-2em, 0ex)
				\curveto(-2em,-4ex)(2em,-1ex)(2em,-5ex)
			}
			\psline[border=3pt](-1em,-5ex)(-1em,5ex)
			\psline[border=3pt]( 1em,-5ex)( 1em,5ex)
			\psframe[fillstyle=solid](-1.5em,2.5ex)(1.5em,-0.75ex)
			\rput[c](0, 0.75ex){$T_*$}
			\rput(0,-4ex){$\cdots$}
			\rput(0, 4ex){$\cdots$}
		\end{pspicture}}}}	\]
	where $T_*$ is any of the~four tangles $T_+$, $T_-$, $T_0$, or $T_h$.
	On the~other hand, each map $\KhBracket{\tau^{T_*}}$ is a~composition
	of the~above equivalences, and so is the~chain map $\alpha$.
	Hence, $\alpha$ and $\KhBracket{\tau^{T_+}}$ agree up to a~chain homotopy,
	which makes the~middle and right squares
	in \eqref{diag:nat-of-dist-tris} commute in the~homotopy category.
		
	The~quantum annular bracket is a~quantum Lefschetz shadow on $\HCom^b(\BBN)$.
	Hence, taking the~degree shifts into acount, we have from
	Proposition~\ref{prop:additivity-of-trian-shadow} that
	\[
		\KhBracket{K^{T_+}}_q = q\KhBracket{K^{T_0}}_q - q^2\KhBracket{K^{T_h}}_q.
	\]
	Likewise, using the~other distinguished triangle we obtain
	\[
		\KhBracket{K^{T_-}}_q = q^{-2}\KhBracket{K^{T_h}}_q - q^{-1}\KhBracket{K^{T_0}}_q.
	\]
	These two equalities imply the~Jones skein relation.
\end{proof}

The~commutativity of \eqref{diag:nat-of-dist-tris} can be also checked directly,
as there are explicit formulas for all the~chain maps, see \cite{Cob}. Although
it requires more work, with this approach one can prove
Theorem~\ref{thm:action-on-qAKh}
for the~homology with integral coefficients, once a~strictly functorial version of
the~construction due to Chen and Khovanov is used. This approach is used in
a~following paper.
Alternatively, one can fix the~signs explicitly as it was done successfully
in \cite{KhHomSL2}.

One can forget the~membrane and work with true annular links,
instead of $\cat{qLinks}(\Ann)$. Although the~map $\qAKh(K^T)$ is \emph{a~priori}
defined only up to an~overall power of $q$, making the~relation
\eqref{eq:skein-relation-for-K^T} problematic, our definition of the~map
uses a~special presentation of the~cobordism: as the~isotopy sliding $T$
along $\tau$, the~opening of $K$. For this particular presentation
the~power of $q$ is well-defined, so that the~Jones relation holds.
However, one must be careful when starting to isotope the~cobordism.
\subsection{Quantum homology of links in a~thickened M\"obius band}
\label{sec:quantum-MB}

Let $\rho\in\mathit{Diff}(\RxI)$ be the~flip along the~interval,
i.e.\ $\rho(x,t) = (-x,t)$. It induces an~endobifunctor $\rho_*$
on $\BBN$. Identifying $(x,1)$ with $(-x,0)$ produces a~M\"obius band $\MB$,
and, by the~argument from the~proof of Theorem~\ref{thm:tr-of-tan},
$\hTr^\oplus(\BBN,\rho_*)$ and $\BN(\MB)$ are equivalent categories.
In an~analogy to the~case of annulus we define
\begin{equation}
	\BNq(\MB) := \qhTr^\oplus(\BBN(\RxI),\rho_*)
\end{equation}

The~reflection $\rho$ induces also an~automorphism of $\CKalg n\!$,
which we denote with the~same symbol.
The~following result is immediate from the~definition of $\FFCK$.

\begin{lemma}
	Given an~$(m,n)$--tangle\/ $T\noic$ write\/ $T^{\textit{flip}}$
	for its reflection along the~vertical axis. Then
	\begin{equation}
		{}_\rho\FFCK(T)_\rho \cong \FFCK(T^{\textit{flip}})
	\end{equation}
	as\/ $(\CKalg m,\CKalg n)$--bimodules.
\end{lemma}

Thence, we can define a~TQFT functor $\FqMKh\colon \BNq(\MB)\to\cat{gMod}(\scalars)$
with the~help of the~pullback shadow $(\FFCK)^*\qHoHom_{\!0}^\rho$.
The~argument from Proposition~\ref{prop:HH-of-CKalg} applied to $\qHoHom^\rho$
implies that
\begin{equation}
	\qHoHom_{\!>0}^\rho(\CKalg n) = 0
\qquad\text{and}\qquad
	\qHoHom_{\!0}^\rho(\CKalg n) \cong \coinv^\rho(\CKalg n_0),
\end{equation}
where the~latter is generated by idempotents corresponding to symmetric cup
diagrams (i.e.\ those fixed by $\rho$). In particular, if we write $c_n$ for
the~collection of $n$ parallel separating curves, each wrapping $\MB$ twice,
and $\gamma$ for the~nonseparating curve, then
\begin{equation}\label{eq:isom-qHH-with-E-for-MB}
	\FqMKh(c_n) \cong \qHoHom^\rho_{\!0}(\CKalg{2n}) \cong \CKalg n_0 \cong V^{\otimes n}
\qquad\text{and}\qquad
	\FqMKh(\gamma \cup c_n) = 0
\end{equation}
where $\CKalg n_0 \cong \qHoHom^\rho_{\!0}(\CKalg{2n})$ takes an~idempotent
$e\in\CKalg n_0$ to $e\otimes\rho(e)\in\CKalg{2n}_0$. In particular, we are
interested only in the~subcategory $\evBNq(\MB)$ generated by those collections
that do not contain $\gamma$. These collections are characterized by the~following
property: if $\mu\subset\MB$ cuts the~band into a~square, then each object from
$\evBNq(\MB)$ intersects $\mu$ in an~even number of points.
The~argument from Lemma~\ref{lem:vTr(BN)=hTr(BN)} adapted to this case shows that
$\qvTr^\oplus(\BBN^{ev},\rho_*) \to \qhTr^\oplus(\BBN^{eb},\rho_*) = \evBNq(\MB)$
is an~equivalence of categories. The~following result is straightforward.

\begin{lemma}
	Let $\gamma\subset\MB$ be a~nonseparating curve. Then there is a~commutative diagram of functors
	\begin{equation}
		\begin{diagps}(0,-0.75ex)(8em,11.25ex)
			\Vtriangle<8em,9ex>[%
				\BN\!_{q^2}(\Ann)`\evBNq(\MB)`\Mod(\scalars);%
				\Phi`\FqAKh`\FqMKh
			]
		\end{diagps}
	\end{equation}
	where $\Phi\colon \BN\!_{q^2}(\Ann) \to \BNq(\MB)$ is induced by the~diffeomorphism $\Ann \approx \MB-\gamma$.
\end{lemma}
\begin{proof}
	The~functor $\Phi$,
	seen as $\hTr^\oplus_{q^2}(\BBN) \to \qhTr^\oplus(\BBN,\rho_*)$,
	adds $n$ vertical lines next to an~$(n,n)$--tangle and sends a~morphism $[T, W]$
	to $[T\cup T^{\mathit{flip}}, W\cup(I\times T^{\mathit{flip}})]$.
	The~difference in powers of $q$ appears, because going
	once through a~membrane in $\hTr^\oplus_{q^2}(\BBN) \cong \BN\!_{q^2}(\Ann)$ corresponds to going twice through
	the~membrane in $\qhTr^\oplus(\BBN,\rho_*)\cong\BNq(\MB)$.
\end{proof}

Each cobordism in $\qvTr^\oplus(\BBN^{ev},\rho_*)$ can be decomposed into
a~composition of cobordisms from the~image of $\Phi$ and projective planes
with a~disk removed, each with a~nontrivial curve in $\MB$ as its boundary.
These can be seen in turn as the~saddle cobordisms between a~trivial and
a~nontrivial curve, with the~trivial curve capped off. They correspond under
the~quotient map $(\RxI)\times I\to \MBxI$ to surfaces $\fntCup\times I$
and $\fntCap\times I$, and we denote them by $S_{\fntCup}$ and $S_{\fntCap}$
respectively.

\begin{lemma}
	$\FqMKh(S_{\fntCap})\colon V\to \scalars$ evaluates\/ $v_+$ and\/ $v_-$
	to\/ $q\noic$ and\/ $1$ respectively, whereas
	$\FqMKh(S_{\fntCup})\colon\scalars\to V$ takes\/ $1$ to\/ $v_+ + q^{-1}v_-$.
\end{lemma}
\begin{proof}
	The~module $V := \qHoHom_{\!0}^\rho(\CKalg 2)$ assigned by $\FqMKh$
	to a~nontrivial curve has canonical generators
	\begin{equation}
		b_0 := \pairofmatchings{10}{10}
	\qquad\text{and}\qquad
		b_1 := \pairofmatchings{01}{01}
	\end{equation}
	which corresponds under \eqref{eq:isom-qHH-with-E-for-MB} to
	$v_+ + q^{-1}v_-$ and $v_-$ respectively.
	The~thesis now follows from direct computations.
	Indeed, the~map $\FqMKh(S_{\fntCup})$ is given by the~sequence
	\begin{align*}
		\Zq = \qHoHom_{\!0}^\rho(\Z)
			&\to \qHoHom_{\!0}^\rho(\Z,
					\FFCK(\fntCup) \CKtimes{2} {}^*\FFCK(\fntCup))\\\
			&\to^\theta \qHoHom_{\!0}^\rho(\CKalg 2,
					{}^*\FFCK(\fntCup) \otimes \FFCK(\fntCup))
			 \to \qHoHom_{\!0}^\rho(\CKalg 2)
	\end{align*}
	which takes $1\in\scalars$ into $b_0$.
	Dually, $S_{\fntCap}$ is assigned the~sequence
	\begin{align*}
		\qHoHom_{\!0}^\rho(\CKalg 2)
			&\to \qHoHom_{\!0}^\rho(\CKalg 2,
					\FFCK(\fntCap) \otimes {}^*\FFCK(\fntCap)) \\
			&\to^\theta \qHoHom_{\!0}^\rho(\Z,
					{}^*\FFCK(\fntCap) \CKtimes{2} \FFCK(\fntCap))
			 \to \qHoHom_{\!0}^\rho(\Z) = \Zq,
	\end{align*}
	which takes $b_0$ and $b_1$ to $q+q^{-1}$ and $1$ respectively.
\end{proof}

The~above is enough to compute $\FqMKh$ on all morphisms in $\BNq(\MB)$.
A~comparison with the~formulas from Section~\ref{sec:khov-hom-mobius}
shows that $\FqMKh$ is a~deformation of the~APS construction.

\begin{theorem}
	Let $\widehat T$ be the~closure in $\MB$ of a~$(2n,2n)$--tangle $T$.
	Then there is an~isomorphism
	\begin{equation}
		\MKh(\widehat T) \cong \HoHom_{\!\bullet}^\rho(\CKalg{2n}, \CKmod T)
	\end{equation}
	natural with respect to the~chain maps associated to link cobordisms.
\end{theorem}

% ==================================================================
%  Section: twisting homology by bimodules
% ==================================================================

\subsection{Annular twistors}

A~family of homotopy equivalences
$\CKmod{\tau^{\underline\epsilon'}} \CKtimes{n} \CKmod{T} \simeq
	\CKmod{T} \CKtimes{m} \CKmod{\tau^{\underline\epsilon}}$
was used in Section~\ref{sec:action} to construct the~action of oriented tangles,
where $T$ is placed vertically on the~membrane.
However, one can also place $T$ horizontal and understand the~homotopy equivalences
as components of a~natural endotransformation $\tau_*$ of the~identity functor on
$\HCom^b(\DiagBimod)$, where $\DiagBimod := \FFCK(\BBN)$ is the~bicategory
of \emph{diagrammatic bimodules}. Therefore, the~pair $(\Id,\tau_*)$
is a~functor from the~endocategory $(\HCom^b(\DiagBimod),\Id)$ to itself,
and we can obtain a~new homology theory by pulling back the~shadow
$(\FFCK)^*\qHoHom_{\!0}$.
This motivates the~following definition.

\begin{definition}
	An~\emph{annular twistor} is a~family $\underline M:=\{M_{\underline\epsilon}\}$ of chain complexes of $(\CKalg n,\CKalg n)$--bimodules, parametrized by finite
	sequences $\underline\epsilon = (\epsilon_1,\dots,\epsilon_n)$ with $\epsilon_i
	\in \{-1,+1\}$, together with natural chain maps
	\begin{equation}\label{def:ann-family-isoms}
		\phi_T\colon M_{\underline\epsilon} \CKtimes{m} \CKmod{T}
			\to \CKmod{T} \CKtimes{n} M_{\underline\epsilon'},
	\end{equation}
	one per oriented $(m,n)$--tangle $T$ with orientation of input and output
	encoded by $\underline\epsilon$ and $\underline\epsilon'$ respectively.
\end{definition}

Although not stated this way, the~chain map \eqref{def:ann-family-isoms}
is a~homotopy equivalence. Up to degree shifts, its inverses can be built
from $\phi_{T^!}$ together with evaluation and coevalution maps.

Let $L$ be the~annular closure of an~$(n,n)$--tangle $T$ with orientation
of endpoints encoded by a~sequence $\underline\epsilon$.
It follow from Proposition~\ref{prop:shadow-from-nat-trans} that
$$
	\qAKh(L, \underline M) := \qHoHom_{\!\bullet}
	(\CKalg n, M_{\underline\epsilon}\CKtimes{n} \CKmod T)
$$
does not depends on the~choice of $T$. We call it the~\emph{annular homology
twisted by\/ $\underline M$}.
When $M_{\underline\epsilon} = \CKmod{\tau^{\underline\epsilon}}$ is the~twistor
described above, $\qAKh(L,\underline M)$ is the~quantum annular homology
of the~satellite of $L$ with companion $\widehat T$.

% ==================================================================
%  Section: Generalized annular homology
% ==================================================================

\subsection{Generalized annular homology}

There is another generalization of the~annular homology, which is very close to
twisting. In this section we work in characteristic two, so that the~Chen--Khovanov
construction is strictly functorial.

Let us fix a~$(1,1)$--tangle $T$ and denote by $T^!$ its mirror image.
Write $V_T:=\qAKhCom(\widehat T)$ and $V_T^* := \qAKhCom(\widehat T^!)$
for the~quantized annular chain complexes of the~annular closures of $T$ and $T^!$.
They form a~dual pair, with evaluation and coevalution maps induced by cobordisms
$T\times\fntCap$ and $T\times\fntCup$ in $(\AxI)\times I$ respectively,
see Appendix~\ref{sec:bicats} and Figure~\ref{fig:ev-for-tangle} for a~picture
of the~evaluation cobordism.

We generalize $\FqAKh\colon \BNq(\Ann) \to \gRep(\Uqsl 2)$ by assigning to essential
circles alternatively $V_T$ and $V_T^*$. An~annulus with essential boundary is then
assigned one of the~evaluation of coevaluation maps, and merging a~trivial circle to
an~essential one is determined by the~cobordisms merging an~unknot to $\widehat T$
or $\widehat T'$. It follows from the~functoriality that this produces
a~well-defined TQFT functor $\FqAKh^T\colon \BN(\Ann) \to \HCom^b(\Uqsl 2)$,
valued in the~homotopy category of representations of $\Uqsl 2$.

\begin{definition}
	Let $T$ be a~$(1,1)$--tangle and choose an~annular link $L$ with diagram $D$.
	The~\emph{$T$--annular homology} $\qAKh^T(L)$ of $L$ is the~homology of
	the~chain complex $\FqAKh^T\KhBracket{D}$.
\end{definition}

It is tempting to express the~above construction using annular twistors.
Indeed, one can define $M_{\underline\epsilon}$, where $\underline\epsilon$
has length $n$, as the~alternating tensor product
$V_T\otimes V_T^*\otimes V_T\otimes\dots$ of $n$ factors.
The~duality between $V_T$ and $V_T^*$ can be used to
construct a~chain map
\[
	\nu_T\colon M_{\underline\epsilon} \CKtimes{m} \CKmod{T}
	\to \CKmod{T} \CKtimes{n} M_{\underline\epsilon'}
\]
for every $(m,n)$--tangle $T$. However, it is not natural.
For instance, seeing $T=\fntCircle$ as a~composition of a~cup
and a~cap, the~composition
\begin{multline*}
	M_0\otimes \CKmod{\fntCup} \CKtimes{2} \CKmod\fntCap
		\to^{\nu_{\fntCup}}
	\CKmod{\fntCup} \CKtimes{2} M_2 \CKtimes{2} \CKmod\fntCap
		\to^{\nu_{\fntCap}}
	\CKmod{\fntCup} \CKtimes{2} \CKmod\fntCap \otimes M_0
\end{multline*}
is the~multiplication by the~quantum dimension of $V_T\otimes V_T^*$,
the~Jones polynomial of $\widehat T \# \widehat T'$.
On the~other hand, $\CKmod\fntCircle = \scalars\{-1\} \oplus
\scalars\{+1\}$, suggesting that $\nu_{\fntCircle} = \id$.
This seams to be the~only issue, because for the~decomposition of a~zigzag
\begin{center}
	\begin{pspicture}(-0.5,-1.5)(2.5,1.5)
		\pscustom[style=diagarc]{
			\moveto(2,-1.5)
			\lineto(2,0)
			\psellipticarc (1.5,0)(0.5,0.75){0}{180}
			\psellipticarcn(0.5,0)(0.5,0.75){360}{180}
			\lineto(0,1.5)
		}
		\psset{style=diaghline}
		\psline(-0.5, 1.5)(2.5, 1.5)
		\psline(-0.5, 0.0)(2.5, 0.0)
		\psline(-0.5,-1.5)(2.5,-1.5)
	\end{pspicture}
\end{center}
the~induced map is a~homotopy equivalence due to the~relation between
the~evaluation and coevaluation morphisms.

\appendix

\section{Background survey}\label{sec:appendix}

The~material presented here is widely known, and the~main goal of this section is to
fix the~notation. Bicategories are treated in the~excellent paper \cite{Bicats}, whereas
\cite{BasicBicats} is a~brief list of basic definitions. The~reader is also referred
to \cite{MonCats}, because many results about monoidal categories immediately
translates to bicategories. 

% ==================================================================
%  Section: Uq(sl2) and dual canonical bases
% ==================================================================

\subsection{Representations of \texorpdfstring{$\Uqsl{2}$}{Uq(sl2)}}
\label{sec:sl2-rep}

As usual we fix a~commutative unital ring $\scalars$ together with
an~invertible element $q$.
By definition, $\Uqsl{2}$ is the~unital associative $\scalars$--algebra
with generators $E$, $F$, $K$, $K^{-1}$ and relations
\begin{align*}
	KE &= q^2EK,	&
	KK^{-1} &= 1 = K^{-1}K, \\
	KF &= q^{-2}FK, &
	K-K^{-1} &= (q-q^{-1})(EF-FE).
\end{align*}
It is a~Hopf algebra with the~comultiplication
$\Delta\colon \Uqsl{2}\to\Uqsl{2}\otimes\Uqsl{2}$,
the~counit $\epsilon\colon\Uqsl{2}\to\scalars$,
and the~antipode $S\colon\Uqsl{2}\to\Uqsl{2}$ defined by
\begin{align}
	  \Delta(E) &=  E\otimes K+1\otimes E, &
	\epsilon(E) &=  0,	   &
	       S(E) &= -EK^{-1},
\\
	  \Delta(F) &=  F\otimes 1+K^{-1}\otimes F, &
	\epsilon(F) &=  0,  &
	       S(F) &= -KF,
\\
	  \Delta(K^{\pm 1}) &= K^{\pm 1}\otimes K^{\pm 1}, &
	\epsilon(K^{\pm 1}) &= 1, &
	       S(K^{\pm 1}) &= K^{\mp 1}.
\end{align}
Using this Hopf algebra structure, we can regard the~category of finite-dimensional
representations of $\Uqsl{2}$ as a~monoidal category with duals. The unit in this
monoidal category is given by the~trivial representation $V_0=\scalars$, on which
$\Uqsl{2}$ acts by multiplication by $\epsilon(X)$ for any $X\in\Uqsl{2}$.

We write $V_1 := \Span_{\scalars}\{v_1,v_{-1}\}$ and $V_1^* :=
\Span_{\scalars}\{v_1^*, v_{-1}^*\}$ for the~fundamental representation and its
dual. We identify both with the~rank two module $V := \Span_{\scalars}\{v_+,v_-\}$
using the~isomorphisms
\begin{equation}\label{eq:isom-V-V1-V1*}
	\setlength\arraycolsep{0pt}
	\left\{\begin{array}{rl}
		v_1   &\mapsto v_+ \\
		v_{-1}&\mapsto v_-
	\end{array}\right.
\qquad\text{and}\qquad
	\left\{\begin{array}{rl}
		v^*_1   &\mapsto v_- \\
		v^*_{-1}&\mapsto q^{-1}v_+
	\end{array}\right..
\end{equation}
This equips $V$ with two actions of $\Uqsl 2$ that differ by signs:
\begin{equation}
	\begin{tabular}{llp{1cm}ll}
		\multicolumn{2}{c}{The~action on $V\cong V_1$} && \multicolumn{2}{c}{The~action on $V\cong V_1^*$}\\
		\cline{1-2} \cline{4-5}
		 $Ev_+ = 0$		& $Ev_- = v_+$		&& $Ev_+ =  0$   & $Ev_- = -v_+$ \\
	     $Fv_+ = v_-$	& $Fv_- = 0$			&& $Fv_+ = -v_-$ & $Fv_- = 0$    \\
		 $Kv_+ = qv_+$	& $Kv_- = q^{-1}v_-$	&& $Kv_+ = qv_+$ & $Kv_- = q^{-1}v_-$
	\end{tabular}
\end{equation}
The~duality between $V_1$ and $V_1^*$ comes with the~evaluation and coevaluation maps
\begin{equation}\label{eq:ev-coev-for-V}
	\begin{tabular}{llp{1cm}l}
		\multicolumn{2}{c}{$\ev\colon V \otimes V \to \scalars$}
		&&
		\multicolumn{1}{c}{$\coev\colon \scalars \to V \otimes V$}
		\\
		\cline{1-2} \cline{4-4}
		$\tensor{v+,v+}\mapsto 0$ & $\tensor{v+,v-}\mapsto q$
		&\rule{0pt}{2.5ex}&
		$1\mapsto \tensor{v+,v-}+q^{-1}\tensor{v-,v+}$
		\\
		$\tensor{v-,v-}\mapsto 0$ & $\tensor{v-,v+}\mapsto 1$
	\end{tabular}
\end{equation}
that intertwine the~action $\Uqsl 2$ if $V\otimes V$ is identified with either $V_1\otimes V_1^*$ or $V_1^*\otimes V_1$.

The~full subcategory of $\Uqsl 2$ generated by tensor powers of $V$ admits
a~graphical representation. Let $\TL$ be the~\emph{Temperly--Lieb category},
the~linear category with objects finite collections of points on a~real line and
morphisms generated by flat loopless tangles, i.e.\ collection of disjoint 
intervals in $\RxI$ with endpoints on the~boundary lines. Composition is defined
by stacking pictures one onto another and trading each closed component for
$q+q^{-1}$. For example,
\begin{equation}
	\psset{yunit=3}
	\begin{centerpict}(-1.25,0)(4.25,1)
		\psframe[style=mysurface,linestyle=none](-0.75,0)(3.75,1)
		\psline[style=diaghline](-0.75,1)(3.75,1)
		\psline[style=diaghline](-0.75,0)(3.75,0)
		\psset{style=diagarc,fillstyle=none}
		\psline(0,0)(0,1)
		\psbezier(1,0)(1,0.5)(3,0.5)(3,1)
		\psellipticarc(2.5,0)(0.5,0.25){0}{180}
		\psellipticarc(1.5,1)(0.5,0.25){180}{360}
	\end{centerpict}
	\circ
	\begin{centerpict}(-1.25,0)(4.25,1)
		\psframe[style=mysurface,linestyle=none](-0.75,0)(3.75,1)
		\psline[style=diaghline](-0.75,1)(3.75,1)
		\psline[style=diaghline](-0.75,0)(3.75,0)
		\psset{style=diagarc,fillstyle=none}
		\psbezier(1,0)(1,0.5)(0,0.5)(0,1)
		\psbezier(2,0)(2,0.5)(1,0.5)(1,1)
		\psellipticarc(2.5,1)(0.5,0.25){180}{360}
	\end{centerpict}
	=
	(q+q^{-1})
	\begin{centerpict}(-1,0)(4.25,1)
		\psframe[style=mysurface,linestyle=none](-0.75,0)(3.75,1)
		\psline[style=diaghline](-0.75,1)(3.75,1)
		\psline[style=diaghline](-0.75,0)(3.75,0)
		\psset{style=diagarc,fillstyle=none}
		\psbezier(1,0)(1,0.5)(0,0.5)(0,1)
		\psbezier(2,0)(2,0.5)(3,0.5)(3,1)
		\psellipticarc(1.5,1)(0.5,0.25){180}{360}
	\end{centerpict}
\end{equation}
There is a~functor $\FTL\colon \TL\to \Rep(\Uqsl 2)$ that assigns $V^{\otimes n}$
to a~collection of $n$ points, whereas caps and cups are sent to the~evaluation
and evaluation homomorphisms. It is known that $\FTL$ is faithful \cite{TLFunctor}.

\subsection{Knots and tangles}
\label{sec:def-tangles}

Let $M$ be an~oriented smooth $3$--manifold. A~proper 1--submanifold $T\subset M$ is called a~\emph{tangle}. We call it a~\emph{link} if it has no boundary, and a~\emph{knot} if in addition it has one component. All tangles and links in this paper are assumed to be oriented unless stated otherwise. An~\emph{isotopy of tangles} $T$ and $T'$ is a~smooth map $\Phi\colon M\times I\to M$ such that each $\Phi_t := \Phi(\blank, t)$ is a~diffeomorphism fixed at the~boundary, $\Phi_0 = \id$, and $\Phi_1(T) = T'$. If $T$ and $T'$ are oriented, then we require that the~orientation is preserved by $\Phi_1$.

Denote by $-T$ the~~tangle $T$ with reversed orientation of all its components. A~\emph{tangle cobordism} from a~tangle $T_0$ to $T_1$ is an~oriented surface $S\subset M\times I$ with boundary %$\partial S = T_0\times\{0\} \cup -T_1\times\{1\} \cup (\partial T_0\times I)$.
$\partial S = -T_0\times\{0\} \cup T_1\times\{1\} \cup (-\partial T_0\times I)$.
We shall consider cobordisms only up to an~isotopy, in which case they form a~category: composition is given by gluing cobordisms, and the~identity morphism on a~tangle $T$ is represented by the~cylinder $T\times I\subset M\times I$.

When $M= \sfce\times I$ is a~thickened surface, then isotopy classes of oriented tangles in $M$ form a~category, with the~product induced by stacking, $M\cup M \cong M$, and tangles with tangle cobordisms form a~2--category.

\begin{notation*}
	We shall write $\sLinks(M)$ for the~set of isotopy classes of oriented links in $M$, and $\Links(M)$ for the~category of oriented links in $M$ and cobordisms between them. Isotopy classes of oriented tangles in a~thickened surface $\sfce\times I$ form a~category $\Tan(F)$ with the~composition induced by stacking, and similarly cobordisms between tangles in $\sfce\times I$ form a~2--category $\TTan(F)$. We write simply $\Tan$ and $\TTan$ when $\sfce=\RxI$.
\end{notation*}

\wrapfigure[r](0pt,0.75ex)<0>{\psset{unit=8mm}\begin{pspicture}(-0.5,-1ex)(2.5,3)
	\psframe[style=mysurface,linestyle=none](-0.5,0)(2.5,3)
	\psset{style=diagarc,fillstyle=none}
	\rput{30}(0.2,1.7){\psellipticarc(0,0)(0.4,0.25){0}{180}}
	\psecurve(1,0)(1.5,1.2)(2,2)(1,3)
	\begin{psset}{border=3\pslinewidth,bordercolor=sfceFillColor}%
		\psecurve(0.5, 4)(0,3)(0.75,1)(2,0.8)(1.3,2)(2,3)(1.7,4)
		\rput{30}(0.2,1.7){\psellipticarc(0,0)(0.4,0.25){180}{360}}
		\psecurve(1.1,-1)(1,0)(1.5,1.2)(2,2)
		\psecurve(1.5,1.2)(2,2)(1,3)(1.1,4)
	\end{psset}
	\psset{style=diaghline}
	\psline[style=diaghline](-0.5,0)(2.5,0)
	\psline[style=diaghline](-0.5,3)(2.5,3)
\end{pspicture}}
Assume $M$ is a~line bundle over a~surface $\sfce$ and consider the~projection $M\to \sfce$ onto the~zero section. It maps a~generic tangle $T$ to an~immersed collection of intervals and circles $\tilde T\subset \sfce$ with only finitely many multiple points, each a~transverse intersection of two arcs. A~\emph{diagram} of $T$ is constructed from $\tilde T$ by breaking one of the~arcs at each double point as follows. When $\sfce$ is oriented, then fibers of $M$ admit a~canonical orientation and we break the~lower arc at each double point of $\tilde T$ (see an~example above for $\sfce =\RxI$). In case $\sfce$ is nonorientable, choose a~minimal collection of curves $\gamma$ that cuts $\sfce$ into an~orientable surface. Then there is a~normal field over $\sfce - \gamma$ and we can construct the~diagram as before. Then Reidemeister moves and planar isotopies relate diagrams between isotopic tangles if a~crossing is switched when moved through $\gamma$:
\begin{equation*}
\psset{unit=1cm}%
\begin{centerpict}(-1,-1)(1,1)
	\psframe[linestyle=none,fillstyle=solid,fillcolor=sfceFillColor](-1,-1)(1,1)
	\begin{psset}{style=diagarc}
		\pscustom{%
			\moveto(-1,-0.5)
			\curveto(-0.3,-0.5)(-0.7,0.5)(0,0.5)
			\lineto(1,0.5)
			\stroke
		}%
		\pscustom{
			\moveto(-1,0.5)
			\curveto(-0.3,0.5)(-0.7,-0.5)(0,-0.5)
			\lineto(1,-0.5)
			\stroke[linecolor=sfceFillColor,linewidth=5\pslinewidth]
		}%
	\end{psset}
	\psline[linecolor=red](0,1)(0,-1)
	\rput[l](0.1,0){\textcolor{red}{$\scriptstyle\gamma$}}
\end{centerpict}
\quad=\quad
\begin{centerpict}(-1,-1)(1,1)
	\psframe[linestyle=none,fillstyle=solid,fillcolor=sfceFillColor](-1,-1)(1,1)
	\begin{psset}{style=diagarc}
		\pscustom{%
			\moveto(1,-0.5)
			\curveto(0.3,-0.5)(0.7,0.5)(0,0.5)
			\lineto(-1,0.5)
			\stroke
		}%
		\pscustom{
			\moveto(1,0.5)
			\curveto(0.3,0.5)(0.7,-0.5)(0,-0.5)
			\lineto(-1,-0.5)
			\stroke[linecolor=sfceFillColor,linewidth=7\pslinewidth]
		}%
	\end{psset}
	\psline[linecolor=red](0,1)(0,-1)
	\rput[r](-0.1,0){\textcolor{red}{$\scriptstyle\gamma$}}
\end{centerpict}
\end{equation*}
This follows, because the~normal field is reversed at points of $\gamma$.

\subsection{Constructions on categories}
\label{sec:categories}

Below we review definitions of certain constructions on categories
that appear throughout the~paper.

\subsubsection*{Additive closure}

We say that a~linear category is \emph{additive} if it has finite
direct sums. Each category $\cat C$ admits the~\emph{additive closure}
$\catAdd{\cat C}\!$, the~smallest additive category containing $\cat C$.
It is constructed by introducing formal direct sums:
\begin{itemize}
	\item an~object of $\catAdd{\cat C}$ is a~finite sequence $(x_1,\dots,x_r)$ with $x_i\in \Ob(\cat C)$, possibly empty,
	\item a~morphism from $(x_1,\dots,x_r)$ to $(y_1,\dots,y_s)$ is a~matrix
	$(f_{ij})$ of morphisms $f_{ij}\in \cat C(x_i,y_j)$, and
	\item composition is defined by the~matrix multiplication rule.
\end{itemize}
There is an~inclusion $\cat C\to \catAdd{\cat C}$, which takes an~object $x$
to the~1-element sequence $(x)$. It is an~equivalence of categories if
$\cat C$ is already additive.

\subsubsection*{Idempotent completion}

An~endomorphism $p\in \cat C(x,x)$ satisfying $p\circ p = p$ is an~\emph{idempotent}.
We say that $p$ \emph{splits} if it decomposes $p = s\circ r$ such that $r\circ s$
is an~identity morphism. In such a~case $r$ is an~epimorphism and its codomain
is called the~\emph{image of\/ $p$}.

A~category is \emph{idempotent complete} if all its idempotents split.
Each category $\cat C$ admits its \emph{idempotent completion} $\Kar(\cat C)$,
also called also the~\emph{Karoubi envelope of\/ $\cat C$}, which is the~smallest
idempotent complete category containing $\cat C$. It is constructed by taking
all idempotents of $\cat C$ as objects, and defining morphisms from $e$ to $e'$
as those morphisms $f$ from $\cat C$ that $e'\circ f\circ e$ is well-defined
and equal to $f$. The~identity morphism on $e$ is given by the~idempotent itself.

\subsubsection*{Formal complexes and homotopy category}

A~\emph{formal complexes} over a~category $\cat C$ is a~sequence of objects and
morphisms from $\cat C$
\begin{equation}
	(C^\bullet, d) = \left( \dots \to C^i \to^{d^i} C^{i+1} \to^{d^{i+1}} C^{i+2} \to \dots \right)
\end{equation}
satisfying $d^{i+1}d^i = 0$ at each place. The~morphisms $d^i$ are called
the~\emph{differential}. We say that $(C^\bullet, d)$ is \emph{bounded}
if $C^i = 0$ except finitely many indices.

A~\emph{formal chain map} from $(C^\bullet, d)$ to $(D^\bullet, d)$ is
a~collection of morphisms $f^i \colon C^i \to D^i$ fitting into a~commuting ladder
\begin{equation}
	\begin{diagps}(-3em,-0.5ex)(14em,10ex)
		\hsquares(0,0)<5em,8ex>|a`a`a`a`a`a`a|[%
			C^i`C^{i+1}`C^{i+2}`
			D^i`D^{i+1}`D^{i+2};
			d`d`
			f^i`f^{i+1}`f^{i+2}`
			d`d]
		\node uld(-3em,8ex)[\dots]	\node urd(14em,8ex)[\dots]
		\node dld(-3em,0ex)[\dots]	\node drd(14em,0ex)[\dots]
		\arrow{->}[uld`ul;]			\arrow{->}[ur`urd;]
		\arrow{->}[dld`dl;]			\arrow{->}[dr`drd;]
	\end{diagps}
\end{equation}
Finally, a~\emph{formal chain homotopy} from $f^\bullet$ to $g^\bullet$,
both chain maps from $(C^\bullet,d_C)$ to ($D^\bullet, d_D)$, is a~collection
of morphisms $h^i\colon C^i\to D^{i-1}$ satisfying
$d^{i-1}_D\circ h^i + h^{i+1}\circ d^i_C = g^i - f^i$. In such case we say
that $f^\bullet$ and $g^\bullet$ are \emph{homotopic}, which we write
$f\sim g$.

Formal complexes (reps.\ bounded formal complexes) and chain maps modulo chain
homotopies constitute the~\emph{homotopy category of complexes} $\HCom(\cat C)$
(resp.\ $\HCom^b(\cat C)$). Isomorphism in these categories are called
\emph{homotopy equivalences} and we usually write $C^\bullet \simeq D^\bullet$
for complexes that are homotopically equivalent.

The~categories $\HCom(\cat C)$ and $\HCom^b(\cat C)$ are \emph{triangulated}
\cite{MethodsOfHA, HomAlg}, which means that they come with a~\emph{homological
degree shift functor} and a~collection of \emph{distinguish triangles} satisfying
certain axioms. The~degree shift functor is usually denoted by $[1]$ and it shifts
a~complex leftwards, negating the~differential at the~same time:
\begin{equation}
	C[1]^i := C^{i+1},
\qquad
	d[1]^i := -d^{i+1}.
\end{equation}
Distinguished triangles are of the~form
\begin{equation}
	C^\bullet
		\to^{f^\bullet}
	D^\bullet
		\to^{\mathit{in}^\bullet}
	\cone^\bullet(f)
		\to^{\mathit{pr}^\bullet}
	C[1]^\bullet,
\end{equation}
where $\cone^\bullet(f)$ stands for the~\emph{mapping cone} of $f\noic$,
the~formal complex
\begin{equation}
	\cone^i(f^\bullet) := C^{i+1}\oplus D^i,
\qquad
	d^i = \begin{pmatrix}
		-d^{i+1}_C & 0 \\
		f^i & d^i_D
	\end{pmatrix}.
\end{equation}
The~morphisms $\mathit{in}^i\colon D^i\to \cone^i(f^\bullet)$
and $\mathit{pr}^i\colon\cone^i(f^\bullet)\to C[1]^i = C^{i+1}$
are the~inclusion and projection respectively.

%% ============================================================================
%%  Section: Bicategories, duals, and adjoints
%% ============================================================================

\subsection{Bicategories}\label{sec:bicats}

%
% Definition: a bicategory
%
A~\emph{bicategory} $\ccat C$ is a~`higher level' analogue of a~category.
It consists of
\begin{itemize}
	\item a~class of objects $\Ob(\ccat C)$,
	\item	a~category $\ccat C(x,y)$ for each pair of objects $(x,y)$, whose objects
	and morphisms are called \emph{1--} and \emph{2--morphisms} respectively
	and represented by single and double arrows,
	\item a~unit $\id_x\in \ccat{C}(x,x)$ for each object $x$,
	\item a~functor $\circ: \ccat{C}(y,z) \times \ccat{C}(x,y) \to \ccat{C}(x,z)$
	for each triple of objects $(x,y,z)$, and
	\item natural isomorphisms
	\begin{equation}\label{eq:can-isoms-for-bicat}
		\mathfrak{a} \colon f\circ (g\circ h) \dblto^\cong (f\circ g)\circ h
			\qquad
		\mathfrak{l} \colon \id_y \circ f \dblto^\cong f
			\qquad
		\mathfrak{r} \colon f \circ \id_x \dblto^\cong f
	\end{equation}
	called \emph{associators} and \emph{unitors}, which satisfy the~pentagon
	and triangle axioms \cite{Bicats}.
\end{itemize}
A~bicategory is called \emph{strict} or a~\emph{2--category} if the~natural
isomorphisms are identities. We call $\ccat C$ a~\emph{locally small} bicategory
when each category $\ccat C(x,y)$ is small, and $\ccat C$ is \emph{small} when
also $\Ob(\ccat C)$ is a~set.

Associators and unitors are often omitted for clarity. According to the~MacLane's
Coherence Theorem \cite[Chapter VII.2]{Cats} there is only one way how to insert
these isomorphisms back when necessary.

%
% Notation: vertical composition, identity 1-morphisms
%
\begin{notation*}
\wrapfigure[r]<2>{\begin{diagps}(-2.5em,-4ex)(2.5em,5.5ex)
	\node dom(-2em,0)[x]
	\node cod( 2em,0)[y]
%	\carrow[angleA= 45,angleB= 135,arm=0.7,linearc=1.2]|a|{->}[dom`cod;\Rnode{f}{\scriptstyle f}]
%	\carrow[angleA=-45,angleB=-135,arm=0.7,linearc=1.2]|b|{->}[dom`cod;\Rnode{g}{\scriptstyle g}]
	\carrow*[arcangle= 50]|a|{->}[dom`cod;\Rnode{f}{\scriptstyle f}]
	\carrow*[arcangle=-50]|b|{->}[dom`cod;\Rnode{g}{\scriptstyle g}]
	\arrow[nodesep=1.5ex]|a{npos=0.55,labelsep=1pt}|{=>}[f`g;\alpha]
\end{diagps}}
	In this paper we denote categories with calligraphic letters $\cat C$, $\cat D$, etc., whereas bold letters $\ccat C$, $\ccat D$, etc.\ are reserved for bicategories. Identity morphisms are written as $\id_x$, and identity 2--morphisms as $\mathbf 1_f$. If $\ccat C$ is a~bicategory, then the~composition in $\ccat C(x,y)$ is denoted by $\ast$ and called \emph{vertical}, whereas $\circ$ is the~\emph{horizontal composition}. These come from the~common convention to draw 1--morphisms horizontally and 2--morphisms vertically, see the~diagram to the~right of a~2--morphism $\alpha\colon f\dblto<1em> g$.
\end{notation*}

%
% Definition: bifunctor
%
Choose bicategories $\ccat C$ and $\ccat D$.
A~\emph{bifunctor} $\ccat F\colon \ccat C\to \ccat D$ consists of a~function
of objects $\Ob(\ccat C) \to \Ob(\ccat D)$ and a~collection of functors
$\ccat C(x,y) \to \ccat D(\ccat Fx,\ccat Fy)$, together with natural 1--morphisms
\begin{equation}
	\mathfrak m\colon \ccat F(g) \circ \ccat F(f) \dblto \ccat F(g\circ f)
\qquad
	\mathfrak i\colon \id_{\ccat Fx} \dblto \ccat F(\id_x)
\end{equation}
satisfying certain coherence axioms. They are called \emph{morphisms of
bicategories} in \cite{Bicats}, whereas the~word \emph{homomorphism} is reserved
for the~case when $\mathfrak m$ and $\mathfrak i$ are invertible. In such case
we say that $\ccat F$ is a~\emph{strong bifunctor}.

%
% Definition: natural transformation
%
Choose bifunctors $\ccat F,\ccat G\colon \ccat C\to \ccat D$.
A~\emph{natural transformation} $\eta\colon \ccat F \to \ccat G$
is a~collection of 1--morphisms $\eta_x\colon \ccat F(x) \to \ccat G(x)$,
one per object $x\in\ccat C$, and 2--morphisms
$\eta_f \colon \ccat G(f) \circ \eta_x \dblto \eta_y \circ \ccat F(f)$,
one per 1--morphisms $f\in\ccat C(x,y)$, such that
\begin{equation}
	\begin{diagps}(0,0)(10em,13ex)
		\square<7em,11ex>{->``->`->}[%
			\ccat F(x)`\ccat G(x)`\ccat F(y)`\ccat G(y);
			\eta_x``\ccat G(f)`\eta_y
		]
		\carrow*[arcangle= 45]{->}[v0`v2;\ccat F(f)]
		\carrow*[arcangle=-45]|b|{->}[v0`v2;\ccat F(g)]
		\morphism[nodesep=1.6em]{=>}|a{npos=0.53}|(6.5em,10.5ex)<-4.5em,-9ex>[`;\eta_f]
		\morphism{<=}(-1em,5.5ex)<2em,0ex>[`;\ccat F(\alpha)]
	\end{diagps}
	=
	\begin{diagps}(-3em,0)(7em,13ex)
		\square<7em,11ex>{->`->``->}[%
			\ccat F(x)`\ccat G(x)`\ccat F(y)`\ccat G(y);
			\eta_x`\ccat F(g)``\eta_y
		]
		\carrow*[arcangle= 45]{->}[v1`v3;\ccat G(f)]
		\carrow*[arcangle=-45]|b|{->}[v1`v3;\ccat G(g)]
		\morphism[nodesep=1.6em]{=>}|b{npos=0.53}|(5em,9.5ex)<-4.5em,-9ex>[`;\eta_g]
		\morphism{<=}(6em,5.5ex)<2em,0ex>[`;\ccat G(\alpha)]
	\end{diagps}
\end{equation}
for every 2--morphism $\alpha\colon f\dblto g$. Moreover, $\eta_f$ must be
coherent with all the~other canonical 2--isomorphisms (associators, unitors,
the~structure 2--isomorphisms of $\ccat F$ and $\ccat G$), see \cite{BasicBicats}.
We say that $\eta$ is \emph{strong} if each $\eta_f$ is invertible.

%
% Definition: modification
%
Finally, let $\eta,\nu\colon \ccat F\to \ccat G$ be two natural transformations.
A~\emph{modification} $\Gamma\colon \eta\to \nu$ is a~collection of 1--morphisms
$\Gamma_x\colon \eta_x\to \nu_x$, such that
$\nu_f\ast (\iid_{\ccat Gf}\circ\Gamma_x)
	= (\Gamma_y\circ\iid_{\ccat Ff})\ast \eta_f$ 
for every 1--morphism $f\colon x\to y$.

\subsubsection*{Duals}

%
% Definition: dual/adjoint morphisms
%
A~bicategory $\ccat{C}$ has \emph{left duals} if each $f\in \ccat{C}(x,y)$ admits ${}^*\!f\in \ccat{C}(y,x)$ together with coevaluation and evaluation 2--morphisms
\begin{align}
\label{def:left-duals}
	\id_y \to/=>/^{\mathit{coev}\,} {}^*\!f\circ f
		\hskip 2em&\hskip 2em
	f\circ {}^*\!f\to/=>/^{\mathit{ev}\,} \id_x\\
\intertext{fitting into commuting triangles}
\label{rel:left-duals}
	\begin{diagps}(-0.5em,-0.5ex)(9.5em,12.5ex)
		\qtriangle<7em,10ex>|a{npos=0.4}`b`*c|{=>`=>`=>}[%
			f`f \circ {}^*\!f \circ f`f;
			\mathbf 1_f\mcircm\mathit{coev}`
			\mathbf 1_f`
			\mathit{ev}\mcircm\mathbf 1_f%
		]
	\end{diagps}
		\hskip 2em&\hskip 2em
	\begin{diagps}(-0.5em,-0.5ex)(9.5em,12.5ex)
		\qtriangle<7em,10ex>|a{npos=0.4}`b`*c|{=>`=>`=>}[%
			{}^*\!f`{}^*\!f \circ f\circ {}^*\!f`{}^*\!f;
			\mathit{coev}\mcircm\mathbf 1_{{}^*\!f}`
			\mathbf 1_{{}^*\!f}`
			\mathbf 1_{{}^*\!f}\mcircm\mathit{ev}%
		]
	\end{diagps}
\end{align}
where for clarity associators and unitors are omitted. The~morphism ${}^*\!f$ is called the~\emph{left dual} to $f$. We define the~\emph{right dual} $f^*$ of $f$ by reversing the~order of the~horizontal composition in \eqref{def:left-duals} and \eqref{rel:left-duals}. If a~dual 1--morphism exists, then it is unique up to an~isomorphism. In particular, dual pairs are preserved by strong bifunctors.

\subsubsection*{Examples}

Small categories, functors, and natural transformations form a~strong bicategory
$\ccat{Cat}$. A~left (resp.\ right) dual to a~functor $F$ is its left (resp.\ 
right) adjoint. Hence, not all 1--morphisms in $\ccat{Cat}$ are dualizable.

Tangles in a~thickened surface $\sfce\times I$ constitute a~bicategory
$\TTan(\sfce)$, which objects are finite collections of points in $\sfce$,
1--morphisms are tangles with endpoints on $\sfce\times\partial I$,
and 2--morphisms are tangle cobordisms. This bicategory has both left and
right duals. Indeed, the~mirror image $T^!$ of a~tangle $T$, obtained by
flipping $\sfce\times I$, is both the~left and right dual of $T$.
The~evaluation 2-morphism is obtained by revolving $T$ in four dimensions
along the~input surface $\sfce\times\{0\}$, i.e.\ it is the~image of the~map
$(p,t,s)\mapsto (p,t\cos(s\pi),t\sin(s\pi))$ with $(p,t)\in T$ and $s\in I$,
suitably normalized (see Fig.~\ref{fig:ev-for-tangle}). The~coevaluation is
defined dually by a~rotation along the~output surface $\sfce\times\{1\}$.

\begin{figure}[h]%
	\centering
	\psset{unit=7mm}%
	\begin{pspicture}(-4,-1)(4,4)
		\def\mytangle[#1]{%
			\psframe[framearc=0.5,fillstyle=solid,#1](-1.5,2)(1.5,1)
			\rput(0,1.5){$T$}
		}%
		\psline(-2,-1)(-0.5,-0.5) \psline(2.5,0.5)(4,1)
		\psline(-3,-1)( 3, 1)
		\psline(-4,-1)(-2.5,-0.5) \psline(0.5,0.5)(2,1)
		\pstilt{18}{\mytangle[linewidth=1.5\pslinewidth]}%
		\pstilt{198}{\mytangle[linewidth=1.5\pslinewidth]}%
		\psline(-1, 3)(-1, 1.5)
		\psline[linestyle=none,border=2\pslinewidth](0,0.2)(0,1)
		\psline( 0, 3)( 0, 0)
		\psline( 1, 3)( 1, 1.5)
		\psline[border=2\pslinewidth](-3.4,-0.8)(-4,-1)(-4, 2)(2, 4)(2, 1)(1.4,0.8)\psline(-1, 3)(-1, 2.5)
		\psline[border=2\pslinewidth](-2.4,-0.8)(-3,-1)(-3, 2)(3, 4)(3, 1)(2.4,0.8)\psline( 0, 3)( 0, 2.5)
		\psline[border=2\pslinewidth](-1.4,-0.8)(-2,-1)(-2, 2)(4, 4)(4, 1)(3.4,0.8)\psline( 1, 3)( 1, 2.5)
		\begin{psset}{linewidth=0.5\pslinewidth,linestyle=dashed,dash=2pt 1pt}
			\rput(-2.78,0){% outer left
				\psbezier{-c}(3.21,0.59)(3.21,1.59)(2.61,2.31)(1.41,1.91)
				\psbezier{c-}(1.41,1.91)(0.21,1.51)(-0.46,0.41)(-0.46,-0.59)
			}
			\rput(-2.79,0){% inner left
				\psbezier{-c}(2.41,0.34)(2.41,0.91)(2.09,1.32)(1.41,1.09)
				\psbezier{c-}(1.41,1.09)(0.73,0.86)(0.34,0.23)(0.34,-0.34)
			}%
		\end{psset}
		\mytangle[]
		\begin{psset}{linewidth=0.5\pslinewidth,linestyle=dashed,dash=2pt 1pt}
			% outer right
			\psbezier{-c}(3.21,0.59)(3.21,1.59)(2.61,2.31)(1.41,1.91)
			\psbezier{c-}(1.41,1.91)(0.21,1.51)(-0.46,0.41)(-0.46,-0.59)
			% inner right
			\psbezier{-c}(2.41,0.34)(2.41,0.91)(2.09,1.32)(1.41,1.09)
			\psbezier{c-}(1.41,1.09)(0.73,0.86)(0.34,0.23)(0.34,-0.34)
		\end{psset}
	\end{pspicture}
	\caption{The~evaluation cobordism $TT^! \protect\dblto \mathbf 1_3$ for a~tangle $T\in\TTan(1,3)$.}%
	\label{fig:ev-for-tangle}%
\end{figure}

Rings, bimodules, and bimodule maps form a~bicategory $\BBimod$,
with horizontal composition given by the~tensor product: $N\circ M := M\utimes{B}N$
for an~$(A,B)$--bimodule $M$ and a~$(B,C)$--bimodule $N$. This formula comes
from interpreting an~$(A,B)$--bimodule $M$ as a~functor
$(\blank)\utimes{A} M \colon \Mod(A)\to \Mod(B)$ between the~categories
of right modules. This bicategory does not have duals. Indeed,
an~$(A,B)$--bimodule $M$ has a~left (resp.\ right) dual if and only if it is
finitely generated and projective as a~right $B$--module (resp.\ left $A$--module).
If so, the~left and right dual modules are given as the~modules of right $B$--linear
and left $A$--linear morphisms respectively:
\[
	\ldual M := \Hom_B(M,B)
\qquad
	M^* := \Hom_A(M,A)
\]
For this reason we usually restrict either to $\RRep$ (bimodules with left duals) or $\ccat{Birep}$ (bimodules with both left and right duals).

\def\authordata#1#2#3{\begingroup
	\def\\{,}%
	\bigskip
	\noindent\small\textit{#1:} \textsc{#2}\newline\texttt{email: #3}
\endgroup}

\authordata{A.B.}{%
	University of Z\"urich\\
	Institute of Mathematics\\
	8057 Z\"urich\\
	Switzerland
}{anna@math.uzh.ch}

\authordata{K.K.P.}{%
	ETH Z\"urich\\
	Institute for Theoretical Studies\\
	8096 Z\"urich\\
	Switzerland
}{krzysztof.putyra@eth-its.ethz.ch}

\authordata{S.M.W.}{%
	Syracuse University\\
	Department of Mathematics\\
	Syracuse, NY 13244
}{smwehrli@syr.edu}

\end{document}